\documentclass[a4paper, 11pt, english]{article}
\usepackage{graphicx}
\usepackage{stmaryrd}
\usepackage[a4paper]{geometry}
\usepackage{amsmath,amsfonts,amssymb,amsthm,epsfig,epstopdf,url,array}
\numberwithin{equation}{section}
\usepackage{rotating}

\usepackage[colorlinks=true,citecolor=red,linkcolor=blue,pdfpagetransition=Blinds,backref=page]{hyperref}
\usepackage{xcolor}
\definecolor{citesky}{RGB}{70,130,180}
\usepackage{aliascnt}
\usepackage{cleveref}
\usepackage{orcidlink}

\theoremstyle{plain}

\newtheorem{proposition}{Proposition}[section]
\newaliascnt{theo}{proposition}
\newtheorem{theo}[theo]{Theorem}
\aliascntresetthe{theo}

\newaliascnt{lemma}{proposition}
\newtheorem{lemma}[lemma]{Lemma}
\aliascntresetthe{lemma}

\newaliascnt{defprop}{proposition}

\aliascntresetthe{defprop}

\newaliascnt{cor}{proposition}
\newtheorem{cor}[cor]{Corollary}
\aliascntresetthe{cor}

\newaliascnt{rmk}{proposition}
\newtheorem{rmk}[rmk]{Remark}
\aliascntresetthe{rmk}

\theoremstyle{definition}

\newaliascnt{defi}{proposition}
\newtheorem{defi}[defi]{Definition}
\aliascntresetthe{defi}

\newaliascnt{app}{proposition}

\aliascntresetthe{app}

\newaliascnt{claim}{proposition}

\aliascntresetthe{claim}

\newaliascnt{ass}{proposition}

\aliascntresetthe{ass}

\crefname{proposition}{Proposition}{Propositions}
\Crefname{proposition}{Proposition}{Propositions}

\crefname{theo}{Theorem}{Theorems}
\Crefname{theo}{Theorem}{Theorems}

\crefname{lemma}{Lemma}{Lemmas}
\Crefname{lemma}{Lemma}{Lemmas}

\crefname{defprop}{Definition--Proposition}{Definition--Propositions}
\Crefname{defprop}{Definition--Proposition}{Definition--Propositions}

\crefname{cor}{Corollary}{Corollaries}
\Crefname{cor}{Corollary}{Corollaries}

\crefname{rmk}{Remark}{Remarks}
\Crefname{rmk}{Remark}{Remarks}

\crefname{defi}{Definition}{Definitions}
\Crefname{defi}{Definition}{Definitions}

\crefname{app}{Application}{Applications}
\Crefname{app}{Application}{Applications}

\crefname{claim}{Claim}{Claims}
\Crefname{claim}{Claim}{Claims}

\crefname{ass}{Assumption}{Assumptions}
\Crefname{ass}{Assumption}{Assumptions}

\usepackage{enumitem}

\usepackage{comment}
\Crefname{paragraph}{Section}{Sections}
\usepackage{fancyhdr}
\usepackage{cases}
\usepackage{mathrsfs}

\usepackage{pgfplots}

\usepackage{fullpage}
\newcommand\inter[1]{\left\llbracket #1\right\rrbracket}

\newcommand{\TT}{\mathbb T^3}

\DeclareMathOperator{\diver}{div}

\providecommand{\keywords}[1]{\noindent {\textit{Keywords:}} #1}

\def\dt{\textnormal{d}t}
\def\d{\textnormal{d}}

\newcommand{\vertiii}[1]{{\left\vert\kern-0.25ex\left\vert\kern-0.25ex\left\vert #1 
		\right\vert\kern-0.25ex\right\vert\kern-0.25ex\right\vert}}

\makeatletter
\let\original@addcontentsline\addcontentsline
\newcommand{\dummy@addcontentsline}[3]{}
\newcommand{\DeactivateToc}{\let\addcontentsline\dummy@addcontentsline}
\newcommand{\ActivateToc}{\let\addcontentsline\original@addcontentsline}
\makeatother

\newcommand{\norm}[1]{\left\|#1\right\|}
\begin{document}
	
	\title{Global strong solutions of the three-dimensional Cahn--Hilliard--Navier--Stokes system via finite-dimensional scalar forcing\thanks{This work has received support from UNAM-DGAPA-PAPIIT grant IN117525 (Mexico).}}
	\author{V\'ictor Hern\'andez-Santamar\'ia\thanks{V. Hern\'andez-Santamar\'ia is supported by the Project CBF2023-2024-116 of SECIHTHI and by UNAM-DGAPA-PAPIIT grants IA100324 and IN102925 (Mexico).} \and Subrata Majumdar\thanks{Subrata Majumdar is supported by the Institute Postdoctoral Fellowship of TIFR CAM, Bangalore, India} \and Yingying Wu-Zhang\thanks{Y. Wu-Zhang is partially supported by the Fundación Sofía Kovalévskaya-SMM grant (Mexico).}}
	
	\maketitle
\begin{abstract}
In this paper, we establish global existence of strong solutions for the three-dimensional Cahn--Hilliard--Navier--Stokes system on the torus through finite-dimensional scalar modifications of a prescribed external force. For arbitrary initial data in $H^1\times H^2$ and any force that is locally bounded in time, we construct two smooth scalar corrections, one added to the third component of the fluid force and the other to the phase source. These corrections involve three fixed Fourier frequencies in the fluid equation and one in the phase equation, and act only on a set of times of arbitrarily small total measure. The proof combines a quantitative local theory with repeated returns of the solution to a fixed ball around zero. These returns rely on small-time asymptotic limits that balance the different scales of the two equations, together with a saturation procedure carefully adapted to the structure of the coupled system. Fluid transport generates the phase frequencies that the chosen phase forcing cannot produce on its own. As a consequence, we also obtain global strong solutions for the three-dimensional Navier--Stokes equations with a single smooth scalar correction acting on one component of the force at three fixed Fourier frequencies.
\end{abstract}

\footnotesize
\begingroup
\hypersetup{linkcolor=black}
\tableofcontents
\endgroup
\normalsize

\section{Introduction}

\subsection{The system and the question of global existence}

The Cahn--Hilliard--Navier--Stokes (CHNS) system describes the motion of a
binary mixture of incompressible fluids undergoing phase separation.
The Cahn--Hilliard equation describes the redistribution of the two
components, while the Navier--Stokes equation describes the fluid motion.
A typical phase-separation experiment consists in rapidly cooling an
initially homogeneous mixture below its critical temperature. When the
homogeneous state becomes unstable after this quench, concentration
fluctuations grow and lead to the formation of regions occupied
predominantly by one component. The description of this
mechanism originates in the work of Cahn and Hilliard \cite{CH58} and
the analysis of spinodal decomposition by Cahn \cite{Cah61}.

In a fluid mixture, the velocity transports and deforms the concentration
profile, while concentration gradients produce capillary forces in the
momentum equation. The diffuse-interface formulation represents the
interface by a concentration profile varying across a transition
layer. Its evolution is determined by the coupled equations. Examples
of this formulation and its hydrodynamic coupling are discussed in
\cite{LS03,AGG12}.

We consider a bulk mixture with periodic boundary conditions, so that
the analysis concerns the interaction between phase separation and fluid
motion without boundary effects. The spatial domain is the
three-dimensional torus $\mathbb T^3=\mathbb R^3/(2\pi\mathbb Z^3)$.

Phase separation is driven by a double-well potential.
In this work, we use the regular potential
\begin{equation}\label{reg potential}
 F(r)=\frac14(r^2-1)^2.
\end{equation}
Its minima at $\phi=\pm1$ correspond to the pure phases. Regular
double-well potentials are widely used in the analysis and numerical
approximation of phase separation. They retain the two preferred
concentrations, although they do
not enforce the physical range $|\phi|\leq1$. For a more complete
discussion of regular and singular potentials, see, for instance,
\cite{Mir17}. %A coupled system with singular potentials is studied in \cite{Abels2009}.

With constant mobility normalized to one, the concentration
evolves by the Cahn--Hilliard equation
$\partial_t\phi=\Delta\mu+f_\phi$, where
$\mu=-\Delta\phi+F'(\phi)$ is the chemical potential.
The mixture is treated as an incompressible Newtonian fluid with equal
constant densities (normalized to one) and viscosity $\nu>0$.
Including transport of the concentration and the capillary force in the
momentum equation yields the coupled system
\begin{equation}\label{chn}
\begin{cases}
 \partial_tu-\nu\Delta u+(u\cdot\nabla)u+\nabla p
       =\mu\nabla\phi+f_u,
       & t>0,\ x\in\mathbb T^3,\\
 \partial_t\phi+u\cdot\nabla\phi=\Delta\mu+f_\phi,
       & t>0,\ x\in\mathbb T^3,\\
 \mu=-\Delta\phi+\phi^3-\phi, & t>0,\ x\in\mathbb T^3, \\
  \operatorname{div}u=0, 
       & t>0,\ x\in\mathbb T^3,\\
 u(0)=u_0,\qquad \phi(0)=\phi_0,
       & x\in\mathbb T^3.
\end{cases}
\end{equation}
Here $u$, $p$, $\phi$, and $\mu$ denote the velocity, pressure, relative
concentration, and chemical potential, respectively, while $f_u$ is the external fluid force and $f_\phi$
is the phase source. A direct computation using periodicity and
incompressibility gives
\begin{equation*}%\label{eq:mass_intro}
 \frac{d}{dt}\int_{\mathbb T^3}\phi\,dx
       =\int_{\mathbb T^3}f_\phi\,dx.
\end{equation*}
We use zero-average initial data and forces, so the phase mass remains
zero. See \Cref{ff} for further details.

The term $u\cdot\nabla\phi$ transports the concentration along fluid
trajectories and can change its gradient energy. The force
$\mu\nabla\phi$ describes the action of the phase profile on the fluid.
Its stress form follows from
\begin{equation}\label{eq:capillary_stress_intro}
\begin{aligned}
 \mu\nabla\phi
 &=-(\Delta\phi)\nabla\phi+\nabla F(\phi)=-\operatorname{div}(\nabla\phi\otimes\nabla\phi)
   +\nabla\left(\frac12|\nabla\phi|^2+F(\phi)\right),
\end{aligned}
\end{equation}
where $(\nabla\phi\otimes\nabla\phi)_{ij}
=\partial_i\phi\,\partial_j\phi$. Absorbing the gradient into the
pressure leaves the capillary stress associated with concentration
gradients. This stress couples the phase free energy to the momentum
balance (see, e.g., \cite{AMW98,LS03} for more details).

The existence of global weak solutions in two and three dimensions is
established in the classical works \cite{Boyer1999,Abels2009}.
In two dimensions, global well-posedness and long-time behavior for
regular potentials are studied in \cite{GG10a}, while global strong
well-posedness for the logarithmic potential is proved in \cite{GMT19}. A semigroup approach to local existence and uniqueness is developed
in \cite{Dlo22}.

For the unforced three-dimensional system, global strong solutions
are known under additional restrictions. For instance, \cite{ZWH09}
proves global strong existence for sufficiently large viscosity or
for initial states sufficiently close to an absolute minimizer of
the energy. Regularity criteria for
three-dimensional solutions are established, for example, in
\cite{GPGP16,FL19}. For a model with unmatched densities,
\cite{AGG24} proves eventual regularity of global weak solutions:
the phase regularizes at positive times and the velocity at
sufficiently large times.

Further results on long-time dynamics can be found, for instance, in \cite{GG10b,BG14,GT22}. Models with variable density are studied in \cite{GT20}, and models with mass transfer, chemotaxis, or active transport in \cite{LW18,HW24,HW26}.

For sufficiently regular initial data and locally bounded forces, the
system considered here admits a unique strong solution on a short
time interval. Global strong existence for general three-dimensional
data remains open. The following formal calculations
illustrate why the natural a priori estimates do not suffice to
establish global strong existence.

For a concentration profile $\phi$, the Ginzburg--Landau free energy is given by
\begin{equation*}%\label{eq:phase_energy_intro}
 \mathscr F(\phi)
 =\int_{\mathbb T^3}
     \left(\frac12|\nabla\phi|^2+F(\phi)\right)\,dx.
\end{equation*}
The gradient term penalizes spatial variations of the concentration,
and the bulk potential favors the pure phases. This form of the energy
is based on the theory of Cahn and Hilliard \cite{CH58}.
{The chemical potential $\mu$ in \eqref{chn} is its first
variation.}
The total energy of the mixture is the sum of the fluid kinetic energy
and the phase free energy
\begin{equation*}%\label{eq:energy_intro}
 \mathscr E(u,\phi)=
 \frac12\|u\|_{L^2}^2+\mathscr F(\phi).
\end{equation*}
Formally, we test the momentum equation by $u$ and the phase equation
by $\mu$, and integrate by parts. The coupling terms
{$\pm\int_{\mathbb T^3}\mu\,u\cdot\nabla\phi\,dx$}
exchange energy between the two components and
cancel upon addition, yielding
\begin{equation*}%\label{eq:energy_balance_intro}
 \frac{d}{dt}\mathscr E(u,\phi)
  +\nu\|\nabla u\|_{L^2}^2
  +\|\nabla\mu\|_{L^2}^2
 =(f_u,u)_{L^2}+(f_\phi,\mu)_{L^2}.
\end{equation*}
When $f_u=f_\phi=0$, the total energy is nonincreasing. Integration
in time then bounds both the energy and the accumulated dissipation
by the initial energy. Moreover, \Cref{lem:quartic-potential}
shows that the potential energy controls
$\|\phi\|_{L^4}^4$, while the gradient contribution and
Poincar\'e's inequality control $\|\phi\|_{H^1}^2$ in the
zero-average setting. Thus, for a smooth unforced solution on $[0,T]$,
\begin{align}
 &\sup_{0\leq t\leq T}
  \left(\|u(t)\|_{L^2}^2+\|\phi(t)\|_{H^1}^2
                            +\|\phi(t)\|_{L^4}^4\right)\notag\\
 &\qquad+\int_0^T\left(\nu\|\nabla u(t)\|_{L^2}^2
                              +\|\nabla\mu(t)\|_{L^2}^2\right)\,dt
 \leq C\left(1+\mathscr E(u_0,\phi_0)\right),
 \label{eq:energy_bounds_intro}
\end{align}
where $C$ is independent of $T$.
Estimates of this type are fundamental to the construction of
global weak solutions; see, for instance, \cite{Boyer1999,Abels2009} and \Cref{sec:wp_local} below.

For strong solutions, the question is whether the velocity and phase
profile remain bounded in $H^1$ and $H^2$, respectively, on every
finite time interval. Such bounds allow the local existence argument
to be repeated and the solution to be extended. The estimate
\eqref{eq:energy_bounds_intro} alone does not provide this control:
it bounds the time integral of $\|\nabla u\|_{L^2}^2$, but not
its supremum in time. To estimate these higher-order norms in the
unforced case, we test the momentum equation by $-\Delta u$ and
the phase equation by $\Delta^2\phi$. Using \eqref{eq:energy_bounds_intro} and
usual interpolation inequalities, we obtain
\begin{align*}
 &\frac{d}{dt}\left(\|\nabla u\|_{L^2}^2
                         +\|\Delta\phi\|_{L^2}^2\right)
 +c\left(\|\Delta u\|_{L^2}^2
                         +\|\Delta^2\phi\|_{L^2}^2\right)\leq
 C\left(1+\|\nabla u\|_{L^2}^2
                         +\|\Delta\phi\|_{L^2}^2\right)^3,
\end{align*}
where $c>0$ depends on $\nu$ and $C>0$ on $\nu$ and
$\mathscr E(u_0,\phi_0)$ (see \Cref{thm:existence-energy-m1} for precise further details).
Writing $\mathcal Y(t)=\|\nabla u(t)\|_{L^2}^2
+\|\Delta\phi(t)\|_{L^2}^2$, this gives
$\mathcal Y'\leq C(1+\mathcal Y)^3$. A comparison argument
bounds $\mathcal Y$ on a short interval depending on the initial data and $\nu$, but yields no bound for arbitrary times.

This difficulty is already present in the Navier--Stokes equations alone.
Indeed, when $\phi_0=0$ and $f_u=f_\phi=0$, the phase remains zero
and \eqref{chn} reduces to the unforced Navier--Stokes system,
for which the velocity estimate takes the classical form
\begin{equation*}%\label{eq:ns_gradient_intro}
 \frac{d}{dt}\|\nabla u\|_{L^2}^2
 +\nu\|\Delta u\|_{L^2}^2
 \leq C_\nu\|\nabla u\|_{L^2}^6.
\end{equation*}
Thus, even without the phase variable, the energy and gradient
estimates do not suffice to establish global strong existence for
arbitrary data; see, for instance, \cite{Tem95}. We therefore ask

\begin{center}
\begin{minipage}{0.9\textwidth}
\centering\itshape
Can a modification of only a few Fourier modes of a prescribed force ensure global strong existence?
\end{minipage}
\end{center}

\subsection{Finite-dimensional forcing and the main result}

Our answer is affirmative. To state the result precisely,
we denote by $\mathbf L^2_{\mathrm{div},0}$ and $\mathbf H^m$ the spaces
of divergence-free vector fields with zero spatial mean in $L^2$ and
$H^m$, respectively, and by $\dot H^m$ the space of functions in
$H^m(\mathbb T^3)$ with zero spatial mean. We set
\begin{equation*}
 \mathbb V^0:=\mathbf L^2_{\mathrm{div},0}\times\dot H^1,
 \qquad
 \mathbb V^1:=\mathbf H^1\times\dot H^2.
\end{equation*}
Precise definitions and norms are given in \Cref{ff}.

Let $z_0=(u_0,\phi_0)\in\mathbb V^1$ and
$f\in L^\infty_{\mathrm{loc}}(\mathbb R_+;\mathbb V^0)$.
{For $T\in(0,\infty]$, we consider strong solutions
$z=(u,\phi)$ in the class
\[
 z\in C([0,T);\mathbb V^1)
 \cap L^2_{\mathrm{loc}}([0,T);\mathbf H^2\times\dot H^4).
\]
The precise notion of strong solution will be given in
\Cref{sec:wp_local}.}
% TODO: Introduce the notions of solution in Section 3.

%Let $z_0=(u_0,\phi_0)\in\mathbb V^1$ and
%$f\in L^\infty_{\mathrm{loc}}(\mathbb R_+;\mathbb V^0)$. For
%$T\in(0,\infty]$, we call $z=(u,\phi)$ a \emph{strong solution} of
%\eqref{chn} on $[0,T)$ if
%\begin{equation*}%\label{eq:strong_class_intro}
% z\in C([0,T);\mathbb V^1)
%   \cap L^2_{\mathrm{loc}}([0,T);\mathbf H^2\times\dot H^4)
%\end{equation*}
%and \eqref{chn} holds almost everywhere in $(0,T)\times\mathbb T^3$ for some pressure $p$.  
%Indeed, the above regularity and the equations imply
%\[
%\partial_t z\in L^2_{\mathrm{loc}}([0,T);\mathbb V^0),
%\quad
%\mu\in L^2_{\mathrm{loc}}([0,T);H^2(\mathbb T^3)),
%\quad
%p\in L^2_{\mathrm{loc}}([0,T);\dot H^1),
%\]
%where the pressure is normalized to have zero spatial mean, that is $\int_{\mathbb T^3}p(t,x)\,dx=0$.
%

As recalled above, such a solution exists and is
unique on a short time interval, and it can be continued as long as its
$\mathbb V^1$ norm remains bounded. This gives a unique
\emph{maximal strong solution} on $[0,T_{\max}(z_0,f))$, with
$T_{\max}(z_0,f)\in(0,\infty]$. If $T_{\max}(z_0,f)<\infty$, then
\begin{equation*}%\label{eq:blowup_intro}
 \lim_{t\uparrow T_{\max}(z_0,f)}\|z(t)\|_{\mathbb V^1}=\infty,
\end{equation*}
see \Cref{prop:maximal-strong} for a precise statement.
{Since the system contains the three-dimensional Navier--Stokes
equations, it is not known whether $T_{\max}(z_0,f)=\infty$ for
general initial data, even when $f=0$.}\footnote{The recent manuscript by OpenAI
\cite{OpenAI26} announces finite-time blow-up for the
three-dimensional Navier--Stokes equations with a suitably constructed
smooth force. If confirmed, it would show
that global strong existence does not hold for all smooth given forces. Note however, that such construction does not address the global regularity problem for arbitrary smooth initial data without forcing, which remains open.}
Instead, we keep $z_0$ fixed and look for a modification $g$ of $f$ such that $T_{\max}(z_0,g)=\infty$.

The question is therefore how much of $f$ must be changed.
We keep $f$ unchanged outside a fixed finite-dimensional space of
spatial profiles, and modify it only inside this space.
For a finite-dimensional space $\mathcal H\subset\mathbb V^0$ consisting of
smooth spatial profiles, the forces that differ from $f$ only in
$\mathcal H$ form the equivalence class
\begin{align*}
 f\sim_{\mathcal H}g
 &\quad\Longleftrightarrow\quad
 g-f\in L^\infty_{\mathrm{loc}}(\mathbb R_+;\mathcal H),\notag\\
 [f]_{\mathcal H}
 &=f+L^\infty_{\mathrm{loc}}(\mathbb R_+;\mathcal H).
% \label{eq:force_class_intro}
\end{align*}
Both $f$ and $g$ belong to
$L^\infty_{\mathrm{loc}}(\mathbb R_+;\mathbb V^0)$: $f$ is given and $g$ is to be constructed. If $Q_{\mathcal H}$ is the
orthogonal projection in
$\mathbb V^0$ onto $\mathcal H^\perp$, then
\begin{equation*}%\label{eq:force_projection_intro}
 g\in[f]_{\mathcal H}
 \quad\Longleftrightarrow\quad
 Q_{\mathcal H}g(t)=Q_{\mathcal H}f(t)
 \quad\text{for almost every }t\geq0.
\end{equation*}
We now specify the Fourier modes defining
the space of modifications used in the present setting. For a finite set $\mathcal K\subset\mathbb Z^3\setminus\{0\}$, let
\begin{equation*}%\label{eq:scalar_modes_intro}
 \mathcal E_{\mathcal K}
 =\operatorname{span}\{\cos(k\cdot x),\sin(k\cdot x):k\in\mathcal K\}.
\end{equation*}
Writing $e_1,e_2,e_3$ for the standard basis of $\mathbb R^3$, we choose
\begin{equation}\label{eq:mode_sets_intro}
 \mathcal K_u={\{e_2,e_1+e_3,e_2+e_3\}},\qquad
 \mathcal K_\phi=\{e_3\},\qquad
 \mathcal E=\mathcal E_{\mathcal K_u}\times\mathcal E_{\mathcal K_{\phi}}.
\end{equation}
We add a scalar function $\eta_u$ to the third component of the fluid
force and a scalar function $\eta_\phi$ to the phase source, with $(\eta_u,\eta_\phi)$ taking values in $\mathcal E$. The modified CHNS equations are 
\begin{equation}\label{eq:chn_modified_intro}
\begin{cases}
 \partial_tu-\nu\Delta u+(u\cdot\nabla)u+\nabla p
       =\mu\nabla\phi+f_u+\eta_ue_3,
       & t>0,\ x\in\mathbb T^3,\\
 \partial_t\phi+u\cdot\nabla\phi=\Delta\mu+f_\phi+\eta_\phi,
       & t>0,\ x\in\mathbb T^3.
\end{cases}
\end{equation}
Note that here $\eta_ue_3=(0,0,\eta_u)$, that is, the first two components of the
prescribed fluid force are left unchanged. Thus the whole
correction consists of two scalar functions, supported on the four frequencies in $\mathcal K_u\cup\mathcal K_\phi$.

To relate \eqref{eq:chn_modified_intro} to the class
$[f]_{\mathcal H}$, we apply the Leray projection $\Pi$ onto
divergence-free, zero-average vector fields. The pressure
term vanishes, and \eqref{eq:capillary_stress_intro} gives
$\Pi(\mu\nabla\phi)=-\Pi((\Delta\phi)\nabla\phi)$.
Thus the modified system takes the form
\begin{equation}\label{eq:chn_projected_intro}
\begin{cases}
 \partial_tu-\nu\Delta u+\Pi((u\cdot\nabla)u)
       =-\Pi((\Delta\phi)\nabla\phi)+{g_u},\\
 \partial_t\phi+u\cdot\nabla\phi
       =-\Delta^2\phi+\Delta(\phi^3-\phi)+g_\phi,\\
u(0)=u_0,\quad\phi(0)=\phi_0.
\end{cases}
\end{equation}
where the total force $g=(g_u,g_\phi)$ is given by
\begin{equation}\label{eq:modified_force_intro}
 g_u=f_u+\Pi(\eta_ue_3),\qquad g_\phi=f_\phi+\eta_\phi.
\end{equation}
Here $\Pi f_u=f_u$ because the prescribed momentum force $f_u$ is divergence free {and has zero spatial mean}. The operator formulation and the associated functional setting are developed in \Cref{ff}.

Thus the space of modifications appearing in the force class is
\begin{equation}\label{eq:projected_modes_intro}
 \mathcal H
 :=\{(\Pi(\eta_ue_3),\eta_\phi):(\eta_u,\eta_\phi)\in\mathcal E\}.
\end{equation}
For $L^\infty_{\mathrm{loc}}(\mathbb R_+;\mathcal H)\ni \eta:=(\Pi(\eta_ue_3),\eta_\phi)$,
\eqref{eq:modified_force_intro} gives $g\in[f]_{\mathcal H}$.
Conversely, every $g\in[f]_{\mathcal H}$ has this form.

With $\mathcal H$ defined by \eqref{eq:projected_modes_intro}, our main result gives a global strong solution for every initial state $z_0\in\mathbb V^1$ and every
prescribed force $f\in L^\infty_{\mathrm{loc}}(\mathbb R_+;\mathbb V^0)$.
Given any $\rho>0$, the new force can be chosen to agree with the prescribed force $f$ outside a union $I$ of intervals with total Lebesgue measure less than $\rho$. Each bounded time interval meets only finitely many of these intervals, so the modifications do not accumulate at a finite time.

The result reads as follows.
\begin{theo}\label{mnthm}
Let $\mathcal H$ be defined by \eqref{eq:projected_modes_intro},
$z_0=(u_0,\phi_0)\in\mathbb V^1$, and
$f\in L^\infty_{\mathrm{loc}}(\mathbb R_+;\mathbb V^0)$. 
There exists $g\in[f]_{\mathcal H}$ such that $T_{\max}(z_0,g)=\infty$,
that is, the maximal strong solution of \eqref{eq:chn_projected_intro}
with initial state $z_0$ and force $g$ satisfies
\begin{equation}\label{eq:global_solution_intro}
 z=(u,\phi)\in C([0,\infty);\mathbb V^1)
       \cap L^2_{\mathrm{loc}}([0,\infty);{\mathbf H^2\times\dot H^4}).
\end{equation}
Moreover, given any $\rho>0$, the force $g$ can be chosen to satisfy
\begin{equation}\label{eq:force_agreement_intro}
 g(t)=f(t)\qquad\text{for almost every }t\in\mathbb R_+\setminus I,
\end{equation}
where $I\subset\mathbb R_+$ is a measurable set with Lebesgue
measure $|I|<\rho$.
\end{theo}
\begin{rmk}\label{rem:main-force-choice}
Note that we do not assume that the solution corresponding to $f$
is global. It is the force $g$ that we choose so that
$T_{\max}(z_0,g)=\infty$, without requiring additional regularity
of $z_0$ or $f$. We change only the $\mathcal H$-component of the
force, so $Q_{\mathcal H}g=Q_{\mathcal H}f$, and leave the
differential operators and nonlinear terms of \eqref{chn} unchanged.

Although the prescribed force need not be smooth, we show in the
proof that its scalar corrections can be chosen with
$\eta_u,\eta_\phi\in C^\infty([0,\infty)\times\mathbb T^3)$.
They involve only the four fixed frequencies above, and the set of
times on which they act can still have arbitrarily small total
measure. Finally, the theorem does not assert uniqueness of $g$.
Once $g$ is fixed, uniqueness of the strong solution follows from
the local theory.
\end{rmk}
\begin{rmk}\label{rem:ns-consequence-intro}
Taking $\phi_0=f_\phi=0$, we may keep $\eta_\phi=0$ and $\phi\equiv0$
throughout the construction, obtaining global strong solutions for
the three-dimensional Navier--Stokes equations with only the scalar
correction $\eta_u$. We discuss this consequence and possible
extensions in \Cref{sec:further-results}.
\end{rmk}
\subsection{Relation with earlier work}
A classical approach to obtaining global regular
solutions is to modify the fluid equations. An early example is
Leray's work
\cite{Leray34}, which uses a smoothed transport velocity to obtain
globally regular approximations of the Navier--Stokes
equations. Passing to the limit gives global weak solutions for
finite-energy data, whose existence was also established by Hopf
\cite{H51}\footnote{Nevertheless, global weak existence leaves the fundamental question of
uniqueness. For Navier--Stokes, \cite{BV19} proves
nonuniqueness of finite-energy weak solutions on the torus. In the
forced problem on $\mathbb R^3$, \cite{ABC22} constructs distinct
Leray--Hopf solutions with zero initial velocity and a suitable common
external force.}

These regularization ideas have also been developed for CHNS, as in
the family of two-phase flow models studied in \cite{GM15}.
Another approach for the Navier--Stokes equations multiplies
$u\cdot\nabla u$ by a cutoff depending on the norm of the velocity
gradient, limiting the effect of convection when that norm is large
\cite{CKR06}. Such cutoff modifications also give globally
well-posed three-dimensional CHNS
models \cite{TM17}. A Kelvin--Voigt regularization of the momentum
equation gives global strong solutions in two and three dimensions, as shown in \cite{GNTT23}.

In a different approach, an additional force can be chosen as a function of the current state. This allows the force to provide dissipation or to guide the evolution towards a prescribed state or trajectory. Damped and tamed Navier--Stokes equations use this dependence to
obtain global regular solutions under suitable growth conditions
\cite{CJ08,RZ09}. Nonlinear damping has
also been studied for CHNS, for instance \cite{BDM26} proves global weak existence and uniqueness under suitable assumptions on the damping.

The force can also be designed to influence the asymptotic behavior.
For the Camassa--Holm equation, a spatially localized state-dependent
term gives global solutions converging to an equilibrium
\cite{Glass08}. For the three-dimensional Navier--Stokes equations,
\cite{BP21} constructs a global strong solution of an assimilation
system that approaches a prescribed Leray--Hopf weak solution, using
corrections based on coarse observations. The result requires conditions on those
observations, without assuming global strong regularity of the
reference trajectory.\footnote{The same idea appears in other
settings. In stochastic fluid models, state dependence enters through
multiplicative noise and can prevent finite-time blow-up; see, for
instance, \cite{GHV14,BMX25,LC26}. For the two-dimensional
Kuramoto--Sivashinsky equation, the modification in
\cite[Remark~5.3]{LM24} can be written as a force on low Fourier
modes, with a state-dependent frequency cutoff, canceling selected
nonlinear interactions. The use of low-mode forcing is close in
spirit to the present work, although the regularizing mechanism is
different.}

The closest precedents for modifications within a fixed
finite-dimensional forcing space are \cite{S07} and \cite{GHHM18}.
The interpretation of the force in \Cref{mnthm} is inspired by
Shirikyan's work \cite{S07}, where the component of the force outside
a finite-dimensional space is prescribed and the remaining component
is chosen as a control. We express this through the equivalence
classes $[f]_{\mathcal H}$.
For the three-dimensional Navier--Stokes equations,
that work establishes exact controllability in finite-dimensional
projections. In particular, within each prescribed class, one can
choose a force for which a strong solution exists throughout any
given finite time interval (see \cite[Proposition 2.9]{S07}). Note, however, that the sufficient condition in
\cite[Corollary~2.3]{S07} uses a control space containing the first
$N$ eigenfunctions of the Stokes operator, with $N$ sufficiently
large. In \eqref{eq:chn_modified_intro}, only the third component of
the fluid force is corrected and at three fixed frequencies.

For the three-dimensional Euler equations, \cite{GHHM18} starts
from the locally defined evolution and establishes accessibility
by forcing finitely many low modes.
As observed in \cite[Remark 5.33]{GHHM18}, these accessibility
properties can be used to prevent blow-up even in the presence
of a time-independent background force. The remark
states this possibility without giving a full proof.
We take this observation as a starting point and give a complete
construction for CHNS, allowing a prescribed time-dependent force
and corrections acting on a set of times of arbitrarily small
total measure.

The restricted forcing space makes the nonlinear coupling
essential since the directly forced phase profiles
depend only on $x_3$, and the phase nonlinearity preserves this
dependence. We therefore use transport by the fluid to
produce the missing horizontal phase frequencies. The main new step for
CHNS is a mixed small-time limit that balances the different scales of
the two equations to retain this transport together with the cubic
phase term, while keeping the fluid evolution under control. Iterating
the resulting transport displacements generates all phase frequencies
from $e_3$, using only the three prescribed frequencies acting on one
component of the fluid force. These same frequencies also suffice for
fluid saturation.
\subsection{Main ideas and difficulties}\label{sec:main-ideas}
Our construction alternates small-time modifications of the force with
periods in which we let the solution evolve under the prescribed force.
The idea is to bring the solution
close to zero, so that we can apply the local existence theory and
obtain a further interval of existence. The main challenge is to do this using only the two scalar corrections
introduced above. Below we explain the main ideas and difficulties in
turning these repeated returns into a global strong solution.

\begin{itemize}[leftmargin=0pt]
\item{\textbf{Different scales in the fluid and phase equations.}}
Following the scaling and saturation arguments in \cite{GHHM18,Ner21}, we use large
profiles over small time intervals to generate directions that are not
directly forced. For CHNS, the difficulty comes from the coupling:
the velocity transports the phase, while the phase acts back on the
fluid through capillarity. The amplitudes must therefore be chosen together, taking into account their effect on both equations.

We first illustrate the idea with a scalar equation
\begin{equation*}
 \partial_t y+\mathcal L y+y^p=f+\eta,\qquad y(0)=y_0,
\end{equation*}
where $f$ is the prescribed force, $\eta$ is an additive correction,
$p\geq2$ is an integer, and $\mathcal L$ is a linear
differential operator representing, for instance, diffusion.

First, a large force acting over a short interval produces
a finite displacement. Indeed, for a fixed smooth profile $h$ in the forcing space, take
$\eta=\delta^{-1}h$ on $[0,\delta]$ and set
$w_\delta(s)=y(\delta s)$, $0\leq s\leq1$. Then
\begin{equation*}%\label{eq:eq_w_delta}
 \partial_s w_\delta
 =h-\delta\bigl(\mathcal Lw_\delta+w_\delta^p\bigr)
   +\delta f(\delta s),
 \qquad w_\delta(0)=y_0.
\end{equation*}
The formal limiting equation is $\partial_s w=h$, whose
solution is $w(s)=y_0+sh$. We therefore expect
$y(\delta)\to y_0+h$ as $\delta\to0$. On this time scale, the prescribed force and the evolution terms have a vanishing contribution.

The nonlinear evolution can also produce a finite displacement.
Suppose that translations in a smooth direction $q$ can
already be approximated in arbitrarily short time. We first
bring the state close to $y_0+\delta^{-1/p}q$, and then set
$\eta=0$ for time $\delta$. To examine the effect on this interval, let
$y_\delta$ denote the solution starting exactly from
$y_0+\delta^{-1/p}q$ and write
$r_\delta(s)=y_\delta(\delta s)-\delta^{-1/p}q$. A direct computation shows that it satisfies
\begin{equation*}
 \partial_s r_\delta
 +\delta\mathcal Lr_\delta
 +\delta^{1-1/p}\mathcal Lq
 +(q+\delta^{1/p}r_\delta)^p
 =\delta f(\delta s),
 \qquad r_\delta(0)=y_0,
\end{equation*}
and the formal limiting equation is
$\partial_s r=-q^p$, suggesting that
\begin{equation*}
 y_\delta(\delta)-\delta^{-1/p}q
 \longrightarrow y_0-q^p \quad\text{as } \delta \to 0.
\end{equation*}
A final short translation in the opposite direction removes the large profile $\delta^{-1/p}q$, leaving a displacement close to $-q^p$, even when this direction is not directly available in the forcing space. All three steps can be carried out in arbitrarily short time. Earlier uses of nonlinear interactions to generate new directions can be found, for instance, in \cite{AS05,AS06}.

For CHNS, we implement this construction while balancing both equations of the system. The cubic phase nonlinearity suggests the
following amplitude for a smooth phase profile $\psi$
\begin{equation*}
\delta\Delta\bigl((\delta^{-1/3}\psi)^3\bigr)
 =\Delta(\psi^3).
\end{equation*}
On the other hand, for a smooth divergence-free velocity profile
$v$, the nonlinear term in the fluid equation suggests the amplitude
$\delta^{-1/2}$ since
\begin{equation*}
 \delta\Pi\bigl(((\delta^{-1/2}v)\cdot\nabla)(\delta^{-1/2}v)\bigr)
 =\Pi((v\cdot\nabla)v).
\end{equation*}
With these two amplitudes, the transport term in the phase equation contributes only at order $\delta^{1/6}$ and disappears
in the limit. We therefore need a larger velocity profile $\delta^{-2/3}v$, so that transport and the cubic phase term both have a finite effect over time $\delta$. More precisely,
\begin{align}
 \delta\Delta((\delta^{-1/3}\psi)^3)
 &=\Delta(\psi^3),\notag\\
 \delta(\delta^{-2/3}v)\cdot\nabla(\delta^{-1/3}\psi)
 &=v\cdot\nabla\psi.
 \label{eq:mixed_balance_intro}
\end{align}
After removing the large profiles, the expected phase
displacement is $\Delta(\psi^3)-v\cdot\nabla\psi$.

However, this choice also amplifies fluid self-advection, namely its contribution over time $\delta$ is
$\delta^{-1/3}\Pi((v\cdot\nabla)v)$, which in general diverges as
$\delta\to0$.
To eliminate this term, we choose the intermediate velocity
profile $v$ to be a single divergence-free Fourier mode of the form
\begin{equation*}%\label{eq:safe_mode_intro}
 v=a\cos(k\cdot x),\qquad a\cdot k=0.
\end{equation*}
Then $(v\cdot\nabla)v=0$, so this contribution vanishes.
The sine mode has the same property.
By contrast, over time $\delta$ the projected capillary force contributes only
\begin{equation*}
 -\delta\Pi\bigl(\Delta(\delta^{-1/3}\psi)\,
                    \nabla(\delta^{-1/3}\psi)\bigr)
 =-\delta^{1/3}\Pi(\Delta\psi\,\nabla\psi)
 \longrightarrow0 \ \text{ as $\delta\to 0$}.
\end{equation*}
Making this intuition rigorous also requires showing that strong
solutions exist over the required intervals, even with these large
intermediate profiles. These choices allow us to use the transport term
in the phase equation without a divergent contribution in the fluid
equation. By transporting suitable
phase profiles, we can then generate phase directions that are not
directly forced, as we explain in the next point.

Related scaling arguments have also been developed
for a few coupled fluid models, including
the primitive equations in \cite{BGN23} and Boussinesq
flows in \cite[Section 5.3]{GHHM18} and
\cite{NR25b,rissel2025controllability}.
In these models, temperature enters the fluid equation
linearly, which makes that part of the coupling simpler. In CHNS,
both transport and capillarity are nonlinear and, as we have seen, obtaining the phase
displacement in \eqref{eq:mixed_balance_intro} also requires controlling
the effect of the large profiles in the fluid equation. Our goal here is to
combine these displacements, using only corrections in $\mathcal H$,
to obtain $\|z(\tau)\|_{\mathbb V^1}<1$ from any
$z_0\in\mathbb V^1$, with $\tau>0$ arbitrarily small. We call this
property small-time accessibility.

\item{\textbf{Saturation in the fluid and phase equations.}}
The small-time limits above give us directions of motion
that can in turn generate others through the nonlinear terms. The idea
of saturation is to repeat this process and obtain directions whose
combinations approximate any desired displacement, without enlarging
the original forcing space.

For the fluid equation, \cite{N15} shows that saturation in three dimensions can be obtained by forcing a single component at four Fourier frequencies. Here we reduce this number to three while retaining the same restriction. This gives a six-dimensional saturating space that we use to prove the Navier--Stokes result in \Cref{prop:ns-scalar-continuation}. See \Cref{rem:three-frequency-saturation} for further discussion.

For the phase, the initial directions
$\operatorname{span}\{\cos x_3,\sin x_3\}$ do not saturate through
the cubic Cahn--Hilliard nonlinearity itself. Products and spatial
derivatives preserve dependence on $x_3$, so they cannot generate
a dense family of functions on $\mathbb T^3$. In
\cite{hernandez2025global}, for the Cahn--Hilliard equation alone, the phase forcing space contains
frequencies in all spatial directions and saturation follows
from the cubic phase term. The new ingredient here is the
generation of the missing phase directions through transport.

The mixed balance \eqref{eq:mixed_balance_intro} gives the phase displacement
$\Delta(\psi^3)-v\cdot\nabla\psi$.
A second application with zero velocity profile
and phase profile $-\psi$ cancels the cubic term, leaving
$-v\cdot\nabla\psi$. 
For instance, ${v=e_3\cos x_2}$ has no
self-advection and, with $\psi=\sin x_3$, the transport term is
\begin{equation*}
 v\cdot\nabla\psi
 ={\cos x_2\cos x_3
 =\tfrac12\bigl(\cos(x_2+x_3)+\cos(x_2-x_3)\bigr)}.
\end{equation*}
Therefore, transport introduces horizontal dependence
into the phase. The fluid modes at frequencies $e_1+e_3$ and $e_2+e_3$
also permit the generation of purely horizontal phase frequencies.
{Iterating these transport operations with the original
fluid profiles generates all nonzero phase frequencies. For the
velocity, a separate iteration using quadratic interactions yields a
dense family of directions.}

Each generated change is approximated by a finite sequence of
additions taking values in $\mathcal H$.
Density and local stability of the solutions 
then give the small-time accessibility of the unit ball described above.

\item{\textbf{From accessibility to global existence.}}
For CHNS, we quantify the local existence time uniformly for
initial states in the unit ball $B_1$ of $\mathbb V^1$ and forces
bounded by a given constant. Our construction combines this estimate
with small-time accessibility
\begin{equation*}
 z_0
 \xrightarrow[\text{scaling and saturation}]{f+\eta}
 B_1
 \xrightarrow[\text{uniform existence time}]{f}
 \text{new state}
 \xrightarrow[\text{ \ \ return  \ \ }]{f+\eta}
 B_1
 \longrightarrow\cdots .
\end{equation*}
After each return, we remove the correction and let the solution
evolve under the prescribed force $f$ for a fixed fraction of the guaranteed
existence time. These
repeated returns are similar to those used in the proof of
\cite[Lemma~3.7]{GHHM18}.

Since $f$ is locally bounded, the time gained after each return
has a uniform positive lower bound on every bounded time interval.
The returns therefore cannot accumulate in finite time, and we
can continue the strong solution for all positive times.

The correction intervals $I_n$ are chosen with $\sum_n|I_n|<\rho$,
and $g=f$ outside their union. Only finitely many bounded corrections act on any bounded time interval, so $g-f$ is locally bounded.
{We give a more detailed explanation of this construction in
\Cref{sec:global-prescribed-force}.}
\end{itemize}
%%%%%%%%%%%%%%%
\subsection{\texorpdfstring{{Outline of the paper}}{Outline of the paper}}
The paper is organized as follows. Section~\ref{ff} introduces
the functional framework, and Section~\ref{sec:wp_local} establishes
quantitative existence, stability, and continuation results for CHNS. Section~\ref{sec:global-prescribed-force} states the small-time
accessibility result and uses it to prove \Cref{mnthm}.
The proof of accessibility is developed in three blocks:
small-time asymptotic limits in Section~\ref{sec:asymp_limits},
saturation in Section~\ref{sec:saturation}, and the accessibility
construction in Section~\ref{sec:accessibility}.

\section{Functional framework}\label{ff}

We first introduce the Lebesgue spaces on the three--dimensional torus $\mathbb T^3$.
For $p\in[1,\infty]$, we denote by $L^p(\mathbb T^3)$ the usual Lebesgue spaces,
either scalar-- or vector--valued.
Throughout the paper, we do not distinguish between scalar-- and vector--valued
function spaces; the precise meaning is always clear from the context.
{All functions considered in this paper are real valued.}
The associated norms are defined by
\[
\|u\|_{L^p}
:=
\left(\int_{\mathbb T^3} |u(x)|^p\,dx\right)^{1/p},
\quad 1\le p<\infty,
\qquad\text{and}\qquad
\|u\|_{L^\infty}
:=\operatorname*{ess\,sup}_{x\in\mathbb T^3}|u(x)|,
\]
where $|\cdot|$ denotes the Euclidean norm in $\mathbb R^k$, with $k=1$ or $3$.

We now introduce the usual Sobolev spaces on $\mathbb T^3$.
For {every integer $m\geq0$}, we denote by $H^m(\mathbb T^3)$ the standard Sobolev space,
either scalar- or vector-valued, endowed with the norm
\begin{equation*}%\label{Sobolev_norm}
\|u\|_{H^m}
:=
\left(
\sum_{|\alpha|\le m}
\|\partial^\alpha u\|_{L^2}^2
\right)^{1/2},
\end{equation*}
associated with the inner product
\[
(u,v)_{H^m}
:=
\sum_{|\alpha|\le m}
(\partial^\alpha u,\partial^\alpha v)_{L^2},
\]
where $(\cdot,\cdot)_{L^2}$ denotes the usual $L^2(\mathbb T^3)$ inner product.

For $\phi\in L^1(\mathbb T^3)$, we denote its spatial average by $\langle\phi\rangle
 :=\frac{1}{(2\pi)^3}\int_{\mathbb T^3}\phi(x)\,dx$. We write
\begin{equation*}
 L^2_{\mathrm{avg}}(\mathbb T^3)
 :=\{v\in L^2(\mathbb T^3):\langle v\rangle=0\}.
\end{equation*}
For every integer $m\geq0$, we set
\begin{equation*}
 \dot H^m(\mathbb T^3):=H^m_{\mathrm{avg}}(\mathbb T^3)
 :=\{v\in H^m(\mathbb T^3):\langle v\rangle=0\},
\end{equation*}
and endow this space with the homogeneous norm
\begin{equation*}
 \|v\|_{\dot H^m}
 :=\left(\sum_{|\alpha|=m}
                   \|\partial^\alpha v\|_{L^2}^2\right)^{1/2}.
\end{equation*}
In particular, $\dot H^0=L^2_{\mathrm{avg}}$ with the $L^2$ norm.

For $v\in L^2(\mathbb T^3)$, let $\widehat v(k):=\frac{1}{(2\pi)^3}
      \int_{\mathbb T^3}v(x)e^{-ik\cdot x}\,dx$ for $k\in\mathbb Z^3$. {Note that since $v$ is real valued, $\widehat v(-k)=\overline{\widehat v(k)}$.} The norms above have the useful characterizations
\begin{align}
 \|v\|_{H^m}^2
 &\sim\sum_{k\in\mathbb Z^3}
             (1+|k|^2)^m|\widehat v(k)|^2,
 \quad v\in H^m,\notag\\
 \|v\|_{\dot H^m}^2
 &\sim\sum_{k\in\mathbb Z^3\setminus\{0\}}
             |k|^{2m}|\widehat v(k)|^2,
 \quad v\in\dot H^m.
 \label{eq:fourier-sobolev-norms}
\end{align}
where $\sim$ denotes equivalence up to positive multiplicative
constants depending only on $m$. Indeed, since $\widehat v(0)=0$ on $\dot H^m$ and
$|k|^2\leq1+|k|^2\leq2|k|^2$ for $k\neq0$,
\eqref{eq:fourier-sobolev-norms} gives
$\|v\|_{H^m}\sim\|v\|_{\dot H^m}$ on this space. Note that the characterizations in \eqref{eq:fourier-sobolev-norms} extend to all real values of $m$. In this case, for $v\in\dot H^{m+2}(\mathbb T^3)$, it also gives
\begin{equation}\label{eq:periodic-norm-equivalences}
 \|v\|_{H^{m+2}}
 \sim\|\nabla v\|_{H^{m+1}}
 \sim\|\Delta v\|_{H^m} \sim \|\Delta^2 v\|_{H^{m-2}}.
\end{equation}
In view of this, hereinafter we use the $H^m$- and $\dot H^m$-norms interchangeably in
estimates for functions with zero spatial average.

We now introduce the functional spaces relevant for incompressible fluid flows. We define
\begin{align*}
 \mathbf L^2_{\mathrm{div},0}
 &:=\{u\in L^2(\mathbb T^3)^3:
                   \operatorname{div}u=0,\ \langle u\rangle=0\},\\
 \mathbf H^m
 &:=\mathbf L^2_{\mathrm{div},0}\cap H^m(\mathbb T^3)^3,
 \qquad m\geq0.
\end{align*}
In particular, $\mathbf H^0=\mathbf L^2_{\mathrm{div},0}$.
These spaces inherit the Sobolev norms and inner products
defined above.
 Using this notation, we define the product spaces
\begin{equation*}
\mathbb V^0:= \mathbf L^2_{\mathrm{div},0}\times \dot H^1
\end{equation*}
and for $m\in\mathbb N^*$, 
\[
\mathbb V^m
:= \mathbf H^{m} \times \dot H^{m+1}.
\]
The space $\mathbb V^m$ with $m\in\mathbb N^*$ is endowed with the natural product norm
\begin{equation*}%\label{not_Vm}
\|(u,\phi)\|_{\mathbb V^m}
:=
\bigl(
\|u\|_{H^{m}}^2 + \|\phi\|_{\dot H^{m+1}}^2
\bigr)^{1/2},
\end{equation*}
associated with the inner product
\[
\bigl((u,\phi),(v,\psi)\bigr)_{\mathbb V^m}
:=
(u,v)_{H^m} + (\phi,\psi)_{\dot H^{m+1}}.
\]
Similar definitions hold for $\mathbb V^0$.

We define the operator $A_0$ by
\[
A_0 u := - \Pi \Delta u, \qquad \forall u \in \mathbf H^2,
\]
where $\Pi : L^2(\mathbb T^3)^3 \to \mathbf L^2_{\mathrm{div},0}$ denotes the Leray projector.
Then $A_0$ is a positive self--adjoint unbounded operator on
$\mathbf L^2_{\mathrm{div},0}$, with compact inverse $A_0^{-1}$.
Moreover, the norm $\|A_0^{1/2}u\|_{L^2}$ is equivalent to the $H^1$-norm on
$\mathbf H^1$, and the norm $\|A_0 u\|_{L^2}$ is equivalent to the $H^2$-norm
on $\mathbf H^2$.

In order to define the variational setting for the CHNS equations, we introduce the bilinear operators $B_0$, $B_1$ and the coupling operator $R_0$,
together with their associated trilinear forms $b_0$ and $b_1$.

The operator
$B_0:\mathbf H^2\times\mathbf H^2\to\mathbf L^2_{\mathrm{div},0}$
is associated with the trilinear form
\begin{equation*}
 b_0(u,v,w):=\int_{\mathbb T^3}(u\cdot\nabla)v\cdot w\,dx,
\end{equation*}
and is defined by
\begin{equation*}
 (B_0(u,v),w)_{L^2}=b_0(u,v,w),
 \qquad w\in\mathbf H^2.
\end{equation*}

The operator
$B_1:\mathbf H^2\times\dot H^2\to L^2_{\mathrm{avg}}$
is associated with the trilinear form
\begin{equation}\label{def_b1}
 b_1(u,\phi,\psi):=\int_{\mathbb T^3}(u\cdot\nabla)\phi\,\psi\,dx,
\end{equation}
and is defined by
\begin{equation*}%\label{def_B1}
 (B_1(u,\phi),\psi)_{L^2}=b_1(u,\phi,\psi),
 \qquad\psi\in\dot H^2.
\end{equation*}
The forms $b_0$ and $b_1$ extend continuously to $H^1$ arguments.

Finally, the coupling operator
$R_0:L^2(\mathbb T^3)\times{\dot H^3}\to \mathbf L^2_{\mathrm{div},0}$
is defined by
\begin{equation}\label{def_R0}
 (R_0(\xi,\phi),w)_{L^2}
 :=\int_{\mathbb T^3}\xi\,\nabla\phi\cdot w\,dx,
 \qquad w\in\mathbf H^2.
\end{equation}
Equivalently,
\begin{align*}
&B_0(u,v)=\Pi\bigl((u\cdot \nabla )v\bigr), \\
&R_0(\xi,\phi)=\Pi\bigl(\xi\,\nabla\phi\bigr).
\end{align*}
These identities for $B_0$ and $R_0$ are understood
in the sense of distributions.
\begin{lemma}\label{lem-cancellation-B0}
Let $u\in \mathbf H^1$. Then, for all $v,w \in \mathbf H^1$,
\[
b_0(u,v,w) = -\,b_0(u,w,v), \qquad b_0(u,v,v) = 0.
\]
\end{lemma}
\begin{proof}
The identities follow from integration by parts and the divergence-free condition.
\end{proof}
For a force $f=(f_u,f_\phi)$ taking values in $\mathbb V^0$,
system \eqref{chn} has the projected form
\begin{equation}\label{chn-abstract}
\begin{cases}
\partial_t u
+ \nu A_0 u
+ B_0(u,u)
= {-R_0(\Delta\phi,\phi)}+f_u \\
\partial_t \phi + B_1(u,\phi)
= \Delta \mu+f_\phi \\
\mu=-\Delta \phi+\phi^3-\phi \\
u(0)=u_0,\qquad \phi(0)=\phi_0.
\end{cases}
\end{equation}
%here $\widetilde\eta_1=\Pi(\eta_1e_3).$ Then the control $\widetilde{\eta}_1$ takes values in the space
%\begin{equation*}
%\widetilde{\mathcal E}_{\mathcal K_1 } := \operatorname{span} \left\{ (P_k e_3)\cos(x\cdot k),\, (P_k e_3)\sin(x\cdot k) : k \in \mathcal K_1 \right\}, 
%\end{equation*}
%where for any $k \in \mathbb{Z}^3 \setminus \{0\}$,  $P_k$ is the orthogonal projection in $\mathbb{R}^3$ onto $k^\perp$.
The above formulation is obtained by applying the Leray--Helmholtz projection to the velocity equation in \eqref{chn}. {Since $f_u\in\mathbf L^2_{\mathrm{div},0}$, we have $\Pi f_u=f_u$.}  Also note that all other gradient 
terms are eliminated.  More precisely, the coupling term on the right-hand side of the velocity
equation in \eqref{chn} can be decomposed as
\begin{equation}\label{iden_chem_pot}
\mu\nabla \phi=(-\Delta\phi+\phi^3-\phi) \, \nabla\phi
= (-\Delta\phi)\nabla\phi + \nabla F(\phi),
\end{equation}
where $F$ is the regular potential defined in \eqref{reg potential}.
Since $\nabla F(\phi)$ is a gradient field, it can be absorbed into the
pressure term and therefore vanishes after applying the projection.
As a consequence, only the term $(-\Delta\phi)\nabla\phi$ contributes to
the coupling in the projected velocity equation, which is represented by $-R_0(\Delta\phi,\phi)$ in \eqref{chn-abstract}.

For notational convenience, we introduce the vector unknown
\begin{equation}\label{def_z}
z(t) := (u(t),\phi(t)),
\end{equation}
so that system \eqref{chn-abstract} can be written in the abstract form
\begin{equation}\label{abs}
\partial_t z + Lz + B(z) + Q(z) = f,
\end{equation}
where the operators $L$, $B$ and $Q$ are defined componentwise by
\begin{equation}\label{def_ops}
\begin{gathered}
Lz :=
\bigl(
\nu A_0 u,\;
\Delta^2\phi + \Delta\phi
\bigr),
\\
B(z) :=
\bigl(
B_0(u,u),\;
B_1(u,\phi)
\bigr),
\qquad
Q(z) :=
\bigl(
R_0(\Delta\phi,\phi),\;
-\Delta(\phi^3)
\bigr), \quad f:=(f_u, f_\phi).
\end{gathered}
\end{equation}
%
%For this auxiliary equation, the finite-time approximate accessibility
%statement takes the following form in projected variables:
%\begin{theo}
%Let $T>0$ and assume that $z_0:=(u_0,\phi_0)$, and $ z_1:=(u_1,\phi_1) \in \mathbb{V}^1$. For any $\epsilon>0$ there exists $\eta \in L^{\infty}(0,T;\widetilde{\mathcal E}_{\mathcal K_1}\times \mathcal E_{\mathcal K_2})$ such that the solution of \eqref{abs} satisfies 
%   \begin{align*}
%       \norm{z(T)-z_1}_{\mathbb V^1}<\epsilon.
%   \end{align*}
%\end{theo}
\section{Local well-posedness and quantitative results}\label{sec:wp_local}
In this section, we introduce the notions of weak and strong solution of
\eqref{abs} and state the corresponding existence results and related
energy estimates in \Cref{sec:weak-strong-solutions}.
The existence arguments follow the approach of \cite[Sections~3.2--3.3]{Boyer1999},
where a more general CHNS system is studied. These results provide the starting point of our analysis.

We then establish the following results, which will be used repeatedly
in the subsequent analysis:
\begin{enumerate}[label=\textup{(\roman*)}]
\item stability under small perturbations of the initial state and
the force, and uniqueness of strong solutions
(\Cref{thm:stability-strong,cor:unique-strong-common} in
\Cref{sec:strong});
\item a common local existence time for bounded initial states and
forces (\Cref{thm:uniform-local-existence} in
\Cref{sec:uniform-local-existence});
\item maximal strong solutions and a blow-up result, followed
by existence and stability of perturbed solutions on the full interval
of a reference solution
(\Cref{prop:maximal-strong,thm:stability-strong1} in
\Cref{sec:maximal-strong-solutions});
\item the selection of states in $\mathbf H^2\times\dot H^4$
arbitrarily close in $\mathbb V^1$ to a given state along a strong
solution (\Cref{lem:regular-times} in \Cref{sec:regular-times}).
\end{enumerate}
Some of the ingredients are classical (or can be obtained by adaptations of classical results), but we formulate them with the
quantitative dependence required below or in a slightly different form.
One example is the uniqueness of strong solutions. In a usual
uniqueness argument (see, e.g., \cite[Proof of Theorem 2.2, Step 6]{Boyer1999}),
the constants in the estimate for the difference
may depend on both solutions.
For our purposes, we need constants that are
determined by a single reference solution, its force, and the time
horizon. We establish this estimate for small perturbations,
uniformly in the length of the common interval of existence.
Thus the reference trajectory provides the bounds needed to estimate
the perturbation, even when the perturbed solution is only known to
exist on a shorter interval; see \Cref{rem:stability-smallness}.
\subsection{Weak and strong solutions}\label{sec:weak-strong-solutions}
We use the projected formulation \eqref{chn-abstract} to define
solutions, so the pressure does not appear in the definition.\footnote{
The pressure can be recovered in the usual way from the momentum
equation; see \Cref{rem:pressure-recovery} in
\Cref{app:auxiliary-results}.}

\begin{defi}[Weak solutions]\label{def:weak-solutions}
	Let $T>0$, $z_0=(u_0,\phi_0)\in\mathbb V^0$, and
	$f=(f_u,f_\phi)\in L^\infty(0,T;\mathbb V^0)$.
	A pair
	\begin{equation*}
		z=(u,\phi)\in L^\infty(0,T;\mathbb V^0)
		\cap L^2(0,T;\mathbb V^1)
	\end{equation*}
	is a weak solution of \eqref{abs} on $[0,T]$ with initial state
	$z_0$ and force $f$ if
	$\mu=-\Delta\phi+\phi^3-\phi\in L^2(0,T;H^1)$
	and, for every $v\in\mathbf H^1$ and $\psi\in H^1$,
	the identities
	\begin{align*}
		&\frac{\d}{\d t}(u,v)_{L^2}
		+\nu(\nabla u,\nabla v)_{L^2}+b_0(u,u,v)\\
		&\qquad=-\int_{\mathbb T^3}\Delta\phi\,\nabla\phi\cdot v\,dx
		+(f_u,v)_{L^2},\\
		&\frac{\d}{\d t}(\phi,\psi)_{L^2}
		+b_1(u,\phi,\psi)+(\nabla\mu,\nabla\psi)_{L^2}
		=(f_\phi,\psi)_{L^2}
	\end{align*}
	hold in the sense of distributions (in time), with the initial conditions
	\begin{equation*}
		u(0)=u_0,\qquad\phi(0)=\phi_0.
	\end{equation*}
\end{defi}

%{\color{red} We can not write Pointwise notion or distribution notion together as they are different. If there is time regularity in the definition then only we can write pointwise notion as in that case we have absolutely continuous representative, otherwise we only need to write with distribution}
%
%Also both of this papers are considering the time derivative term as $\int_0^T(\partial_t  u,v)_{L^2}\,dt$ with out intigration by parts and takes $v\in C_c^\infty((0,T)\times\mathbb T^3;\mathbb R^3)$ without the initial data in the weak formulation. We need to find the reason here.}
\begin{defi}[Strong solutions]\label{def:strong-solutions}
Let $T>0$, $z_0\in\mathbb V^1$, and
$f\in L^\infty(0,T;\mathbb V^0)$.
A weak solution $z=(u,\phi)$ of \eqref{abs} on $[0,T]$ with
initial state $z_0$ and force $f$ is called a strong solution if
\begin{equation*}
 z\in C([0,T];\mathbb V^1)
       \cap L^2(0,T;\mathbf H^2\times\dot H^4).
\end{equation*}
\end{defi}

\begin{rmk}\label{rem:solution-intervals}
For $T\in(0,\infty]$, a strong solution on $[0,T)$ is understood in
the sense of \Cref{def:strong-solutions} on every $[0,S]$, $0<S<T$.
\end{rmk}
We now state the existence results and energy estimates for weak and
strong solutions. Their proofs
are given in \Cref{app:auxiliary-results}.
\begin{theo}[Existence of weak solutions]\label{thm:existence-energy-m0}
 Let $T>0$, $z_0=(u_0,\phi_0)\in \mathbb V^0$, and
$f=(f_u,f_\phi)\in L^\infty(0,T;\mathbb V^0)$.
Then there exists at least one solution
\[
z=(u,\phi)\in L^\infty(0,T;\mathbb V^0)\cap L^2(0,T;\mathbb V^1)
\]
to \eqref{abs}. Moreover, the following estimate holds for almost every $t\in[0,T]$
\begin{equation}\label{energy-estimate-m0}
\|z(t)\|_{\mathbb V^0}^2
+
\int_0^t \|z(s)\|_{\mathbb V^1}^2\,ds
\le
C\Bigl(1 + \|f\|_{L^\infty(0,T;\mathbb V^0)}^2\Bigr),
\end{equation}
where $C>0$ depends at most on $T$, $u_0$ and $\phi_0$.
\end{theo}

\begin{theo}[Existence of strong solutions]\label{thm:existence-energy-m1}
Let $T>0$, $z_0=(u_0,\phi_0)\in \mathbb V^1$, and
$f=(f_u,f_\phi)\in L^\infty(0,T;\mathbb V^0)$.
Then there exists a time $T_*\in(0,T]$, depending only on
$\|z_0\|_{\mathbb V^1}$ and $\|f\|_{L^\infty(0,T;\mathbb V^0)}$, and at least one solution
\[
z=(u,\phi)\in {C([0,T_*);\mathbb V^1)}\cap {L^2(0,T_*;\mathbf H^2\times\dot H^4)}
\]
to \eqref{abs}. Moreover, for all $t\in[0,T_*)$, one has
\begin{equation}\label{energy-estimate-m1}
\|z(t)\|_{\mathbb V^1}^2
+
{\int_0^t \|z(s)\|_{\mathbf H^2\times\dot H^4}^2\,ds}
\leq
C,%\bigl(\|z_0\|_{\mathbb V^1},\|\eta\|_{L^\infty(0,T;\mathbb V^0)},T_*\bigr),
\end{equation}
where $C>0$ is a constant depending only on
$\|z_0\|_{\mathbb V^1}$, $\|f\|_{L^\infty(0,T;\mathbb V^0)}$, and $T_*$.
\end{theo}
\subsection{Stability and uniqueness of strong solutions}\label{sec:strong}
We now establish a stability estimate that provides continuous dependence on the initial data and the forcing term.
The proof is inspired by the arguments in
\cite[Proposition~1.1]{Ner21}. The result reads as follows.
\begin{theo}[Stability of strong solutions]\label{thm:stability-strong}
Let $\widetilde z =(\widetilde u,\widetilde\phi)$ be a strong solution of
\eqref{abs} on $[0,T]$ corresponding to the initial state
$\widetilde z_0=(\widetilde u_0,\widetilde\phi_0)\in\mathbb V^1$ and the force
$\widetilde f\in L^\infty(0,T;\mathbb V^0)$, and set
\begin{equation}\label{lambda_stability}
 \Lambda:=
 \|\widetilde z\|_{C([0,T];\mathbb V^1)}
 +{\|\widetilde z\|_{L^2(0,T;\mathbf H^2\times\dot H^4)}}
 +\|\widetilde f\|_{L^\infty(0,T;\mathbb V^0)}.
\end{equation}
There exist $\delta_*>0$ and $C>0$, depending only on $T$ and
$\Lambda$, such that, if $z_0=(u_0,\phi_0)\in\mathbb V^1$ and
$f\in L^\infty(0,T;\mathbb V^0)$ satisfy
\begin{equation}\label{eq:stability-smallness}
 d^2:=\|z_0-\widetilde z_0\|_{\mathbb V^1}^2
       +{\|f-\widetilde f\|_{L^2(0,T;\mathbb V^0)}^2}
 <\delta_*^2
\end{equation}
and $z=(u,\phi)$ is a strong solution of \eqref{abs} with initial
state $z_0$ and force $f$ on $[0,S]$, where $0<S\le T$, then
\begin{align}
&\sup_{t\in[0,S]}
\Bigl(
\|u(t)-\widetilde u(t)\|_{H^1}^2
+
\|\phi(t)-\widetilde\phi(t)\|_{\dot H^2}^2
\Bigr) \notag\\
&\qquad
+{\int_0^S}
\Bigl( 
\|u(t)-\widetilde u(t)\|_{H^2}^2
+
\|\phi(t)-\widetilde\phi(t)\|_{\dot H^4}^2
\Bigr)\,dt \notag\\
&\le
C\Bigl(
\|u_0-\widetilde u_0\|_{H^1}^2
+
\|\phi_0-\widetilde\phi_0\|_{\dot H^2}^2
+
{\|f-\widetilde f\|_{L^2(0,T;\mathbb V^0)}^2}
\Bigr).
\label{eq:strong-stability-estimate}
\end{align}
The constants $\delta_*$ and $C$ are independent of $S$.
\end{theo}
%{\color{red}In this theorem we are working on the common interval $[0,S]$ of existence. So this is a question is that in this case why do we need the smallness condition of the data? As I understand the smallness condition is needed to prove that new trajectory is also well-posed upto the time of existence for the main solution????
%
%{\color{blue} Indeed, we cannot avoid the smallness assumption. This is due to the argument since we only have one trajectory and we build another one close to it. I have explanied this further in Remark \ref{rem:stability-smallness} }
%
%Also in this result did we prove that $S$ can be $T$ for some suitable $\delta$?}
%
%{\color{blue}This is not the general case. It only happens if both exist until $T$ or by doing all the maximal stuff in the next section. I have rephrased slightly the theorem to reflect this. You can also check the new corollary \ref{thm:stability-strong1} inspired in what you wrote.}
\begin{proof}
Let $(u,\phi)$ and $(\widetilde u,\widetilde\phi)$ be two strong solutions of \eqref{abs}
corresponding to the initial data $z_0$ and $\widetilde z_0$, and forcing terms $f$ and $\widetilde f$, respectively, {on the common interval $[0,S]$}.
Define the differences
\begin{equation*}
w:=u-\widetilde u,\qquad \psi:=\phi-\widetilde\phi,\qquad \bar f:=f-\widetilde f.
\end{equation*}
Then $(w,\psi)$ satisfy the system
\begin{equation}\label{chnp2}
\begin{cases}
\partial_t w + \nu A_0 w + B_0(u,w) + B_0(w,\widetilde u)
= -R_0(\Delta\phi,\psi) - R_0(\Delta\psi,\widetilde\phi)
+ \bar f_u
& {0<t<S},\ x\in\mathbb T^3, \\
\partial_t \psi + \Delta^2\psi + \Delta\psi
+ B_1(w,\widetilde\phi) + B_1(u,\psi)
= \Delta(\phi^3)-\Delta(\widetilde\phi^3)
+ \bar f_\phi
& {0<t<S},\ x\in\mathbb T^3, \\
% \operatorname{div} w = 0
% & 0<t<T,\ x\in\mathbb T^3, \\
w(0)=u_0-\widetilde u_0, \qquad
\psi(0)=\phi_0-\widetilde\phi_0,
& x\in\mathbb T^3.
\end{cases}
\end{equation}
For readability, we divide the rest of the proof into three steps.

\smallskip
--\textit{Step 1. Estimates for $w$.} Multiplying the first equation in \eqref{chnp2} by $A_0w$ and integrating the result over $\mathbb{T}^3$, we have 
\begin{align} \notag
    \frac{1}{2}&\frac{\d}{\d t}\norm{\nabla w(t)}_{L^2}^2 
    + \nu\norm{\Delta w(t)}_{L^2}^2
    \\ \notag
    &= b_0(u,w,\Delta w) 
    + b_0(w,\widetilde u,\Delta w)
    + (R_0(\Delta\phi,\psi),\Delta w)_{L^2}
    + (R_0(\Delta\psi,\widetilde \phi),\Delta w)_{L^2}-(\bar f_u,\Delta w) \\ \label{eq:est_init_w}
    &=: \sum_{i=1}^{5}J_i.
\end{align}
Let us estimate each term $J_i$, $i\in\inter{1,5}$. For the first one, using
H\"older inequality together with \eqref{eq:u-L6} and \eqref{eq:grad-u-L3} in \Cref{lem:u-L6-grad-u-L3}, we obtain
\begin{align*}
|J_1|
&=|b_0(u,w,\Delta w)|
\le \|u\|_{L^6}\,\|\nabla w\|_{L^3}\|\Delta w\|_{L^2} \\
&\le C\|\nabla u\|_{L^2}
        \|\nabla w\|_{L^2}^{1/2}\|\Delta w\|_{L^2}^{3/2}.
\end{align*}
Then, by Young's inequality, for any $\epsilon\in(0,1)$,
\begin{equation*}
|J_1|
\le \epsilon \|\Delta w\|_{L^2}^{2}
+ C_{\epsilon}\|\nabla u\|_{L^2}^{4}\,
  \|\nabla w\|_{L^2}^{2}.
\end{equation*}
For the second term, using H\"older inequality together with
\eqref{eq:u-L6} and \eqref{eq:grad-u-L3} in \Cref{lem:u-L6-grad-u-L3}, we obtain
\begin{align*}
|b_0(w,\widetilde u,\Delta w)|
&\le \|w\|_{L^6} \|\nabla \widetilde u\|_{L^3}\,\|\Delta w\|_{L^2} \\
&\le C \|\nabla w\|_{L^2} \|\nabla \widetilde u\|_{L^2}^{1/2} \|\Delta \widetilde u\|_{L^2}^{1/2}
        \|\Delta w\|_{L^2}.
\end{align*}
Hence, by Young's inequality, for any $\epsilon\in(0,1)$,
\begin{align*}
|J_2|=|b_0(w,\widetilde u,\Delta w)|
\le \epsilon \|\Delta w\|_{L^2}^2
+ C_\epsilon 
\|\nabla \widetilde u\|_{L^2} \|\Delta \widetilde u\|_{L^2} 
\|\nabla w\|_{L^2}^{2}.
\end{align*}
Moving to the third term, using H\"older inequality and estimate \eqref{eq:grad-infty} in \Cref{lem:interp-estimates}, we obtain
\begin{align*}
|J_3|
&=\bigl|(R_0(\Delta\phi,\psi),\Delta w)_{L^2}\bigr|
\le \|\Delta\phi\|_{L^2}\,\|\nabla\psi\|_{L^\infty}\,\|\Delta w\|_{L^2} \\
&\le C\,\|\Delta\phi\|_{L^2}\,\|\nabla\psi\|_{L^2}^{1/3}\,\|\Delta^2\psi\|_{L^2}^{2/3}\,\|\Delta w\|_{L^2}.
\end{align*}
Consequently, by Young's inequality, for any $\epsilon,\delta\in(0,1)$,
\begin{align*}
|J_3|
\le \epsilon\|\Delta w\|_{L^2}^2
+\delta\|\Delta^2\psi\|_{L^2}^2
+ C_{\epsilon,\delta}\,\|\Delta\phi\|_{L^2}^6\,\|\nabla\psi\|^2_{L^2}.
\end{align*}
Finally, using Poincar\'e inequality, we deduce that
\begin{align*}
|J_3|
\le \epsilon\|\Delta w\|_{L^2}^2
+\delta\|\Delta^2\psi\|_{L^2}^2
+ C_{\epsilon,\delta}\,\|\Delta\phi\|_{L^2}^6\,\|\Delta\psi\|^2_{L^2}.
\end{align*}
For the fourth term, using \eqref{def_R0}, H\"older inequality and
estimate \eqref{eq:grad-infty} in \Cref{lem:interp-estimates}, we obtain
\begin{align*}
|J_4|
&=\bigl|(R_0(\Delta\psi,\widetilde\phi),\Delta w)_{L^2}\bigr|
\le \|\Delta\psi\|_{L^2}\,\|\nabla\widetilde\phi\|_{L^\infty}\,\|\Delta w\|_{L^2} \\
&\le C\,\|\Delta\psi\|_{L^2}\,
        \|\nabla\widetilde\phi\|_{L^2}^{1/3}\,\|\Delta^2\widetilde\phi\|_{L^2}^{2/3}\,
        \|\Delta w\|_{L^2}.
\end{align*}
Consequently, by Young's inequality, for any $\epsilon\in(0,1)$,
\begin{align*}
|J_4|
&\le \epsilon\|\Delta w\|_{L^2}^2
+ C_{\epsilon}\,\|\nabla\widetilde\phi\|_{L^2}^{2/3}\,\|\Delta^2\widetilde\phi\|_{L^2}^{4/3}\,
\|\Delta\psi\|_{L^2}^2.
\end{align*}
For the last term, by Cauchy--Schwarz and Young's inequalities, for any
$\epsilon\in(0,1)$,
\begin{align*}
|J_5|=\bigl|(\bar f_u,\Delta w)_{L^2}\bigr|
\le \|\bar f_u\|_{L^2}\,\|\Delta w\|_{L^2}
\le \epsilon\|\Delta w\|_{L^2}^2 + C_\epsilon\|\bar f_u\|_{L^2}^2.
\end{align*}
Collecting the estimates for $J_i$, $i\in\inter{1,5}$, and using them in \eqref{eq:est_init_w}, we can choose
$\epsilon\in(0,1)$ small enough (depending only on $\nu$), to deduce that for any
$\delta\in(0,1)$
\begin{align}\label{eq:w-H1-diff-ineq}
\frac12\frac{\d}{\d t}\|\nabla w\|_{L^2}^2
+ \frac{\nu}{2}\|\Delta w\|_{L^2}^2
\le\,
& C\|\nabla u\|_{L^2}^{4}\,\|\nabla w\|_{L^2}^2 + C\|\nabla \widetilde u\|_{L^2}\,\|\Delta \widetilde u\|_{L^2}
\,\|\nabla w\|_{L^2}^2 \notag\\
&+ C\Bigl(\|\Delta\phi\|_{L^2}^{6}
+\|\nabla\widetilde\phi\|_{L^2}^{2/3}\,\|\Delta^2\widetilde\phi\|_{L^2}^{4/3}\Bigr)
\,\|\Delta\psi\|_{L^2}^2 \notag\\
&+ \delta\|\Delta^2\psi\|_{L^2}^2
+ C\|\bar f_u\|_{L^2}^2.
\end{align}

\smallskip
--\textit{Step 2. Estimates for $\psi$.} We now derive the corresponding estimate for $\psi$. To this end, we multiply the
second equation in \eqref{chnp2} by $\Delta^2\psi$ and integrate over
$\mathbb T^3$. This yields
\begin{align}\notag
\frac12\frac{\d}{\d t}\|\Delta\psi(t)\|_{L^2}^2
+\|\Delta^2\psi(t)\|_{L^2}^2
&=
-b_1(w,\widetilde\phi,\Delta^2\psi)
-b_1(u,\psi,\Delta^2\psi) \\
&\quad
+(\Delta(\phi^3)-\Delta(\widetilde\phi^3),\Delta^2\psi)_{L^2}
+(\bar f_\phi,\Delta^2\psi)_{L^2}
-(\Delta\psi,\Delta^2\psi)_{L^2} \notag\\
&=: \sum_{i=1}^{5} Q_i. \label{eq:iden-init-Dpsi}
\end{align}

Let us estimate each of the $Q_i$, $i\in\inter{1,5}$. For the first term, using H\"older inequality together with
\eqref{eq:u-L6} and \eqref{eq:grad-u-L3} in \Cref{lem:u-L6-grad-u-L3}, we obtain
\begin{align*}
|Q_1|
&=|b_1(w,\widetilde\phi,\Delta^2\psi)|
\le \|w\|_{L^6}\,\|\nabla\widetilde\phi\|_{L^3}\,\|\Delta^2\psi\|_{L^2} \\
&\le C\,\|\nabla w\|_{L^2}\,
\|\nabla\widetilde\phi\|_{L^2}^{1/2}\|\Delta\widetilde\phi\|_{L^2}^{1/2}
\,\|\Delta^2\psi\|_{L^2}. 
\end{align*}
Consequently, by Young's inequality,
\begin{align}
|Q_1|
&\le \delta\|\Delta^2\psi\|_{L^2}^2
+ C_\delta\,
\|\nabla\widetilde\phi\|_{L^2}\|\Delta\widetilde\phi\|_{L^2}\,
\|\nabla w\|_{L^2}^{2}. \label{eq:J6-est}
\end{align}
For the second term, by H\"older's inequality together with \eqref{eq:u-L6} in
\Cref{lem:u-L6-grad-u-L3}, and Poincar\'e inequality, we obtain
\begin{align*}
|Q_2|
&=|b_1(u,\psi,\Delta^2\psi)|
\le \|u\|_{L^6}\|\nabla\psi\|_{L^3}\|\Delta^2\psi\|_{L^2} \\
&\le C\,\|\nabla u\|_{L^2}\,
        \|\nabla\psi\|_{L^2}^{1/2}\|\Delta\psi\|_{L^2}^{1/2}\,
        \|\Delta^2\psi\|_{L^2} \\
&\le C\,\|\nabla u\|_{L^2}\,
        \|\Delta\psi\|_{L^2}\,
        \|\Delta^2\psi\|_{L^2}.
\end{align*}
Consequently, by Young's inequality,
\begin{align}
|Q_2|
\le \delta\|\Delta^2\psi\|_{L^2}^2
+ C_\delta\|\nabla u\|_{L^2}^2\|\Delta\psi\|_{L^2}^2.
\label{eq:Q2-est}
\end{align}
For the third term, we observe that $\phi^3-\widetilde\phi^3=\psi G$ where $G:=\phi^2+\phi\widetilde\phi+\widetilde\phi^2$. Then
\begin{equation*}
Q_3=(\Delta(\phi^3)-\Delta(\widetilde\phi^3),\Delta^2\psi)_{L^2}
=(\Delta(\psi G),\Delta^2\psi)_{L^2},
\end{equation*}
and therefore
\begin{equation*}
|Q_3|
\le \|\Delta(\psi G)\|_{L^2}\|\Delta^2\psi\|_{L^2}
\le \|\psi G\|_{H^2}\|\Delta^2\psi\|_{L^2}.
\end{equation*}
Since $H^2(\mathbb T^3)$ is an algebra, and using that $\psi$ has zero mean, we obtain
\begin{equation*}
\|\psi G\|_{H^2}
\le C\|\psi\|_{H^2}\|G\|_{H^2}
\le C\|\Delta\psi\|_{L^2}\|G\|_{H^2}.
\end{equation*}
Using again that $H^2(\mathbb T^3)$ is an algebra, and that $\phi$ and $\widetilde\phi$ have zero mean, we get
\begin{align*}
\|G\|_{H^2}
&\le C\Bigl(\|\phi\|_{H^2}^2
+\|\phi\|_{H^2}\|\widetilde\phi\|_{H^2}
+\|\widetilde\phi\|_{H^2}^2\Bigr)
\le C\Bigl(\|\Delta\phi\|_{L^2}^2+\|\Delta\widetilde\phi\|_{L^2}^2\Bigr).
\end{align*}
Consequently,
\begin{equation*}
|Q_3|
\le C\Bigl(\|\Delta\phi\|_{L^2}^2+\|\Delta\widetilde\phi\|_{L^2}^2\Bigr)
\|\Delta\psi\|_{L^2}\|\Delta^2\psi\|_{L^2}.
\end{equation*}
Finally, by Young's inequality,
\begin{equation}
|Q_3|
\le \delta\|\Delta^2\psi\|_{L^2}^2
+ C_\delta\Bigl(\|\Delta\phi\|_{L^2}^4+\|\Delta\widetilde\phi\|_{L^2}^4\Bigr)
\|\Delta\psi\|_{L^2}^2.
\label{eq:Q3-est}
\end{equation}
For the fourth term, by Cauchy--Schwarz and Young's inequalities,
\begin{align}
|Q_4|
&=\bigl|(\bar f_\phi,\Delta^2\psi)_{L^2}\bigr|
\le \|\bar f_\phi\|_{L^2}\,\|\Delta^2\psi\|_{L^2}
\le \delta\|\Delta^2\psi\|_{L^2}^2 + C_\delta\|\bar f_\phi\|_{L^2}^2.
\label{eq:Q4-est}
\end{align}
and, similarly, 
\begin{align}
|Q_5|
&=\bigl|(\Delta\psi,\Delta^2\psi)_{L^2}\bigr|
\le \|\Delta\psi\|_{L^2}\,\|\Delta^2\psi\|_{L^2}
\le \delta\|\Delta^2\psi\|_{L^2}^2 + C_\delta\|\Delta\psi\|_{L^2}^2.
\label{eq:Q5-est}
\end{align}
Putting estimates \eqref{eq:J6-est}--\eqref{eq:Q5-est} into
\eqref{eq:iden-init-Dpsi}, we deduce that, choosing $\delta\in(0,1)$ sufficiently small,
\begin{align}\label{eq:psi-H2-diff-ineq-abs}
\frac12\frac{\d}{\dt}\|\Delta\psi(t)\|_{L^2}^2
+\frac12\|\Delta^2\psi(t)\|_{L^2}^2
\le\,
& C\,\|\nabla\widetilde\phi(t)\|_{L^2}\,\|\Delta\widetilde\phi(t)\|_{L^2}
\,\|\nabla w(t)\|_{L^2}^2 \notag\\
&+ C\Bigl(1+\|\nabla u(t)\|_{L^2}^2
+\|\Delta\phi(t)\|_{L^2}^4+\|\Delta\widetilde\phi(t)\|_{L^2}^4\Bigr)\,
\|\Delta\psi(t)\|_{L^2}^2 \notag\\
&+ C\,\|\bar f_\phi(t)\|_{L^2}^2.
\end{align}

\smallskip
-- \textit{Step 3. Collecting the estimates and conclusion}. We now sum \eqref{eq:w-H1-diff-ineq} and \eqref{eq:psi-H2-diff-ineq-abs}. 
Decreasing (if necessary) $\delta\in(0,1)$, the term
$\delta\|\Delta^2\psi(t)\|_{L^2}^2$ on the right-hand side of
\eqref{eq:w-H1-diff-ineq} is absorbed into the left-hand side of
\eqref{eq:psi-H2-diff-ineq-abs}. More precisely, setting
\begin{equation*}
E(t):=\|\nabla w(t)\|_{L^2}^2+\|\Delta\psi(t)\|_{L^2}^2,
\end{equation*}
we obtain
\begin{align}\notag
\frac12&\frac{\d}{\dt}E(t)
+\frac{\nu}{2}\|\Delta w(t)\|_{L^2}^2
+\frac14\|\Delta^2\psi(t)\|_{L^2}^2
\\ \notag
& \leq C\Bigl(\|\nabla u(t)\|_{L^2}^{4}
+\|\nabla\widetilde u(t)\|_{L^2}\|\Delta\widetilde u(t)\|_{L^2}
+\|\nabla\widetilde\phi(t)\|_{L^2}\|\Delta\widetilde\phi(t)\|_{L^2}\Bigr)
\|\nabla w(t)\|_{L^2}^2 \\ \notag
&\quad + C\Bigl(1+\|\nabla u(t)\|_{L^2}^{2}
+\|\Delta\phi(t)\|_{L^2}^{6}
+\|\nabla\widetilde\phi(t)\|_{L^2}^{2/3}\|\Delta^2\widetilde\phi(t)\|_{L^2}^{4/3} \\
&\qquad\quad  
+\|\Delta\phi(t)\|_{L^2}^{4}
+\|\Delta\widetilde\phi(t)\|_{L^2}^{4}\Bigr)
\|\Delta\psi(t)\|_{L^2}^{2}  \label{eq:combined-diff-ineq}
+ C\bigl(\|\bar f_u(t)\|_{L^2}^2+\|\bar f_\phi(t)\|_{L^2}^2\bigr).
\end{align}
By \eqref{eq:periodic-norm-equivalences}, there is a constant $C'\ge1$ such that $E(0)\le C'd^2$. Let us choose $\delta_*\le(2C')^{-1/2}$, so $E(0)<1/2$.
Since $z,\widetilde z\in C([0,S];\mathbb V^1)$, the function $E$
is continuous. Set
$\tau:=\inf\{t\in(0,S]:E(t)=1\}\wedge S$.
Since $E(0)<1/2$, we have $\tau>0$ and
$E(t)\le1$ for every $t\in[0,\tau]$.
The definition of $\Lambda$ and
\eqref{eq:periodic-norm-equivalences} give
\begin{equation*}
 \sup_{t\in[0,T]}
 \bigl(\|\nabla\widetilde u(t)\|_{L^2}^2
       +\|\Delta\widetilde\phi(t)\|_{L^2}^2\bigr)
 +\int_0^T
 \bigl(\|\Delta\widetilde u(t)\|_{L^2}^2
       +\|\Delta^2\widetilde\phi(t)\|_{L^2}^2\bigr)\,dt
 \le C(1+\Lambda^2) =: K_{\Lambda}.
\end{equation*}
Since $u=\widetilde u+w$ and $\phi=\widetilde\phi+\psi$, increasing the
constant in $K_\Lambda$ if necessary and using that $E(t)\le1$ for all $t\in[0,\tau]$, we also have
\begin{equation*}
 \sup_{t\in[0,\tau]}
 \bigl(\|\nabla u(t)\|_{L^2}^2+\|\Delta\phi(t)\|_{L^2}^2\bigr)
 \le K_\Lambda.
\end{equation*}
Using these bounds, we now estimate the coefficients in \eqref{eq:combined-diff-ineq} on {$[0,\tau]$}. First,
\begin{equation*}
{\int_0^{\tau}}\|\nabla u(t)\|_{L^2}^{4}\,dt
\le
T{\sup_{t\in[0,\tau]}}\|\nabla u(t)\|_{L^2}^{4}
\le T {K_\Lambda^2},
\end{equation*}
and similarly
\begin{equation*}
{\int_0^{\tau}}\|\nabla u(t)\|_{L^2}^{2}\,dt
\le
T{\sup_{t\in[0,\tau]}}\|\nabla u(t)\|_{L^2}^{2}
\le T {K_\Lambda}.
\end{equation*}
Next, by Cauchy--Schwarz inequality,
\begin{align*}
\int_0^{T}\|\nabla\widetilde u(t)\|_{L^2}\|\Delta\widetilde u(t)\|_{L^2}\,dt
&\le
\left(\int_0^{T}\|\nabla\widetilde u(t)\|_{L^2}^{2}\,dt\right)^{1/2}
\left(\int_0^{T}\|\Delta\widetilde u(t)\|_{L^2}^{2}\,dt\right)^{1/2} \\
&\le
T^{1/2}
\left(\sup_{t\in[0,T]}\|\nabla\widetilde u(t)\|_{L^2}^{2}\right)^{1/2}
\left(\int_0^{T}\|\Delta\widetilde u(t)\|_{L^2}^{2}\,dt\right)^{1/2} \\
&\le T^{1/2}\,{K_\Lambda}.
\end{align*}
Similarly, using Poincar\'e inequality,
\begin{align*}
\int_0^{T}\|\nabla\widetilde\phi(t)\|_{L^2}\|\Delta\widetilde\phi(t)\|_{L^2}\,dt
&\le
T\sup_{t\in[0,T]}\|\nabla\widetilde\phi(t)\|_{L^2}\|\Delta\widetilde\phi(t)\|_{L^2} \\
&\le
C\,T\sup_{t\in[0,T]}\|\Delta\widetilde\phi(t)\|_{L^2}^{2}
\le C\,T\,{K_\Lambda}.
\end{align*}
Moreover,
\begin{equation*}
{\int_0^{\tau}}\|\Delta\phi(t)\|_{L^2}^{6}\,dt
\le
T{\sup_{t\in[0,\tau]}}\|\Delta\phi(t)\|_{L^2}^{6}
\le T {K_\Lambda^3},
\end{equation*}
and
\begin{equation*}
{\int_0^{\tau}}\|\Delta\phi(t)\|_{L^2}^{4}\,dt
+\int_0^{T}\|\Delta\widetilde\phi(t)\|_{L^2}^{4}\,dt
\le
T{\sup_{t\in[0,\tau]}}\|\Delta\phi(t)\|_{L^2}^{4}
+
T\sup_{t\in[0,T]}\|\Delta\widetilde\phi(t)\|_{L^2}^{4}
\le 2T {K_\Lambda^2}.
\end{equation*}
Finally, since $\|\nabla\widetilde\phi\|_{L^2}\in L^\infty(0,T)$ and
$\|\Delta^2\widetilde\phi\|_{L^2}\in L^2(0,T)$, and since
$L^2(0,T)\hookrightarrow L^{4/3}(0,T)$, we have
\begin{align*}
\int_0^{T}\|\nabla\widetilde\phi(t)\|_{L^2}^{2/3}\|\Delta^2\widetilde\phi(t)\|_{L^2}^{4/3}\,dt
&\le
\left(\sup_{t\in[0,T]}\|\nabla\widetilde\phi(t)\|_{L^2}^{2/3}\right)
\int_0^{T}\|\Delta^2\widetilde\phi(t)\|_{L^2}^{4/3}\,dt \\
&\le
T^{1/3}
\left(\sup_{t\in[0,T]}\|\nabla\widetilde\phi(t)\|_{L^2}^{2/3}\right)
\left(\int_0^{T}\|\Delta^2\widetilde\phi(t)\|_{L^2}^{2}\,dt\right)^{2/3} \\
&\le C\,T^{1/3}\,{K_\Lambda}.
\end{align*}
Collecting the above estimates, the coefficients in
\eqref{eq:combined-diff-ineq} on $[0,\tau]$ are bounded by a
nonnegative function $\Gamma\in L^1(0,T)$ satisfying
\begin{equation*}
 \int_0^T\Gamma(t)\,dt\le C(T,\Lambda).
\end{equation*}
Consequently, \eqref{eq:combined-diff-ineq} gives
\begin{equation}\label{eq:stability-reference-diff-ineq}
 \frac{\d}{\dt}E(t)
 +c\bigl(\|\Delta w(t)\|_{L^2}^2
         +\|\Delta^2\psi(t)\|_{L^2}^2\bigr)
 \le\Gamma(t)E(t)+C\|\bar f(t)\|_{\mathbb V^0}^2,
 \qquad 0<t<\tau,
\end{equation}
for some $c>0$. Applying Gr\"onwall's inequality, we obtain
\begin{align}\label{eq:bound_casi_sup_E}
 E(t)
 &\le
 \exp\left(\int_0^t\Gamma(s)\,ds\right)
 \left(E(0)+C\int_0^t\|\bar f(s)\|_{\mathbb V^0}^2\,ds\right)\le {C(T,\Lambda)d^2},
 \qquad {t\in[0,\tau]}.
\end{align}
Decreasing $\delta_*$ if necessary, we may assume
$C(T,\Lambda)\delta_*^2<1/2$. If $E$ attains $1$ on $[0,S]$,
then $E(\tau)=1$ by the definition of $\tau$, whereas the estimate
above gives $E(\tau)<1/2$, a contradiction. Thus $\tau=S$ and
$E(t)<1$ on $[0,S]$, 
and
\eqref{eq:stability-reference-diff-ineq} and the above estimate hold on the whole interval.

Integrating \eqref{eq:stability-reference-diff-ineq} over $[0,S]$
and using $E(S)\ge0$, we obtain
\begin{align*}
 c\int_0^S
 \bigl(\|\Delta w(t)\|_{L^2}^2
       +\|\Delta^2\psi(t)\|_{L^2}^2\bigr)\,dt
 &\le E(0)+\int_0^S\Gamma(t)E(t)\,dt
                  +C\int_0^S\|\bar f(t)\|_{\mathbb V^0}^2\,dt\\
 &\le C(T,\Lambda)d^2.
\end{align*}
This, together with \eqref{eq:bound_casi_sup_E} and
\eqref{eq:periodic-norm-equivalences}, yields the desired estimate with $C=C(T,\Lambda)$, independently of $S$.
\end{proof}
\begin{rmk}\label{rem:stability-smallness}
In a usual proof of uniqueness, the coefficients in the estimate
for the difference may depend on both solutions.
In \Cref{thm:stability-strong}, the constants depend only on $T$
and $\Lambda$, uniformly in $S$.
Condition \eqref{eq:stability-smallness} is used after
\eqref{eq:combined-diff-ineq} to keep $E\le1$ and bound the
$L^1(0,S)$ norms of the coefficients in terms of $T$ and $\Lambda$.
\end{rmk}
The previous stability estimate immediately yields uniqueness of strong solutions on their interval of existence.
\begin{proposition}[Uniqueness of strong solutions]\label{cor:unique-strong-common}
Let $T>0$ and let $z$ be a strong solution of \eqref{abs}
on $[0,T)$ with initial state $z_0\in\mathbb V^1$ and force
$f\in L^\infty(0,T;\mathbb V^0)$.
Then \eqref{abs} admits at most one strong solution in
$C([0,T);\mathbb V^1)
\cap L^2_{\mathrm{loc}}([0,T);\mathbf H^2\times\dot H^4)$
with initial state $z_0$ and force $f$.
\end{proposition}
\begin{proof}
Let $z$ and $\widetilde z$ be two strong solutions of \eqref{abs}
on ${[0,T)}$ with the same initial data $z_0$ and forcing term $f$.
Fix $\tau\in(0,{T})$. Since
$\widetilde z\in C([0,\tau];\mathbb V^1)\cap L^2(0,\tau;\mathbf H^2\times\dot H^4)$,
{the quantity $\Lambda$ defined by} \eqref{lambda_stability} in \Cref{thm:stability-strong} with $T=\tau$
is finite. The initial data and forcing terms coincide so $d=0$.
Applying \Cref{thm:stability-strong} with $T=S=\tau$, we obtain
$z(t)=\widetilde z(t)$ for every $t\in[0,\tau]$.
Since $\tau<T$ is arbitrary, the two solutions coincide on
${[0,T)}$.
\end{proof}
\subsection{Uniform local existence}
\label{sec:uniform-local-existence}
The following result is a quantitative version of
\Cref{thm:existence-energy-m1}. The comparison argument in
\Cref{app:strong-existence}, with explicit constants for bounded initial
states and forces, gives a common local existence time.
\begin{theo}[Uniform local existence]
\label{thm:uniform-local-existence}
%\label{thm:global-uniform-existence}
For every $R>0$ and $M\geq 0$, there exists $\theta_{R,M}\in(0,1/4]$,
depending only on $R$, $M$, and $\nu$ such that,
for $z_0\in\mathbb V^1$ and $f\in L^\infty(0,1;\mathbb V^0)$
satisfying
\begin{equation}\label{eq:uniform-local-bounds}
 \|z_0\|_{\mathbb V^1}\leq R,\qquad
 \|f\|_{L^\infty(0,1;\mathbb V^0)}\leq M,
\end{equation}
system \eqref{abs} admits a unique strong solution
\begin{equation}\label{eq:uniform-local-regularity}
 z\in C([0,2\theta_{R,M}];\mathbb V^1)
       \cap L^2(0,2\theta_{R,M};\mathbf H^2\times\dot H^4).
\end{equation}
Moreover,
\begin{equation}\label{eq:uniform-local-estimate}
 \sup_{t\in[0,2\theta_{R,M}]}\|z(t)\|_{\mathbb V^1}^2
 +\int_0^{2\theta_{R,M}}\|z(t)\|_{\mathbf H^2\times\dot H^4}^2\,dt
 \leq C(R,M),
\end{equation}
where $C(R,M)>0$ depends only on $R$, $M$, and $\nu$.
The time $\theta_{R,M}$ can be chosen so that it does not increase
when either $R$ or $M$ is increased.
\end{theo}
\begin{proof}
We follow the proof of \Cref{thm:existence-energy-m1}
with $T=1$ and keep track of the constants for data satisfying
\eqref{eq:uniform-local-bounds}.
Let us set
\begin{align*}
 E(t)&:=\|u(t)\|_{L^2}^2+\|\nabla\phi(t)\|_{L^2}^2
                   +2\int_{\mathbb T^3}F(\phi(t))\,dx,\\
 \mathcal Y(t)&:=\|\nabla u(t)\|_{L^2}^2
                   +\|\Delta\phi(t)\|_{L^2}^2,\\
 \mathcal D(t)&:=\|\Delta u(t)\|_{L^2}^2
                   +\|\Delta^2\phi(t)\|_{L^2}^2.
\end{align*}
By \Cref{lem:quartic-potential}, \eqref{eq:phi-infty}, and
\eqref{eq:periodic-norm-equivalences}, the bound
$\|z_0\|_{\mathbb V^1}\leq R$ in
\eqref{eq:uniform-local-bounds} gives
\begin{equation*}
 E(0)+\mathcal Y(0)
 \leq C_0(1+R^2+R^4)=:C_0(R),
\end{equation*}
for a fixed constant $C_0\geq1$.
By \eqref{iden_energy_f}, we have
\begin{equation*}
 E'\leq C_1E+C_1M^2.
\end{equation*}
Here $C_0$ and $C_1>0$ are independent of $z_0$, $f$, $R$,
and $M$. Gronwall's inequality gives
$E(t)\leq e^{C_1t}(C_0(R)+C_1M^2t)$.
Since $F\geq0$, we have
\begin{equation}\label{eq:uniform-local-low-energy}
 \|u(t)\|_{L^2}^2+\|\nabla\phi(t)\|_{L^2}^2
 \leq E(t)\leq e^{C_1}(C_0(R)+C_1M^2)=:K_{R,M}.
\end{equation}
Following the proof of \Cref{thm:existence-energy-m1}, adding \eqref{eq:est_nabla-ut} and \eqref{est_delta_phi},
absorbing the phase term, and using
\eqref{eq:uniform-local-low-energy}, we obtain
\begin{align}
 \mathcal Y'+c\mathcal D
 &\leq C\bigl((1+K_{R,M})\mathcal Y^3
       +(1+K_{R,M}^3+K_{R,M}^4)\mathcal Y+M^2\bigr)\notag\\
 &\leq C_2(1+K_{R,M}^4+M^2)(1+\mathcal Y)^3,
 \label{eq:uniform-local-differential}
\end{align}
where the constants $c>0$ and $C_2>0$ are independent of $z_0$, $f$, $R$,
and $M$. Set $C_2(1+K_{R,M}^4+M^2)=:C_{R,M}$.
Consider the equation
\begin{equation*}
 q'=C_{R,M}(1+q)^3,\qquad q(0)=C_0(R).
\end{equation*}
Its solution is
\begin{equation*}
 q(t)=\bigl((1+C_0(R))^{-2}-2C_{R,M}t\bigr)^{-1/2}-1,
 \qquad 0\leq t<T_q,
\end{equation*}
where $T_q:=\frac{1}{2C_{R,M}(1+C_0(R))^2}$. Since $\mathcal D\geq0$ and $\mathcal Y(0)\leq q(0)$,
the comparison principle applied to
\eqref{eq:uniform-local-differential} gives
$\mathcal Y(t)\leq q(t)$ on their common interval of existence. Hence
\begin{equation}\label{eq:uniform-local-inverse-bound}
 (1+\mathcal Y(t))^{-2}
 \geq(1+q(t))^{-2}
 =(1+C_0(R))^{-2}-2C_{R,M}t.
\end{equation}
Choose
\begin{equation}\label{eq:uniform-local-time-choice}
 h_{R,M}:=\min\{1,T_q/2\},
 \qquad \theta_{R,M}:=\frac{h_{R,M}}4.
\end{equation}
Then $h_{R,M}<T_q$, and for $0\leq t\leq h_{R,M}$,
\eqref{eq:uniform-local-inverse-bound} gives
\begin{equation*}
 (1+\mathcal Y(t))^{-2}\geq\tfrac12(1+C_0(R))^{-2},
 \qquad 1+\mathcal Y(t)\leq\sqrt2(1+C_0(R)).
\end{equation*}
Integrating \eqref{eq:uniform-local-differential}, we obtain
\begin{align*}
 c\int_0^t\mathcal D(r)\,dr
 &\leq\mathcal Y(0) + C_{R,M}\int_0^t(1+\mathcal Y(r))^3\,dr\notag\\
 &\leq C_0(R)+2^{3/2}C_{R,M}h_{R,M}(1+C_0(R))^3.
% \label{eq:uniform-local-dissipation}
\end{align*}
Since $h_{R,M}<T_q$, the existence time in the proof of
\Cref{thm:existence-energy-m1}, with $T=1$, can be taken equal
to $h_{R,M}$, and the solution satisfies
\begin{equation*}
z\in C([0,h_{R,M});\mathbb V^1)
\cap L^2(0,h_{R,M};\mathbf H^2\times\dot H^4).
\end{equation*}
By \Cref{cor:unique-strong-common}, this solution
is unique.
Since $2\theta_{R,M}=h_{R,M}/2$, restriction to
$[0,2\theta_{R,M}]$ gives \eqref{eq:uniform-local-regularity}.
The bounds above and \eqref{eq:periodic-norm-equivalences}
yield \eqref{eq:uniform-local-estimate}.
Finally, if either $R$ or $M$ is increased, the quantities
$C_0(R)$, $K_{R,M}$, and $C_{R,M}$ do not decrease.
Thus $T_q$ does not increase, and
\eqref{eq:uniform-local-time-choice} gives the same property
for $\theta_{R,M}$.
\end{proof}
\subsection{Maximal strong solutions and some consequences}
\label{sec:maximal-strong-solutions}
A strong solution of \eqref{abs} on $[0,T)$ is called maximal if
it admits no strong extension, with the same force, to an interval
$[0,T')$ with $T'>T$.
The uniform existence time in \Cref{thm:uniform-local-existence},
together with uniqueness in \Cref{cor:unique-strong-common}, yields
the following result.
\begin{cor}[Blow-up of strong solutions]\label{prop:maximal-strong}
Let $z_0\in\mathbb V^1$ and
$f\in L^\infty_{\mathrm{loc}}(\mathbb R_+;\mathbb V^0)$.
Then there exist $T_{\max}=T_{\max}(z_0,f)\in(0,\infty]$ and a unique maximal strong solution
\begin{equation*}
 z\in C([0,T_{\max});\mathbb V^1)
       \cap L^2_{\mathrm{loc}}([0,T_{\max});\mathbf H^2\times\dot H^4)
\end{equation*}
of \eqref{abs} with initial state $z_0$ and force $f$.
Moreover, if $T_{\max}<\infty$, then
\begin{equation}\label{eq:maximal-strong-blowup}
 \lim_{t\uparrow T_{\max}}\|z(t)\|_{\mathbb V^1}=\infty.
\end{equation}
We omit the arguments of $T_{\max}$ when no confusion can arise.
\end{cor}
\begin{proof}
The proof is standard; see, for instance,
\cite[Section~5, Corollary~5.2]{Tao13}.
We include it for completeness.
Define the set of admissible times
\begin{equation*}
 \begin{aligned}
 \mathscr A:=\Bigl\{T>0:\ &\eqref{abs} \text{ with force } f
 \text{ admits a solution } z_T \text{ with } z_T(0)=z_0
 \text{ and}\\
 &z_T\in C([0,T);\mathbb V^1)
 \cap L^2_{\mathrm{loc}}([0,T);\mathbf H^2\times\dot H^4)\Bigr\},
 \end{aligned}
\end{equation*}
and set $T_{\max}:=\sup\mathscr A$.
Applying \Cref{thm:uniform-local-existence} to $z_0$ and
$f|_{(0,1)}$, we see that $\mathscr A\neq\emptyset$, so that
$T_{\max}\in(0,\infty]$.
Let $T,T'\in\mathscr A$ with $T<T'$. Since $z_T$ and the restriction of
$z_{T'}$ to $[0,T)$ have the same initial state and force,
\Cref{cor:unique-strong-common} gives $z_T=z_{T'}$ on $[0,T)$.
Therefore, setting $z(t):=z_T(t)$ for any $T\in\mathscr A$ with $T>t$
(such $T$ exists since $t<T_{\max}=\sup\mathscr A$), we obtain a function
on $[0,T_{\max})$ which, by the previous identity, does not depend on the
choice of $T$.
For every $S\in(0,T_{\max})$, choose $T\in\mathscr A$ with $S<T$.
Since $z=z_T$ on $[0,S]$, $z$ has the asserted regularity and solves
\eqref{abs} with initial state $z_0$ and force $f$.
Uniqueness follows once again from \Cref{cor:unique-strong-common},
and maximality from the definition of $T_{\max}$.
It remains to prove \eqref{eq:maximal-strong-blowup}.

Suppose that $T_{\max}<\infty$ and that
\eqref{eq:maximal-strong-blowup} does not hold.
There exist $R>0$ and a sequence $t_n\uparrow T_{\max}$ such that
\begin{equation*}
 \|z(t_n)\|_{\mathbb V^1}\leq R,\qquad n\geq1.
\end{equation*}
For every $n\geq1$, define the translated force
\begin{equation*}
 f_n(s):=f(t_n+s),\qquad 0<s<1.
\end{equation*}
Since $f\in L^\infty_{\mathrm{loc}}(\mathbb R_+;\mathbb V^0)$,
\begin{equation*}
 \|f_n\|_{L^\infty(0,1;\mathbb V^0)}
 \leq \|f\|_{L^\infty(0,T_{\max}+1;\mathbb V^0)}
 =:M<\infty.
\end{equation*}
By \Cref{thm:uniform-local-existence}, the initial state $z(t_n)$
and force $f_n$ give a strong solution
\begin{equation*}
 v_n\in C([0,2\theta_{R,M}];\mathbb V^1)
       \cap L^2(0,2\theta_{R,M};\mathbf H^2\times\dot H^4),
 \qquad v_n(0)=z(t_n).
\end{equation*}
The time $\theta_{R,M}$ is independent of $n$.
The functions $v_n$ and $z(t_n+\cdot)$ have the same initial state
and solve \eqref{abs} with force $f_n$ on their common interval.
By uniqueness of strong solutions on their common
interval of existence (\Cref{cor:unique-strong-common}), we have
\begin{equation*}
 v_n(s)=z(t_n+s),\qquad
 0\leq s<\min\{2\theta_{R,M},T_{\max}-t_n\}.
\end{equation*}
Choose $n$ sufficiently large so that
$t_n+2\theta_{R,M}>T_{\max}$.
The equality above then gives
$v_n(t-t_n)=z(t)$ for every $t\in[t_n,T_{\max})$.
Define
\begin{equation*}
 Z(t):=
 \begin{cases}
 z(t),&0\leq t\leq t_n,\\
 v_n(t-t_n),&t_n\leq t\leq t_n+2\theta_{R,M}.
 \end{cases}
\end{equation*}
Since $f_n(t-t_n)=f(t)$ for almost every
$t\in(t_n,t_n+2\theta_{R,M})$, $Z$ is a strong solution
of \eqref{abs} with initial state
$z_0$ and force $f$, satisfying
\begin{equation*}
 Z\in C([0,t_n+2\theta_{R,M}];\mathbb V^1)
       \cap L^2(0,t_n+2\theta_{R,M};\mathbf H^2\times\dot H^4),
\end{equation*}
and $Z(t)=z(t)$ for every $t\in[0,T_{\max})$.
Hence $t_n+2\theta_{R,M}\in\mathscr A$. This contradicts the
definition of $T_{\max}$, since $t_n+2\theta_{R,M}>T_{\max}$,
and proves \eqref{eq:maximal-strong-blowup}.
\end{proof}

\Cref{prop:maximal-strong} allows us to strengthen \Cref{thm:stability-strong}, more precisely,
under \eqref{eq:stability-smallness}, the perturbed solution
exists on $[0,T]$ and satisfies
\eqref{eq:strong-stability-estimate} with $S=T$.
\begin{cor}[Existence and stability under small perturbations]
	\label{thm:stability-strong1}
	Let $T>0$, and let $\widetilde z=(\widetilde u,\widetilde\phi)$ be a
	strong solution of \eqref{abs} on $[0,T]$, corresponding to
	the initial state $\widetilde z_0\in\mathbb V^1$ and the force
	$\widetilde f\in L^\infty(0,T;\mathbb V^0)$. Set
	\[
	\Lambda:=
	\|\widetilde z\|_{C([0,T];\mathbb V^1)}
	+\|\widetilde z\|_{L^2(0,T;\mathbf H^2\times\dot H^4)}
	+\|\widetilde f\|_{L^\infty(0,T;\mathbb V^0)}.
	\]
	Then there exist $\delta_*>0$ and $C>0$, depending only
	on $T$ and $\Lambda$, such that, for any $z_0\in\mathbb V^1$
	and $f\in L^\infty(0,T;\mathbb V^0)$ satisfying
	\[
	\|z_0-\widetilde z_0\|_{\mathbb V^1}^2
	+\|f-\widetilde f\|_{L^2(0,T;\mathbb V^0)}^2
	<\delta_*^2,
	\]
	system \eqref{abs} admits a unique strong solution
	$
	z\in C([0,T];\mathbb V^1)
	\cap L^2(0,T;\mathbf H^2\times\dot H^4).
	$
	Moreover,
	\begin{align}
		&\|z-\widetilde z\|_{C([0,T];\mathbb V^1)}^2
		+\|z-\widetilde z\|_{L^2(0,T;\mathbf H^2\times\dot H^4)}^2
		\notag\leq  C\Bigl(
		\|z_0-\widetilde z_0\|_{\mathbb V^1}^2
		+\|f-\widetilde f\|_{L^2(0,T;\mathbb V^0)}^2
		\Bigr).
		\label{eq:stability-full-interval}
	\end{align}
\end{cor}

%\begin{cor}[Existence and stability under small perturbations]\label{thm:stability-strong1}
%Let $\widetilde z$ be the reference solution in
%\Cref{thm:stability-strong} on $[0,T]$, and let $\delta_*>0$
%and $C>0$ be the constants given there.
%If $z_0\in\mathbb V^1$ and
%$f\in L^\infty(0,T;\mathbb V^0)$ satisfy
%\eqref{eq:stability-smallness}, then \eqref{abs} with initial
%state $z_0$ and force $f$ admits a unique strong solution
%$z$ on $[0,T]$ satisfying \eqref{eq:strong-stability-estimate}
%with $S=T$.
%\end{cor}
%{\color{red}this statement is difficult to follow?? we will cite this every places when we need it. so it will be difficult to read going back to previous right? can we write this with details?
%
%
%}
\begin{proof}
Let $\delta_*>0$ and $C>0$ be given by \Cref{thm:stability-strong},
with $d$ as in \eqref{eq:stability-smallness}.
Extend $f$ to an element of
$L^\infty_{\mathrm{loc}}(\mathbb R_+;\mathbb V^0)$, still denoted
by $f$. By \Cref{prop:maximal-strong}, there is a unique maximal
strong solution $z$ of \eqref{abs} with initial state $z_0$ and
force $f$, defined on $[0,T_{\max})$.
Suppose that $T_{\max}\le T$. For every $S\in(0,T_{\max})$,
\Cref{thm:stability-strong} gives
\begin{equation*}
 \sup_{t\in[0,S]}\|z(t)-\widetilde z(t)\|_{\mathbb V^1}^2
 \le C d^2.
\end{equation*}
Since $C$ is independent of $S$, we obtain
\begin{equation*}
 \sup_{0\le t<T_{\max}}\|z(t)\|_{\mathbb V^1}
 \le \|\widetilde z\|_{C([0,T];\mathbb V^1)}
      +C^{1/2}d<\infty,
\end{equation*}
contradicting \eqref{eq:maximal-strong-blowup}.
Thus $T_{\max}>T$. Applying \Cref{thm:stability-strong} with
$S=T$ gives the desired result.
\end{proof}
\Cref{prop:maximal-strong} also allows us to restart a maximal strong
solution from a state reached at a later time.
For $z_0\in\mathbb V^1$ and
$f\in L^\infty_{\mathrm{loc}}(\mathbb R_+;\mathbb V^0)$,
we denote the corresponding maximal strong solution by
$\mathcal S(z_0,f):=z$.
We write
\begin{equation*}
 \mathcal S_t(z_0,f):=z(t),\qquad t\in[0,T_{\max}).
\end{equation*}
Taking $z(s)$, $s\in[0,T_{\max})$, as the initial state at time zero
requires the translated force
\begin{equation*}
 f_s(r):=f(s+r),\qquad s,r\ge0.
\end{equation*}
The new solution is $z(s+\cdot)$, and its maximal time is exactly
$T_{\max}-s$, the time remaining after $s$.
For a time-dependent force, the
semigroup property takes the following form.
\begin{cor}\label{cor:strong-restart}
For every $s\in[0,T_{\max})$,
\begin{equation}\label{eq:strong-restart-maximal-time}
 T_{\max}\bigl(\mathcal S_s(z_0,f),f_s\bigr)=T_{\max}-s.
\end{equation}
Moreover, for $0\le t<T_{\max}-s$,
\begin{equation}\label{eq:strong-restart-identity}
 \mathcal S_{s+t}(z_0,f)
 =\mathcal S_t\bigl(\mathcal S_s(z_0,f),f_s\bigr).
\end{equation}
\end{cor}
\begin{proof}
Set $w:=\mathcal S(z(s),f_s)$ and
$T_s:=T_{\max}(z(s),f_s)$.
The function $t\mapsto z(s+t)$ is a strong solution of \eqref{abs}
with initial state $z(s)$ and force $f_s$ on $[0,T_{\max}-s)$.
By the definition of the maximal time and
\Cref{cor:unique-strong-common},
\begin{equation*}
 T_s\ge T_{\max}-s,\qquad
 w(t)=z(s+t),\quad 0\le t<T_{\max}-s.
\end{equation*}
If $T_s>T_{\max}-s$, then $T_{\max}<\infty$, and continuity of $w$
gives
\begin{equation*}
 \lim_{t\uparrow T_{\max}}\|z(t)\|_{\mathbb V^1}
 =\|w(T_{\max}-s)\|_{\mathbb V^1}<\infty,
\end{equation*}
contradicting \eqref{eq:maximal-strong-blowup}.
Thus $T_s=T_{\max}-s$, and the two identities follow.
\end{proof}
\subsection{Regular times for strong solutions}
\label{sec:regular-times}
Following the arguments used in
\cite[Section~4.1, proof of Lemma~4.2]{Tem95} and
\cite[proof of Theorem~2.3, p.~781]{PetcuWirosoetisno2005},
we use the spatial
regularity of a strong solution on $[0,T]$ and its continuity in
$\mathbb V^1$ to select states in
$\mathbf H^2\times\dot H^4$ arbitrarily close in $\mathbb V^1$
to the state at a given time.
%
%\begin{lemma}\label{lem:regular-times}
%Let $z=(u,\phi)$ be a maximal strong solution of \eqref{abs} on
%$[0,T_{\max})$, where $T_{\max}\in(0,\infty]$. For every
%$t_0\in[0,T_{\max})$ and every $\epsilon,\sigma>0$, there exists
%$t_*\in(t_0,T_{\max})$, with $0<t_*-t_0<\sigma$, such that
%\begin{equation*}
% z(t_*)\in\mathbf H^2\times\dot H^4,
% \qquad
% \|z(t_*)-z(t_0)\|_{\mathbb V^1}<\epsilon.
%\end{equation*}
%\end{lemma}
%
%\begin{proof}
%Applying \eqref{eq:energy-m1-local} on the local existence intervals
%of the maximal solution in \Cref{prop:maximal-strong},
%we have $z(t)\in\mathbf H^2\times\dot H^4$
%for almost every $t\in(0,T_{\max})$.
%By continuity, there exists $h\in(0,\sigma)$ such that
%$t_0+h<T_{\max}$ and
%$\|z(t)-z(t_0)\|_{\mathbb V^1}<\epsilon$ for
%$t\in[t_0,t_0+h]$.
%Since $z\in L^2(t_0,t_0+h;\mathbf H^2\times\dot H^4)$,
%we may choose $t_*\in(t_0,t_0+h)$ such that
%\begin{equation*}
% \|z(t_*)\|_{\mathbf H^2\times\dot H^4}^2
% \leq \frac1h\int_{t_0}^{t_0+h}
%       \|z(t)\|_{\mathbf H^2\times\dot H^4}^2\,dt=:A<\infty.
%\end{equation*}
%Indeed, otherwise $\|z(t)\|_{\mathbf H^2\times\dot H^4}^2>A$
%for almost every $t\in(t_0,t_0+h)$, and integration would give
%\begin{equation*}
% \int_{t_0}^{t_0+h}\|z(t)\|_{\mathbf H^2\times\dot H^4}^2\,dt
% >hA
% =\int_{t_0}^{t_0+h}\|z(t)\|_{\mathbf H^2\times\dot H^4}^2\,dt,
%\end{equation*}
%a contradiction.
%This $t_*$ has the required properties.
%\end{proof}
%}
\begin{lemma}\label{lem:regular-times}
Let $T>0$ and let
$z=(u,\phi)\in C([0,T];\mathbb V^1)
\cap L^2(0,T;\mathbf H^2\times\dot H^4)$
be a strong solution of \eqref{abs}.
For every $t_0\in(0,T]$ and every $\epsilon,\sigma>0$, there exists
$t_*\in(\max\{0,t_0-\sigma\},t_0)$ such that
\begin{equation*}
 z(t_*)\in\mathbf H^2\times\dot H^4,
 \qquad
 \|z(t_*)-z(t_0)\|_{\mathbb V^1}<\epsilon.
\end{equation*}
Similarly, for every $t_0\in[0,T)$ and every
$\epsilon,\sigma>0$, there exists
$t_*\in(t_0,\min\{T,t_0+\sigma\})$
satisfying the same properties.
\end{lemma}
\begin{proof}
Since $z\in L^2(0,T;\mathbf H^2\times\dot H^4)$,
we have $z(t)\in\mathbf H^2\times\dot H^4$ for almost every
$t\in(0,T)$. By continuity, there exists $h\in(0,\sigma)$ such that
\begin{equation*}
 \|z(t)-z(t_0)\|_{\mathbb V^1}<\epsilon,
 \qquad t\in[0,T],\quad |t-t_0|<h.
\end{equation*}
If $t_0>0$, the interval $(\max\{0,t_0-h\},t_0)$
has positive length, since $h>0$ and $t_0>0$.
As $z(t)\in\mathbf H^2\times\dot H^4$ for almost
every $t\in(0,T)$, this interval contains a time
$t_*$ such that
$z(t_*)\in\mathbf H^2\times\dot H^4$.
This choice satisfies the required properties
with $t_*<t_0$.
If $t_0<T$, the same argument on
$(t_0,\min\{T,t_0+h\})$ gives a choice
with $t_*>t_0$.
\end{proof}

\section{From the local theory to global strong solutions}\label{sec:global-prescribed-force}
In this section we prove \Cref{mnthm}, combining the local theory of
\Cref{sec:wp_local} with the accessibility property described in
\Cref{sec:main-ideas}. Let $z_0\in\mathbb V^1$ and
$f\in L^\infty_{\mathrm{loc}}(\mathbb R_+;\mathbb V^0)$ be given.

Recall that we look for the force in the form $g=f+\eta$, with a correction $\eta$
taking values in $\mathcal H$. We therefore consider the system
\begin{equation}\label{eq:chn-accessibility}
\begin{cases}
 \partial_tu+\nu A_0u+B_0(u,u)
 =-R_0(\Delta\phi,\phi)+f_u+\Pi(\eta_ue_3),\\
 \partial_t\phi+\Delta^2\phi+\Delta\phi+B_1(u,\phi)
 =\Delta(\phi^3)+f_\phi+\eta_\phi,\\
 u(0)=u_0,\qquad \phi(0)=\phi_0,
\end{cases}
\end{equation}
where $f=(f_u,f_\phi)$ is given and the scalar corrections
$(\eta_u,\eta_\phi)$ take values in $\mathcal E$, defined by
\eqref{eq:mode_sets_intro}. Writing
\begin{equation*}
 \eta(t):=\bigl(\Pi(\eta_u(t)e_3),\eta_\phi(t)\bigr)\in\mathcal H,
\end{equation*}
with $\mathcal H$ as in \eqref{eq:projected_modes_intro}, the system
takes the abstract form \eqref{abs} with force $f+\eta$.

When $\eta\equiv0$, let $z$ be the maximal strong solution of
\eqref{abs} with initial state $z_0$ and force $f$. By
\Cref{prop:maximal-strong}, {$z$} may have a finite
maximal time $T_{\max}$, in which case its $\mathbb V^1$-norm tends
to infinity. For the sake of exposition, we recall the idea behind the proof of \Cref{prop:maximal-strong}. For $R>0$ and $M\ge0$,
\Cref{thm:uniform-local-existence} gives a strong solution on
$[0,2\theta_{R,M}]$ for every initial state $\widehat z$ and force satisfying
\begin{equation*}
 \|\widehat z\|_{\mathbb V^1}\le R,\qquad
 \|f\|_{L^\infty(0,1;\mathbb V^0)}\le M,
\end{equation*}
and this time does not depend on the particular state.
If $T_{\max}<\infty$, we can choose $M$ so that
$\|f\|_{L^\infty(0,T_{\max}+1;\mathbb V^0)}\le M$.
Now suppose that, for some $R>0$, there are times $t_n\uparrow T_{\max}$ such that
\begin{equation*}
 \|z(t_n)\|_{\mathbb V^1}\le R,\qquad n\ge1.
\end{equation*}
We can then apply \Cref{thm:uniform-local-existence} with initial state $z(t_n)$
and translated force $f_{t_n}$, obtaining a strong solution on
$[0,2\theta_{R,M}]$. By uniqueness, this solution agrees with
$z(t_n+\cdot)$ on their common interval of existence and therefore continues
$z$ up to $t_n+2\theta_{R,M}$. Since $\theta_{R,M}$ is independent
of $n$, this extends $z$ beyond $T_{\max}$ for sufficiently large $n$.

With $\eta\equiv0$, however, the usual local existence theory does not provide any estimate
ensuring that the states $z(t_n)$ (used as new initial data) remain bounded by a constant $R$ independent of $n$. 
Since $\theta_{R,M}$ may decrease as $R$ grows,
applying the local theorem naively along a solution whose norm grows
may only produce summable existence times.

Our main idea in the proof of \Cref{mnthm} is to use
$\eta$ to obtain the bound on the new initial states that is
missing when $\eta\equiv0$. Suppose that, from any state, we can choose $\eta$ to bring the norm
below one in arbitrarily short time. After each such return we set
$\eta=0$ and apply \Cref{thm:uniform-local-existence} with $R=1$. Since $f$ is locally bounded, for each
bounded time interval there is a bound $M$ for the force,
so the solution evolves under $f$ for at least $\theta_{1,M}>0$
between consecutive corrections, and the returns cannot accumulate
in finite time. Moreover, since each return can be made in
arbitrarily short time, the correction intervals can be chosen with
total length below any given $\rho>0$.

In this regard, the property we need is the following adaptation of the accessibility sets of
\cite[Definition~3.3]{GHHM18} to a prescribed time-dependent force.
\begin{defi}[Small-time accessibility]\label{def:small-time-accessibility}
Let $\mathscr B:=\{z\in\mathbb V^1:\|z\|_{\mathbb V^1}<1\}$.
We say that $\mathscr B$ is accessible in small time with
corrections in $\mathcal H$ if, for every $z_0\in\mathbb V^1$,
every $f\in L^\infty_{\mathrm{loc}}(\mathbb R_+;\mathbb V^0)$,
and every $\sigma>0$, there exist $\tau\in(0,\sigma)$ and
$\eta\in L^\infty(0,\tau;\mathcal H)$ such that
system \eqref{eq:chn-accessibility} has a strong solution $z$
on $[0,\tau]$ satisfying
\begin{equation}\label{eq:global-small-time-access}
 \|z(\tau)\|_{\mathbb V^1}<1.
\end{equation}
\end{defi}
Note that the
existence of a strong solution throughout the correction interval is
part of the property, since the corrections may be large.

The following result is key in our construction.
\begin{theo}\label{thm:global-small-time-control}
For $\mathcal H$ defined by \eqref{eq:projected_modes_intro},
the unit ball $\mathscr B$ is accessible in small time with
corrections in $\mathcal H$.
\end{theo}

The difficulty in \Cref{thm:global-small-time-control} comes from the
restrictions on $\mathcal H$: the corrections involve only four
frequencies, the fluid correction acts only on the third component (before the Leray projection), and the phase profiles depend only on
$x_3$. The proof, postponed to
\Cref{sec:accessibility}, uses the asymptotic limits of
\Cref{sec:asymp_limits} and the saturation properties of
\Cref{sec:saturation} to generate the additional directions through
nonlinear interactions.

For the moment, we use \Cref{thm:global-small-time-control} to prove
\Cref{mnthm}. The main idea is to construct
consecutive intervals of the form
\begin{equation*}
 0=a_0<b_0<a_1<b_1<a_2<\cdots.
\end{equation*}
On $[a_n,b_n)$ we use a correction to reach $\mathscr B$.
We then let the solution evolve under $f$ on $[b_n,a_{n+1})$
for a fixed fraction of the local existence time. Finally,
\Cref{cor:unique-strong-common,prop:maximal-strong} show that these evolutions agree with the maximal strong solution for the resulting force $g$. Since $a_n\to\infty$, this solution is global.
\begin{proof}[Proof of \Cref{mnthm}]
Let $z_0$ and $f$ satisfy the hypotheses of \Cref{mnthm}, and
fix $\rho>0$. Set
\begin{equation}\label{eq:global-time-budget}
 \sigma_n:=\rho\,2^{-n-2},\, n\geq 0 \, \,
 \text{ and thus } \sum_{n=0}^{\infty}\sigma_n=\frac{\rho}{2}.
\end{equation}
For $M\ge0$, define $\theta_M:=\theta_{1,M}$, with
$\theta_{R,M}$ as in \Cref{thm:uniform-local-existence}.

\textit{--Step 1. Construction of the intervals and the solutions.}
We construct the times $a_n$ and states $z_n$ recursively,
starting from $a_0:=0$ and the given initial state $z_0$.

{Given $a_n\geq0$ and $z_n\in\mathbb V^1$,} the translated force $f_{a_n}$ satisfies
the same hypotheses as $f$. By
\Cref{thm:global-small-time-control}, there exist
$\tau_n\in(0,\sigma_n)$ and
{$\bar\eta^n\in L^\infty(0,\tau_n;\mathcal H)$ such that}
the strong solution $\bar z^{(n)}$ on $[0,\tau_n]$, with
{$\bar z^{(n)}(0)=z_n$} and force $f_{a_n}+\bar\eta^n$, satisfies
$\|\bar z^{(n)}(\tau_n)\|_{\mathbb V^1}<1$.
{Set $b_n:=a_n+\tau_n$.}

Let $\delta_n>0$ and $C_n\ge1$ be given by
\Cref{thm:stability-strong1}, applied to the reference solution
$\bar z^{(n)}$ with force $f_{a_n}+\bar\eta^n$
{(increasing $C_n$ if necessary)}.
{By density of $C_c^\infty((0,\tau_n);\mathcal H)$ in
$L^2(0,\tau_n;\mathcal H)$, we can choose
$\eta^n\in C_c^\infty((0,\tau_n);\mathcal H)$ such that}
\begin{equation}\label{eq:global-smooth-control-choice}
 d_n:=\|\eta^n-\bar\eta^n\|_{L^2(0,\tau_n;\mathbb V^0)}
 <\min\left\{\delta_n,
 \frac{1-\|\bar z^{(n)}(\tau_n)\|_{\mathbb V^1}}
      {2C_n^{1/2}}\right\}.
\end{equation}
{By \eqref{eq:projected_modes_intro}, we write
$\eta^n=(\Pi(\eta_u^ne_3),\eta_\phi^n)$, with
$(\eta_u^n,\eta_\phi^n)\in C_c^\infty((0,\tau_n);\mathcal E)$.}
By \Cref{thm:stability-strong1}, the force
$f_{a_n}+\eta^n$ gives a unique strong solution $z^{(n)}$
on $[0,\tau_n]$, {with $z^{(n)}(0)=z_n$,} satisfying
\begin{equation}\label{eq:global-smooth-control-stability}
 \|z^{(n)}-\bar z^{(n)}\|_{C([0,\tau_n];\mathbb V^1)}^2
 \le C_n d_n^2.
\end{equation}
{Consequently, setting $\widehat z_n:=z^{(n)}(\tau_n)$,
we obtain
\begin{align}\label{eq:global-return-state}
 \|\widehat z_n\|_{\mathbb V^1}
 &\le \|z^{(n)}(\tau_n)-\bar z^{(n)}(\tau_n)\|_{\mathbb V^1}
       +\|\bar z^{(n)}(\tau_n)\|_{\mathbb V^1}\notag\\
 &\le C_n^{1/2}d_n+\|\bar z^{(n)}(\tau_n)\|_{\mathbb V^1}\notag\\
 &<\frac{1-\|\bar z^{(n)}(\tau_n)\|_{\mathbb V^1}}{2}
       +\|\bar z^{(n)}(\tau_n)\|_{\mathbb V^1}
 =\frac{1+\|\bar z^{(n)}(\tau_n)\|_{\mathbb V^1}}{2}<1.
\end{align}
Thus the smooth correction $\eta^n$ still brings the state into
$\mathscr B$.}
Set
\begin{equation}\label{eq:global-force-bound}
 M_n:=\|f\|_{L^\infty(0,b_n+1;\mathbb V^0)}.
\end{equation}
Since
\begin{equation*}
 \|\widehat z_n\|_{\mathbb V^1}<1,\qquad
 \|f_{b_n}\|_{L^\infty(0,1;\mathbb V^0)}\le M_n,
\end{equation*}
\Cref{thm:uniform-local-existence} gives a strong solution
$w^{(n)}$ {with $w^{(n)}(0)=\widehat z_n$} and force $f_{b_n}$
such that
\begin{equation}\label{eq:global-uniform-local-existence}
 w^{(n)}\in C([0,2\theta_{M_n}];\mathbb V^1)
       \cap L^2(0,2\theta_{M_n};\mathbf H^2\times\dot H^4).
\end{equation}
We define
\begin{equation}\label{eq:global-next-stage}
 a_{n+1}:=b_n+\theta_{M_n},
 \qquad z_{n+1}:=w^{(n)}(\theta_{M_n}).
\end{equation}
{Thus each step has the form
\begin{equation*}
 (a_n,z_n)
 \xrightarrow[\tau_n]{\,f_{a_n}+\eta^n\,}
 (b_n,\widehat z_n)
 \xrightarrow[\theta_{M_n}]{\,f_{b_n}\,}
 (a_{n+1},z_{n+1}).
\end{equation*}
In Step~4 we show that the states $z_n$ defined by this recursion
coincide with the values of the maximal solution at the times $a_n$.}
In particular, $z_{n+1}\in\mathbb V^1$, so the construction
can be repeated {to obtain} a sequence of times such that
\begin{equation}\label{eq:global-ordered-times}
 0=a_0<b_0<a_1<b_1<a_2<b_2<\cdots.
\end{equation}

\smallskip

\textit{--Step 2. The times tend to infinity.}
Suppose that $a_n\uparrow A<\infty$. {Since $(a_n)$ is
strictly increasing and converges to $A$,
\eqref{eq:global-ordered-times} gives $b_n<a_{n+1}<A$ for every
$n\ge0$. Hence $(0,b_n+1)\subset(0,A+1)$, and}
\begin{equation*}
 {M_n=\|f\|_{L^\infty(0,b_n+1;\mathbb V^0)}}\le
 \|f\|_{L^\infty(0,A+1;\mathbb V^0)}=:M_A<\infty.
\end{equation*}
{Since the function $M\mapsto\theta_M$ is nonincreasing,}
\eqref{eq:global-next-stage} {gives}
\begin{equation*}
 a_{n+1}-a_n
 =\tau_n+\theta_{M_n}
 \ge\theta_{M_A}>0.
\end{equation*}
Summing these inequalities and using that $a_0=0$, we obtain
 $a_n=\sum_{j=0}^{n-1}(a_{j+1}-a_j)\ge n\theta_{M_A},
 \, n\ge1.$
This contradicts $a_n<A$ for all $n\ge0$. Therefore
\begin{equation}\label{eq:global-time-divergence}
 a_n\longrightarrow\infty.
\end{equation}
{Consequently, for every $L>0$, only finitely many
of the intervals $[a_n,b_n)$ intersect $[0,L]$, since $a_n>L$
for all sufficiently large $n$.}

\smallskip

\textit{--Step 3. Definition of the force.}
Extend $\eta_u^n$, $\eta_\phi^n$, and $\eta^n$ by zero to
$\mathbb R$. Set
\begin{equation}\label{eq:global-correction-set}
 I_n:=[a_n,b_n),\qquad I:=\bigcup_{n=0}^{\infty}I_n,
\end{equation}
and define, for $t\ge0$,
\begin{equation*}
 \eta_u(t):=\sum_{n=0}^{\infty}\eta_u^n(t-a_n),
 \qquad
 \eta_\phi(t):=\sum_{n=0}^{\infty}\eta_\phi^n(t-a_n).
\end{equation*}
Then set
\begin{equation}\label{eq:global-concatenated-force}
 \eta(t):=(\Pi(\eta_u(t)e_3),\eta_\phi(t))
       =\sum_{n=0}^{\infty}\eta^n(t-a_n),
 \qquad g(t):=f(t)+\eta(t).
\end{equation}
Fix $T>0$. By \eqref{eq:global-time-divergence}, there exists
$N$ such that $a_N>T$. For $n\ge N$ and $0\le t\le T$,
we have $t-a_n<0$, so
$\eta_u^n(t-a_n)=\eta_\phi^n(t-a_n)=0$. Hence
\begin{equation*}
 (\eta_u(t),\eta_\phi(t))
 =\sum_{n=0}^{N-1}
   (\eta_u^n(t-a_n),\eta_\phi^n(t-a_n)),
 \qquad 0\le t\le T.
\end{equation*}
Since $(\eta_u^n,\eta_\phi^n)\in
C_c^\infty((0,\tau_n);\mathcal E)$, its zero extension belongs
to $C^\infty(\mathbb R;\mathcal E)$. The sum above is finite,
so $(\eta_u,\eta_\phi)\in C^\infty([0,T];\mathcal E)$.
As $T>0$ was arbitrary, we conclude that
$(\eta_u,\eta_\phi)\in C^\infty([0,\infty);\mathcal E)$.
Since $\mathcal E$ is spanned by finitely many trigonometric
functions, we obtain $\eta_u,\eta_\phi
 \in C^\infty([0,\infty)\times\mathbb T^3).$
 Moreover,
$\eta\in C^\infty([0,\infty);\mathcal H)$, 
 $\eta_u(t)=\eta_\phi(t)=0$
  for every $t\ge0$ outside $I$.

For this choice of $N$, we also have
\begin{equation}\label{eq:global-local-force-bound}
 \|\eta\|_{L^\infty(0,T;\mathcal H)}
 \le\max_{0\le n<N}
       \|\eta^n\|_{L^\infty(0,\tau_n;\mathcal H)}<\infty.
\end{equation}
Thus $\eta\in L^\infty_{\mathrm{loc}}(\mathbb R_+;\mathcal H)$
and $g\in[f]_{\mathcal H}$. Moreover, the intervals $I_n$
are disjoint, so
\begin{equation}\label{eq:global-measure-bound}
 |I|=\sum_{n=0}^{\infty}\tau_n
 \le\sum_{n=0}^{\infty}\sigma_n
 =\frac{\rho}{2}<\rho,
\end{equation}
and $g=f$ almost everywhere outside $I$.

\textit{--Step 4. Global existence.}
Let $z$ be the maximal strong solution with initial state
$z_0$ and force $g$, and write
$T_{\max}:=T_{\max}(z_0,g)$. Such a solution exists by
\Cref{prop:maximal-strong}. We show by induction that
\begin{equation}\label{eq:global-stage-values}
 T_{\max}>a_n,\qquad z(a_n)=z_n,\qquad\text{for every }n\ge0.
\end{equation}
Since $a_0=0$ and $z(0)=z_0$, the assertion holds for $n=0$.
Assume that it holds for some $n\ge0$. We prove it for $n+1$.
For almost every $r\in(0,\tau_n)$,
\begin{equation*}
 g(a_n+r)=f_{a_n}(r)+\eta^n(r).
\end{equation*}
By {\eqref{eq:strong-restart-identity}} and
\Cref{cor:unique-strong-common},
\begin{equation}\label{eq:global-controlled-restart}
 z(a_n+r)=z^{(n)}(r),
 \qquad
 0\le r<\min\{\tau_n,T_{\max}-a_n\}.
\end{equation}
If $T_{\max}\le b_n$, then $T_{\max}-a_n\in(0,\tau_n]$.
By \eqref{eq:global-controlled-restart},
$z(t)=z^{(n)}(t-a_n)$ for $a_n\le t<T_{\max}$.
Since $z^{(n)}\in C([0,\tau_n];\mathbb V^1)$, we obtain
\begin{align*}
 \lim_{t\uparrow T_{\max}}\|z(t)\|_{\mathbb V^1}
 &=\|z^{(n)}(T_{\max}-a_n)\|_{\mathbb V^1}\le\|z^{(n)}\|_{C([0,\tau_n];\mathbb V^1)}<\infty,
\end{align*}
contradicting \eqref{eq:maximal-strong-blowup}. Thus
$T_{\max}>b_n$, and continuity in
\eqref{eq:global-controlled-restart} gives
$z(b_n)=\widehat z_n$.

On the next interval,
$g(b_n+r)=f_{b_n}(r)$ for almost every
$r\in(0,\theta_{M_n})$.
By {\eqref{eq:strong-restart-identity}} and
\Cref{cor:unique-strong-common},
\begin{equation}\label{eq:global-free-restart}
 z(b_n+r)=w^{(n)}(r),
 \qquad
 0\le r<\min\{\theta_{M_n},T_{\max}-b_n\}.
\end{equation}
If $T_{\max}\le a_{n+1}$, then
$T_{\max}-b_n\in(0,\theta_{M_n}]$.
By \eqref{eq:global-free-restart},
$z(t)=w^{(n)}(t-b_n)$ for $b_n\le t<T_{\max}$.
Since $w^{(n)}\in C([0,\theta_{M_n}];\mathbb V^1)$, we obtain
\begin{align*}
 \lim_{t\uparrow T_{\max}}\|z(t)\|_{\mathbb V^1}
 &=\|w^{(n)}(T_{\max}-b_n)\|_{\mathbb V^1}\le\|w^{(n)}\|_{C([0,\theta_{M_n}];\mathbb V^1)}<\infty,
\end{align*}
again contradicting \eqref{eq:maximal-strong-blowup}.
Consequently, $T_{\max}>a_{n+1}$ and
$z(a_{n+1})=w^{(n)}(\theta_{M_n})=z_{n+1}$.
This proves \eqref{eq:global-stage-values}.

Since $a_n\to\infty$, we have $T_{\max}=\infty$.
By {\Cref{prop:maximal-strong}}, the solution has the regularity
\eqref{eq:global_solution_intro} and is unique for the fixed
force $g$. The properties proved in Step 3 give
\eqref{eq:force_agreement_intro} and complete the proof.
\end{proof}
\section{Small-time asymptotic limits}\label{sec:asymp_limits}
In this section, we establish two small-time asymptotic limits for strong solutions
of the CHNS system. These limits show how large profiles can
produce finite displacements through the nonlinear terms in arbitrarily
short time. Together with the saturation properties of
\Cref{sec:saturation}, they will allow us to generate directions beyond
those directly available in $\mathcal H$ and prove small-time
accessibility in \Cref{sec:accessibility}.

\subsection{A shifted formulation}
We begin by introducing a version of CHNS with time-independent shifts and additional source terms, and define its solution map. For time-independent functions $\xi=(\xi_u,\xi_\phi)$ and
a source term $\eta=(\eta_u,\eta_\phi)$\footnote{{For simplicity of notation,
$\eta_u$ here and throughout this section denotes a vector-valued
source, rather than the scalar correction in
\eqref{eq:chn-accessibility}, and we state the results for these more general sources. Clearly, the projected scalar corrections considered there are
included as a particular case.}}, we consider the following equation
\begin{align*}
\begin{cases}
\partial_t u+\nu A_0(u+\xi_u)+B_0(u+\xi_u,u+\xi_u)=-R_0(\Delta(\phi+\xi_\phi),\phi+\xi_\phi)+f_u+{\eta_u} \\
\partial_t \phi + \Delta^2(\phi+\xi_\phi)+\Delta(\phi+\xi_\phi)+B_1(u+\xi_u,\phi+\xi_\phi)=\Delta((\phi+\xi_\phi)^3)+ f_\phi+{\eta_\phi}, \\
u(0)=u_0, \quad \phi(0)=\phi_0
\end{cases}
\end{align*}
or, equivalently, using notation \eqref{def_z}--\eqref{def_ops},
\begin{equation}\label{ext}
\partial_t z + L(z+\xi) + B(z+\xi) + Q(z+\xi) = f+\eta, \qquad
z(0)=z_0.
\end{equation}
We extend the notation of \Cref{cor:strong-restart} by writing
$\mathcal S(z_0,\xi,f+\eta):=z$ for the maximal strong solution of
\eqref{ext}. We denote its maximal existence time by $\mathcal T$ and write
\begin{equation*}
 \mathcal T:=\mathcal T(z_0,\xi,f+\eta)\in(0,\infty],\qquad
 \mathcal S_t(z_0,\xi,f+\eta):=z(t),\quad 0\le t<\mathcal T.
\end{equation*}
For $\xi=0$, we recover the notation of \Cref{sec:maximal-strong-solutions}, namely
\begin{equation*}
\begin{aligned}
 \mathcal S_t(z_0,0,f+\eta)&=\mathcal S_t(z_0,f+\eta),\\
 \mathcal T(z_0,0,f+\eta)&=T_{\max}(z_0,f+\eta).
\end{aligned}
\end{equation*}
The following result relates systems
\eqref{ext} and \eqref{abs} and some quantities and properties of interest.
\begin{proposition}\label{wp1}
Let $z_0\in\mathbb V^1$, $\xi_u\in\mathbf H^2$, $\xi_\phi\in\dot H^4$, and
$f, \eta\in L^\infty_{\mathrm{loc}}(\mathbb R_+;\mathbb V^0)$.
Then there exists a maximal time
$\mathcal T=\mathcal T(z_0,\xi,f+\eta)\in(0,\infty]$ and a unique
strong solution
$z\in C([0,\mathcal T);\mathbb V^1)
\cap L^2_{\mathrm{loc}}([0,\mathcal T);{\mathbf H^2\times\dot H^4})$
of \eqref{ext} on $[0,\mathcal T)$.
The maximal times satisfy
\begin{equation}\label{eq:shift_maximal_time}
 \mathcal T(z_0,\xi,f+\eta)
 =T_{\max}(z_0+\xi,f+\eta).
\end{equation}
Moreover, the associated solution map satisfies
\begin{equation}\label{eq:shift_iden}
\mathcal S_t(z_0,\xi,f+\eta)=\mathcal S_t(z_0+\xi,0,f+\eta)-\xi,
\qquad t\in[0,\mathcal T).
\end{equation}
For $s\in[0,\mathcal T)$, set
%\begin{equation*}
$\eta_s(r):=\eta(s+r)$ for $r\ge0$.
%\end{equation*}
Then
\begin{equation}\label{eq:restart_maximal_time}
 \mathcal T\bigl(\mathcal S_s(z_0,\xi,f+\eta),\xi,f_s+\eta_s\bigr)
 =\mathcal T-s.
\end{equation}
Moreover, for $0\le t<\mathcal T-s$,
\begin{equation}\label{eq:restart_identity}
 \mathcal S_{s+t}(z_0,\xi,f+\eta)
 =\mathcal S_t\bigl(\mathcal S_s(z_0,\xi,f+\eta),\xi,f_s+\eta_s\bigr).
\end{equation}
\end{proposition}
\begin{proof}
We introduce the change of variable $z':=z+\xi$.
Since $\xi$ is independent of time, \eqref{ext} becomes
\eqref{abs} with initial state $z_0+\xi$ and force $f+\eta$.
By \Cref{prop:maximal-strong}, this problem has a unique maximal
strong solution $z'$ on $[0,T_{\max}(z_0+\xi,f+\eta))$.
Set $z:=z'-\xi$. Note that any (strong) extension of $z$ would, after addition of $\xi$,
give a strong extension of $z'$, contradicting maximality.
Thus \eqref{eq:shift_maximal_time} holds and $z\in C([0,\mathcal T);\mathbb V^1)
\cap L^2_{\mathrm{loc}}([0,\mathcal T);{\mathbf H^2\times\dot H^4})$. Uniqueness follows from \Cref{cor:unique-strong-common}
after the same change of variables, and
\eqref{eq:shift_iden} follows from $z(t)=z'(t)-\xi$.

Fix $s\in[0,\mathcal T)$ and define
\begin{equation*}
 \mathcal T_s:=\mathcal T(z(s),\xi,f_s+\eta_s),
 \qquad
 w(r):=\mathcal S_r(z(s),\xi,f_s+\eta_s),
 \quad 0\le r<\mathcal T_s.
 \end{equation*}
The functions $r\mapsto z(s+r)$ and $w$ solve \eqref{ext}
with initial state $z(s)$, shift $\xi$, and force $f_s+\eta_s$.
By uniqueness,
%\begin{equation*}
 $w(r)=z(s+r),
 \, 0\le r<\min\{\mathcal T_s,\mathcal T-s\}.$
%\end{equation*}
If $\mathcal T_s<\mathcal T-s$, continuity of $z$ gives
\begin{equation*}
 \lim_{r\uparrow\mathcal T_s}\|w(r)+\xi\|_{\mathbb V^1}
 =\|z(s+\mathcal T_s)+\xi\|_{\mathbb V^1}<\infty.
\end{equation*}
By \eqref{eq:shift_maximal_time}, $w+\xi$ is a maximal strong
solution of \eqref{abs} with maximal time $\mathcal T_s$.
The above limit contradicts \eqref{eq:maximal-strong-blowup}.
Hence $\mathcal T_s\ge\mathcal T-s$.

If $\mathcal T<\infty$ and $\mathcal T_s>\mathcal T-s$,
continuity of $w$ gives
\begin{equation*}
 \lim_{t\uparrow\mathcal T}\|z(t)+\xi\|_{\mathbb V^1}
 =\|w(\mathcal T-s)+\xi\|_{\mathbb V^1}<\infty,
\end{equation*}
contradicting \eqref{eq:maximal-strong-blowup} for $z'=z+\xi$.
Therefore $\mathcal T_s=\mathcal T-s$, which proves
\eqref{eq:restart_maximal_time}.
The identity \eqref{eq:restart_identity} follows from
$w(r)=z(s+r)$ on this interval.
\end{proof}
\begin{proposition}
	\label{cont}
	Let $z_0\in\mathbb V^1$, $\xi_u\in\mathbf H^2$, $\xi_\phi\in\dot H^4$, and
	$f,\eta\in L^\infty_{\mathrm{loc}}(\mathbb R_+;\mathbb V^0) $.
	Suppose that $z=\mathcal S(\,z_0,\xi,f+\eta)$ $\in C([0,\mathcal T);\mathbb V^1)
	\cap L^2_{\mathrm{loc}}([0,\mathcal T);{\mathbf H^2\times\dot H^4})$ is the
	corresponding unique maximal strong solution of \eqref{ext}, and let
	$0<T<\mathcal T(z_0,\xi,f+\eta) $ . Set
	\[
	\Lambda:=
	\|z\|_{C([0,T];\mathbb V^1)}
	+\|z\|_{L^2(0,T;\mathbf H^2\times\dot H^4)}
	+\|f+\eta\|_{L^\infty(0,T;\mathbb V^0)}
	+\|\xi_u\|_{H^2}+\|\xi_\phi\|_{\dot H^4}.
	\]
	Then there exist $\delta>0$ and $C>0$, depending only on
	$T$ and $\Lambda$, such that, for every
	$\widehat z_0\in\mathbb V^1$ and
	$\widehat f\in L^\infty_{\mathrm{loc}}(\mathbb R_+;\mathbb V^0)$
	satisfying
	\begin{equation}\label{eq:shifted-perturbation-smallness}
		\|z_0-\widehat z_0\|_{\mathbb V^1}
		+\|f-\widehat f\|_{L^2(0,T;\mathbb V^0)}
		<\delta,
	\end{equation}
	the corresponding maximal strong solution
	$\widehat z=\mathcal S(\widehat z_0,\xi,\widehat f+\eta)$
	satisfies $
	\mathcal T(\widehat z_0,\xi,\widehat f+\eta)>T
	$, that means $\widehat z \in C([0, T];\mathbb V^1)
	\cap L^2(0, T;{\mathbf H^2\times\dot H^4}).$ Moreover,
	\begin{equation}\label{eq:cont_stability}
		\begin{aligned}
			&\|z-\widehat z\|_{C([0,T];\mathbb V^1)}
			+\|z-\widehat z\|_{L^2(0,T;\mathbf H^2\times\dot H^4)}\leq C\Bigl(
			\|z_0-\widehat z_0\|_{\mathbb V^1}
			+\|f-\widehat f\|_{L^2(0,T;\mathbb V^0)}
			\Bigr).
		\end{aligned}
	\end{equation}
\end{proposition}
\begin{proof}
	Set $z':=z+\xi$. By \Cref{wp1}, $z'$ is a strong
	solution of \eqref{abs} on $[0,T]$, corresponding to
	the initial state $z_0+\xi$ and the force $f+\eta$.
	Moreover, since $\xi_u\in\mathbf H^2$ and $\xi_\phi\in\dot H^4$, we have
	\[
	\|z'\|_{C([0,T];\mathbb V^1)}
	+\|z'\|_{L^2(0,T;\mathbf H^2\times\dot H^4)}
	+\|f+\eta\|_{L^\infty(0,T;\mathbb V^0)}
	\le C(1+T^{1/2})\Lambda.
	\]
	By \Cref{thm:stability-strong1}, there exist
	$\delta,C>0$, depending only on $T$ and $\Lambda$,
	such that, whenever \eqref{eq:shifted-perturbation-smallness}
	holds, the equation \eqref{abs} with initial state
	$\widehat z_0+\xi$ and force $\widehat f+\eta$
	admits a unique strong solution $\widehat z'$
	on $[0,T]$, satisfying
	\begin{equation}\label{eq:cont_stability_unshifted}
	\begin{aligned}
		&\|z'-\widehat z'\|_{C([0,T];\mathbb V^1)}
		+\|z'-\widehat z'\|_{L^2(0,T;\mathbf H^2\times\dot H^4)}\le C\Bigl(
		\|z_0-\widehat z_0\|_{\mathbb V^1}
		+\|f-\widehat f\|_{L^2(0,T;\mathbb V^0)}
		\Bigr).
	\end{aligned}
	\end{equation}
	By {\eqref{eq:cont_stability_unshifted}}, we have
	\begin{align*}
		\|\widehat z'\|_{C([0,T];\mathbb V^1)}
		&\le \|z'\|_{C([0,T];\mathbb V^1)}
		+\|z'-\widehat z'\|_{C([0,T];\mathbb V^1)}\\
		&\le \|z'\|_{C([0,T];\mathbb V^1)}
		+C\Bigl(\|z_0-\widehat z_0\|_{\mathbb V^1}
		+\|f-\widehat f\|_{L^2(0,T;\mathbb V^0)}\Bigr)\\
		&\le \|z'\|_{C([0,T];\mathbb V^1)}+C\delta<\infty,
	\end{align*}
	which, by \Cref{cor:unique-strong-common} and
	\eqref{eq:maximal-strong-blowup}, implies $
	T_{\max}(\widehat z_0+\xi,\widehat f+\eta)>T.
	$ Hence, by \eqref{eq:shift_maximal_time},
	$\mathcal T(\widehat z_0,\xi,\widehat f+\eta)>T$. Therefore $\widehat z \in C([0, T];\mathbb V^1)
	\cap L^2(0, T;{\mathbf H^2\times\dot H^4}).$
	Finally, $\widehat z=\widehat z'-\xi$ and
	$z=z'-\xi$, so that
	$z-\widehat z=z'-\widehat z'$.
	This proves \eqref{eq:cont_stability}.
\end{proof}
\subsection{First asymptotic limit}
For any fixed $m\in\mathbb N$, let $\pi_{u}$ and $\pi_{\phi}$
be the projections onto the velocity and phase components of
$\mathbb V^m$, respectively, so that
$\pi_{u}(\xi)=\xi_u$ and $\pi_{\phi}(\xi)=\xi_\phi$ for
$\xi=(\xi_u,\xi_\phi)\in\mathbb V^m$.
%
%{\color{magenta}We set $M_f:=\|f\|_{L^\infty(0,2;\mathbb V^0)}$.}
%{\color{red}
%	At the beginning of the section, $\eta=(\Pi(\eta_u e_3),\eta_\phi)$,
%	whereas the theorem assumes $\eta\in\mathbb V^5$.
%	These assumptions are not equivalent.  Should we  allow
%	general $\eta\in\mathbb V^5$ in the shifted formulation or we just take $(\Pi(\eta_u e_3),\eta_\phi)$ type of $\eta$? Let me know if my question is correct? or I making any mistake??
%}
\begin{theo}\label{thm:asymp-fluid}
Let $z_0=(u_0,\phi_0)\in \mathbf H^2\times\dot H^4$, $ f=(f_u, f_\phi) \in L^\infty_{\mathrm{loc}}(\mathbb R_+;\mathbb V^0)$,
$\eta=(\eta_u,\eta_\phi)\in\mathbf H^2\times\dot H^4$, and
$\zeta_u\in\mathbf H^3$. Set $\zeta=(\zeta_u,0)$. There exists $\delta_0>0$
such that for all $\delta\in(0,\delta_0)$, the solution
\begin{equation*}
t\mapsto \mathcal S_t(z_0,\delta^{-1/2}\zeta,f+\delta^{-1}\eta)
\end{equation*}
is well defined on the time interval $[0,\delta]$. Moreover, the following
convergence holds in $\mathbb V^1$
\begin{equation}\label{asymp_phi}
\mathcal S_\delta(z_0,\delta^{-1/2}\zeta,f+\delta^{-1}\eta)
\longrightarrow z_0+ \eta - (B_0({\zeta_u},{\zeta_u}),0)
\qquad\text{as } \delta\to0^+.
\end{equation}
\end{theo}
\begin{proof}
For clarity, we divide the proof into five steps. In what follows, $\epsilon\in(0,1)$ denotes a small parameter that will be fixed later in the computations and comes from employing Young's inequality.

\smallskip
\textit{-- Step 1. A change of variable.}
Let us take any $\delta\in(0,1)$. Since
$\delta^{-1/2}\zeta\in\mathbf H^2\times\dot H^4$,
\Cref{wp1} gives a maximal strong solution with
\begin{equation*}
 \mathcal T_\delta
 :=\mathcal T\bigl(z_0,\delta^{-1/2}\zeta,f+\delta^{-1}\eta\bigr)
 =T_{\max}\bigl(z_0+\delta^{-1/2}\zeta,f+\delta^{-1}\eta\bigr).
\end{equation*}
We define
\begin{equation*}
w(t)=z(\delta t)-q(t),
\end{equation*}
where
\begin{equation*}
z(t)=\mathcal S_t\bigl(z_0,\delta^{-1/2}\zeta,{f+\delta^{-1}\eta}\bigr),
\qquad
q(t)=\bigl({q_{u}}(t),{q_{\phi}}(t)\bigr),
\end{equation*}
with
\begin{equation*}
{q_{u}}(t)=u_0+t\bigl({\eta_u}-B_0({\zeta_{u}},{\zeta_{u}})\bigr),
\qquad
{q_{\phi}}(t)=\phi_0+t\,{\eta_\phi}.
\end{equation*}
We also set
\begin{equation*}
U:=v+{q_{u}}+\delta^{-1/2}{\zeta_{u}},
\qquad
\Phi:=\varphi+{q_{\phi}},
\end{equation*}
where $(v,\varphi)=({\pi_{u}} w,{\pi_{\phi}} w)$. In this case, $(v,\varphi)$ satisfies
the following system in $(0,\delta^{-1}{\mathcal T_\delta})\times\mathbb T^3$:
\begin{equation}\label{chnp_asymp}
\begin{cases}
\partial_t v+\delta\nu A_0U+\delta B_0(U,U)
+\delta R_0(\Delta\Phi,\Phi)
=
B_0({\zeta_{u}},{\zeta_{u}})+\delta f_u(\delta t),
\\[0.3em]
\partial_t\varphi+\delta\Delta^2\Phi+\delta\Delta\Phi
+\delta B_1(U,\Phi)
=
\delta\Delta(\Phi^3){+\delta f_\phi(\delta t)},
\end{cases}
\end{equation}
with initial data
\begin{equation*}
v(0)=0,
\qquad
\varphi(0)=0
\qquad\text{in } \mathbb T^3.
\end{equation*}
Our main goal is to show that $w(1)\to0$ in $\mathbb V^1$ as
$\delta\to0^+$.
The algebra property of $H^2$ gives
\begin{equation} \label{defr}
 \sup_{t\in[0,2]}
 \left(\|q_u(t)\|_{H^2}+\|q_\phi(t)\|_{H^4}\right)
 \leq C\left(\norm{z_0}_{\mathbf H^2\times\dot H^4}
       +{\norm{\eta}_{\mathbf H^2\times\dot H^4}}
       +{\norm{\zeta_u}_{H^3}^2}\right)=:r.
\end{equation}
%{\color{red}Here the algebra property of $H^4$ is not needed in \eqref{defr} right? Since $q_\phi(t)=\phi_0+t\,{\color{teal}\eta_\phi}$ is linear in ${\color{teal}\eta_\phi}$, it suffices that $\phi_0,{\color{teal}\eta_\phi}\in\dot H^4$, while $q_u(t)=u_0+t\bigl({\color{teal}\eta_u}-B_0(\zeta_u,\zeta_u)\bigr)$ requires only the algebra property of $H^2$ to get $\|B_0(\zeta_u,\zeta_u)\|_{H^2}\lesssim\|\zeta_u\|_{H^3}^2$, right?? Also do we need $\mathbb V^5$ of $\zeta$ here in this proof?}

\smallskip

\textit{-- Step 2. $H^1$-estimate for $v$.}
Multiplying the first equation of \eqref{chnp_asymp} by $A_0v$, integrating over $\TT$ and then performing integration by parts, we get
\begin{equation}\label{est2_asymp}
	\begin{split}
		\frac{1}{2}&\frac{\d}{\d t}\norm{\nabla v}^2_{L^2}   +  \delta\nu \norm{\Delta v}^2_{L^2}  
		\\ 
		& = -\delta \nu (\Delta({q_{u}}+\delta^{-1/2}{\zeta_{u}}),\Delta v)_{L^2} +\delta b_0(v,v, \Delta v) + \delta b_0({q_{u}},v,\Delta v) +\delta^{1/2} b_0({\zeta_{u}},v,\Delta v) \\
		&\quad  + \delta b_0(v,{q_{u}}, \Delta v) +\delta b_0({q_{u}},{q_{u}},\Delta v) +\delta^{1/2}b_0(v,{\zeta_{u}},\Delta v) +\delta^{1/2} b_0({\zeta_{u}},{q_{u}}, \Delta v)   \\
		&\quad + \delta^{1/2}b_0({q_{u}},{\zeta_{u}},\Delta v) + \delta \left(R_0(\Delta \varphi,\varphi),\Delta v\right)_{L^2} + \delta \left(R_0(\Delta \varphi,{q_{\phi}}), \Delta v\right)_{L^2} \\
		&\quad +\delta \left(R_0(\Delta {q_{\phi}},\varphi), \Delta v\right)_{L^2} + \delta \left(R_0(\Delta {q_{\phi}},{q_{\phi}}), \Delta v\right)_{L^2}
\\
&\quad {-\delta(f_u(\delta t),\Delta v)_{L^2}}.
\end{split}
\end{equation}
To bound the first term, we estimate the two  contributions separately. By
Cauchy--Schwarz and Young's inequalities, together with integration by parts
for the term involving ${\zeta_u}$, we have
\begin{equation}\label{G1_H1_v}
\begin{split}
 \delta\nu\left|(\Delta{q_u},\Delta v)_{L^2}\right|
 &\leq\epsilon\delta\nu\|\Delta v\|_{L^2}^2
       +C\delta\|{q_u}\|_{H^2}^2,\\
 \delta^{1/2}\nu
 \left|(\Delta{\zeta_u},\Delta v)_{L^2}\right|
 &=\delta^{1/2}\nu
 \left|(\nabla\Delta{\zeta_u},\nabla v)_{L^2}\right|
 \\
 &\leq C\delta^{1/2}
 \left(\|{\zeta_u}\|_{H^3}^2
             +\|\nabla v\|_{L^2}^2\right).
\end{split}
\end{equation}
For the second term, we recall the estimate obtained in \eqref{est_B0}
\begin{equation}\label{G2_H1_v}
	\delta \left| b_0(v,v, \Delta v) \right| \leq \epsilon\delta\nu \|\Delta v\|_{L^2}^2 + C\delta \| \nabla v\|_{L^2}^6
\end{equation}
For the third and fourth terms, H\"older's and Young's inequalities
and $H^2\hookrightarrow L^\infty$ give
\begin{align}
 \delta|b_0(q_u,v,\Delta v)|
 &\leq \epsilon\delta\nu\|\Delta v\|_{L^2}^2
       +C\delta\|q_u\|_{H^2}^2\|\nabla v\|_{L^2}^2,\notag\\
 \delta^{1/2}|b_0(\zeta_u,v,\Delta v)|
 &\leq \epsilon\delta\nu\|\Delta v\|_{L^2}^2
       +C\|\zeta_u\|_{H^2}^2\|\nabla v\|_{L^2}^2.\notag\\
\intertext{For the fifth term, we use \Cref{lem:u-L6-grad-u-L3}.
Since $q_u$ has zero mean, \eqref{eq:periodic-norm-equivalences}
gives $\|\nabla q_u\|_{L^2}\leq C\|\Delta q_u\|_{L^2}$. Thus}
 \delta|b_0(v,q_u,\Delta v)|
 &\leq C\delta\|\nabla v\|_{L^2}
       \|\nabla q_u\|_{L^2}^{1/2}\|\Delta q_u\|_{L^2}^{1/2}
       \|\Delta v\|_{L^2}\notag\\
 &\leq C\delta\|\nabla v\|_{L^2}\|\Delta q_u\|_{L^2}
       \|\Delta v\|_{L^2}\notag\\
 &\leq \epsilon\delta\nu\|\Delta v\|_{L^2}^2
       +C\delta\|\nabla v\|_{L^2}^2\|q_u\|_{H^2}^2.
 \label{eq:asymp-fluid-fifth-term}\\
\intertext{The same argument gives the estimate for the sixth term}
 \delta|b_0(q_u,q_u,\Delta v)|
 &\leq \epsilon\delta\nu\|\Delta v\|_{L^2}^2
       +C\delta\|q_u\|_{H^2}^4.
 \label{eq:asymp-fluid-sixth-term}\\
\intertext{For the seventh term, we apply the embedding
$H^2\hookrightarrow L^\infty$ to $\nabla\zeta_u$ and use
Young's inequality to get}
 \delta^{1/2}|b_0(v,\zeta_u,\Delta v)|
\notag &\leq \epsilon\delta\nu\|\Delta v\|_{L^2}^2
       +C\|v\|_{L^2}^2\|\zeta_u\|_{H^3}^2.
% \label{G3_H1_v}
\end{align}
For the eighth term, we use \eqref{eq:ibp_vector_transport}
and \eqref{eq:leibniz_transport_vector}. Applying H\"older's inequality, the embedding $H^2\hookrightarrow L^\infty$,
and Young's inequality, we obtain
\begin{align}\label{ocho_H1_v}
 \delta^{1/2}|b_0(\zeta_u,q_u,\Delta v)|
 &=\delta^{1/2}\left|
   \int_{\TT}\bigl(\nabla\zeta_u\cdot\nabla q_u
          +(\zeta_u\cdot\nabla)\nabla q_u\bigr):\nabla v\,dx
   \right|\notag\\
 &\leq C\delta^{1/2}
 \bigl(\|\nabla\zeta_u\|_{L^\infty}\|\nabla q_u\|_{L^2}
       +\|\zeta_u\|_{L^\infty}\|\nabla^2q_u\|_{L^2}\bigr)
       \|\nabla v\|_{L^2}\notag\\
 &\leq C\delta^{1/2}\|\zeta_u\|_{H^3}\|q_u\|_{H^2}
                        \|\nabla v\|_{L^2}\notag\\
 &\leq C\delta^{1/2}\|\nabla v\|_{L^2}^2
       +C\delta^{1/2}\|\zeta_u\|_{H^3}^2\|q_u\|_{H^2}^2.
\end{align}
For the ninth term, the same identities and inequalities give
\begin{align}\label{nueve_H1_v}
 \delta^{1/2}|b_0(q_u,\zeta_u,\Delta v)|
 &\leq C\delta^{1/2}
 \bigl(\|\nabla q_u\|_{L^2}\|\nabla\zeta_u\|_{L^\infty}
       +\|q_u\|_{L^\infty}\|\nabla^2\zeta_u\|_{L^2}\bigr)
       \|\nabla v\|_{L^2}\notag\\
 &\leq C\delta^{1/2}\|q_u\|_{H^2}\|\zeta_u\|_{H^3}
                        \|\nabla v\|_{L^2}\notag\\
 &\leq C\delta^{1/2}\|\nabla v\|_{L^2}^2
       +C\delta^{1/2}\|q_u\|_{H^2}^2\|\zeta_u\|_{H^3}^2.
\end{align}
For the tenth term, we recall the estimate obtained in \eqref{est_R0_Delta-u}
\begin{equation}\label{G4_H1_v}
\begin{split}
	\delta \left| \left(R_0(\Delta \varphi,\varphi),\Delta v\right)_{L^2} \right| &\leq \epsilon\delta\nu \|\Delta v\|_{L^2}^2 + \epsilon\delta\|\Delta^2\varphi\|^2_{L^2} + C\delta \|\nabla\varphi\|_{L^2}^2\|\Delta\varphi\|_{L^2}^6 \\
	&\leq \epsilon\delta\nu \|\Delta v\|_{L^2}^2 + \epsilon\delta\|\Delta^2\varphi\|^2_{L^2} + C\delta\|\Delta\varphi\|_{L^2}^8
\end{split}
\end{equation}
The {three remaining capillary} terms can be bounded directly since ${q_{\phi}\in\dot H^4}$. More precisely,
\begin{equation}\label{G5_H1_v}
\begin{split}
	\delta \left| \left(R_0(\Delta \varphi,{q_{\phi}}), \Delta v\right)_{L^2} \right|
	& \leq \epsilon\delta\nu \|\Delta v \|^2_{L^2} + C\delta \|\Delta\varphi\|^2_{L^2}\|{q_{\phi}}\|^2_{H^3} \\
	\delta \left| \left(R_0(\Delta {q_{\phi}},\varphi), \Delta v\right)_{L^2} \right|
	& \leq \epsilon\delta\nu \|\Delta v \|^2_{L^2} + C\delta \|{q_{\phi}}\|^2_{H^4}\|\nabla\varphi\|^2_{L^2} \\
	\delta \left| \left(R_0(\Delta {q_{\phi}},{q_{\phi}}), \Delta v\right)_{L^2} \right|
	&\leq \epsilon\delta\nu \|\Delta v \|^2_{L^2} + C\delta \|{q_{\phi}}\|^4_{H^4}
\end{split}
\end{equation}
{Finally, for almost every $t\in(0,\min\{2,\delta^{-1}{\mathcal T_\delta}\})$,
Cauchy--Schwarz and Young's inequalities give
\begin{equation}\label{eq:asymp-background-force-u}
 \delta\bigl|(f_u(\delta t),\Delta v)_{L^2}\bigr|
 \leq\epsilon\delta\nu\|\Delta v\|_{L^2}^2
     +C_{\epsilon,\nu}\delta M_f^2.
\end{equation}
where  $M_f:=\|f\|_{L^\infty(0,2;\mathbb V^0)}$.}

\smallskip
\textit{-- Step 3. $H^2$-estimate for $\varphi$.} Multiplying the second equation in \eqref{chnp_asymp} by $\Delta^2 \varphi$, integrating over in $\TT$ and using integration by parts, we get
\begin{equation}\label{eq:estphi2_asymp}
	\begin{split}
		\frac{1}{2}&\frac{\d}{\d t}\norm{\Delta \varphi}^2_{L^2}   +  \delta \norm{\Delta^2 \varphi}^2_{L^2}  
		\\ 
		& = - \delta (\Delta^2{q_{\phi}},\Delta^2 \varphi)_{L^2} - \delta(\Delta \varphi,\Delta^2 \varphi)_{L^2} - \delta(\Delta {q_{\phi}},\Delta^2 \varphi)_{L^2} - \delta b_1(v,\varphi, \Delta^2 \varphi)
		\\
		& \quad - \delta b_1(v,{q_{\phi}}, \Delta^2 \varphi) - \delta b_1({q_{u}},\varphi, \Delta^2 \varphi) - \delta b_1({q_{u}},{q_{\phi}}, \Delta^2 \varphi) - \delta^{1/2} b_1({\zeta_{u}},\varphi, \Delta^2 \varphi)
		\\
		& \quad - \delta^{1/2} b_1({\zeta_{u}},{q_{\phi}}, \Delta^2 \varphi) + \delta \left(\Delta(\varphi+{q_{\phi}})^3,\Delta^2 \varphi\right)_{L^2}
\\
&\quad {+\delta(f_\phi(\delta t),\Delta^2\varphi)_{L^2}}.
\end{split}
\end{equation}
The first three terms can be bounded directly using Cauchy-Schwarz and Young inequalities 
\begin{equation}\label{G1_H2_phi}
\begin{split}
	\delta \left| (\Delta^2{q_{\phi}},\Delta^2 \varphi)_{L^2} \right| 
	&\leq \epsilon\delta\|\Delta^2 \varphi\|_{L^2}^2 + C\delta \|\Delta^2 {q_{\phi}} \|_{L^2}^2
	\leq \epsilon\delta\|\Delta^2 \varphi\|_{L^2}^2 + C\delta \|{q_{\phi}} \|_{H^4}^2 \\
	\delta \left| (\Delta \varphi,\Delta^2 \varphi)_{L^2} \right| 
	&\leq \epsilon\delta\|\Delta^2 \varphi\|_{L^2}^2 + C\delta \|\Delta \varphi \|_{L^2}^2 \\
	\delta \left| (\Delta {q_{\phi}},\Delta^2 \varphi)_{L^2} \right| 
	&\leq \epsilon\delta\|\Delta^2 \varphi\|_{L^2}^2 + C\delta \|\Delta {q_{\phi}} \|_{L^2}^2
	\leq \epsilon\delta\|\Delta^2 \varphi\|_{L^2}^2 + C\delta \|{q_{\phi}} \|_{H^2}^2.
\end{split}
\end{equation}
For the fourth term, using H\"older inequality, \Cref{lem:u-L6-grad-u-L3} and Young inequalities yields to
\begin{align}\label{G2_H2_phi}
\begin{split}
	\delta \left| b_1(v,\varphi, \Delta^2 \varphi) \right | 
	&\leq {\delta}\|v\|_{L^6}\|\nabla \varphi\|_{L^3}\|\Delta^2\varphi\|_{L^2}
	\leq C{\delta} \|\nabla v\|_{L^2}\|\Delta\varphi\|_{L^2}\|\Delta^2\varphi\|_{L^2} \\
	&\leq \epsilon\delta \|\Delta^2\varphi\|_{L^2}^2 + C\delta \|\nabla v\|_{L^2}^2\|\Delta\varphi\|_{L^2}^2 \\
	&\leq \epsilon\delta \|\Delta^2\varphi\|_{L^2}^2 + C\delta \|\nabla v\|_{L^2}^4+C\delta\|\Delta\varphi\|_{L^2}^4.
\end{split}
\end{align}
{
For the fifth to eighth terms, we use Cauchy--Schwarz and Young's
inequalities to obtain
}
\begin{equation}\label{G3_H2_phi}
\begin{split}
	\delta \left| b_1(v,{q_{\phi}}, \Delta^2 \varphi) \right|
	&\leq \epsilon\delta \| \Delta^2 \varphi \|_{L^2}^2 + C\delta \|v\|_{L^2}^2\|{q_{\phi}}\|_{H^3}^2\\
	\delta \left| b_1({q_{u}},\varphi, \Delta^2 \varphi) \right| 
	&\leq \epsilon\delta \| \Delta^2 \varphi \|_{L^2}^2 + C\delta \|{q_{u}}\|_{H^2}^2\|\nabla\varphi\|_{L^2}^2\\
	\delta \left| b_1({q_{u}},{q_{\phi}}, \Delta^2 \varphi) \right|
	&\leq \epsilon\delta \| \Delta^2 \varphi \|_{L^2}^2 + C\delta \|{q_{u}}\|_{H^2}^2\|{q_{\phi}}\|_{H^3}^2\\
	\delta^{1/2} \left| b_1({\zeta_{u}},\varphi, \Delta^2 \varphi) \right|
	&\leq \epsilon\delta \| \Delta^2 \varphi \|_{L^2}^2 + C \|{\zeta_{u}}\|_{H^2}^2 \|\nabla\varphi\|_{L^2}^2. 
\end{split}
\end{equation}
For the ninth term, by integration by parts and using Cauchy-Schwarz inequality, we get
\begin{align*}
\delta^{1/2}
\left|b_1({\zeta_{u}},{q_{\phi}},\Delta^2\varphi)\right| &= \delta^{1/2}\left|  \int_{\TT} \Delta ({\zeta_{u}}\cdot \nabla {q_{\phi}})\Delta \varphi \right|
\\
& \leq \delta^{1/2}\|\Delta({\zeta_{u}}\cdot \nabla {q_{\phi}})\|_{L^2}\|\Delta\varphi\|_{L^2} \leq \delta^{1/2}\|{\zeta_{u}}\cdot \nabla {q_{\phi}}\|_{H^2}\|\Delta\varphi\|_{L^2}. 
\end{align*}
Since $H^2$ is algebra, $\|{\zeta_{u}}\cdot \nabla {q_{\phi}}\|_{H^2}\leq \|{\zeta_{u}}\|_{H^2}\|{q_{\phi}}\|_{H^3}$. Therefore, we have
\begin{align}\label{G8_H2_phi}
\delta^{1/2}\left|b_1({\zeta_{u}},{q_{\phi}},\Delta^2\varphi)\right|\leq
C\delta^{1/2}
\norm{{\zeta_{u}}}_{H^2}^2
\norm{{q_{\phi}}}_{H^3}^2
+
C\delta^{1/2}
\norm{\Delta\varphi}_{L^2}^2.
\end{align}
For the cubic term, let us recall that $\Phi=\varphi+{q_{\phi}}$ and using the identity
\begin{equation*}
	\Delta(\Phi^3)=6\Phi|\nabla\Phi|^2+3\Phi^2\Delta\Phi,
\end{equation*}
we have
\begin{equation}\label{estDeltaPhi}
	\delta |(\Delta(\Phi^3),\Delta^2\varphi )_{L^2}| \leq 6\delta \int_{\TT}  |\Phi||\nabla\Phi|^2|\Delta^2\varphi| +
	3\delta \int_{\TT} |\Phi|^2|\Delta\Phi||\Delta^2\varphi|.
\end{equation}
To estimate the first term on the right-hand side of \eqref{estDeltaPhi}, we use H\"older inequality, \eqref{eq:paso_interp}, Poincaré and  Young's inequalities to get
\begin{align*}
	6\delta \int_{\TT}  |\Phi||\nabla\Phi|^2|\Delta^2\varphi| 
	&\leq C{\delta}\|\Phi\|_{L^\infty} \|\nabla\Phi\|^2_{L^4}\|\Delta^2\varphi\|_{L^2} \\
	&\leq C{\delta}\|\Phi\|_{L^\infty}\|\Phi\|^2_{H^2} \|\Delta^2\varphi\|_{L^2} \\
	&\leq {\frac{\epsilon}{2}}\delta\|\Delta^2\varphi\|_{L^2}^2 + C\delta\left(\|\Phi \|_{H^2}^3 \right)^2 \\
	&\leq {\frac{\epsilon}{2}}\delta\|\Delta^2\varphi\|_{L^2}^2 + C\delta\|\varphi \|_{H^2}^6 + C\delta\|{q_{\phi}} \|_{H^2}^6.
\end{align*}
To estimate the second term on the right-hand side of \eqref{estDeltaPhi}, we use H\"older inequality, Young's inequalities and \cref{lem:interp-estimates} to get
\begin{align*}
	3\delta \int_{\TT} |\Phi|^2|\Delta\Phi||\Delta^2\varphi| 
	&\leq C{\delta}\|\Phi\|_{L^\infty}^2{\|\Delta\Phi\|_{L^2}}\|\Delta^2\varphi\|_{L^2} \\
	&\leq C{\delta}\|\Phi\|_{H^2}^3\|\Delta^2\varphi\|_{L^2} \\
	&\leq {\frac{\epsilon}{2}}\delta \|\Delta^2\varphi\|_{L^2}^2 + C\delta\|\varphi \|_{H^2}^6 + C\delta\|{q_{\phi}} \|_{H^2}^6.
\end{align*}
Adding these two bounds on the right-hand side of
\eqref{estDeltaPhi}, we obtain
\begin{equation}\label{G4_H2_phi}
	\delta |(\Delta(\Phi^3),\Delta^2\varphi )_{L^2}| \leq \epsilon\delta \|\Delta^2\varphi\|_{L^2}^2 + C\delta\|\varphi \|_{H^2}^6 + C\delta\|{q_{\phi}} \|_{H^2}^6.
\end{equation}
{Finally, for almost every $t\in(0,\min\{2,\delta^{-1}{\mathcal T_\delta}\})$,
Cauchy--Schwarz and Young's inequalities give
\begin{equation}\label{eq:asymp-background-force-phi}
 \delta\bigl|(f_\phi(\delta t),\Delta^2\varphi)_{L^2}\bigr|
 \leq\epsilon\delta\|\Delta^2\varphi\|_{L^2}^2
     +C_\epsilon\delta M_f^2.
\end{equation}}

 \textit{-- Step 4. Combining the estimates.} 
% First, by combining \eqref{est1_asymp} with \eqref{eq:estphi_asymp} and using the estimates \eqref{G1_L2_v}--\eqref{G5_L2_v} together with \eqref{G1_H1_phi}--\eqref{G4_H1_phi}, we get 
%\begin{equation} \label{eq:estL2H1}
%\begin{split}
%	\frac{1}{2}\frac{\d}{\d t}\norm{v}^2_{L^2} &+ \frac{1}{2}\frac{\d}{\d t}\norm{\nabla \varphi}^2_{L^2}   
%	+ \delta\nu \norm{\nabla v}^2_{L^2} 
%	+ \delta \norm{\nabla \Delta \varphi}^2_{L^2} \\ 
%	&\leq C(1+\delta^{1/2}) + C(1+\delta^{1/2})\|v\|_{L^2}^2  +C(1+\delta^{1/2})\|\varphi\|_{H^1}^2  + C\delta\|\varphi\|_{H^2}^6
%\end{split}
%\end{equation}
%where $C>0$ may depend on $r$, defined in \eqref{defr}, but not on $\delta>0$.
%
Combining \eqref{est2_asymp} with \eqref{eq:estphi2_asymp}
and using \eqref{G1_H1_v}--\eqref{eq:asymp-background-force-u}
and \eqref{G1_H2_phi}--\eqref{eq:asymp-background-force-phi},
we choose $\epsilon>0$ sufficiently small and absorb the
corresponding dissipation terms to obtain
\begin{equation} \label{eq:estH1H2}
\begin{split}
	\frac{1}{2} \frac{\d}{\d t}\norm{\nabla v}^2_{L^2}  &+ \frac{1}{2}\frac{\d}{\d t}\norm{\Delta \varphi}^2_{L^2}  +  {\frac{\delta\nu}{2}} \norm{\Delta v}^2_{L^2}   +  {\frac{\delta}{2}} \norm{\Delta^2 \varphi}^2_{L^2} \\
	&\leq C(\delta+\delta^{1/2}) + C(1+\delta+\delta^{1/2})\|\nabla v\|_{L^2}^2 + C(1+\delta+\delta^{1/2})\|\Delta \varphi\|_{L^2}^2 \\
	&\quad + C\delta (\|v\|_{H^1}^4 +\|v\|_{H^1}^6 + \|\varphi\|_{H^2}^4 + \|\varphi\|_{H^2}^6 + \|\varphi\|_{H^2}^8 ),
\end{split}
\end{equation}
for some constant $C>0$ that may depend on $r$ {and $M_f$} but not on $\delta\in(0,1)$. Here we have used repeatedly \eqref{eq:periodic-norm-equivalences} and the  fact $H^{s_1}\hookrightarrow H^{s_0}$ for $s_1>s_0$.

Finally, using that $\delta\in(0,1)$, we obtain from \eqref{eq:estH1H2}
\begin{align*}
	\frac{\d}{\d t}\left(\norm{\nabla v}^2_{L^2}+\norm{\Delta \varphi}^2_{L^2}\right) \leq  C\left(\norm{\nabla v}^2_{L^2}+\norm{\Delta \varphi}_{L^2}^2\right)+C\delta^{1/2} \left(1+\norm{v}_{H^1}^6+\norm{\varphi}^8_{H^2}\right). 
\end{align*}
By Gr\"onwall's inequality and \eqref{eq:periodic-norm-equivalences}, we have the following estimate $\text{ for }t<  \min\{2, \delta^{-1}{\mathcal T_\delta}\}$
\begin{align}\label{grn}
	\norm{v(t)}^2_{H^1}+\norm{\varphi(t)}^2_{H^2}\leq C\delta^{1/2}t +C \delta^{1/2}\int_{0}^{t}\left(\norm{v(s)}^6_{H^1}+\norm{\varphi(s)}_{H^2}^8\right) \, ds.
\end{align}

\textit{-- Step 5. Asymptotic limit and conclusion.}
Choose $\delta_0\in(0,1)$ such that $3C\delta_0^{1/2}<1$,
where $C$ is the constant in \eqref{grn}, and fix
$\delta\in(0,\delta_0)$. Define
\begin{equation*}
\begin{aligned}
 \gamma^\delta:=\sup\bigl\{
 t\in(0,\delta^{-1}{\mathcal T_\delta}):\;&
 \|v(s)\|_{H^1}^2+\|\varphi(s)\|_{H^2}^2<1
 \text{ for every }s\in[0,t]\bigr\}.
\end{aligned}
\end{equation*}
Since $v(0)=0$, $\varphi(0)=0$, and $w$ is continuous in
$\mathbb V^1$, we have
$0<\gamma^\delta{\leq}\delta^{-1}{\mathcal T_\delta}$.
Moreover,
$\|v(t)\|_{H^1}^2+\|\varphi(t)\|_{H^2}^2<1$
for every $t\in[0,\gamma^\delta)$.

If $\gamma^\delta=+\infty,$ then $\delta^{-1}\mathcal T_\delta> 1$ follows trivially. Let us consider the case $\gamma^\delta<+\infty.$ 
In order to prove $\delta^{-1}{\mathcal T_\delta}>1,$ it is enough to show that, for all $\delta\in (0,\delta_0),$  $\gamma^\delta> 1.$ We prove it by contradiction. If possible, let us assume that there exists $\delta\in (0,\delta_0)$ such that $\gamma^{\delta}\leq 1.$ For $0\leq t<\gamma^\delta$,
\eqref{grn} gives
\begin{equation}\label{eq:asymp-first-exit-estimate}
	\|v(t)\|_{H^1}^2+\|\varphi(t)\|_{H^2}^2
	\leq C\delta^{1/2}t+2C\delta^{1/2}t
	=3C\delta^{1/2}t<1.
\end{equation}
If $\gamma^\delta=\delta^{-1}\mathcal T_\delta$, then
$z(s)+\delta^{-1/2}\zeta=w(s/\delta)+q(s/\delta)+\delta^{-1/2}\zeta$
remains bounded in $\mathbb V^1$ for $0\leq s<\mathcal T_\delta$,
by \eqref{defr} and \eqref{eq:asymp-first-exit-estimate}.
By the definition of $\mathcal T_\delta$, this contradicts the blow-up
alternative in \Cref{prop:maximal-strong}.
Hence $\gamma^\delta<\delta^{-1}{\mathcal T_\delta}$, and
continuity of $v$ and $\varphi$ and the definition of $\gamma^\delta$ imply
\begin{equation*}
	1=\|v(\gamma^\delta)\|_{H^1}^2
	+\|\varphi(\gamma^\delta)\|_{H^2}^2
	\leq3C\delta^{1/2}\gamma^\delta
	\leq3C\delta_0^{1/2}<1,
\end{equation*}
a contradiction. Thus we can conclude that, there exists a $\delta_0>0$ such that when $0<\delta<\delta_0,$ $\gamma^\delta> 1$ and thus $\delta^{-1}{\mathcal T_\delta}> 1.$
In particular, the
solution $\mathcal S\bigl(z_0,\delta^{-1/2}\zeta,{f+\delta^{-1}\eta}\bigr)$ is defined on $[0,\delta]$. Applying \eqref{grn}
at $t=1$, we obtain
\begin{equation*}
	\|w(1)\|_{\mathbb V^1}^2
	\leq\|v(1)\|_{H^1}^2+\|\varphi(1)\|_{H^2}^2
	\leq3C\delta^{1/2}\longrightarrow0
	\qquad\text{as }\delta\to0^+.
\end{equation*}
Finally,
\begin{align*}
	\mathcal S_\delta(z_0,\delta^{-1/2}\zeta,
	{f+\delta^{-1}\eta})
	&=q(1)+w(1)\longrightarrow
	z_0+\eta-\bigl(B_0({\zeta_{u}},
	{\zeta_{u}}),0\bigr)
	\qquad\text{as }\delta\to0^+,
\end{align*}
which proves \eqref{asymp_phi}.
\end{proof}

In a similar manner, we establish another asymptotic limit that will play a crucial role in proving our desired result.
\subsection{Second asymptotic limit}

{\begin{theo}\label{thm:asymp-mixed}
Let $z_0=(u_0,\phi_0)\in{\mathbf H^2\times\dot H^4}$, $f=(f_u, f_\phi) \in L^\infty_{\mathrm{loc}}(\mathbb R_+;\mathbb V^0),$
$\zeta_u\in\mathbf H^4$, and $\xi_\phi\in\dot H^8$.
Set $\zeta=(\zeta_u,0)$ and $\xi=(0,\xi_\phi)$, and assume that
\begin{equation}\label{hyp_mix}
B_0({\zeta_{u}}, \zeta_{u})=0.
\end{equation}
Then there exists $\delta_0>0$ such that for all $\delta\in(0,\delta_0)$, the solution
\begin{equation*}
\mathcal S_t\bigl(z_0,\delta^{-2/3}\zeta+\delta^{-1/3}\xi,f\bigr)
\end{equation*}
is well defined on the time interval $[0,\delta]$. Moreover, the following
convergence holds in $\mathbb V^1$:
\begin{equation}\label{asymp_mixed}
\mathcal S_\delta\bigl(z_0,\delta^{-2/3}\zeta+\delta^{-1/3}\xi,{f}\bigr)
\longrightarrow
z_0+\bigl(0,\Delta({\xi_\phi}^3)-B_1({\zeta_{u}},{\xi_\phi})\bigr)
\qquad\text{as } \delta\to0^+.
\end{equation}
\end{theo}
\begin{proof}
We divide the proof into five steps.

\smallskip
\noindent
\textit{Step 1. A change of variable.}
Let $\delta\in(0,1)$ be arbitrary and let us define
\begin{equation*}
w(t)=z(\delta t)-q(t),
\end{equation*}
where
\begin{equation*}
z(t)=\mathcal S_t\bigl(z_0,\delta^{-2/3}\zeta+\delta^{-1/3}\xi,{f}\bigr),
\qquad
q(t)=\bigl({q_{u}}(t),{q_{\phi}}(t)\bigr),
\end{equation*}
with
\begin{equation*}
{q_{u}}(t)=u_0,
\qquad
{q_{\phi}}(t)=\phi_0+t\bigl(\Delta({\xi_\phi}^3)-B_1({\zeta_{u}},{\xi_\phi})\bigr).
\end{equation*}
Since $\delta^{-2/3}\zeta+\delta^{-1/3}\xi\in \mathbf H^2\times\dot H^4$,
\Cref{wp1} gives a maximal strong solution with
\begin{equation*}
 \mathcal T_\delta
 :=\mathcal T\bigl(z_0,\delta^{-2/3}\zeta+\delta^{-1/3}\xi,f\bigr)
 =T_{\max}\bigl(z_0+\delta^{-2/3}\zeta+\delta^{-1/3}\xi,f\bigr).
\end{equation*}
The functions $z$ and $w$ are defined on
$[0,\mathcal T_\delta)$ and $[0,\delta^{-1}\mathcal T_\delta)$,
respectively.

We also set
\begin{equation*}
U:=v+{q_{u}}+\delta^{-2/3}{\zeta_{u}},
\qquad
\Phi:=\varphi+{q_{\phi}}+\delta^{-1/3}{\xi_\phi},
\end{equation*}
where $(v,\varphi)=({\pi_{u}} w,{\pi_{\phi}} w)$. In this case, $(v,\varphi)$ satisfies
the following system in $(0,\delta^{-1}{\mathcal T_\delta})\times\mathbb T^3$
\begin{equation}\label{chnp_asymp_mix}
\begin{cases}
\partial_t v+\delta\nu A_0U+\delta B_0(U,U)
+\delta R_0(\Delta\Phi,\Phi)
=
{\delta f_u(\delta t)},
\\[0.3em]
\partial_t\varphi+\delta\Delta^2\Phi+\delta\Delta\Phi
+\delta B_1(U,\Phi)
=
\delta\Delta(\Phi^3)
-\Delta({\xi_\phi}^3)+B_1({\zeta_{u}},{\xi_\phi}){+\delta f_\phi(\delta t)},
\end{cases}
\end{equation}
with initial data
\begin{equation*}
v(0)=0,
\qquad
\varphi(0)=0
\qquad\text{in } \mathbb T^3.
\end{equation*}
 
Our main goal is to show that $w(1)\to0$ in $\mathbb V^1$ as
$\delta\to0^+$.

\smallskip
\textit{-- Step 2. $H^1$-estimate for $v$.}
{
The algebra properties of $H^4$ and $H^6$, together with Young's
inequality, give
\begin{equation} \label{defr_mix}
\begin{aligned}
 &\sup_{t\in[0,2]}
     \left(\|q_u(t)\|_{H^2}+\|q_\phi(t)\|_{H^4}\right)\\
 &\qquad\leq C\left(\|z_0\|_{\mathbf H^2\times\dot H^4}
       +\|\zeta_u\|_{H^4}^2
       +\|\xi_\phi\|_{H^8}^2
       +\|\xi_\phi\|_{H^8}^3\right)=:r.
\end{aligned}
\end{equation}
}
Multiplying the first equation of \eqref{chnp_asymp_mix} by $A_0v$ in $L^2(\TT)$, integrating by parts and using the cancellation property in \Cref{lem-cancellation-B0} we get
\begin{align}\label{est2_asymp_mix}
%\begin{split}
\notag\frac{1}{2}&\frac{\d}{\d t}\norm{\nabla v}^2_{L^2}   +  \delta\nu \norm{\Delta v}^2_{L^2}  
\\ 
\notag& = -\delta \nu (\Delta({q_{u}}+\delta^{-2/3}{\zeta_{u}}),\Delta v)_{L^2} + \delta b_0(v,v,\Delta v) + \delta b_0(v,{q_{u}}, \Delta v) 
\\
\notag& \quad + \delta b_0({q_{u}},v,\Delta v) + \delta b_0({q_{u}},{q_{u}}, \Delta v) + \delta^{1/3} b_0({q_{u}},{\zeta_{u}}, \Delta v) + \delta^{1/3} b_0({\zeta_{u}},{q_{u}}, \Delta v)
+ \delta^{1/3}b_0({\zeta_{u}},v,\Delta v) \\ 
\notag& \quad + \delta^{1/3} b_0(v,{\zeta_{u}}, \Delta v) +\delta^{-1/3}b_0({\zeta_{u}},{\zeta_{u}},\Delta v) + \delta \left(R_0(\Delta \varphi,\varphi),\Delta v\right)_{L^2} + \delta \left(R_0(\Delta \varphi,{q_{\phi}}),\Delta v\right)_{L^2}
\\
\notag& \quad  + \delta \left(R_0(\Delta {q_{\phi}},\varphi),\Delta v\right)_{L^2} + \delta \left(R_0(\Delta {q_{\phi}},{q_{\phi}}),\Delta v\right)_{L^2} + \delta^{2/3}\left(R_0(\Delta \varphi,{\xi_\phi}),\Delta v\right)_{L^2}
\\
\notag& \quad  + \delta^{2/3}\left(R_0(\Delta {q_{\phi}},{\xi_\phi}), \Delta v\right)_{L^2} + \delta^{2/3}\left(R_0(\Delta {\xi_\phi},\varphi),\Delta v\right)_{L^2} + \delta^{2/3}\left(R_0(\Delta {\xi_\phi},{q_{\phi}}),\Delta v\right)_{L^2}
 \\
& \quad  + \delta^{1/3}\left(R_0(\Delta {\xi_\phi},{\xi_\phi}),\Delta v\right)_{L^2}
 {-\delta(f_u(\delta t),\Delta v)_{L^2}}.
%\end{split}
\end{align}
To bound the first term, we estimate the two contributions separately.
By Cauchy--Schwarz and Young's inequalities, together with integration
by parts for the term involving $\zeta_u$, we have
\begin{align}\label{G1_H1_v_mix}
 \delta\nu|(\Delta q_u,\Delta v)_{L^2}|
 &\leq\epsilon\delta\nu\|\Delta v\|_{L^2}^2
       +C\delta\|q_u\|_{H^2}^2,\notag\\
 \delta^{1/3}\nu|(\Delta\zeta_u,\Delta v)_{L^2}|
 &=\delta^{1/3}\nu|(\nabla\Delta\zeta_u,\nabla v)_{L^2}|
       \notag\\
 &\leq C\delta^{1/3}
       \bigl(\|\zeta_u\|_{H^3}^2+\|\nabla v\|_{L^2}^2\bigr).
\end{align}
The second to fifth terms are estimated as in \eqref{G2_H1_v}
and \eqref{eq:asymp-fluid-fifth-term}--\eqref{eq:asymp-fluid-sixth-term},
using $H^2\hookrightarrow L^\infty$ for $b_0(q_u,v,\Delta v)$.
For brevity, we only present the final estimates
\begin{align}\label{grupo1_v_b0}
\begin{split}
	\delta \left| b_0(v,v, \Delta v) \right| &\leq \epsilon\delta\nu \|\Delta v\|_{L^2}^2 + C\delta \| \nabla v\|_{L^2}^6 \\
\delta \left| b_0(v,{q_{u}}, \Delta v) \right|
	&\leq \epsilon \delta \nu \|\Delta  v\|^2_{L^2}+{C\delta \|\nabla v\|_{L^2}^2\|q_u\|_{H^2}^2} \\
	\delta \left| b_0({q_{u}},v,\Delta v) \right| 
	&\leq \epsilon \delta \nu \|\Delta  v\|^2_{L^2}+C\delta \|{q_{u}}\|_{H^2}^2\|\nabla v\|_{L^2}^2 \\
	\delta \left| b_0({q_{u}},{q_{u}},\Delta v) \right| 
	& \leq \epsilon \delta \nu \|\Delta  v\|^2_{L^2}+{C\delta \|q_u\|_{H^2}^4}.
	\end{split}
\end{align}
{
For the sixth and seventh terms, we apply \eqref{ocho_H1_v}
and \eqref{nueve_H1_v}, with $\delta^{1/3}$ in place of
$\delta^{1/2}$, to obtain
\begin{align}\label{grupo2_v_b0}
\begin{split}
 \delta^{1/3}|b_0(\zeta_u,q_u,\Delta v)|
 &\leq C\delta^{1/3}\|\nabla v\|_{L^2}^2
       +C\delta^{1/3}\|\zeta_u\|_{H^3}^2\|q_u\|_{H^2}^2,\\
 \delta^{1/3}|b_0(q_u,\zeta_u,\Delta v)|
 &\leq C\delta^{1/3}\|\nabla v\|_{L^2}^2
       +C\delta^{1/3}\|q_u\|_{H^2}^2\|\zeta_u\|_{H^3}^2.
\end{split}
\end{align}
}
For the eight and nine terms, we use {\Cref{lem:transport_fluid_first_argument,lem:transport_fluid_second_argument}, which are} based on suitable integration by parts. The final estimates are
\begin{align}\label{grupo3_v_b0}
\begin{split}
\delta^{1/3} \left| b_0({\zeta_{u}},v, \Delta v) \right| &\leq C\delta^{1/3}\norm{{\zeta_{u}}}_{H^3} \norm{\nabla v}_{L^2}^2\\
\delta^{1/3} \left| b_0(v,{\zeta_{u}}, \Delta v) \right| &\leq C\delta^{1/3}\|{\zeta_{u}}\|_{H^4} \|\nabla v\|^2_{L^2}.
\end{split}
\end{align} 
The tenth term 
\begin{align}\label{cancelacion_v_b0}
\delta^{-1/3}b_0({\zeta_{u}},{\zeta_{u}},\Delta v)=0
\end{align}
due to hypothesis \eqref{hyp_mix}. Notice that if this is not the case, this term will blow up as $\delta\to 0^+$.

The eleventh to fourteenth  is the same as in \eqref{G4_H1_v} and \eqref{G5_H1_v}
\begin{align}\label{G4_H1_v_mix}
\begin{split}
	\delta \left| \left(R_0(\Delta \varphi,\varphi),\Delta v\right)_{L^2} \right| &\leq \epsilon\delta\nu \|\Delta v\|_{L^2}^2 + \epsilon\delta{\|\Delta^2\varphi\|^2_{L^2}} + C\delta\|\Delta\varphi\|_{L^2}^8 \\
	\delta \left| \left(R_0(\Delta \varphi,{q_{\phi}}), \Delta v\right)_{L^2} \right|
	& \leq \epsilon\delta\nu \|\Delta v \|^2_{L^2} + C\delta \|\Delta\varphi\|^2_{L^2}\|{q_{\phi}}\|^2_{H^3} \\
	\delta \left| \left(R_0(\Delta {q_{\phi}},\varphi), \Delta v\right)_{L^2} \right|
	& \leq \epsilon\delta\nu \|\Delta v \|^2_{L^2} + C\delta \|{q_{\phi}}\|^2_{H^4}\|\nabla\varphi\|^2_{L^2} \\
	\delta \left| \left(R_0(\Delta {q_{\phi}},{q_{\phi}}), \Delta v\right)_{L^2} \right|
	&\leq \epsilon\delta\nu \|\Delta v \|^2_{L^2} + C\delta \|{q_{\phi}}\|^4_{H^4}.
\end{split}
\end{align}
For the terms fifteen through eighteen, we apply H\"older and Young inequalities and use the fact that {$q_\phi,\xi_\phi\in\dot H^4$} to obtain
\begin{equation}\label{grupo5_v_R0}
\begin{split}
\delta^{2/3}\left|\left(R_0(\Delta \varphi,{\xi_\phi}),\Delta v\right)_{L^2}\right| & \leq \epsilon\delta \nu \|\Delta v\|_{L^2}^2+C\delta^{1/3}\|{\xi_\phi}\|_{H^3}^2\|\Delta \varphi\|_{L^2}^2
\\
\delta^{2/3}\left|\left(R_0(\Delta {q_{\phi}},{\xi_\phi}),\Delta v\right)_{L^2}\right| & \leq \epsilon\delta \nu \|\Delta v\|_{L^2}^2+C\delta^{1/3}\|{q_{\phi}}\|_{H^2}^2\|{\xi_\phi}\|_{H^3}^2
\\
\delta^{2/3}\left|\left(R_0(\Delta {\xi_\phi},\varphi),\Delta v\right)_{L^2}\right| & \leq \epsilon\delta \nu \|\Delta v\|_{L^2}^2+C\delta^{1/3}\|{\xi_\phi}\|_{H^4}^2\|\nabla \varphi\|_{L^2}^2
\\
\delta^{2/3}\left|\left(R_0(\Delta {\xi_\phi},{q_{\phi}}),\Delta v\right)_{L^2}\right| & \leq \epsilon\delta \nu \|\Delta v\|_{L^2}^2+C\delta^{1/3}\|{\xi_\phi}\|_{H^2}^2\|{q_{\phi}}\|_{H^3}^2.
\end{split}
\end{equation}
For the remaining capillary term, we argue as in \eqref{G8_H2_phi}, that is, we integrate by parts to get
\begin{equation}\label{grupo6_v_R0}
\begin{split}
\delta^{1/3}\left|\left(R_0(\Delta {\xi_\phi},{\xi_\phi}),\Delta v\right)_{L^2}\right|&=\delta^{1/3}\left|\sum_{j=1}^{3}\int_{\TT}\Delta(\Delta{\xi_\phi}\partial_j {\xi_\phi}) v_j\right|
\\
&\leq C\delta^{1/3}\|v\|_{L^2}^2+C\delta^{1/3}\|{\xi_\phi}\|_{H^4}^4.
\end{split}
\end{equation}
Finally, set $M_f:=\|f\|_{L^\infty(0,2;\mathbb V^0)}$.
For almost every $t\in(0,\min\{2,\delta^{-1}\mathcal T_\delta\})$,
we have $\delta t<2$. Cauchy--Schwarz and Young's inequalities give,
for every $\epsilon>0$,
\begin{align}\label{eq:asymp-mixed-background-force-u}
 \delta\bigl|(f_u(\delta t),\Delta v)_{L^2}\bigr|
 &\leq\delta\|f_u(\delta t)\|_{L^2}\|\Delta v\|_{L^2}\notag\\
 &\leq\epsilon\delta\nu\|\Delta v\|_{L^2}^2
       +C_{\epsilon,\nu}\delta M_f^2.
\end{align}

%\medskip
\textit{-- Step 3. $H^2$-estimate for $\varphi$.} Multiplying the second equation in {\eqref{chnp_asymp_mix}} by $\Delta^2 \varphi$ in $L^2(\TT)$ and using integration by parts, we get
\begin{align}\label{eq:estphi2_asymp_mix}
	\notag	\frac{1}{2}&\frac{\d}{\d t}\norm{\Delta \varphi}^2_{L^2}   +  \delta \norm{\Delta^2 \varphi}^2_{L^2}  
		\\ 
	\notag	& = - \delta (\Delta^2{q_{\phi}},\Delta^2 \varphi)_{L^2} - \delta^{2/3}(\Delta^2{\xi_\phi},\Delta^2\varphi)_{L^2} - \delta(\Delta \varphi,\Delta^2 \varphi)_{L^2} - \delta(\Delta {q_{\phi}},\Delta^2 \varphi)_{L^2} - \delta^{2/3}(\Delta {\xi_\phi},\Delta^2 \varphi)_{L^2}
		\\ 
	\notag	& \quad - \delta b_1(v,\varphi, \Delta^2 \varphi) - \delta b_1(v,{q_{\phi}}, \Delta^2 \varphi) - \delta b_1({q_{u}},\varphi, \Delta^2 \varphi) - \delta b_1({q_{u}},{q_{\phi}}, \Delta^2 \varphi)
		\\
	\notag	& \quad - \delta^{2/3}b_1(v,{\xi_\phi},\Delta^2\varphi) - \delta^{2/3}b_1({q_{u}},{\xi_\phi},\Delta^2\varphi)   - \delta^{1/3} b_1({\zeta_{u}},{q_{\phi}}, \Delta^2 \varphi) - \delta^{1/3} b_1({\zeta_{u}},\varphi, \Delta^2 \varphi)
		\\
	\notag	& \quad  
		 + \left(\delta \Delta(\Phi^3)
-\Delta({\xi_\phi}^3),\Delta^2 \varphi\right)_{L^2}
\\
&\quad {+\delta(f_\phi(\delta t),\Delta^2\varphi)_{L^2}}.
\end{align}
As before, we will estimate each of the terms in the right-hand side of \eqref{eq:estphi2_asymp_mix}. The first five terms can be estimated directly by Cauchy-Schwarz and Young inequalities together with the fact that {$q_\phi,\xi_\phi\in\dot H^4$}. More precisely,
\begin{equation}\label{G1_H1_phi2_mix}
\begin{split}
\delta \left|(\Delta^2{q_{\phi}},\Delta^2 \varphi)_{L^2}\right| & \leq \epsilon\delta\|\Delta^2 \varphi\|_{L^2}^2 + C\delta\|{q_{\phi}}\|_{H^4}^2 
\\ 
\delta^{2/3} \left|(\Delta^2{\xi_\phi},\Delta^2 \varphi)_{L^2}\right| & \leq \epsilon\delta \|\Delta^2\varphi\|^2_{L^2}+C\delta^{1/3}\|{\xi_\phi}\|_{H^4}^2
\\
\delta\left|(\Delta \varphi,\Delta^2 \varphi)_{L^2}\right| & \leq \epsilon\delta \|\Delta^2\varphi\|_{L^2}^2+C\delta \|\Delta \varphi\|_{L^2}^2 
\\
\delta\left|(\Delta {q_{\phi}},\Delta^2 \varphi)_{L^2}\right| & \leq \epsilon\delta\|\Delta^2 \varphi\|_{L^2}^2 + C\delta\|{q_{\phi}}\|_{H^2}^2 
\\
\delta^{2/3} \left|(\Delta {\xi_\phi},\Delta^2 \varphi)_{L^2}\right| & \leq \epsilon\delta\|\Delta^2\varphi\|^2_{L^2}+C\delta^{1/3}\|{\xi_\phi}\|_{H^2}^2.
\end{split}
\end{equation}
The sixth through ninth terms can be estimated exactly as we have done in \eqref{G2_H2_phi} and \eqref{G3_H2_phi} in \Cref{thm:asymp-fluid}, namely, 
\begin{align}\label{G2_H2_phi2_mix}
\begin{split}
	\delta \left| b_1(v,\varphi, \Delta^2 \varphi) \right | 
	&\leq \epsilon\delta \|\Delta^2\varphi\|_{L^2}^2 + C\delta \|\nabla v\|_{L^2}^4+C\delta\|\Delta\varphi\|_{L^2}^4 \\
	\delta \left| b_1(v,{q_{\phi}}, \Delta^2 \varphi) \right|
	&\leq \epsilon\delta \| \Delta^2 \varphi \|_{L^2}^2 + C\delta \|v\|_{L^2}^2\|{q_{\phi}}\|_{H^3}^2\\
	\delta \left| b_1({q_{u}},\varphi, \Delta^2 \varphi) \right| 
	&\leq \epsilon\delta \| \Delta^2 \varphi \|_{L^2}^2 + C\delta \|{q_{u}}\|_{H^2}^2\|\nabla\varphi\|_{L^2}^2\\
	\delta \left| b_1({q_{u}},{q_{\phi}}, \Delta^2 \varphi) \right|
	&\leq \epsilon\delta \| \Delta^2 \varphi \|_{L^2}^2 + C\delta \|{q_{u}}\|_{H^2}^2\|{q_{\phi}}\|_{H^3}^2.
\end{split}
\end{align}
For terms ten and eleven, we use H\"older and Young inequalities to readily obtain
\begin{align}\label{G3_H2_phi_mix}
\begin{split}
\delta^{2/3}\left|b_1(v,{\xi_\phi},\Delta^2\varphi)\right| &\leq \epsilon\delta \|\Delta^2\varphi\|_{L^2}^2+C\delta^{1/3}\|{\xi_\phi}\|_{H^3}^2\|v\|_{L^2}^2
\\
\delta^{2/3}\left|b_1({q_{u}},{\xi_\phi},\Delta^2\varphi)\right|&\leq \epsilon\delta \|\Delta^2\varphi\|_{L^2}^2+C\delta^{1/3}\|{q_{u}}\|_{H^2}^2\|{\xi_\phi}\|_{H^3}^2.
\end{split}
\end{align}
For the twelfth term, we use integration by parts and H\"older and Young inequalities to get
\begin{align}\label{G4_H2_phi_mix}
\begin{split}
\delta^{1/3} \left|b_1({\zeta_{u}},{q_{\phi}}, \Delta^2 \varphi)\right|&=\delta^{1/3}\left|\int_{\TT}\Delta({\zeta_{u}}\cdot\nabla {q_{\phi}})\Delta \varphi\right| \leq \delta^{1/3}\|{\zeta_{u}}\cdot\nabla {q_{\phi}}\|_{H^2}\|\Delta\varphi\|_{L^2} \\
& \leq {C}\delta^{1/3}\|{\zeta_{u}}\|_{H^2}\|{q_{\phi}}\|_{H^3}\|\Delta \varphi\|_{L^2}
\\
&\leq C\delta^{1/3}\|\Delta\varphi\|_{L^2}^2+C\delta^{1/3}\|{\zeta_{u}}\|^2_{H^2}\|{q_{\phi}}\|^2_{H^3}.
\end{split}
\end{align}
For the thirteenth term, we use {\Cref{lem:transport_phase_high_order}} to directly obtain
\begin{equation}\label{G5_H2_phi_mix}
\delta^{1/3}
\left|
b_1({\zeta_{u}},\varphi,\Delta^2\varphi)
\right|
\leq
C\delta^{1/3}
\norm{{\zeta_{u}}}_{H^4}
\norm{\Delta \varphi}_{L^2}^2.
\end{equation}
For the {cubic} term we argue as follows. Using that
\begin{align*}
\Phi^3=(\varphi+&{q_{\phi}}+\delta^{-1/3}{\xi_\phi})^3\\
&=(\varphi+{q_{\phi}})^3+3\delta^{-1/3}(\varphi+{q_{\phi}})^2{\xi_\phi}+3\delta^{-2/3}(\varphi+{q_{\phi}}){\xi_\phi}^2+\delta^{-1}{\xi_\phi}^3.
\end{align*}
and we see that
\begin{align}\label{G6_H2_phi_mix}
\begin{split}
\left(\delta \Delta(\Phi^3)-\Delta({\xi_\phi}^3),\Delta^2 \varphi\right)_{L^2}&= \delta \left(\Delta(\varphi+{q_{\phi}})^3,\Delta^2\varphi\right)_{L^2}+3\delta^{2/3}\left(\Delta\left((\varphi+{q_{\phi}})^2{\xi_\phi}\right),\Delta^2 \varphi\right)_{L^2} \\
&\quad +3\delta^{1/3}\left(\Delta\left((\varphi+{q_{\phi}}){\xi_\phi}^2\right),\Delta^2\varphi\right)_{L^2}=: \sum_{i=1}^{3}M_i.
\end{split}
\end{align}
We proceed to estimate each term $M_i$. For the first one, we have 
\begin{align}\label{est_M1}
\begin{split}
	|M_1|&=\delta \left|\left(\Delta(\varphi+{q_{\phi}})^3,\Delta^2 \varphi\right)_{L^2}\right|
	\\
	& \leq {C}\delta \|\Delta^2\varphi\|_{L^2}\|\varphi+{q_{\phi}}\|_{H^2}^3 \\
	& \leq  \epsilon\delta \|\Delta^2\varphi\|^2_{L^2}+C\delta\|\varphi\|_{H^2}^6+C\delta \|{q_{\phi}}\|_{H^2}^6.
	\end{split}
\end{align}
Similarly, for $M_2$ we get
\begin{align}\label{est_M2}
\begin{split}
	|M_2|&=3\delta^{2/3} \left|\left(\Delta\left((\varphi+{q_{\phi}})^2{\xi_\phi}\right),\Delta^2 \varphi\right)_{L^2}\right|
	\\
	& \leq 3{C}\delta^{2/3} \|\Delta^2\varphi\|_{L^2}\|\varphi+{q_{\phi}}\|_{H^2}^2\|{\xi_\phi}\|_{H^2} \\
	& \leq  \epsilon\delta \|\Delta^2\varphi\|^2_{L^2}+C\delta^{1/3}\|\varphi\|_{H^2}^4\|{\xi_\phi}\|_{H^2}^2+C\delta^{1/3} \|{q_{\phi}}\|_{H^2}^4\|{\xi_\phi}\|_{{H^2}}^{2}.
	\end{split}
\end{align}
For $M_3$, applying \cite[Lemma~B.1]{hernandez2025global} with $k=2$
and using \eqref{eq:periodic-norm-equivalences}, we obtain
\begin{align}\label{est_M3}
M_3\leq C\delta^{1/3}\|\xi_\phi\|_{H^8}^2
\left(\|\Delta\varphi\|_{L^2}^2+\|q_\phi\|_{H^4}^2\right).
\end{align}
For the forcing term, Cauchy--Schwarz, Poincar\'e, and Young's
inequalities give, for almost every
$t\in(0,\min\{2,\delta^{-1}\mathcal T_\delta\})$,
\begin{align}\label{eq:asymp-mixed-background-force-phi}
 \delta\bigl|(f_\phi(\delta t),\Delta^2\varphi)_{L^2}\bigr|
 &\leq\delta\|f_\phi(\delta t)\|_{L^2}
                      \|\Delta^2\varphi\|_{L^2}\notag\\
 &\leq\epsilon\delta\|\Delta^2\varphi\|_{L^2}^2
       +C_\epsilon\delta M_f^2.
\end{align}

\textit{-- Step 4. Combining the estimates.} 
{
By combining \eqref{est2_asymp_mix} with \eqref{eq:estphi2_asymp_mix}
and using the estimates
\eqref{G1_H1_v_mix}--\eqref{eq:asymp-mixed-background-force-u}
together with
\eqref{G1_H1_phi2_mix}--\eqref{eq:asymp-mixed-background-force-phi},
and choosing $\epsilon$ sufficiently small, we obtain
}
\begin{equation} \label{eq:estH1H2_mix}
\begin{split}
	\frac{1}{2} \frac{\d}{\d t}\norm{\nabla v}^2_{L^2}  &+ \frac{1}{2}\frac{\d}{\d t}\norm{\Delta \varphi}^2_{L^2}  +  \frac{\delta\nu}{2} \norm{\Delta v}^2_{L^2}   +  \frac{\delta}{2} \norm{\Delta^2 \varphi}^2_{L^2} \\
	&\leq C(\delta+\delta^{1/3}) + C(\delta+\delta^{1/3})\|\nabla v\|_{L^2}^2 + C(\delta+\delta^{1/3})\|\Delta \varphi\|_{L^2}^2 \\
	&\quad + C\delta (\|v\|_{H^1}^4 +\|v\|_{H^1}^6 + \|\varphi\|_{H^2}^4 + \|\varphi\|_{H^2}^6 + \|\varphi\|_{H^2}^8 )+C\delta^{1/3}\|\varphi\|_{H^2}^4,
\end{split}
\end{equation}
for some constant $C>0$ that may depend on $r$ {and $M_f$} but not on $\delta\in(0,1)$. Here we have used repeatedly {\eqref{eq:periodic-norm-equivalences}}. Finally, using that $\delta\in(0,1)$, we obtain from \eqref{eq:estH1H2_mix}
\begin{align*}
	\frac{\d}{\d t}\left(\norm{\nabla v}^2_{L^2}+\norm{\Delta \varphi}^2_{L^2}\right) \leq  C\left(\norm{\nabla v}^2_{L^2}+\norm{\Delta \varphi}_{L^2}^2\right)+C\delta^{1/3} \left(1+\norm{v}_{H^1}^6+\norm{\varphi}^8_{H^2}\right). 
\end{align*}
By Gr\"onwall's inequality and {\eqref{eq:periodic-norm-equivalences}}, we have the following estimate $\text{for }t<  \min\{2, \delta^{-1}{\mathcal T_\delta}\}$
\begin{align}\label{grn_mix}
	\norm{v(t)}^2_{H^1}+\norm{\varphi(t)}^2_{H^2}\leq C\delta^{1/3}t +C \delta^{1/3}\int_{0}^{t}\left(\norm{v(s)}^6_{H^1}+\norm{\varphi(s)}_{H^2}^8\right) \, ds.
\end{align}

\textit{-- Step 5. Asymptotic limit and conclusion.} The last step is totally analogous to that appearing in the proof of \Cref{thm:asymp-fluid}, we just need to use \eqref{grn_mix} instead of \eqref{grn} and follow the same arguments. For brevity, we skip the details. This ends the proof.
\end{proof}

}

\section{Saturating subspaces}\label{sec:saturation}
We now study the algebraic properties of the nonlinear interactions
appearing in \Cref{thm:asymp-fluid,thm:asymp-mixed}.
Starting from a fixed finite set of Fourier modes, we prove that
iterating the fluid and transport nonlinearities generates
a dense subspace of $\mathbb V^1$.
The construction of forcing terms using these properties
is carried out in the next section.

For $l\in\mathbb Z^3\setminus\{0\}$, we denote
\begin{equation*}
\mathcal C_l
:=
\operatorname{span}\{\cos(l\cdot x),\sin(l\cdot x)\}.
\end{equation*}
Notice that $\mathcal C_l=\mathcal C_{-l}$. If
$\mathcal K\subset\mathbb Z^3\setminus\{0\}$ is finite, we set
\begin{equation*}
\mathcal C_{\mathcal K}
:=
\operatorname{span}
\left\{
\cos(k\cdot x),\,
\sin(k\cdot x)
:\ k\in\mathcal K
\right\}.
\end{equation*}
We take
\begin{equation}\label{K_sets}
\mathcal K_1
=
{\{e_2,e_1+e_3,e_2+e_3\}},
\qquad
\mathcal K_2
=
\{e_3\}.
\end{equation}
For $k\in\mathbb Z^3\setminus\{0\}$, let $P_k$ denote the orthogonal
projection of $\mathbb R^3$ onto $k^\perp$  given by
$
P_k a=a-\frac{a\cdot k}{|k|^2}k,\, a\in\mathbb R^3.
$
 We also introduce the family
\begin{equation}\label{eq:G1_def}
\mathscr G_1
:=
\left\{
(P_k e_3)\cos(x\cdot k),\,
(P_k e_3)\sin(x\cdot k)
:
k\in\mathcal K_1
\right\}.
\end{equation}
We define the initial controlled space by
\begin{equation}\label{eq:H1_initial}
\mathcal H_1
:=
\operatorname{span}\mathscr G_1
\subset \mathbf H^5,
\end{equation}
and
\begin{equation}\label{eq:H2_initial}
\mathcal H_2
:=
\mathcal C_{\mathcal K_2}
=
\mathcal C_{e_3}
=
\operatorname{span}\{\cos x_3,\sin x_3\}
\subset \dot H^6.
\end{equation}
Thus,
\begin{equation}\label{h}
\mathcal H
:=
\mathcal H_1\times\mathcal H_2
\subset \mathbb V^5.
\end{equation}
Notice that, since
$\Pi(e_3\cos(k\cdot x))=P_ke_3\cos(k\cdot x)$
and
$\Pi(e_3\sin(k\cdot x))=P_ke_3\sin(k\cdot x)$,
we have
\[
\mathcal H
=\{(\Pi(\eta_ue_3),\eta_\phi):
(\eta_u,\eta_\phi)\in
\mathcal C_{\mathcal K_1}\times\mathcal C_{\mathcal K_2}\}.
\]
Thus the projected velocity controls correspond to
forces acting only in the third component of the
original momentum equation of \eqref{eq:chn_modified_intro}.

Let $\mathcal E=\mathcal E_1\times\mathcal E_2$ be a finite-dimensional
subspace of $\mathbb V^5$. We first define the velocity saturation by
\begin{equation}\label{eq:F1_def}
\mathcal F_1(\mathcal E)
:=
\operatorname{span}
\left\{
\zeta_0-\sum_{i=1}^{m}B_0(\zeta_i,\zeta_i)
:\ 
\zeta_i\in\mathcal E_1,\ i\in\llbracket 0,m\rrbracket,\ m\ge 1
\right\}.
\end{equation}
We define the phase saturation by
\begin{equation}\label{eq:F2_mixed_def}
\mathcal F_2(\mathcal E)
:=
\operatorname{span}
\left\{
\xi_0+\sum_{j=1}^{m}B_1(\zeta_j,\xi_j)
:\ 
\xi_j\in\mathcal E_2,\ j\in\llbracket 0,m\rrbracket,\ 
\zeta_j\in\mathscr G_1,\ j\in\llbracket 1,m\rrbracket,\ 
m\ge 1
\right\}.
\end{equation}
We denote by $\mathcal F(\mathcal E)$ the product
\begin{equation*}
\mathcal F(\mathcal E)
:=
\mathcal F_1(\mathcal E)\times\mathcal F_2(\mathcal E).
\end{equation*}
Starting from $\mathcal H(0)=\mathcal H$, we define a sequence
$\{\mathcal H(j)\}_{j\ge0}$ of finite-dimensional subspaces of
$\mathbf H^5\times\dot H^6$ by
\begin{equation}\label{eq:Hj_def}
	\mathcal H(j)=\mathcal F(\mathcal H(j-1)),
	\qquad j\ge1.
\end{equation}
By definition, $\mathcal E_i\subset\mathcal F_i(\mathcal E)$
for $i=1,2$, and hence
$\mathcal H(j-1)\subset\mathcal H(j)$ for every $j\ge1$.
Equivalently, writing
$\mathcal H(j)=\mathcal H_1(j)\times\mathcal H_2(j)$,
we have
\begin{equation*}
	\mathcal H_1(j)=\mathcal F_1(\mathcal H(j-1)),
	\qquad
	\mathcal H_2(j)=\mathcal F_2(\mathcal H(j-1)).
\end{equation*}
\begin{rmk}\label{rem:trigonometric-saturation}
The components of every element of $\mathcal H(j)$ are trigonometric
polynomials and hence smooth. Indeed, this holds for $\mathcal H(0)$
and follows by induction from \eqref{eq:Hj_def}, since differentiation,
multiplication, and the Leray projector preserve trigonometric polynomials.
\end{rmk}
We then set
\begin{equation}\label{eq:Hinf_def}
\mathcal H(\infty)
:=
\bigcup_{j=0}^{\infty}\mathcal H(j).
\end{equation}
We also write
\begin{equation*}
\mathcal H_i(\infty)
:=
\bigcup_{j=0}^{\infty}\mathcal H_i(j),
\qquad i=1,2.
\end{equation*}

\begin{lemma}\label{lem:safe_modes}
For every $\zeta\in\mathscr G_1$, we have $B_0(\zeta,\zeta)=0$.
\end{lemma}

\begin{proof}
Let $\zeta=(P_k e_3)\cos(x\cdot k)$ for some $k\in\mathcal K_1$. The case of the sine mode is treated in the same way. Since $P_k e_3\in k^\perp$, we have
\begin{equation*}
(P_k e_3)\cdot k=0.
\end{equation*}
Therefore,
\begin{equation*}
(\zeta\cdot\nabla)\zeta
=
\bigl((P_k e_3)\cdot k\bigr)
(P_k e_3)\cos(x\cdot k)(-\sin(x\cdot k))
=
0.
\end{equation*}
Hence $B_0(\zeta,\zeta)=0$.
\end{proof}
%%%%%%%%%%%%%%%%
Given a finite-dimensional subspace $E\subset\dot H^6$, we define the following subspace
\begin{equation}\label{def_transport_space}
\mathcal T(E)
:=
E+\operatorname{span}
\left\{
B_1(\zeta,\xi):\ \zeta\in\mathscr G_1,\ \xi\in E
\right\}.
\end{equation}

\begin{lemma}\label{lem:transport_modes}
Let $E\subset \dot H^6$ be a finite-dimensional subspace. Assume that, for some
$l\in\mathbb Z^3\setminus\{0\}$, we have $\mathcal C_l\subset E$. Let
$k\in\mathcal K_1$ be such that $(P_k e_3)\cdot l\neq0$. Then $l+k$ and
$l-k$ are nonzero, and
\begin{equation*}
\mathcal C_{l+k}\subset\mathcal T(E),
\qquad
\mathcal C_{l-k}\subset\mathcal T(E).
\end{equation*}
\end{lemma}
\begin{proof}
	Set $a_k := P_k e_3$. Since $a_k \cdot k = 0$, the assumption $a_k \cdot l \neq 0$ implies that $l \neq k$ and $l \neq -k$. Consequently, $l - k \neq 0$ and $l + k \neq 0$.
	
	By definition, $\mathcal T(E)$ contains $B_1(\zeta, \xi)$ for all choices of $\zeta \in \mathscr G_1$ and $\xi \in E$. Since $\mathcal C_l = \operatorname{span}\{\cos(l\cdot x), \sin(l\cdot x)\} \subset E$ and $a_k\cos(k\cdot x), a_k\sin(k\cdot x) \in \mathscr G_1$, evaluating $B_1$ on these generator pairs yields:
	\begin{align*}
		B_1(a_k\cos(k\cdot x), \sin(l\cdot x)) &= (a_k \cdot l)\cos(k\cdot x)\cos(l\cdot x) \in \mathcal T(E), \\
		B_1(a_k\sin(k\cdot x), \sin(l\cdot x)) &= (a_k \cdot l)\sin(k\cdot x)\cos(l\cdot x) \in \mathcal T(E), \\
		B_1(a_k\cos(k\cdot x), \cos(l\cdot x)) &= -(a_k \cdot l)\cos(k\cdot x)\sin(l\cdot x) \in \mathcal T(E), \\
		B_1(a_k\sin(k\cdot x), \cos(l\cdot x)) &= -(a_k \cdot l)\sin(k\cdot x)\sin(l\cdot x) \in \mathcal T(E).
	\end{align*}
	Because $a_k \cdot l \neq 0$, we can rescale by $(a_k \cdot l)^{-1}$. Thus, the linear subspace $\mathcal T(E)$ contains the four product terms:
	\begin{equation*}
		\cos(k\cdot x)\cos(l\cdot x), \quad \sin(k\cdot x)\cos(l\cdot x), \quad \cos(k\cdot x)\sin(l\cdot x), \quad \sin(k\cdot x)\sin(l\cdot x).
	\end{equation*}
	Finally, using the trigonometric product identities:
	\begin{align*}
		\cos((l+k)\cdot x) &= \cos(k\cdot x)\cos(l\cdot x) - \sin(k\cdot x)\sin(l\cdot x), \\
		\cos((l-k)\cdot x) &= \cos(k\cdot x)\cos(l\cdot x) + \sin(k\cdot x)\sin(l\cdot x), \\
		\sin((l+k)\cdot x) &= \sin(k\cdot x)\cos(l\cdot x) + \cos(k\cdot x)\sin(l\cdot x), \\
		\sin((l-k)\cdot x) &= \sin(l\cdot x)\cos(k\cdot x) - \cos(l\cdot x)\sin(k\cdot x),
	\end{align*}
	we observe that $\cos((l+k)\cdot x)$, $\sin((l+k)\cdot x)$, $\cos((l-k)\cdot x)$, and $\sin((l-k)\cdot x)$ are all linear combinations of the four product terms in $\mathcal T(E)$. Therefore, $\mathcal C_{l+k} \subset \mathcal T(E)$ and $\mathcal C_{l-k} \subset \mathcal T(E)$.
\end{proof}

\begin{lemma}\label{lem:mixed_lattice}
Let $\mathcal K_1$ be as in \eqref{K_sets}. Let $E_0=\mathcal C_{e_3}$ and
define recursively
\begin{equation*}
E_{n+1}:=\mathcal T(E_n),
\qquad n\ge0.
\end{equation*}
Then, for every $l\in\mathbb Z^3\setminus\{0\}$, there exists $n\ge0$ such
that $\mathcal C_l\subset E_n.$
\end{lemma}

\begin{proof}
By \Cref{lem:transport_modes}, we shall use the following rule: if
$\mathcal C_l\subset E_n$ for some $n\ge0$ and $k\in\mathcal K_1$ satisfies
$(P_k e_3)\cdot l\neq0$, then
\begin{equation}\label{eq:propagation_rule}
\mathcal C_{l+k}\subset E_{n+1},
\qquad
\mathcal C_{l-k}\subset E_{n+1}.
\end{equation}
We now apply this rule in several steps.

\smallskip
\textit{-- Step 1. Generation of the line $l_1=0$, $l_3=1$.}
The first observation is that, as long as the vertical component of the
frequency is nonzero, the direction $e_2$ allows us to move in the
second variable without changing the other two components.

Starting from $\mathcal C_{e_3}\subset E_0$, we move in the
$e_2$-direction. We prove by induction that, for every $m\in\mathbb Z$,
there exists $n_m\ge0$ such that
\begin{equation}\label{eq:e2_line}
\mathcal C_{(0,m,1)}\subset E_{n_m}.
\end{equation}
The case $m=0$ is the initial inclusion. If
$\mathcal C_{(0,m,1)}\subset E_{n_m}$, then, using $P_{e_2}e_3=e_3$,
we have
\begin{equation*}
(P_{e_2}e_3)\cdot(0,m,1)=1\neq0.
\end{equation*}
Therefore, by \eqref{eq:propagation_rule} with $k=e_2$,
\begin{equation*}
\mathcal C_{(0,m+1,1)}\subset E_{n_m+1},
\qquad
\mathcal C_{(0,m-1,1)}\subset E_{n_m+1}.
\end{equation*}
This proves that, for every $m\in\mathbb Z$, the block
$\mathcal C_{(0,m,1)}$ is generated.

%\smallskip
%\textit{-- Step 2. Generation of one representative on each vertical level.}
%The previous step only moves frequencies inside a fixed vertical level. To
%reach other vertical levels, we use the diagonal direction $e_1+e_3$, which
%changes the first and third components at the same time.
%
%Since $P_{e_1+e_3}e_3=\frac12(-e_1+e_3)$, we have
%\begin{equation*}
%(P_{e_1+e_3}e_3)\cdot l
%=
%\frac12(l_3-l_1).
%\end{equation*}
%We prove that, for every $r\in\mathbb Z$, there exists $n_r\ge0$ such that
%\begin{equation}\label{eq:vertical_representatives}
%\mathcal C_{(r-1,0,r)}\subset E_{n_r}.
%\end{equation}
%The case $r=1$ is the initial inclusion $\mathcal C_{e_3}\subset E_0$.
%
%First consider $r\ge1$. Suppose that
%$\mathcal C_{(r-1,0,r)}\subset E_{n_r}$. Then
%\begin{equation*}
%(P_{e_1+e_3}e_3)\cdot(r-1,0,r)
%=
%\frac12\bigl(r-(r-1)\bigr)
%=
%\frac12
%\neq0.
%\end{equation*}
%Hence, by \eqref{eq:propagation_rule} with $k=e_1+e_3$,
%%\begin{equation*}
%$\mathcal C_{(r,0,r+1)}\subset E_{n_r+1}.$
%%\end{equation*}
%This proves \eqref{eq:vertical_representatives} for every $r\ge1$.
%
%{
%Now consider $r\le1$. Suppose again that
%$\mathcal C_{(r-1,0,r)}\subset E_{n_r}$. The same computation gives
%\begin{equation*}
%(P_{e_1+e_3}e_3)\cdot(r-1,0,r)
%=
%\frac12
%\neq0.
%\end{equation*}
%Therefore, by \eqref{eq:propagation_rule} with $k=e_1+e_3$,
%%\begin{equation*}
%$\mathcal C_{(r-2,0,r-1)}\subset E_{n_r+1}.$
%%\end{equation*}
%This proves \eqref{eq:vertical_representatives} for every $r\le1$.}

\smallskip
\textit{-- Step 2. Generation of one horizontal line on each vertical level.}
The previous step only moves frequencies inside the vertical level
$l_3=1$. To reach other vertical levels, we use the diagonal direction
$e_1+e_3$, which changes the first and third components at the same time.

Fix $n\in\mathbb Z$. We prove that, for every $r\in\mathbb Z$,
there exists $N_{r,n}\ge0$ such that
\begin{equation}\label{eq:vertical_representatives}
\mathcal C_{(r-1,n,r)}\subset E_{N_{r,n}}.
\end{equation}
The case $r=1$ follows from \eqref{eq:e2_line}.
Suppose that \eqref{eq:vertical_representatives} holds
for some $r\in\mathbb Z$. Since
\[
(P_{e_1+e_3}e_3)\cdot(r-1,n,r)=\frac12\neq0,
\]
applying \eqref{eq:propagation_rule} with $k=e_1+e_3$ gives
\[
\mathcal C_{(r,n,r+1)}\subset E_{N_{r,n}+1},
\qquad
\mathcal C_{(r-2,n,r-1)}\subset E_{N_{r,n}+1}.
\]
Thus the property propagates from $r$ to both $r+1$ and $r-1$.
Starting from $r=1$, induction in both directions proves
\eqref{eq:vertical_representatives} for every $r,n\in\mathbb Z$.

%
%Now consider $r\le1$. Suppose again that
%$\mathcal C_{(r-1,0,r)}\subset E_{n_r}$. The same computation gives
%\begin{equation*}
%	(P_{e_1+e_3}e_3)\cdot(r-1,0,r)
%	=
%	\frac12
%	\neq0.
%\end{equation*}
%Therefore, by \eqref{eq:propagation_rule} with $k=e_1+e_3$,
%\begin{equation*}
%	\mathcal C_{(r-2,0,r-1)}\subset E_{n_r+1}.
%\end{equation*}
%This proves \eqref{eq:vertical_representatives1} for every $r\le1$.}

\smallskip
\textit{-- Step 3. Generation of all nonzero vertical levels.}
We first generate the frequencies $(m,0,r)$ on each level
$r\neq0$. We then use $e_2$ to vary the second component.

Let $r\in\mathbb Z\setminus\{0\}$ and $m\in\mathbb Z$.
By Step 2, with vertical level $m+1$ and second component
$m+1-r$, there exists $N\ge0$ such that
$
\mathcal C_{(m,m+1-r,m+1)}\subset E_N.
$
We now move in the direction $e_2+e_3$. At each frequency along
this sequence,
\[
(P_{e_2+e_3}e_3)\cdot(m,m+1-r+j,m+1+j)
=\frac12r\neq0,\qquad j\in\mathbb Z.
\]
Moreover,
\[
(m,m+1-r,m+1)+(r-m-1)(e_2+e_3)=(m,0,r).
\]
Thus, applying \eqref{eq:propagation_rule} $|r-m-1|$ times,
with $k=e_2+e_3$ and the corresponding sign, gives
\[
\mathcal C_{(m,0,r)}\subset E_{N+|r-m-1|}.
\]
Now fix such an $m\in\mathbb Z$ and move in the $e_2$-direction. If
$\mathcal C_{(m,n,r)}\subset E_N$, then
\begin{equation*}
(P_{e_2}e_3)\cdot(m,n,r)=r\neq0.
\end{equation*}
Therefore, by \eqref{eq:propagation_rule} with $k=e_2$,
%\begin{equation*}
$\mathcal C_{(m,n+1,r)}\subset E_{N+1},
\,
\mathcal C_{(m,n-1,r)}\subset E_{N+1}.$
%\end{equation*}
Consequently, for every $r\in\mathbb Z\setminus\{0\}$ and every
$m,n\in\mathbb Z$, there exists $N\ge0$ such that
$\mathcal C_{(m,n,r)}\subset E_N.$

\medskip
\textit{-- Step 4. Generation of the zero vertical level.}
The only level not covered by Step 3 is $l_3=0$. We reach this level from the
already generated level $l_3=1$, using the diagonal directions $e_1+e_3$ and
$e_2+e_3$.

Let $(m,n,0)\in\mathbb Z^3\setminus\{0\}$. If $m\neq0$, then by Step 3
$\mathcal C_{(m+1,n,1)}$ has already been generated. Taking $k=e_1+e_3$, we
have 
%\begin{equation*}
$(m+1,n,1)-(e_1+e_3)=(m,n,0).$
%\end{equation*}
Moreover,
\begin{equation*}
(P_{e_1+e_3}e_3)\cdot(m+1,n,1)
=
\frac12\bigl(1-(m+1)\bigr)
=
-\frac{m}{2}
\neq0.
\end{equation*}
Thus \eqref{eq:propagation_rule} gives that $\mathcal C_{(m,n,0)}$ is
generated.

It remains to consider the case $m=0$ and $n\neq0$. By Step 3, the block
$\mathcal C_{(0,n+1,1)}$ has already been generated. Taking $k=e_2+e_3$, we
have $(0,n+1,1)-(e_2+e_3)=(0,n,0).$ 
Since $P_{e_2+e_3}e_3 = \frac12(-e_2+e_3)$, we obtain
\begin{equation*}
(P_{e_2+e_3}e_3)\cdot(0,n+1,1)
=
\frac12\bigl(1-(n+1)\bigr)
=
-\frac{n}{2}
\neq0.
\end{equation*}
Hence \eqref{eq:propagation_rule} gives that $\mathcal C_{(0,n,0)}$ is
generated.

\smallskip
\textit{-- Step 5. Conclusion.}
By Step 3, all frequencies with nonzero vertical component have been
generated: for every $l=(l_1,l_2,l_3)\in\mathbb Z^3$ with $l_3\neq0$, there
exists $n\ge0$ such that $\mathcal C_l\subset E_n$. By Step 4, the same
conclusion holds for every nonzero frequency with $l_3=0$. Therefore, for
every $l\in\mathbb Z^3\setminus\{0\}$, there exists $n\ge0$ such that
$\mathcal C_l\subset E_n$. This proves the claim.
\end{proof}
\begin{defi}
A finite-dimensional subspace $\mathcal H$ of $\mathbb V^5$ is saturating in
$\mathbb V^1$ if $\mathcal H(\infty)$ is dense in $\mathbb V^1$.
\end{defi}
\begin{proposition}\label{prop_sat_phase_mixed}
The space $\mathcal H_2(\infty)$ is dense in $\dot H^2$.
\end{proposition}
\begin{proof}
	We first show that the auxiliary spaces $E_n$ introduced in
	\Cref{lem:mixed_lattice} coincide with the phase saturating
	sequence. Indeed,
	$E_0=\mathcal C_{e_3}=\mathcal H_2(0)$.
	Moreover, by \eqref{def_transport_space} and
	\eqref{eq:F2_mixed_def}, we have
	\[
	\mathcal F_2(\mathcal H(n))
	=\mathcal T(\mathcal H_2(n)).
	\]
	Thus, if $E_n=\mathcal H_2(n)$ for some $n\ge0$, then
	\[
	E_{n+1}=\mathcal T(E_n)
	=\mathcal F_2(\mathcal H(n))
	=\mathcal H_2(n+1).
	\]
	By induction, $E_n=\mathcal H_2(n)$ for every $n\ge0$.
	
	By \Cref{lem:mixed_lattice}, for every
	$l\in\mathbb Z^3\setminus\{0\}$, there exists $n\ge0$
	such that $\mathcal C_l\subset E_n$. Therefore,
	\[
	\mathcal C_l\subset\mathcal H_2(\infty),
	\qquad l\in\mathbb Z^3\setminus\{0\}.
	\]
	Since the linear span of these trigonometric modes
$	\{	\cos(l\cdot x), \,
		\sin(l\cdot x), \,
		l\in\mathbb Z^3\setminus\{0\}\}$
	is dense in $\dot H^2(\mathbb T^3)$, the conclusion follows.
\end{proof}
%
%\begin{proof}
%We first observe that the auxiliary spaces $E_n$ introduced in
%\Cref{lem:mixed_lattice} are contained in the phase saturating sequence.
%Indeed,
%%\begin{equation*}
%$E_0=\mathcal C_{e_3}=\mathcal H_2(0).$
%%\end{equation*}
%Moreover, assume that $E_n\subset\mathcal H_2(n)$ for some $n\ge0$. Then, by
%the definition of $\mathcal T(E_n)$ and the definition
%\eqref{eq:F2_mixed_def} of $\mathcal F_2$, we have
%\begin{equation*}
%E_{n+1}
%=
%\mathcal T(E_n)
%\subset
%\mathcal F_2(\mathcal H(n))
%=
%\mathcal H_2(n+1).
%\end{equation*}
%Thus $E_n\subset\mathcal H_2(n)$ for every $n\ge0$. {\color{red}subset or equal?}
%
%By \Cref{lem:mixed_lattice}, for every $l\in\mathbb Z^3\setminus\{0\}$,
%there exists $n\ge0$ such that $\mathcal C_l\subset E_n$. Therefore,
%\begin{equation*}
%\mathcal C_l\subset \mathcal H_2(\infty),
%\qquad
%l\in\mathbb Z^3\setminus\{0\}.
%\end{equation*}
%Equivalently, $\mathcal H_2(\infty)$ contains all trigonometric modes
%\begin{equation*}
%\cos(l\cdot x),
%\qquad
%\sin(l\cdot x),
%\qquad
%l\in\mathbb Z^3\setminus\{0\}.
%\end{equation*}
%Their linear span is dense in $\dot H^2(\mathbb T^3)$.
%\end{proof}
%Using\cite[Proposition 4.7]{N15}, we have that $\mathcal H_1(\infty)$ is dense in $\mathbf H^1$.
\begin{proposition}\label{prop_sat_phase}
	The space $\mathcal H_1(\infty)$ is dense in $\mathbf H^1.$
	\end{proposition}
			\begin{proof}
			Since $B_0$ is bilinear, we have
			\[
			B_0(v,w)+B_0(w,v)
			=B_0(v+w,v+w)-B_0(v,v)-B_0(w,w).
			\]
			Thus, by the definition of $\mathcal F_1$,
			\begin{equation}\label{eq:polarization-fluid}
				B_0(v,w)+B_0(w,v)\in\mathcal H_1(j+1),
				\qquad v,w\in\mathcal H_1(j).
			\end{equation}
			Let $m,n\in\mathbb Z^3\setminus\{0\}$ with $m-n\neq0$,
			and let $a\in m^\perp$, $b\in n^\perp$.
			Using trigonometric identities and the Leray projector,
			we obtain
			\begin{align}
				&P_{m-n}\bigl((a\cdot n)b-(b\cdot m)a\bigr)
				\cos((m-n)\cdot x)\notag\\
				&\quad=B_0(a\cos(m\cdot x),b\sin(n\cdot x))
				+B_0(b\sin(n\cdot x),a\cos(m\cdot x))\notag\\
				&\qquad-B_0(b\cos(n\cdot x),a\sin(m\cdot x))
				-B_0(a\sin(m\cdot x),b\cos(n\cdot x)).
				\label{eq:fluid-cos-interaction}
			\end{align}
			Similarly,
			\begin{align}
				&P_{m-n}\bigl((a\cdot n)b-(b\cdot m)a\bigr)
				\sin((m-n)\cdot x)\notag\\
				&\quad=B_0(a\sin(m\cdot x),b\sin(n\cdot x))
				+B_0(b\sin(n\cdot x),a\sin(m\cdot x))\notag\\
				&\qquad+B_0(a\cos(m\cdot x),b\cos(n\cdot x))
				+B_0(b\cos(n\cdot x),a\cos(m\cdot x)).
				\label{eq:fluid-sin-interaction}
			\end{align}
			In particular, if the sine and cosine modes associated
			with $(m,a)$ and $(n,b)$ belong to $\mathcal H_1(j)$,
			then the two functions on the left-hand sides of
			\eqref{eq:fluid-cos-interaction} and
			\eqref{eq:fluid-sin-interaction} belong to
			$\mathcal H_1(j+1)$.

			\smallskip
			\textit{-- Step 1. Generation of the frequency $e_1$.}
			Set
			\[
			a:=P_{e_1+e_3}e_3=\frac12(-e_1+e_3),\qquad
			b:=P_{e_2+e_3}e_3=\frac12(-e_2+e_3),\qquad
			p:=e_1+e_2+e_3.
			\]
			Taking $m=e_2+e_3$, $n=e_2$, and velocity directions
			$b,e_3$ in \eqref{eq:fluid-cos-interaction} and
			\eqref{eq:fluid-sin-interaction}, we have
			\[
			P_{e_3}\bigl((b\cdot e_2)e_3
			 -(e_3\cdot(e_2+e_3))b\bigr)
			=P_{e_3}\left(-\frac12e_3-b\right)
			=\frac12e_2.
			\]
			Thus
			$e_2\cos x_3,e_2\sin x_3\in\mathcal H_1(1)$.
			Next take $m=e_1+e_3$, $n=-e_2$, and velocity directions
			$a,e_3$. The profiles at frequency $-e_2$ belong to
			$\mathcal H_1(0)$ since $\cos(-x_2)=\cos x_2$ and
			$\sin(-x_2)=-\sin x_2$. Since $a\cdot p=0$,
			\[
			P_p\bigl((a\cdot(-e_2))e_3
			 -(e_3\cdot(e_1+e_3))a\bigr)
			=P_p(-a)=-a.
			\]
			Therefore,
			$a\cos(p\cdot x),a\sin(p\cdot x)\in\mathcal H_1(1)$.

			We now use \eqref{eq:fluid-cos-interaction} and
			\eqref{eq:fluid-sin-interaction} with
			$m=e_1+e_3$, $n=e_3$, and velocity directions $a,e_2$.
			The coefficient is
			\[
			P_{e_1}\bigl((a\cdot e_3)e_2
			 -(e_2\cdot(e_1+e_3))a\bigr)
			=\frac12e_2.
			\]
			With $m=p$, $n=e_2+e_3$, and velocity directions $a,b$,
			the coefficient is
			\[
			P_{e_1}\bigl((a\cdot(e_2+e_3))b-(b\cdot p)a\bigr)
			=\frac12b,
			\]
			because $a\cdot(e_2+e_3)=1/2$ and $b\cdot p=0$.
			Since
			$e_3=2(\frac12e_2)+4(\frac12b)$, linearity gives
			\[
			e_2\cos x_1,\quad e_2\sin x_1,\quad
			e_3\cos x_1,\quad e_3\sin x_1
			\in\mathcal H_1(2).
			\]

			\smallskip
			\textit{-- Step 2. Generation of the frequencies $e_3$ and $e_2$.}
			Take $m=e_1+e_3$, $n=e_1$, and velocity directions
			$a,e_3$. The required profiles belong to $\mathcal H_1(2)$
			by Step 1 and the inclusion
			$\mathcal H_1(0)\subset\mathcal H_1(2)$. We have
			\[
			P_{e_3}\bigl((a\cdot e_1)e_3
			 -(e_3\cdot(e_1+e_3))a\bigr)
			=P_{e_3}\left(-\frac12e_3-a\right)
			=\frac12e_1.
			\]
			By \eqref{eq:fluid-cos-interaction} and
			\eqref{eq:fluid-sin-interaction},
$
			e_1\cos x_3,\, e_1\sin x_3
			\in\mathcal H_1(3).
			$
			The profiles $e_2\cos x_3,e_2\sin x_3$ were obtained in Step 1.

			Finally, take $m=e_2+e_3$, $n=e_3$, and velocity directions
			$b,e_1$. These profiles belong to $\mathcal H_1(3)$, and
			\[
			P_{e_2}\bigl((b\cdot e_3)e_1
			 -(e_1\cdot(e_2+e_3))b\bigr)
			=\frac12e_1.
			\]
			Hence
			$e_1\cos x_2,e_1\sin x_2\in\mathcal H_1(4)$.
			The profiles $e_3\cos x_2,e_3\sin x_2$ belong to
			$\mathcal H_1(0)$.

			\smallskip
			\textit{-- Step 3. Conclusion.}
		Let $\mathcal G$ be the span of all divergence-free
			sine and cosine modes associated with the frequencies
			$e_1,e_2,e_3$. By Steps 1 and 2,
			$
			\mathcal G\subset\mathcal H_1({4}).
			$
			Since $\{e_1,e_2,e_3\}$ generates $\mathbb Z^3$,
			\cite[Theorem~4.1]{N15} implies that $\mathcal G$
			is saturating in $\mathbf H^1$ 	(see \Cref{ner1}).
			By the monotonicity of the saturation operation,
			$
			\mathcal G(\infty)\subset\mathcal H_1(\infty).
			$
			Therefore $\mathcal H_1(\infty)$ is dense in
			$\mathbf H^1$.
			\end{proof}
	\begin{proposition}\label{prop_sat}
The space $\mathcal H$ defined in \eqref{h} is saturating in $\mathbb V^1$.
\end{proposition}
\begin{proof}
Since both saturating sequences are non-decreasing,
we have
\[
\mathcal H(\infty)
=\mathcal H_1(\infty)\times\mathcal H_2(\infty).
\] It is enough to prove that $\mathcal H_1(\infty)$ and $\mathcal H_2(\infty)$
are dense in $\mathbf H^1$ and $\dot H^2$, respectively.
The density of $\mathcal H_2(\infty)$ in $\dot H^2$ follows from
\Cref{prop_sat_phase_mixed} and the density of $\mathcal H_1(\infty)$ in $\mathbf H^1$ follows from
\Cref{prop_sat_phase} .  Therefore $\mathcal H(\infty)$ is dense in
$\mathbb V^1$.
\end{proof}
\begin{rmk}\label{rem:three-frequency-saturation}
The saturation property in \Cref{prop_sat_phase} is independent of the
phase equation and gives a six-dimensional saturating space for the
three-dimensional Navier--Stokes equations. In
\cite[Proposition~4.7]{N15}, saturation with a force of the form
$\Pi(\eta_u e_3)$ is obtained using four Fourier frequencies.
Our construction reduces this number to three, while retaining the
same restriction to one fixed component of the force. The
six-dimensional example in \cite[Proposition~4.9]{N15} also uses three
frequencies, but without the restriction to a single component on the
force.

For CHNS, \Cref{prop_sat_phase_mixed} shows that the fixed fluid profiles
in $\mathscr G_1$ generate all phase frequencies from the single phase
frequency $e_3$ through the transport nonlinearity. Thus the same three frequencies for the fluid
suffice both for fluid saturation and for generating the missing phase
directions.
\end{rmk}
\section{{Accessibility in small-time}}\label{sec:accessibility}
In this section we prove \Cref{thm:global-small-time-control}.
We combine the asymptotic limits of \Cref{sec:asymp_limits} with the
saturation properties of \Cref{sec:saturation} to construct corrections
in $\mathcal H$ that bring any initial state $z_0\in\mathbb V^1$ into
the unit ball in arbitrarily small time, while keeping $f$ prescribed.

We follow the construction developed for the Cahn--Hilliard equation
in \cite{hernandez2025global}, with notation and several arguments
inspired by \cite{Ner21}. The separate construction of fluid and phase
displacements is also inspired by the treatment of the Boussinesq
equations in \cite[Section~5.3]{GHHM18}.
The proof combines the following displacements, which are constructed
in the subsections below.
\begin{itemize}
\item \emph{Adding an available direction.}
For $h\in\mathcal H$, apply $\eta=\delta^{-1}h$ for a time
$\delta$: the state moves approximately from $z_0$ to $z_0+h$
as $\delta\to0^+$. This construction is given in
\Cref{lem:direct_construction}.

\item \emph{Moving the fluid.}
For an available velocity profile $v$, the convective nonlinear term gives the schematic move
\begin{equation*}
 (u_0,\phi_0)\longrightarrow
 (u_0-B_0(v,v),\phi_0).
\end{equation*}
We add $\delta^{-1/2}v$, let the system evolve under $f$ for
a time $\delta$, and remove the added profile.
The first asymptotic limit identifies the remaining velocity
displacement; see \Cref{second_limit}.

\item \emph{Using the fluid to move the phase.}
For available phase and fluid profiles $\rho$ and $\zeta$ with
$B_0(\zeta,\zeta)=0$, the cubic nonlinearity and the transport term
in the phase equation can be combined schematically as follows
\begin{equation*}
 \begin{aligned}
 (u_0,\phi_0)
 &\xrightarrow{\ \Delta(\rho^3)\ }
 (u_0,\phi_0+\Delta(\rho^3)) \xrightarrow{\ -\Delta(\rho^3)+B_1(\zeta,\rho)\ }
 (u_0,\phi_0+B_1(\zeta,\rho)).
 \end{aligned}
\end{equation*}
The first move produces a cubic phase displacement.
The second uses the velocity profile $\zeta$ to transport the
phase profile $-\rho$. After removing the auxiliary profiles,
this second move adds $-\Delta(\rho^3)+B_1(\zeta,\rho)$.
Its cubic contribution cancels the first displacement, leaving
only the phase displacement generated by transport.
These steps are made precise in
\Cref{lem:cubic_generation,lem:mixed_generation,lem:transport_generation}.

\item \emph{Iteration and combination.}
Iterating the fluid and phase moves gives the displacements of
\Cref{second_limit,small1}, with corrections still taking values
in $\mathcal H$. Their combination gives access to the unit ball,
see the approximation argument in \Cref{sec:unit-ball-accessibility-proof}.
\end{itemize}
\subsection{Direct displacements}
We first consider a constant correction.
\begin{lemma}\label{lem:direct_construction}
Let $h\in\mathcal H$ and
$f\in L^\infty_{\mathrm{loc}}(\mathbb R_+;\mathbb V^0)$.
For any $\epsilon,\sigma>0$ and any $z_0\in\mathbf H^2\times\dot H^4$,
there exist $\tau\in(0,\sigma)$ and
$\eta\in L^\infty(0,\tau;\mathcal H)$ such that
\begin{equation}\label{eq:direct_maneuver}
\norm{\mathcal S_\tau(z_0,0,f+\eta)-(z_0+h)}_{\mathbb V^1}<\epsilon.
\end{equation}
\end{lemma}
\begin{proof}
By \Cref{thm:asymp-fluid}, with $\zeta=0$ and $\eta=h$,
for every sufficiently small $\tau>0$ the solution is well
defined on $[0,\tau]$, and
\begin{equation*}
\mathcal S_\tau(z_0,0,f+\tau^{-1}h)
\longrightarrow z_0+h
\qquad\text{in }\mathbb V^1
\qquad\text{as }\tau\to0^+.
\end{equation*}
We can thus choose $\tau\in(0,\sigma)$ sufficiently small
so that the constant control $\eta(t)=\tau^{-1}h$ satisfies
\eqref{eq:direct_maneuver}.
\end{proof}

\subsection{Velocity displacements}
We now use the saturation properties and the fluid nonlinearity to generate velocity directions beyond those obtained by direct displacements. We begin with presenting a couple of auxiliary lemmas.
\begin{lemma}\label{lem:velocity_fourier_generators}
For every $N\ge0$, the space $\mathcal H_1(N)$ is spanned by finitely
many modes of the form
\begin{equation*}
 a_j\cos(k_j\cdot x),\qquad a_j\sin(k_j\cdot x),
 \qquad k_j\in\mathbb Z^3\setminus\{0\},\quad a_j\cdot k_j=0.
\end{equation*}
\end{lemma}

\begin{proof}
The case $N=0$ follows from \eqref{eq:G1_def}. Assume the assertion
holds for $E:=\mathcal H_1(N)$. Let
\begin{equation*}
 v_c=a\cos(k\cdot x),\quad v_s=a\sin(k\cdot x),\qquad
 w_c=b\cos(l\cdot x),\quad w_s=b\sin(l\cdot x)
\end{equation*}
belong to $E$. For $k\neq\pm l$, set
$p:=P_{k+l}\bigl((a\cdot l)b+(b\cdot k)a\bigr)$. Since
$a\cdot k=b\cdot l=0$, we have
\begin{align}
 B_0(v_s+w_s,v_s+w_s)-B_0(v_c+w_c,v_c+w_c)
 &=p\sin((k+l)\cdot x),\notag\\
 B_0(v_c+w_s,v_c+w_s)+B_0(v_s+w_c,v_s+w_c)
 &=p\cos((k+l)\cdot x).
 \label{eq:velocity-mode-separation}
\end{align}
By \eqref{eq:F1_def}, the two modes belong to $\mathcal H_1(N+1)$.
Replacing $l$ by $-l$ also gives the sine and cosine modes at $k-l$.
For $k=\pm l$, all the quadratic terms vanish.

Let $r_1,\ldots,r_m$ be the sine and cosine functions spanning $E$.
Since $B_0(r_j,r_j)=0$, for $w=\sum_{j=1}^{m}c_jr_j$ we have
\begin{align}\label{eq:velocity-quadratic-expansion}
 B_0(w,w)
 &=\sum_{i<j}c_ic_j
   \bigl(B_0(r_i,r_j)+B_0(r_j,r_i)\bigr)=\sum_{i<j}c_ic_j B_0(r_i+r_j,r_i+r_j).
\end{align}
By \eqref{eq:F1_def}, $\mathcal H_1(N+1)
 =E+\operatorname{span}
 \left\{B_0(r_i+r_j,r_i+r_j):1\le i<j\le m\right\}$. For $k\neq\pm l$, put
$q:=P_{k-l}\bigl((a\cdot l)b-(b\cdot k)a\bigr)$.
Using \eqref{eq:velocity-mode-separation} with $l$ and $-l$ gives
\begin{align*}
 B_0(v_c+w_c,v_c+w_c)
 &=-\tfrac12p\sin((k+l)\cdot x)+\tfrac12q\sin((k-l)\cdot x),\\
 B_0(v_s+w_s,v_s+w_s)
 &=\tfrac12p\sin((k+l)\cdot x)+\tfrac12q\sin((k-l)\cdot x),\\
 B_0(v_c+w_s,v_c+w_s)
 &=\tfrac12p\cos((k+l)\cdot x)+\tfrac12q\cos((k-l)\cdot x),\\
 B_0(v_s+w_c,v_s+w_c)
 &=\tfrac12p\cos((k+l)\cdot x)-\tfrac12q\cos((k-l)\cdot x).
\end{align*}
Each $B_0(r_i+r_j,r_i+r_j)$ has one of these four
forms, or vanishes when $k=\pm l$. Therefore,
$\mathcal H_1(N+1)$ is spanned by $E$ and the finitely many sine
and cosine functions obtained in \eqref{eq:velocity-mode-separation} (with $l$ and $-l$). Since $p\cdot(k+l)=q\cdot(k-l)=0$, the assertion follows for $N+1$.
\end{proof}
\begin{lemma}\label{lem:velocity_representation}
	For every $N\ge1$ and every $\eta\in\mathcal H_1(N)$, there exist
	$d\ge0$ and $\zeta_0,\ldots,\zeta_d\in\mathcal H_1(N-1)$ such that
	\begin{equation*}%\label{eq:velocity_representation}
		\eta= \zeta_0-\sum_{i=1}^{d}B_0(\zeta_i,\zeta_i).
	\end{equation*}
\end{lemma}
\begin{proof}
Set $E=\mathcal H_1(N-1)$. By
	\Cref{lem:velocity_fourier_generators}, $E$ is spanned by finitely
	many divergence-free Fourier modes of the form
	\begin{equation}\label{eq:velocity_fourier_basis}
		v_j=a_j\cos(k_j\cdot x)
		\quad\text{or}\quad
		v_j=a_j\sin(k_j\cdot x),
		\qquad a_j\cdot k_j=0.
	\end{equation}
	For every mode $v_j$ in \eqref{eq:velocity_fourier_basis}, we have
$B_0(v_j,v_j)=0.$ 
We claim that, for every $\zeta\in E$, there exist finitely many
	$w_\ell\in E$ such that
	\begin{equation}\label{eq:velocity_sign_change}
		B_0(\zeta,\zeta)
		=
		-\sum_{\ell}B_0(w_\ell,w_\ell).
	\end{equation}
	To see this, write $\zeta=\sum_{j=1}^{q}c_jv_j$
	and apply \eqref{eq:velocity-quadratic-expansion}.
	For every pair such that $c_ic_j\neq0$, set
$w_{ij}
		=
		\sqrt{|c_ic_j|}
		\bigl(v_i-\operatorname{sgn}(c_ic_j)v_j\bigr).$ 
	Then $w_{ij}\in E$ and
	$B_0(w_{ij},w_{ij})
		=
		-c_ic_j
		\bigl(B_0(v_i,v_j)+B_0(v_j,v_i)\bigr).$ 
	Summing over $i<j$ proves \eqref{eq:velocity_sign_change}.
	
	Now, by \eqref{eq:F1_def}, any
	$\eta\in\mathcal H_1(N)$ can be written as
	\begin{equation*}
		\eta
		=
		\sum_{r=1}^{p}
		\alpha_r
		\left(
		\zeta_{r,0}
		-
		\sum_{j=1}^{m_r}
		B_0(\zeta_{r,j},\zeta_{r,j})
		\right),
	\end{equation*}
	where $\alpha_r\in\mathbb R$ and $\zeta_{r,j}\in E$. Since $E$ is
	a vector space,
$\zeta_0:=\sum_{r=1}^{p}\alpha_r\zeta_{r,0}\in E.$ 
	If $\alpha_r\ge0$, then
	\begin{equation*}
		\alpha_rB_0(\zeta_{r,j},\zeta_{r,j})
		=
		B_0(\sqrt{\alpha_r}\,\zeta_{r,j},
		\sqrt{\alpha_r}\,\zeta_{r,j}).
	\end{equation*}
	If $\alpha_r<0$, set
	$v=\sqrt{|\alpha_r|}\,\zeta_{r,j}\in E$.
	By \eqref{eq:velocity_sign_change}, there exist
	finitely many $w_\ell\in E$ such that
	\begin{equation*}
	-\alpha_r B_0(\zeta_{r,j},\zeta_{r,j})
	=B_0(v,v)
	=-\sum_\ell B_0(w_\ell,w_\ell).
		\end{equation*}
	Therefore all the quadratic terms can
	be collected with the same negative sign, and we obtain
	\begin{equation*}
		\eta
		=
		\zeta_0-\sum_{i=1}^{d}B_0(\zeta_i,\zeta_i),
	\end{equation*}
	for some $\zeta_0,\ldots,\zeta_d\in E=\mathcal H_1(N-1)$. This ends the proof.
\end{proof}
\begin{rmk}\label{ner1}
By \Cref{lem:velocity_representation}, for every
$N\ge1$, each element of $\mathcal H_1(N)$
admits a representation of the form
\[
\zeta_0-\sum_{i=1}^d B_0(\zeta_i,\zeta_i),
\qquad \zeta_i\in\mathcal H_1(N-1).
\]
Consequently, on the spaces $\mathcal H_1(N-1)$,
our saturation operation $\mathcal F_1$
coincides with the representation-based
saturation operation used in \cite{N15}.
\end{rmk}
For convenience, we recall the notation and the shift identity from \Cref{wp1}. For $s,r\ge0$, we write
\begin{equation}\label{eq:control-force-translations}
 f_s(r):=f(s+r),\qquad \eta_s(r):=\eta(s+r).
\end{equation}
Writing $\mathcal T=\mathcal T(z_0,\xi,f+\eta)$, we have
\begin{equation}\label{eq:control-shift-identity}
 \mathcal S_t(z_0,\xi,f+\eta)
 =\mathcal S_t(z_0+\xi,0,f+\eta)-\xi,
 \qquad 0\le t<\mathcal T.
\end{equation}
The next lemma will help us to concatenate controls in the constructions below.
\begin{lemma}\label{lem:concatenation}
Let $z_0\in\mathbb V^1$ and
$f\in L^\infty_{\mathrm{loc}}(\mathbb R_+;\mathbb V^0)$.
Let $\tau_1,\ldots,\tau_m>0$, and set
$s_0=0$ and $s_j=\tau_1+\cdots+\tau_j$.
On the $j$-th interval $[s_{j-1},s_j]$, write
$t=s_{j-1}+r$, with $0\le r\le\tau_j$.
For $j=1,\ldots,m$, let
$\eta^j\in L^\infty(0,\tau_j;\mathcal H)$ and suppose that
\begin{equation*}
 z^{(j)}\in C([0,\tau_j];\mathbb V^1)
 \cap L^2(0,\tau_j;\mathbf H^2\times\dot H^4)
\end{equation*}
solves \eqref{abs} with force
$r\mapsto f(s_{j-1}+r)+\eta^j(r)$, with initial conditions
\begin{equation*}
 z^{(1)}(0)=z_0,\qquad
 z^{(j)}(0)=z^{(j-1)}(\tau_{j-1}),\quad j=2,\ldots,m.
\end{equation*}
Set $z_j:=z^{(j)}(\tau_j)$ for $j=1,\ldots,m$, and define
$\eta$ almost everywhere on $\mathbb R_+$ by
\begin{equation*}
 \eta(t)=
 \begin{cases}
  \eta^j(t-s_{j-1}),&s_{j-1}<t<s_j,\quad j=1,\ldots,m,\\
  0,&t\ge s_m.
 \end{cases}
\end{equation*}
Then $\eta\in L^\infty(\mathbb R_+;\mathcal H)$,
$\mathcal T(z_0,0,f+\eta)>s_m$, and
\begin{equation}\label{eq:concatenated-regularity}
 \mathcal S_\cdot(z_0,0,f+\eta)
 \in C([0,s_m];\mathbb V^1)
 \cap L^2(0,s_m;\mathbf H^2\times\dot H^4).
\end{equation}
Moreover,
\begin{equation}\label{eq:concatenated-solutions}
 \mathcal S_t(z_0,0,f+\eta)=z^{(j)}(t-s_{j-1}),
 \qquad s_{j-1}\le t\le s_j,\quad j=1,\ldots,m.
\end{equation}
\end{lemma}
\begin{proof}
	Since $m$ is finite and each $\eta^j \in L^\infty(0,\tau_j;\mathcal H)$, we have $\eta\in L^\infty(\mathbb R_+;\mathcal H)$.
	Let $\mathcal T:=\mathcal T(z_0,0,f+\eta)$ and
$z(t):=\mathcal S_t(z_0,0,f+\eta)$ be given by \Cref{wp1}.
We prove by induction that $\mathcal T>s_j$ and $z(s_j)=z_j$
for $j=0,\ldots,m$.
For $j=0$, $s_0=0<\mathcal T$ and $z(s_0)=z_0$.

Let $j\in\{1,\ldots,m\}$ and suppose that
$s_{j-1}<\mathcal T$ and $z(s_{j-1})=z_{j-1}$.
Since $\eta_{s_{j-1}}=\eta^j$ almost everywhere on $(0,\tau_j)$,
\eqref{eq:restart_maximal_time} and \eqref{eq:restart_identity},
applied with $\xi=0$, together with \Cref{cor:unique-strong-common}, give
\begin{equation*}
 \begin{aligned}
 z(s_{j-1}+r)
 &=\mathcal S_r(z_{j-1},0,f_{s_{j-1}}+\eta_{s_{j-1}})\\
 &=z^{(j)}(r),
 \qquad 0\le r<\min\{\tau_j,\mathcal T-s_{j-1}\}.
 \end{aligned}
\end{equation*}
If $\mathcal T\le s_j$, then
$0<\mathcal T-s_{j-1}\le\tau_j$.
In particular, $z^{(j)}(\mathcal T-s_{j-1})$ is defined.
Since $z^{(j)}\in C([0,\tau_j];\mathbb V^1)$ and
$z(s_{j-1}+r)=z^{(j)}(r)$ for
$0\le r<\mathcal T-s_{j-1}$, we obtain
\begin{align*}
 \lim_{r\uparrow \mathcal T-s_{j-1}}
 \|z(s_{j-1}+r)\|_{\mathbb V^1}
 &=\lim_{r\uparrow \mathcal T-s_{j-1}}
   \|z^{(j)}(r)\|_{\mathbb V^1}=\|z^{(j)}(\mathcal T-s_{j-1})\|_{\mathbb V^1}<\infty.
\end{align*}
Since $s_{j-1}+r\uparrow\mathcal T$, this contradicts
\eqref{eq:maximal-strong-blowup} by \eqref{eq:shift_maximal_time}
with $\xi=0$.
Thus $\mathcal T>s_j$, and \eqref{eq:concatenated-solutions}
follows by continuity.
In particular, $z(s_j)=z_j$, which completes the induction.
The regularity \eqref{eq:concatenated-regularity} follows from
$s_m<\mathcal T$ and \Cref{wp1}.
\end{proof}
\begin{proposition}\label{second_limit}
	Let $\mathcal H$ be the finite-dimensional space defined in \eqref{h}.
{Let $f=(f_u,f_\phi)\in L^\infty_{\mathrm{loc}}(\mathbb R_+;\mathbb V^0)$.}
For any $\epsilon, \sigma >0$, $z_0 \in \mathbb V^1$, and $\rho \in \mathcal H_1(\infty)$, there {exist} $\tau \in (0,\sigma)$ and a control $\eta\in L^{\infty}(0,\tau;\mathcal H)$ such that
	\begin{align}\label{limit_phi2}
		\norm{{\mathcal S_{\tau}(z_0, 0, f+\eta)} - z_0-\widehat \rho}_{\mathbb V^1} <\epsilon,
	\end{align}
	where $\widehat \rho:=(\rho,0)\in \mathbb V^1.$
\end{proposition}

\begin{proof}
	We first prove the result when $z_0\in{\mathbf H^2\times\dot H^4}$. We proceed by induction on $N\in\mathbb N$. For a given $N\in\mathbb N$, we state the following claim $\mathcal P(N)$:
	\begin{quote}
		\textit{For every force $f=(f_u,f_\phi)\in L^\infty_{\mathrm{loc}}(\mathbb R_+;\mathbb V^0)$, $z_0\in{\mathbf H^2\times\dot H^4}$, $\rho\in\mathcal H_1(N)$, and $\epsilon,\sigma>0$, there exist $\tau\in(0,\sigma)$ and $\eta\in L^\infty(0,\tau;\mathcal H)$ such that \eqref{limit_phi2} holds {and $\mathcal S_\tau(z_0,0,f+\eta)\in\mathbf H^2\times\dot H^4$}.}
	\end{quote}
	The base case $N=0$ follows from the limit \eqref{asymp_phi} with $\zeta=0$ and $\eta=\widehat\rho$. Indeed,
	\begin{align*}
		\mathcal S_\delta(z_0,0,f+\delta^{-1}\widehat\rho)
		\longrightarrow z_0+\widehat\rho
		\quad\text{in }\mathbb V^1
	\end{align*}
	as $\delta\to0^+$. {Choose $\delta\in(0,\sigma)$
sufficiently small and set $\eta(t)=\delta^{-1}\widehat\rho$ for
$0<t<\delta$, so that the corresponding solution $S(z_0,0,f+\eta)$ is well-defined on
$[0,\delta]$ and
\begin{equation*}
 \|\mathcal S_\delta(z_0,0,f+\eta)-z_0-\widehat\rho\|_{\mathbb V^1}
 <\frac{\epsilon}{2}.
\end{equation*}
Keeping $\eta$ fixed, write
$z(t)=\mathcal S_t(z_0,0,f+\eta)$ for $0\leq t\leq\delta$.
By \Cref{wp1}, $z\in C([0,\delta];\mathbb V^1)
\cap L^2(0,\delta;\mathbf H^2\times\dot H^4)$.
By \Cref{lem:regular-times}, we may choose $\tau\in(0,\delta)$,
arbitrarily close to $\delta$, such that
\begin{equation*}
 z(\tau)\in\mathbf H^2\times\dot H^4,\qquad
 \|z(\tau)-z(\delta)\|_{\mathbb V^1}<\frac{\epsilon}{2}.
\end{equation*}
Consequently,
\begin{align*}
 \|z(\tau)-z_0-\widehat\rho\|_{\mathbb V^1}
 &\leq \|z(\tau)-z(\delta)\|_{\mathbb V^1}
       +\|z(\delta)-z_0-\widehat\rho\|_{\mathbb V^1}\\
 &<\frac{\epsilon}{2}+\frac{\epsilon}{2}=\epsilon.
\end{align*}
Set $\eta^\tau:=\eta|_{(0,\tau)}$.
By \Cref{cor:unique-strong-common}, one has
$\mathcal S_\tau(z_0,0,f+\eta^\tau)=z(\tau).$
Thus $\mathcal P(0)$ holds with control $\eta^\tau$.}
	
	Assume that $\mathcal P(N-1)$ holds. We first show that for every $\zeta\in\mathcal H_1(N-1)$, $z_0\in{\mathbf H^2\times\dot H^4}$, and $\epsilon,\sigma>0$, there exist $\tau\in(0,\sigma)$ and $\eta\in L^\infty(0,\tau;\mathcal H)$ such that
	\begin{equation}\label{eq:quadratic_generation_velocity}
		\norm{\mathcal S_\tau(z_0,0,f+\eta)-z_0+B(\widehat\zeta)}_{\mathbb V^1}<\epsilon,
	\end{equation}
	{and
\begin{equation}\label{eq:quadratic_generation_velocity_regularity}
 \mathcal S_\tau(z_0,0,f+\eta)\in\mathbf H^2\times\dot H^4,
\end{equation}}
	where $\widehat\zeta:=(\zeta,0)$.
	Indeed, by \eqref{asymp_phi},
	\begin{align*}
		\mathcal S_\delta(z_0,\delta^{-1/2}\widehat\zeta,f)
		\longrightarrow z_0-B(\widehat\zeta)
		\quad\text{in }\mathbb V^1 \text{ as } \delta\to0^+.
	\end{align*}
Hence, using \eqref{eq:control-shift-identity} we have
$	\mathcal S_\delta(z_0+\delta^{-1/2}\widehat\zeta,0,f)
		=
		\mathcal S_\delta(z_0,\delta^{-1/2}\widehat\zeta,f)
		+\delta^{-1/2}\widehat\zeta,$ 
	{and we first choose $\delta\in(0,\sigma/3)$ and set
$\alpha:=\delta^{-1/2}$ so that the corresponding solution is well-defined
on $[0,\delta]$ and
\begin{equation*}
 \norm{\mathcal S_\delta(z_0+\alpha\widehat\zeta,0,f)
       -z_0-\alpha\widehat\zeta+B(\widehat\zeta)}_{\mathbb V^1}
 <\frac{\epsilon}{6}.
\end{equation*}
Keeping $\alpha$ fixed, write
$y(t):=\mathcal S_t(z_0+\alpha\widehat\zeta,0,f)$ for $0\leq t\leq\delta$.
By \Cref{wp1}, $y\in C([0,\delta];\mathbb V^1)
\cap L^2(0,\delta;\mathbf H^2\times\dot H^4)$.
By \Cref{lem:regular-times}, we may choose $\tau_2\in(0,\delta)$,
sufficiently close to $\delta$, such that
\begin{equation*}
 \widetilde z:=y(\tau_2)\in\mathbf H^2\times\dot H^4,\qquad
 \|y(\tau_2)-y(\delta)\|_{\mathbb V^1}<\frac{\epsilon}{6}.
\end{equation*}
Consequently,
\begin{align}\label{eq:quad_asymptotic}
 \|\widetilde z-z_0-\alpha\widehat\zeta+B(\widehat\zeta)\|_{\mathbb V^1}
 &\leq \|\widetilde z-y(\delta)\|_{\mathbb V^1}
       +\|y(\delta)-z_0-\alpha\widehat\zeta+B(\widehat\zeta)\|_{\mathbb V^1}
       \notag\\
 &<\frac{\epsilon}{6}+\frac{\epsilon}{6}=\frac{\epsilon}{3}.
\end{align}
}
	Since $-\alpha\zeta\in\mathcal H_1(N-1)$, the induction hypothesis, applied to $\widetilde z$ with $f_{\tau_2}$ (using notation \eqref{eq:control-force-translations}) in place of $f$, gives $\tau_3\in(0,\sigma/3)$ and $\eta^3\in L^\infty(0,\tau_3;\mathcal H)$ such that
	\begin{equation}\label{eq:quad_remove}
		\norm{\mathcal S_{\tau_3}(\widetilde z,0,f_{\tau_2}+\eta^3)-\widetilde z+\alpha\widehat\zeta}_{\mathbb V^1}
		<\frac{\epsilon}{3}.
	\end{equation}
	Combining \eqref{eq:quad_asymptotic} and \eqref{eq:quad_remove}, we have
	\begin{equation}\label{eq:quad_reference}
		\norm{\mathcal S_{\tau_3}(\mathcal S_{\tau_2}(z_0+\alpha\widehat\zeta,0,f),0,f_{\tau_2}+\eta^3)-z_0+B(\widehat\zeta)}_{\mathbb V^1}
		<\frac{2\epsilon}{3}.
	\end{equation}
	We now fix $\tau_2,\tau_3$ and $\eta^3$.

Let $r_2,r_3>0$ and $C_2,C_3\ge1$ be given by \Cref{cont}
and \eqref{eq:cont_stability}, applied with $\xi=0$ to the
solution from $z_0+\alpha\widehat\zeta$ with force $f$
on $[0,\tau_2]$, and to the solution from $\mathcal S_{\tau_2}(z_0+\alpha\widehat\zeta,0,f)$ with force
$f_{\tau_2}+\eta^3$ on $[0,\tau_3]$, respectively. Choose
\begin{equation*}
 0<\gamma<
 \min\left\{\frac{r_2}{2},\frac{r_3}{4C_2},
                  \frac{\epsilon}{9C_2C_3}\right\}.
\end{equation*}
Since the translation of $f$ is continuous, we have
\begin{equation*}
 \|f_s-f\|_{L^2(0,\tau_2+\tau_3;\mathbb V^0)}
 \longrightarrow0\qquad\text{as }s\to0^+.
\end{equation*}
Since $\alpha\zeta\in\mathcal H_1(N-1)$, we can apply the induction
hypothesis with a sufficiently small time horizon to obtain
$\tau_1\in(0,\sigma/3)$ and
$\eta^1\in L^\infty(0,\tau_1;\mathcal H)$ satisfying
\begin{equation*}%\label{eq:quad_prepare}
 \norm{\mathcal S_{\tau_1}(z_0,0,f+\eta^1)
                   -z_0-\alpha\widehat\zeta}_{\mathbb V^1}<\gamma,
\text{ and }
 \|f_{\tau_1}-f\|_{L^2(0,\tau_2+\tau_3;\mathbb V^0)}<\gamma.
\end{equation*}
Set $z^1:=\mathcal S_{\tau_1}(z_0,0,f+\eta^1)$.
Since $2\gamma<r_2$,  \, $\mathcal S(z_0+\alpha\widehat\zeta,0,f)$ is well defined on $[0,\tau_2]$, \Cref{cont} gives
$\mathcal T(z^1,0,f_{\tau_1})>\tau_2$.
Set $z^2:=\mathcal S_{\tau_2}(z^1,0,f_{\tau_1})$.
By \eqref{eq:cont_stability},
\begin{align*}
 \norm{z^2-\mathcal S_{\tau_2}(z_0+\alpha\widehat\zeta,0,f)}_{\mathbb V^1}
 &\le C_2\Bigl(
     \norm{z^1-z_0-\alpha\widehat\zeta}_{\mathbb V^1}
     +\|f_{\tau_1}-f\|_{L^2(0,\tau_2;\mathbb V^0)}
     \Bigr)<2C_2\gamma.
\end{align*}
Since $(2C_2+1)\gamma<r_3$, \Cref{cont} gives
$\mathcal T(z^2,0,f_{\tau_1+\tau_2}+\eta^3)>\tau_3$,
and \eqref{eq:cont_stability} gives
\begin{align}\label{eq:quad_final_stability}
 &\norm{\mathcal S_{\tau_3}(z^2,0,f_{\tau_1+\tau_2}+\eta^3)
       -\mathcal S_{\tau_3}(\mathcal S_{\tau_2}(z_0+\alpha\widehat\zeta,0,f),0,f_{\tau_2}+\eta^3)}
             _{\mathbb V^1}\notag\\
 &\quad\le C_3\Bigl(
       \norm{z^2-\mathcal S_{\tau_2}(z_0+\alpha\widehat\zeta,0,f)}_{\mathbb V^1}
       +\|f_{\tau_1+\tau_2}-f_{\tau_2}\|
                    _{L^2(0,\tau_3;\mathbb V^0)}\Bigr)\notag\\
 &\quad<C_3(2C_2+1)\gamma
       \le3C_2C_3\gamma<\frac{\epsilon}{3}.
\end{align}
Note that the last norm equals
$\|f_{\tau_1}-f\|_{L^2(\tau_2,\tau_2+\tau_3;\mathbb V^0)}$
by a change of variables, and is therefore less than $\gamma$.

Set $T:=\tau_1+\tau_2+\tau_3<\sigma$ and define
\begin{equation*}
 \eta(t)=
 \begin{cases}
  \eta^1(t),&t\in(0,\tau_1),\\
  0,&t\in(\tau_1,\tau_1+\tau_2),\\
  \eta^3(t-\tau_1-\tau_2),&t\in(\tau_1+\tau_2,T),\\
  0,&t\ge T.
 \end{cases}
\end{equation*}
Consider the three solutions
\begin{align*}
 z^{(1)}(r)&:=\mathcal S_r(z_0,0,f+\eta^1),
 &&0\le r\le\tau_1,\\
 z^{(2)}(r)&:=\mathcal S_r(z^1,0,f_{\tau_1}),
 &&0\le r\le\tau_2,\\
 z^{(3)}(r)&:=\mathcal S_r(z^2,0,f_{\tau_1+\tau_2}+\eta^3),
 &&0\le r\le\tau_3.
\end{align*}
By \Cref{wp1},
$z^{(j)}\in C([0,\tau_j];\mathbb V^1)
\cap L^2(0,\tau_j;\mathbf H^2\times\dot H^4)$, $j=1,2,3$.
Applying \Cref{lem:concatenation} to
$z^{(1)},z^{(2)},z^{(3)}$ and the controls $\eta^1,0,\eta^3$
gives $\mathcal T(z_0,0,f+\eta)>T$.
Set $z(t):=\mathcal S_t(z_0,0,f+\eta)$ for $0\le t\le T$.
Then
\begin{equation}\label{eq:quad_concatenation}
 z(T)=\mathcal S_T(z_0,0,f+\eta)
 =\mathcal S_{\tau_3}(z^2,0,f_{\tau_1+\tau_2}+\eta^3).
\end{equation}
By \eqref{eq:concatenated-regularity},
$z\in C([0,T];\mathbb V^1)\cap L^2(0,T;\mathbf H^2\times\dot H^4)$.
Combining \eqref{eq:quad_concatenation},
\eqref{eq:quad_final_stability}, and \eqref{eq:quad_reference}, we obtain
$e_T:=\|z(T)-z_0+B(\widehat\zeta)\|_{\mathbb V^1}<\epsilon$.
Keeping $\eta$ fixed, by continuity and \Cref{lem:regular-times},
we may choose $\tau\in(\tau_1+\tau_2,T)$, sufficiently close to $T$, such that
\begin{equation*}
 z(\tau)\in\mathbf H^2\times\dot H^4,\qquad
 \|z(\tau)-z(T)\|_{\mathbb V^1}<\frac{\epsilon-e_T}{2},
\end{equation*}
which further implies the following
\begin{align*}
 \|z(\tau)-z_0+B(\widehat\zeta)\|_{\mathbb V^1}
 &\leq \|z(\tau)-z(T)\|_{\mathbb V^1}
       +\|z(T)-z_0+B(\widehat\zeta)\|_{\mathbb V^1}\\
 &<\frac{\epsilon-e_T}{2}+e_T<\epsilon.
\end{align*}
Set $\eta^\tau:=\eta|_{(0,\tau)}\in L^\infty(0,\tau;\mathcal H)$.
By \Cref{cor:unique-strong-common},
 $ \mathcal S_\tau(z_0,0,f+\eta^\tau)=z(\tau).$ 
Thus \eqref{eq:quadratic_generation_velocity} and
\eqref{eq:quadratic_generation_velocity_regularity} hold with
control $\eta^\tau$.
	Let now $\rho\in\mathcal H_1(N)$. Thanks to \Cref{lem:velocity_representation}, there exist $\zeta_0,\ldots,\zeta_d\in\mathcal H_1(N-1)$ such that
	\begin{equation}\label{eq:eta_decomp}
		\rho=\zeta_0-\sum_{i=1}^{d}B_0(\zeta_i,\zeta_i).
	\end{equation}
	Set $\widetilde \epsilon=\epsilon/(d+1)$. By the induction hypothesis, there exist $\tau_0\in(0,\sigma/(d+1))$ and $\eta^0\in L^\infty(0,\tau_0;\mathcal H)$ such that, setting
$Z_0=\mathcal S_{\tau_0}(z_0,0,f+\eta^0),$
	we have
%	\begin{equation}%\label{eq:first_velocity_direction}
	$	\norm{Z_0-z_0-\widehat\zeta_0}_{\mathbb V^1}<\widetilde \epsilon.$
%	\end{equation}

{The induction hypothesis also gives
$Z_0\in\mathbf H^2\times\dot H^4$.}
For $i=1,\ldots,d$, suppose that $Z_{i-1}\in{\mathbf H^2\times\dot H^4}$ and set
$s_{i-1}:=\tau_0+\cdots+\tau_{i-1}$.
Since $f_{s_{i-1}}$ satisfies the same hypotheses as $f$,
\eqref{eq:quadratic_generation_velocity} {and
\eqref{eq:quadratic_generation_velocity_regularity}} give
$\tau_i\in(0,\sigma/(d+1))$ and
$\eta^i\in L^\infty(0,\tau_i;\mathcal H)$ such that, setting
$Z_i=\mathcal S_{\tau_i}(Z_{i-1},0,f_{s_{i-1}}+\eta^i),$
we have
%	\begin{equation}\label{eq:quadratic_iterations}
$\norm{Z_i-Z_{i-1}+B(\widehat\zeta_i)}_{\mathbb V^1}<\widetilde\epsilon.$
%	\end{equation}
{By \eqref{eq:quadratic_generation_velocity_regularity},
$Z_i\in\mathbf H^2\times\dot H^4$, so the same assertions apply
at the next step.}

	Therefore,
	\begin{align*}
		&\norm{Z_d-z_0-\widehat\zeta_0+\sum_{i=1}^{d}B(\widehat\zeta_i)}_{\mathbb V^1}\leq
		\norm{Z_0-z_0-\widehat\zeta_0}_{\mathbb V^1}
		+\sum_{i=1}^{d}
		\norm{Z_i-Z_{i-1}+B(\widehat\zeta_i)}_{\mathbb V^1}
		<\epsilon.
	\end{align*}
	By \eqref{eq:eta_decomp},
%	\begin{equation*}
	$	\widehat\rho=\widehat\zeta_0-\sum_{i=1}^{d}B(\widehat\zeta_i).$
%	\end{equation*}
	{Set $\tau:=\tau_0+\cdots+\tau_d<\sigma$ and define
\begin{align*}
 \eta(r)&=\eta^0(r), \quad 0<r<\tau_0,\\
 \eta(s_{i-1}+r)&=\eta^i(r), \quad 0<r<\tau_i, \quad i=1,\ldots,d,
\end{align*}
and $\eta(t)=0$ for $t\ge\tau$.
Applying \Cref{lem:concatenation} to the
$d+1$ solutions defining $Z_0,\ldots,Z_d$ gives
$\mathcal T(z_0,0,f+\eta)>\tau$ and
\begin{equation*}
 \mathcal S_\tau(z_0,0,f+\eta)=Z_d\in\mathbf H^2\times\dot H^4.
\end{equation*}
This proves $\mathcal P(N)$.}
	Thus, by induction, \eqref{limit_phi2} is proved for every $z_0\in{\mathbf H^2\times\dot H^4}$ and every $\rho\in\mathcal H_1(\infty)$.

	We finally consider an arbitrary $z_0\in\mathbb V^1$. Let
$z(t):=\mathcal S_t(z_0,0,f).$
	Since the solution is well-defined on some interval $[0,T)$ and $z\in C([0,T);\mathbb V^1),$
	{by continuity we may choose
$b\in(0,\min\{T,\sigma/2\})$ such that}
	\begin{equation}\label{eq:continuity_initial}
		\norm{z(t)-z_0}_{\mathbb V^1}<\frac{\epsilon}{2},
		\qquad t\in{[0,b]}.
	\end{equation}
	{By \Cref{wp1}, $z\in L^2(0,b;\mathbf H^2\times\dot H^4)$.
Applying \Cref{lem:regular-times} on $[0,b]$ with $t_0=0$,
we obtain $\delta\in(0,b)$ such that
$z(\delta)\in\mathbf H^2\times\dot H^4$.}
	Set $	\widehat z:=z(\delta)=\mathcal S_\delta(z_0,0,f).$
	By \eqref{eq:continuity_initial}
	\begin{equation}\label{eq:regularize_velocity}
		\widehat z\in{\mathbf H^2\times\dot H^4},
	\text{ and }
		\norm{\widehat z-z_0}_{\mathbb V^1}<\frac{\epsilon}{2}.
	\end{equation}
	Applying the result already proved {for initial data in
$\mathbf H^2\times\dot H^4$} to $\widehat z$, with $\epsilon/2$ and $\sigma/2$ and $f_\delta$ in place of $f$, we find $\tau_1\in(0,\sigma/2)$ and
	$\eta^1\in L^\infty(0,\tau_1;\mathcal H)$ such that
	\begin{equation}\label{eq:regularized_control}
		\norm{\mathcal S_{\tau_1}(\widehat z,0,f_\delta+\eta^1)-\widehat z-\widehat\rho}_{\mathbb V^1}
		<\frac{\epsilon}{2}.
	\end{equation}
Let $\tau:=\delta+\tau_1<\sigma$ and define
\begin{equation*}
 \eta(t)=
 \begin{cases}
  0,&t\in(0,\delta),\\
  \eta^1(t-\delta),&t\in(\delta,\tau),\\
  0,&t\ge\tau.
 \end{cases}
\end{equation*}
Applying \Cref{lem:concatenation} with controls
$0,\eta^1$ to the solutions
$r\mapsto\mathcal S_r(z_0,0,f)$ on $[0,\delta]$ and
$r\mapsto\mathcal S_r(\widehat z,0,f_\delta+\eta^1)$ on
$[0,\tau_1]$ gives $\mathcal T(z_0,0,f+\eta)>\tau$ and
$\mathcal S_\tau(z_0,0,f+\eta)
 =\mathcal S_{\tau_1}(\widehat z,0,f_\delta+\eta^1).$
 Therefore, using \eqref{eq:regularize_velocity} and
	\eqref{eq:regularized_control}, we obtain
	\begin{align*}
		\norm{\mathcal S_\tau(z_0,0,f+\eta)-z_0-\widehat\rho}_{\mathbb V^1}
		&\leq
		\norm{\mathcal S_{\tau_1}(\widehat z,0,f_\delta+\eta^1)
			-\widehat z-\widehat\rho}_{\mathbb V^1}
		+\norm{\widehat z-z_0}_{\mathbb V^1}\\
		&<\frac{\epsilon}{2}+\frac{\epsilon}{2}
		=\epsilon.
	\end{align*}
	This proves the result for every $z_0\in\mathbb V^1$.
	%
%	We finally consider an arbitrary $z_0\in\mathbb V^1$. By the regularizing property of the system and the continuity of the strong solution in $\mathbb V^1$, there exists $t_0\in(0,\sigma/2)$ sufficiently small such that
%	\begin{equation}\label{eq:regularize_velocity}
%		\widehat z:=\mathcal S_{t_0}(z_0,0,0)\in\mathbb V^4,
%		\qquad
%		\norm{\widehat z-z_0}_{\mathbb V^1}<\frac{\epsilon}{2}.
%	\end{equation}
%	Applying the result already proved for $\mathbb V^4$ initial data to $\widehat z$, with $\epsilon/2$ and $\sigma/2$, we find $\tau_1\in(0,\sigma/2)$ and a control
%	$f^1\in L^\infty(0,\tau_1;\mathcal H)$ such that
%	\begin{equation}\label{eq:regularized_control}
%		\norm{\mathcal S_{\tau_1}(\widehat z,0,f^1)-\widehat z-\widehat\eta}_{\mathbb V^1}
%		<\frac{\epsilon}{2}.
%	\end{equation}
%	Let $\tau=t_0+\tau_1<\sigma$ and define
%	\begin{equation*}
%		f(t)=
%		\begin{cases}
%			0, & t\in[0,t_0],\\
%			f^1(t-t_0), & t\in(t_0,t_0+\tau_1].
%		\end{cases}
%	\end{equation*}
%	Then combining \eqref{eq:regularize_velocity} and \eqref{eq:regularized_control} we have
%	\begin{align*}
%		\norm{\mathcal S_{\tau}(z_0,0,f)-z_0-\widehat\eta}_{\mathbb V^1}
%		&\leq
%		\norm{\mathcal S_{\tau_1}(\widehat z,0,f^1)-\widehat z-\widehat\eta}_{\mathbb V^1}
%		+\norm{\widehat z-z_0}_{\mathbb V^1}\\
%		&<\frac{\epsilon}{2}+\frac{\epsilon}{2}
%		=\epsilon.
%	\end{align*}
%	This proves the result for every $z_0\in\mathbb V^1$.
\end{proof}
\subsection{{Phase displacements}}
The main goal of this subsection is to prove the following result.
\begin{proposition}\label{small1}
Let $\mathcal H$ be the finite-dimensional space defined in \eqref{h}.
{Let $f=(f_u,f_\phi)\in L^\infty_{\mathrm{loc}}(\mathbb R_+;\mathbb V^0)$.}
For any
$\epsilon,\sigma>0$, $z_0\in\mathbb V^1$, and
$\xi\in\mathcal H_2(\infty)$, there exist $\tau\in(0,\sigma)$ and a control
$ \eta\in L^\infty(0,\tau;\mathcal H)$ such that
\begin{equation}\label{limit_phi1}
\norm{\mathcal S_\tau(z_0,0,{f+\eta})-(z_0+\widehat\xi)}_{\mathbb V^1}
<\epsilon,
\end{equation}
where $\widehat\xi:=(0,\xi)\in\mathbb V^1$.
\end{proposition}
To make the argument easier to follow, we split the proof into several
auxiliary lemmas, each one corresponding to one of the main ingredients of the
construction.

For $\rho\in\mathcal H_2(\infty)$ and $\zeta\in\mathscr G_1$, we set
{
\begin{equation*}
 \widehat\rho:=(0,\rho),\qquad \widehat\zeta:=(\zeta,0),
\end{equation*}
}
and define
\begin{equation*}
\mathcal Q(\widehat\rho):=(0,\Delta(\rho^3)),
\qquad
\mathcal M(\widehat\zeta,\widehat\rho):=(0,B_1(\zeta,\rho)).
\end{equation*}

For $N\ge0$, we denote by $\mathcal P_N$ the following property:
\begin{equation}\tag{$\mathcal P_N$}\label{propPN}
\begin{aligned}
&\text{For every force } f
 \text{ satisfying the hypotheses of \Cref{small1},}\\
&\text{and any }\epsilon,\sigma>0, \
 z_0\in{\mathbf H^2\times\dot H^4}, \
 \xi\in\mathcal H_2(N), \text{there exist }\tau\in(0,\sigma)\text{ and }\\
&
 \eta\in L^\infty(0,\tau;\mathcal H)\text{ such that } \eqref{limit_phi1} \text{ holds and } {\mathcal S_\tau(z_0,0,f+\eta)\in\mathbf H^2\times\dot H^4.}
\end{aligned}
\end{equation}
The proof of \Cref{small1} will be obtained by proving $\mathcal P_N$ by
induction on $N$. The next lemmas are auxiliary constructions that will be used in
the induction step.
\begin{lemma}\label{lem:cubic_generation}
	Let $N\ge1$ and assume that $\mathcal P_{N-1}$ holds. Let
	$\rho\in\mathcal H_2(N-1)$. Let $f$ satisfy the
hypotheses of \Cref{small1}.
	Then, for any $\epsilon,\sigma>0$ and any
	$z_0\in{\mathbf H^2\times\dot H^4}$, there exist $\tau\in(0,\sigma)$ and
	$ \eta\in L^\infty(0,\tau;\mathcal H)$ such that
	\begin{equation}\label{eq:cubic_generated}
		\left\|\mathcal S_\tau(z_0,0,f+\eta)-z_0-\mathcal Q(\widehat\rho)
		\right\|_{\mathbb V^1}<\epsilon,
	\end{equation}
	and
	\begin{equation}\label{eq:cubic_generated_regularity}
	 \mathcal S_\tau(z_0,0,f+\eta)\in\mathbf H^2\times\dot H^4.
	\end{equation}
\end{lemma}
\begin{proof}
	Let $\rho\in\mathcal H_2(N-1)$. We show that the direction
	$\mathcal Q(\widehat\rho)$ can be generated in arbitrarily small time.
	This will be obtained from the asymptotic limit \eqref{asymp_mixed} by
	taking the velocity direction equal to zero. More precisely, since
	$B_0(0,0)=0$, applying \eqref{asymp_mixed} with $\zeta=0$ and phase
	direction $\rho$ gives
	\begin{equation*}%\label{eq:cubic_from_mixed}
		\mathcal S_\delta(z_0,\delta^{-1/3}\widehat\rho,f)
		\longrightarrow z_0+\mathcal Q(\widehat\rho)
		\quad\text{in }\mathbb V^1,\quad\delta\to0^+.
	\end{equation*}
	On the other hand, by \eqref{eq:control-shift-identity}, applied with $\xi=\delta^{-1/3}\widehat\rho$ and $\eta=0$,
	we have
	\begin{equation*}%\label{eq:cubic_shift_identity}
		\mathcal S_\delta(z_0+\delta^{-1/3}\widehat\rho,0,f)
		=
		\mathcal S_\delta(z_0,\delta^{-1/3}\widehat\rho,f)
		+\delta^{-1/3}\widehat\rho.
	\end{equation*}
	Thus, we can choose $\delta\in(0,\sigma/3)$ and set
	$\beta:=\delta^{-1/3}$ so that the corresponding solution $\mathcal S(z_0+\delta^{-1/3}\widehat\rho,0,f)$ is well-defined
	on $[0,\delta]$ and
	\begin{equation*}
	 \|\mathcal S_\delta(z_0+\beta\widehat\rho,0,f)
	      -z_0-\mathcal Q(\widehat\rho)-\beta\widehat\rho\|_{\mathbb V^1}
	 <\frac{\epsilon}{8}.
	\end{equation*}
	Keeping $\beta$ fixed, write
	$y(t):=\mathcal S_t(z_0+\beta\widehat\rho,0,f)$ for $0\leq t\leq\delta$.
	By \Cref{wp1}, $y\in C([0,\delta];\mathbb V^1)
	\cap L^2(0,\delta;\mathbf H^2\times\dot H^4)$.
	By \Cref{lem:regular-times}, we can choose $\tau_2\in(0,\delta)$,
	sufficiently close to $\delta$, such that
	\begin{equation*}
	 \widetilde z:=y(\tau_2)\in\mathbf H^2\times\dot H^4,\qquad
	 \|y(\tau_2)-y(\delta)\|_{\mathbb V^1}<\frac{\epsilon}{8}.
	\end{equation*}
	Consequently,
	\begin{align}\label{eq:cubic_free_step}
	 \|\widetilde z-z_0-\mathcal Q(\widehat\rho)-\beta\widehat\rho\|_{\mathbb V^1}&\leq \|\widetilde z-y(\delta)\|_{\mathbb V^1}
	   +\|y(\delta)-z_0-\mathcal Q(\widehat\rho)-\beta\widehat\rho\|_{\mathbb V^1}
	   \notag\\
	 &<\frac{\epsilon}{8}+\frac{\epsilon}{8}=\frac{\epsilon}{4}.
	\end{align}
	We now keep $\tau_2$ fixed.
	Since $-\beta\rho\in\mathcal H_2(N-1)$, the induction hypothesis
	$\mathcal P_{N-1}$, applied to the initial condition
	$\widetilde z$ and the phase
	direction $-\beta\rho$ with $f_{\tau_2}$ in place of $f$, gives $\tau_3\in(0,\sigma/3)$ and
	$\eta^3\in L^\infty(0,\tau_3;\mathcal H)$ such that
	\begin{equation}\label{eq:remove_cubic_shift}
		\left\|
		\mathcal S_{\tau_3}
		\left(
		\widetilde z,
		0,f_{\tau_2}+\eta^3
		\right)
		-
		\widetilde z
		+\beta\widehat\rho
		\right\|_{\mathbb V^1}
		<\frac{\epsilon}{4}.
	\end{equation}
	Combining \eqref{eq:cubic_free_step} and
	\eqref{eq:remove_cubic_shift}, we obtain
	\begin{align}
		\notag
		&\left\|
		\mathcal S_{\tau_3}
		\left(
		\widetilde z,
		0,f_{\tau_2}+\eta^3
		\right)
		-z_0-\mathcal Q(\widehat\rho)
		\right\|_{\mathbb V^1}\\ \label{eq:cubic_ideal_final}
		&\qquad \le
		\left\|
		\mathcal S_{\tau_3}
		\left(
		\widetilde z,
		0,f_{\tau_2}+\eta^3
		\right)
		-
		\widetilde z
		+\beta\widehat\rho
		\right\|_{\mathbb V^1}+
		\left\|
		\widetilde z
		-z_0-\mathcal Q(\widehat\rho)-\beta\widehat\rho
		\right\|_{\mathbb V^1}
		<\frac{\epsilon}{2}.
	\end{align}
By \Cref{cont} and \eqref{eq:cont_stability}, applied to the
	solution $t\mapsto\mathcal S_t(\widetilde z,0,f_{\tau_2}+\eta^3)$
	on $[0,\tau_3]$, there exist $r_3>0$ and $C_3\ge1$ such that,
	for $s\ge0$, whenever
	\begin{equation*}
		\|v_0-\widetilde z\|_{\mathbb V^1}
		+\|f_{\tau_2+s}-f_{\tau_2}\|_{L^2(0,\tau_3;\mathbb V^0)}
		<r_3,
	\end{equation*}
	the solution $\mathcal S(v_0,0,f_{\tau_2+s}+\eta^3)$ is well-defined on
	$[0,\tau_3]$ and
	\begin{equation}\label{eq:cubic_stability3}
		\begin{aligned}
		&\left\|
		\mathcal S_{\tau_3}(v_0,0,f_{\tau_2+s}+\eta^3)
		-\mathcal S_{\tau_3}(\widetilde z,0,f_{\tau_2}+\eta^3)
		\right\|_{\mathbb V^1}\\
		&\quad\le C_3\Bigl(
		\|v_0-\widetilde z\|_{\mathbb V^1}
		+\|f_{\tau_2+s}-f_{\tau_2}\|_{L^2(0,\tau_3;\mathbb V^0)}
		\Bigr).
		\end{aligned}
	\end{equation}
Set
	\begin{equation}\label{eq:cubic_gamma}
		\gamma:=
		\min\left\{
		\frac{r_3}{2},\frac{\epsilon}{8C_3}
		\right\}.
	\end{equation}
We now apply \Cref{cont} and \eqref{eq:cont_stability} to the
	solution $t\mapsto\mathcal S_t(z_0+\beta\widehat\rho,0,f)$
	on $[0,\tau_2]$. There exist $r_2>0$ and $C_2\ge1$ such that,
	for $s\ge0$, whenever
	\begin{equation*}
	 \|v_0-(z_0+\beta\widehat\rho)\|_{\mathbb V^1}
	 +\|f_s-f\|_{L^2(0,\tau_2;\mathbb V^0)}<r_2,
	\end{equation*}
	the solution $\mathcal S(v_0,0,f_s)$ is well-defined on $[0,\tau_2]$
	and
	\begin{equation}\label{eq:cubic_stability2}
		\begin{aligned}
		\left\|
		\mathcal S_{\tau_2}(v_0,0,f_s)
		-\widetilde z
		\right\|_{\mathbb V^1} \le C_2\Bigl(
		\|v_0-(z_0+\beta\widehat\rho)\|_{\mathbb V^1}
		+\|f_s-f\|_{L^2(0,\tau_2;\mathbb V^0)}
		\Bigr).
		\end{aligned}
	\end{equation}
Choose $0<\alpha<
		\min\left\{
		\frac{r_2}{2},\frac{\gamma}{2C_2}
		\right\}.$
	Since $f\in L^2_{\mathrm{loc}}(\mathbb R_+;\mathbb V^0)$,
	\begin{equation*}
	 \|f_s-f\|_{L^2(0,\tau_2+\tau_3;\mathbb V^0)}
	 \longrightarrow0\qquad\text{as }s\to0^+.
	\end{equation*}
Since $\beta\rho\in\mathcal H_2(N-1)$, the induction hypothesis
	$\mathcal P_{N-1}$, applied to the initial condition $z_0$ and phase
	direction $\beta\rho$ with a sufficiently small time horizon, gives $\tau_1\in(0,\sigma/3)$ and
	$ \eta^1\in L^\infty(0,\tau_1;\mathcal H)$ such that
	\begin{equation}\label{eq:prepare_cubic}
		\left\|
		\mathcal S_{\tau_1}(z_0,0,f+\eta^1)
		-z_0-\beta\widehat\rho
		\right\|_{\mathbb V^1}
		<\alpha
	\text{ and }
	 \|f_{\tau_1}-f\|_{L^2(0,\tau_2+\tau_3;\mathbb V^0)}
	 <\alpha.
	\end{equation}
Set $z^1:=\mathcal S_{\tau_1}(z_0,0,f+\eta^1)$.
	Since $2\alpha<r_2$, \Cref{cont} gives
	$\mathcal T(z^1,0,f_{\tau_1})>\tau_2$.
	Set $z^2:=\mathcal S_{\tau_2}(z^1,0,f_{\tau_1})$. By \eqref{eq:cubic_stability2} and \eqref{eq:prepare_cubic},
	\begin{align}\label{eq:cubic_after_free}
	%	\notag
		&\left\|
		z^2
		-\widetilde z
		\right\|_{\mathbb V^1}\le
		C_2\Bigl(
		\left\|
		z^1
		-z_0-\beta\widehat\rho
		\right\|_{\mathbb V^1}
		+\|f_{\tau_1}-f\|_{L^2(0,\tau_2;\mathbb V^0)}
		\Bigr)
		<2C_2\alpha<\gamma.
	\end{align}
By a change of variables,
	\begin{equation*}
	 \|f_{\tau_1+\tau_2}-f_{\tau_2}\|_{L^2(0,\tau_3;\mathbb V^0)}
	 =\|f_{\tau_1}-f\|_{L^2(\tau_2,\tau_2+\tau_3;\mathbb V^0)}
	 <\alpha,
	\end{equation*}
and in view of \eqref{eq:cubic_gamma}, the last inequality \eqref{eq:cubic_after_free} also gives
	\begin{align*}
		&\left\|
		z^2
		-\widetilde z
		\right\|_{\mathbb V^1}+
		\|f_{\tau_1+\tau_2}-f_{\tau_2}\|_{L^2(0,\tau_3;\mathbb V^0)}
		<\gamma+\alpha<2\gamma\le r_3.
	\end{align*}
Therefore, by \Cref{cont} and
	\eqref{eq:cubic_stability3},
	$\mathcal T(z^2,0,f_{\tau_1+\tau_2}+\eta^3)>\tau_3$, and
	\begin{align}\label{eq:cubic_last_stability}
		\notag
		&\left\|
		\mathcal S_{\tau_3}
		\left(
		z^2,
		0,f_{\tau_1+\tau_2}+\eta^3
		\right)
		-
		\mathcal S_{\tau_3}
		\left(
		\widetilde z,
		0,f_{\tau_2}+\eta^3
		\right)
		\right\|_{\mathbb V^1}\\
		&\qquad \le
		C_3\Bigl(
		\left\|
		z^2
		-\widetilde z
		\right\|_{\mathbb V^1}+
		\|f_{\tau_1+\tau_2}-f_{\tau_2}\|_{L^2(0,\tau_3;\mathbb V^0)}
		\Bigr)
		<C_3(\gamma+\alpha)<2C_3\gamma\le\frac{\epsilon}{4}.
	\end{align}
{Finally, set $T:=\tau_1+\tau_2+\tau_3<\sigma$ and define
	\begin{equation*}
		\eta(t)=
		\begin{cases}
			\eta^1(t),&t\in(0,\tau_1),\\
			0,&t\in(\tau_1,\tau_1+\tau_2),\\
			\eta^3(t-\tau_1-\tau_2),&t\in(\tau_1+\tau_2,T),\\
			0,&t\ge T.
		\end{cases}
	\end{equation*}
	Applying \Cref{lem:concatenation} to the
	solutions from $z_0$, $z^1$, and $z^2$ with controls
	$\eta^1$, $0$, and $\eta^3$, respectively, gives
	$\mathcal T(z_0,0,f+\eta)>T$ and
$\mathcal S_T(z_0,0,f+\eta)
		=\mathcal S_{\tau_3}(z^2,0,f_{\tau_1+\tau_2}+\eta^3).$
	Thus \eqref{eq:cubic_ideal_final} and \eqref{eq:cubic_last_stability}
	give}
	\begin{align*}
		&\left\|
		\mathcal S_{T}(z_0,0,f+\eta)
		-z_0-\mathcal Q(\widehat\rho)
		\right\|_{\mathbb V^1}\\
		&\le
		\left\|
		\mathcal S_{\tau_3}
		\left(
		z^2,
		0,f_{\tau_1+\tau_2}+\eta^3
		\right)
		-
		\mathcal S_{\tau_3}
		\left(
		\widetilde z,
		0,f_{\tau_2}+\eta^3
		\right)
		\right\|_{\mathbb V^1}\\
		&\quad+
		\left\|
		\mathcal S_{\tau_3}
		\left(
		\widetilde z,
		0,f_{\tau_2}+\eta^3
		\right)
		-z_0-\mathcal Q(\widehat\rho)
		\right\|_{\mathbb V^1}\\
		&<\frac{\epsilon}{4}+\frac{\epsilon}{2}<\epsilon.
	\end{align*}
	{Set $z(t):=\mathcal S_t(z_0,0,f+\eta)$ for $0\le t\le T$.
	By \eqref{eq:concatenated-regularity}, $z\in C([0,T];\mathbb V^1)
	\cap L^2(0,T;\mathbf H^2\times\dot H^4)$.
	Keeping $\eta$ fixed, by continuity and \Cref{lem:regular-times},
	we may choose $\tau\in(\tau_1+\tau_2,T)$, sufficiently close to $T$, such that
	\begin{equation*}
	 z(\tau)\in\mathbf H^2\times\dot H^4,\qquad
	 \|z(\tau)-z(T)\|_{\mathbb V^1}<\frac{\epsilon}{4}.
	\end{equation*}
	Hence,
	\begin{align*}
	 \|z(\tau)-z_0-\mathcal Q(\widehat\rho)\|_{\mathbb V^1}
	 &\leq \|z(\tau)-z(T)\|_{\mathbb V^1}
	       +\|z(T)-z_0-\mathcal Q(\widehat\rho)\|_{\mathbb V^1}\\
	 &<\frac{\epsilon}{4}+\frac{3\epsilon}{4}=\epsilon.
	\end{align*}
	Set $\eta^\tau:=\eta|_{(0,\tau)}\in L^\infty(0,\tau;\mathcal H)$.
	By \Cref{cor:unique-strong-common},
$ \mathcal S_\tau(z_0,0,f+\eta^\tau)=z(\tau).$ 
	Thus \eqref{eq:cubic_generated} and \eqref{eq:cubic_generated_regularity}
	hold with control $\eta^\tau$.}
	This completes the proof.
\end{proof}
The next construction produces a combination of the cubic term and the transport term.
\begin{lemma}\label{lem:mixed_generation}
	Let $N\ge1$ and assume that $\mathcal P_{N-1}$ holds. Let
	$\zeta\in\mathscr G_1$ and $\rho\in\mathcal H_2(N-1)$. Let $f$
	satisfy the hypotheses of \Cref{small1}.
	Then, for any
	$\epsilon,\sigma>0$ and any $z_0\in{\mathbf H^2\times\dot H^4}$, there exist
	$\tau\in(0,\sigma)$ and $ \eta\in L^\infty(0,\tau;\mathcal H)$ such that
	\begin{equation*}%\label{eq:mixed_goal}
		\left\|
		\mathcal S_\tau(z_0,0,f+\eta)
		-
		z_0
		+
		\mathcal Q(\widehat\rho)
		-
		\mathcal M(\widehat\zeta,\widehat\rho)
		\right\|_{\mathbb V^1}
		<\epsilon,
	\end{equation*}
	{and
	\begin{equation}\label{eq:mixed_generated_regularity}
	 \mathcal S_\tau(z_0,0,f+\eta)\in\mathbf H^2\times\dot H^4.
	\end{equation}}
\end{lemma}
\begin{proof}
	Let $\zeta\in\mathscr G_1$ and $\rho\in\mathcal H_2(N-1)$. We first show
	that the mixed direction
	\begin{equation*}
		-\mathcal Q(\widehat\rho)+\mathcal M(\widehat\zeta,\widehat\rho)
	\end{equation*}
	can be generated in arbitrarily small time. Let $\epsilon,\sigma>0$ and
	$z_0\in{\mathbf H^2\times\dot H^4}$ be given.

	Since $\zeta\in\mathscr G_1$, \Cref{lem:safe_modes} gives
	$B_0(\zeta,\zeta)=0$. Therefore, we can apply the asymptotic limit
	\eqref{asymp_mixed} with velocity direction $\zeta$ and phase direction
	$-\rho$. Since
	\begin{equation*}
		\Delta((-\rho)^3)-B_1(\zeta,-\rho)
		=
		-\Delta(\rho^3)+B_1(\zeta,\rho),
	\end{equation*}
	we obtain
	\begin{equation*}%\label{eq:mixed_limit}
		\mathcal S_\delta
		\bigl(
		z_0,
		\delta^{-2/3}\widehat\zeta-\delta^{-1/3}\widehat\rho,
		f
		\bigr)
		\longrightarrow
		z_0
		-\mathcal Q(\widehat\rho)
		+\mathcal M(\widehat\zeta,\widehat\rho)
	\end{equation*}
	in $\mathbb V^1$ as $\delta\to0^+$.

	On the other hand, by identity {\eqref{eq:control-shift-identity}}, we
	have
	\begin{equation*}%\label{eq:mixed_shift_identity}
		\mathcal S_\delta
		\bigl(
		z_0+\delta^{-2/3}\widehat\zeta-\delta^{-1/3}\widehat\rho,
		0,
		f
		\bigr)
		=
		\mathcal S_\delta
		\bigl(
		z_0,
		\delta^{-2/3}\widehat\zeta-\delta^{-1/3}\widehat\rho,
		f
		\bigr)
		+
		\delta^{-2/3}\widehat\zeta
		-
		\delta^{-1/3}\widehat\rho .
	\end{equation*}
	Thus, we may choose $\delta\in(0,\sigma/5)$ and set
	$\alpha:=\delta^{-2/3}$ and $\beta:=\delta^{-1/3}$ so that the
	corresponding solution is well-defined on $[0,\delta]$ and
	\begin{equation*}
	 \left\|\mathcal S_\delta
	 (z_0+\alpha\widehat\zeta-\beta\widehat\rho,0,f)
	 -z_0-\alpha\widehat\zeta+\beta\widehat\rho
	 +\mathcal Q(\widehat\rho)
	 -\mathcal M(\widehat\zeta,\widehat\rho)\right\|_{\mathbb V^1}
	 <\frac{\epsilon}{12}.
	\end{equation*}
	Keeping $\alpha$ and $\beta$ fixed, write
	$y_3(t):=\mathcal S_t(z_0+\alpha\widehat\zeta-\beta\widehat\rho,0,f)$
	for $0\leq t\leq\delta$.
	By \Cref{wp1}, $y_3\in C([0,\delta];\mathbb V^1)
	\cap L^2(0,\delta;\mathbf H^2\times\dot H^4)$.
	By \Cref{lem:regular-times}, we may choose $\tau_3\in(0,\delta)$,
	sufficiently close to $\delta$, such that
	\begin{equation*}
	 y_3(\tau_3)\in\mathbf H^2\times\dot H^4 \quad \text{ and }
	 \|y_3(\tau_3)-y_3(\delta)\|_{\mathbb V^1}<\frac{\epsilon}{12}.
	\end{equation*}
	Consequently,
	\begin{align}\label{eq:mixed_free_step}
	 &\|\mathcal S_{\tau_3}(z_0+\alpha\widehat\zeta-\beta\widehat\rho,0,f)
	       -z_0-\alpha\widehat\zeta+\beta\widehat\rho
	       +\mathcal Q(\widehat\rho)-\mathcal M(\widehat\zeta,\widehat\rho)\|_{\mathbb V^1}
	       \notag\\
	 &\quad\leq \|y_3(\tau_3)-y_3(\delta)\|_{\mathbb V^1}
	       +\|y_3(\delta)-z_0-\alpha\widehat\zeta+\beta\widehat\rho
	       +\mathcal Q(\widehat\rho)-\mathcal M(\widehat\zeta,\widehat\rho)\|_{\mathbb V^1}
	       \notag\\
	 &\quad<\frac{\epsilon}{12}+\frac{\epsilon}{12}=\frac{\epsilon}{6}.
	\end{align}
	We now keep $\tau_3$, $\alpha$, and $\beta$ fixed.
	Since $\mathcal H_1=\operatorname{span}\mathscr G_1$ and
	$\zeta\in\mathscr G_1$, we have $-\alpha\zeta\in\mathcal H_1$, hence
	$-\alpha\widehat\zeta\in\mathcal H$. Therefore, the direct construction
	\eqref{eq:direct_maneuver}, applied to the initial condition
	$ \mathcal S_{\tau_3}
		\bigl(
		z_0+\alpha\widehat\zeta-\beta\widehat\rho,
		0,
		f
		\bigr),$
	with prescribed force $f_{\tau_3}$ {(using
	\eqref{eq:control-force-translations}), gives $b\in(0,\sigma/5)$ and
	$\eta^4\in L^\infty(0,b;\mathcal H)$ such that
	\begin{equation*}
	 \left\|\mathcal S_b
	 \left(\mathcal S_{\tau_3}
	 (z_0+\alpha\widehat\zeta-\beta\widehat\rho,0,f),
	 0,f_{\tau_3}+\eta^4\right)
	 -\mathcal S_{\tau_3}
	 (z_0+\alpha\widehat\zeta-\beta\widehat\rho,0,f)
	 +\alpha\widehat\zeta\right\|_{\mathbb V^1}
	 <\frac{\epsilon}{12}.
	\end{equation*}
	Keeping $\eta^4$ fixed, write
	\begin{equation*}
	 y_4(t):=\mathcal S_t
	 \bigl(y_3(\tau_3),0,f_{\tau_3}+\eta^4\bigr),\qquad 0\leq t\leq b.
	\end{equation*}
	By \Cref{wp1}, $y_4\in C([0,b];\mathbb V^1)\cap
	L^2(0,b;\mathbf H^2\times\dot H^4)$.
	By continuity and \Cref{lem:regular-times}, we may choose
	$\tau_4\in(0,b)$, sufficiently close to $b$, such that
	\begin{equation*}
	 y_4(\tau_4)\in\mathbf H^2\times\dot H^4\quad \text{ and } \quad
	 \|y_4(\tau_4)-y_4(b)\|_{\mathbb V^1}<\frac{\epsilon}{12}.
	\end{equation*}
	Consequently,
	\begin{align}\label{eq:rm1}
	 &\left\|\mathcal S_{\tau_4}\left(
	       \mathcal S_{\tau_3}(z_0+\alpha\widehat\zeta-\beta\widehat\rho,0,f),
	       0,f_{\tau_3}+\eta^4\right)-\mathcal S_{\tau_3}
	       (z_0+\alpha\widehat\zeta-\beta\widehat\rho,0,f)
	       +\alpha\widehat\zeta\right\|_{\mathbb V^1}\notag\\
	 &\quad\leq \|y_4(\tau_4)-y_4(b)\|_{\mathbb V^1}
	       +\|y_4(b)-y_4(0)+\alpha\widehat\zeta\|_{\mathbb V^1}
	       \notag\\
	 &\quad<\frac{\epsilon}{12}+\frac{\epsilon}{12}=\frac{\epsilon}{6}.
	\end{align}
	By \Cref{cor:unique-strong-common},
	\begin{equation*}
	 \mathcal S_{\tau_4}\bigl(y_3(\tau_3),0,
	 f_{\tau_3}+\eta^4|_{(0,\tau_4)}\bigr)=y_4(\tau_4).
	\end{equation*}
	We continue to denote by $\eta^4$ its restriction to $(0,\tau_4)$.}
	Since $\beta\rho\in\mathcal H_2(N-1)$, the induction hypothesis
	$\mathcal P_{N-1}$, applied to this initial condition and the phase
	direction $\beta\rho$ with $f_{\tau_3+\tau_4}$ in place of $f$, gives $\tau_5\in(0,\sigma/5)$ and
	$\eta^5\in L^\infty(0,\tau_5;\mathcal H)$ such that
	\begin{align}\label{eq:remove_mixed_phase}
		&\Bigl\|
		\mathcal S_{\tau_5}
		\left(
		\mathcal S_{\tau_4}
		\left(
		\mathcal S_{\tau_3}
		\bigl(
		z_0+\alpha\widehat\zeta-\beta\widehat\rho,
		0,
		f
		\bigr),
		0,
		f_{\tau_3}+\eta^4
		\right),
		0,
		f_{\tau_3+\tau_4}+\eta^5
		\right)
		\\
		\notag
		&\qquad
		-
		\mathcal S_{\tau_4}
		\left(
		\mathcal S_{\tau_3}
		\bigl(
		z_0+\alpha\widehat\zeta-\beta\widehat\rho,
		0,
		f
		\bigr),
		0,
		f_{\tau_3}+\eta^4
		\right)
		-
		\beta\widehat\rho
		\Bigr\|_{\mathbb V^1}
		<\frac{\epsilon}{6}.
	\end{align}
%	\begin{align}
%		%\label{eq:remove_mixed_phase}
%		\notag
%		&\left\|
%		\mathcal S_{\tau_5}
%		\left(
%		\mathcal S_{\tau_4}
%		\left(
%		\mathcal S_{\tau_3}
%		\bigl(
%		z_0+\alpha\widehat\zeta-\beta\widehat\rho,
%		0,
%		0
%		\bigr),
%		0,
%		f^4
%		\right),
%		0,
%		f^5
%		\right)
%		\\
%		\notag
%		&\qquad\qquad
%		-
%		\mathcal S_{\tau_4}
%		\left(
%		\mathcal S_{\tau_3}
%		\bigl(
%		z_0+\alpha\widehat\zeta-\beta\widehat\rho,
%		0,
%		0
%		\bigr),
%		0,
%		f^4
%		\right)
%		-
%		\beta\widehat\rho
%		\right\|_{\mathbb V^1}
%		<\frac{\epsilon}{6}.
%	\end{align}
	Combining \eqref{eq:mixed_free_step},
	\eqref{eq:rm1}, and
	\eqref{eq:remove_mixed_phase}, we obtain
\begin{align}\label{eq:mixed_ideal_final}
	&\left\|
	\mathcal S_{\tau_5}
	\left(
	\mathcal S_{\tau_4}
	\left(
	\mathcal S_{\tau_3}
	\bigl(
	z_0+\alpha\widehat\zeta-\beta\widehat\rho,
	0,
	f
	\bigr),
	0,
	f_{\tau_3}+\eta^4
	\right),
	0,
	f_{\tau_3+\tau_4}+\eta^5
	\right)
	-
	z_0
	+
	\mathcal Q(\widehat\rho)
	-
	\mathcal M(\widehat\zeta,\widehat\rho)
	\right\|_{\mathbb V^1}
	<\frac{\epsilon}{2}.
\end{align}
We next fix the stability neighborhoods associated with the three reference trajectories
\[
\begin{split}
	&t\mapsto
	\mathcal S_t
	\bigl(
	z_0+\alpha\widehat\zeta-\beta\widehat\rho,0,f
	\bigr),\, t\in[0,\tau_3],\\
	&t\mapsto
	\mathcal S_t
	\left(
	\mathcal S_{\tau_3}
	\bigl(
	z_0+\alpha\widehat\zeta-\beta\widehat\rho,0,f
	\bigr),
	0,f_{\tau_3}+\eta^4
	\right),\, t\in[0,\tau_4],\\
	&t\mapsto
	\mathcal S_t
	\left(
	\mathcal S_{\tau_4}
	\left(
	\mathcal S_{\tau_3}
	\bigl(
	z_0+\alpha\widehat\zeta-\beta\widehat\rho,0,f
	\bigr),
	0,f_{\tau_3}+\eta^4
	\right),
	0,f_{\tau_3+\tau_4}+\eta^5
	\right),\,  t\in[0,\tau_5].
\end{split}
\]
For each of these three reference solutions, write
$\widetilde z(t)=\mathcal S_t(\widetilde z_0,0,\widetilde f)$ on $[0,\tau_j]$,
with $j=3,4,5$, respectively.
By {\Cref{cont} and \eqref{eq:cont_stability}}, there exist
$r_j>0$ and $C_j\ge1$ such that, for
$\widehat z_0\in\mathbb V^1$ and
$\widehat f\in L^\infty_{\mathrm{loc}}(\mathbb R_+;\mathbb V^0)$,
the condition
\begin{equation}\label{eq:mixed-stability-condition}
 \|\widehat z_0-\widetilde z_0\|_{\mathbb V^1}
 +\|\widehat f-\widetilde f\|_{L^2(0,\tau_j;\mathbb V^0)}
 <r_j
\end{equation}
guarantees that $\widehat z(t):=\mathcal S_t(\widehat z_0,0,\widehat f)$
exists on $[0,\tau_j]$ and satisfies
\begin{align}\label{eq:mixed-stability-estimate}
 &\sup_{t\in[0,\tau_j]}
 \|\widehat z(t)-\widetilde z(t)\|_{\mathbb V^1}\le C_j\Bigl(
 \|\widehat z_0-\widetilde z_0\|_{\mathbb V^1}
 +\|\widehat f-\widetilde f\|_{L^2(0,\tau_j;\mathbb V^0)}
 \Bigr).
\end{align}
All the constants $r_j,C_j$ are fixed at this point.
%Here, $r_j$ and $C_j$ may depend on the corresponding fixed reference trajectory, time interval, and control. In particular, all the constants $r_j,C_j$, $j=3,4,5$, are fixed before choosing the approximation parameter below.
%We next fix the stability neighborhoods corresponding to the three solutions constructed above. 
Choose
	\begin{equation}\label{eq:mixed_gamma}
		0<\gamma<
		\min\left\{
		\frac{r_3}{2},
		\frac{r_4}{3C_3},
		\frac{r_5}{4C_4C_3},
		\frac{\epsilon}{8C_5C_4C_3}
		\right\}.
	\end{equation}
	By continuity of translations in $L^2$, choose $\sigma_0>0$ such that
	\begin{equation*}
	 \|f_s-f\|_{L^2(0,\tau_3+\tau_4+\tau_5;\mathbb V^0)}<\gamma
	 \qquad\text{for every }0<s<\sigma_0.
	\end{equation*}
	Since $-\beta\rho\in\mathcal H_2(N-1)$, the induction hypothesis
	$\mathcal P_{N-1}$, applied to the initial condition
	$z_0+\alpha\widehat\zeta$ and the phase direction $-\beta\rho$, gives
	$ \tau_2\in(0,\min\{\sigma/5,\sigma_0/2\})$ and
	$\eta^2\in L^\infty(0,\tau_2;\mathcal H)$ such that
	\begin{equation}\label{eq:prepare_mixed_phase}
		\left\|
		\mathcal S_{\tau_2}
		\bigl(
		z_0+\alpha\widehat\zeta,
		0,
		f+\eta^2
		\bigr)
		-
		z_0
		-
		\alpha\widehat\zeta
		+
		\beta\widehat\rho
		\right\|_{\mathbb V^1}
		<\frac{\gamma}{2}.
	\end{equation}
	We now fix $\tau_2$ and $\eta^2$. {Applying
	\Cref{cont} and \eqref{eq:cont_stability} to}
	\begin{equation*}
	 t\mapsto\mathcal S_t(z_0+\alpha\widehat\zeta,0,f+\eta^2),
	 \qquad t\in[0,\tau_2],
	\end{equation*}
	gives $r_2>0$ and $C_2\ge1$ for which
	\eqref{eq:mixed-stability-condition}--\eqref{eq:mixed-stability-estimate}
	hold with $j=2$. Choose
	\begin{equation}\label{eq:mixed_eta}
		0<\kappa<
		\min\left\{
		\frac{r_2}{2},\frac{\gamma}{4C_2}
		\right\}.
	\end{equation}
	Choose $\sigma_1\in(0,\min\{\sigma/5,\sigma_0/2\})$ such that
	\begin{equation*}
	 \|f_s-f\|_{L^2(0,\tau_2;\mathbb V^0)}<\kappa
	 \qquad\text{for every }0<s<\sigma_1.
	\end{equation*}
	Since $\alpha\zeta\in\mathcal H_1$, hence
	$\alpha\widehat\zeta\in\mathcal H$, the direct construction
	\eqref{eq:direct_maneuver} with prescribed force $f$
	and time bound $\sigma_1$ gives
	$ \tau_1\in(0,\sigma_1)$ and
	$\eta^1\in L^\infty(0,\tau_1;\mathcal H)$ such that
	\begin{equation}\label{eq:prepare_mixed_velocity}
		\left\|
		\mathcal S_{\tau_1}(z_0,0,f+\eta^1)
		-
		z_0
		-
		\alpha\widehat\zeta
		\right\|_{\mathbb V^1}
		<\kappa.
	\end{equation}
	{Define $\eta^0$ by
	\begin{equation*}
		\eta^0(t)=
		\begin{cases}
			\eta^1(t), & t\in(0,\tau_1),\\
			\eta^2(t-\tau_1), & t\in(\tau_1,\tau_1+\tau_2),\\
			0, & t\ge\tau_1+\tau_2.
		\end{cases}
	\end{equation*}}
	By \eqref{eq:mixed_eta}, \eqref{eq:prepare_mixed_velocity},
	and the choice of $\sigma_1$,
	\begin{align*}
	 &\|\mathcal S_{\tau_1}(z_0,0,f+\eta^1)
	       -z_0-\alpha\widehat\zeta\|_{\mathbb V^1}+\|f_{\tau_1}-f\|_{L^2(0,\tau_2;\mathbb V^0)}
	 <2\kappa<r_2.
	\end{align*}
	Hence \eqref{eq:mixed-stability-condition} holds with $j=2$,
	$\widehat z_0=\mathcal S_{\tau_1}(z_0,0,f+\eta^1)$ and
	$\widehat f=f_{\tau_1}+\eta^2$.
	The solution exists on $[0,\tau_2]$.
	{By \Cref{lem:concatenation} with controls
	$\eta^1,\eta^2$,}
	\eqref{eq:mixed-stability-estimate} and
	\eqref{eq:prepare_mixed_phase} give
	\begin{align}\notag
		\notag
		&\left\|
		\mathcal S_{\tau_1+\tau_2}(z_0,0,f+\eta^0)
		-
		z_0
		-
		\alpha\widehat\zeta
		+
		\beta\widehat\rho
		\right\|_{\mathbb V^1}
		\\
		\notag
		&\le
		\left\|
		\mathcal S_{\tau_2}
		\bigl(
		\mathcal S_{\tau_1}(z_0,0,f+\eta^1),
		0,
		f_{\tau_1}+\eta^2
		\bigr)
		-
		\mathcal S_{\tau_2}
		\bigl(
		z_0+\alpha\widehat\zeta,
		0,
		f+\eta^2
		\bigr)
		\right\|_{\mathbb V^1}
		\\ \notag
		&\quad + 
		\left\|
		\mathcal S_{\tau_2}
		\bigl(
		z_0+\alpha\widehat\zeta,
		0,
		f+\eta^2
		\bigr)
		-
		z_0
		-
		\alpha\widehat\zeta
		+
		\beta\widehat\rho
		\right\|_{\mathbb V^1}
		\\ \label{eq:prepare_mixed_total}
		&<
		C_2\Bigl(\kappa+
		\|f_{\tau_1}-f\|_{L^2(0,\tau_2;\mathbb V^0)}\Bigr)
		+\frac{\gamma}{2}
		<2C_2\kappa+\frac{\gamma}{2}
		<\gamma.
	\end{align}
	Define $\widetilde\eta$ by
	\begin{equation*}
		\widetilde\eta(t)=
		\begin{cases}
			\eta^0(t), & t\in(0,\tau_1+\tau_2),\\
			0, & t\ge\tau_1+\tau_2.
		\end{cases}
	\end{equation*}
Set $s:=\tau_1+\tau_2<\sigma_0$. For $a,b\ge0$ with
	$a+b\le\tau_3+\tau_4+\tau_5$, a change of variables gives
	\begin{equation}\label{eq:mixed-force-translation-bound}
	 \|f_{s+a}-f_a\|_{L^2(0,b;\mathbb V^0)}
	 =\|f_s-f\|_{L^2(a,a+b;\mathbb V^0)}<\gamma.
	\end{equation}
By \eqref{eq:prepare_mixed_total} and
	\eqref{eq:mixed-force-translation-bound},
	\begin{align*}
	 &\|\mathcal S_s(z_0,0,f+\eta^0)
	       -z_0-\alpha\widehat\zeta+\beta\widehat\rho\|_{\mathbb V^1}+\|f_s-f\|_{L^2(0,\tau_3;\mathbb V^0)}
	 <2\gamma<r_3.
	\end{align*}
	By \Cref{lem:concatenation} with controls
	$\eta^0,0$, \eqref{eq:mixed-stability-estimate} with $j=3$ gives
	\begin{align}\label{eq:mixed_after_free}
		\left\|
		\mathcal S_{\tau_1+\tau_2+\tau_3}(z_0,0,f+\widetilde\eta)
		-
		\mathcal S_{\tau_3}
		\bigl(
		z_0+\alpha\widehat\zeta-\beta\widehat\rho,
		0,
		f
		\bigr)
		\right\|_{\mathbb V^1}
		<
		C_3(\gamma+\gamma)=2C_3\gamma.
	\end{align}

	For $j=4$, \eqref{eq:mixed_after_free} and
	\eqref{eq:mixed-force-translation-bound} give $2C_3\gamma+\gamma\le3C_3\gamma<r_4$. Thus the solution from the state at time $s+\tau_3$ with force
	$f_{s+\tau_3}+\eta^4$ exists on $[0,\tau_4]$.
	By \eqref{eq:mixed-stability-estimate},

	\begin{align}\notag
	&\left\|
		\mathcal S_{\tau_4}
		\left(
		\mathcal S_{\tau_1+\tau_2+\tau_3}(z_0,0,f+\widetilde\eta),
		0,
		f_{s+\tau_3}+\eta^4
		\right)
		-
		\mathcal S_{\tau_4}
		\left(
		\mathcal S_{\tau_3}
		\bigl(
		z_0+\alpha\widehat\zeta-\beta\widehat\rho,
		0,
		f
		\bigr),
		0,
		f_{\tau_3}+\eta^4
		\right)
		\right\|_{\mathbb V^1}
		\\ \label{eq:mixed_after_velocity_removal}
		&\qquad <
		C_4(2C_3\gamma+\gamma)\le3C_4C_3\gamma.
	\end{align}
Similarly, \eqref{eq:mixed_after_velocity_removal} and
	\eqref{eq:mixed-force-translation-bound} give
	$3C_4C_3\gamma+\gamma\le4C_4C_3\gamma<r_5.$
	Applying \eqref{eq:mixed-stability-estimate} with $j=5$,
	initial state at time $s+\tau_3+\tau_4$, and force
	$f_{s+\tau_3+\tau_4}+\eta^5$ yields
		\begin{align}\notag 
			&\left\|
		\mathcal S_{\tau_5}
		\left(
		\mathcal S_{\tau_4}
		\left(
		\mathcal S_{\tau_1+\tau_2+\tau_3}(z_0,0,f+\widetilde\eta),
		0,
		f_{s+\tau_3}+\eta^4
		\right),
		0,
		f_{s+\tau_3+\tau_4}+\eta^5
		\right)
		\right.
		\\
		\notag
		&\qquad\left.
		-
		\mathcal S_{\tau_5}
		\left(
		\mathcal S_{\tau_4}
		\left(
		\mathcal S_{\tau_3}
		\bigl(
		z_0+\alpha\widehat\zeta-\beta\widehat\rho,
		0,
		f
		\bigr),
		0,
		f_{\tau_3}+\eta^4
		\right),
		0,
		f_{\tau_3+\tau_4}+\eta^5
		\right)
		\right\|_{\mathbb V^1}
		\\ \label{eq:mixed_final_stability}
		&<
		C_5(3C_4C_3\gamma+\gamma)
		\le4C_5C_4C_3\gamma
		<\frac{\epsilon}{2}.
	\end{align}
	Finally, set
	$T:=\tau_1+\tau_2+\tau_3+\tau_4+\tau_5<\sigma$
	and define $\eta$ by
	\begin{equation*}
		\eta(t)=
		\begin{cases}
			\eta^1(t), & t\in(0,\tau_1),\\
			\eta^2(t-\tau_1), & t\in(\tau_1,\tau_1+\tau_2),\\
			0, & t\in(\tau_1+\tau_2,\tau_1+\tau_2+\tau_3),\\
			\eta^4(t-\tau_1-\tau_2-\tau_3),
			& t\in(\tau_1+\tau_2+\tau_3,\tau_1+\tau_2+\tau_3+\tau_4),\\
			\eta^5(t-\tau_1-\tau_2-\tau_3-\tau_4),
			& t\in(\tau_1+\tau_2+\tau_3+\tau_4,T),\\
			0, & t\ge T.
		\end{cases}
	\end{equation*}
	By \Cref{lem:concatenation} with controls
	$\eta^1,\eta^2,0,\eta^4,\eta^5$, we have
	$\mathcal T(z_0,0,f+\eta)>T$.
	Then \eqref{eq:concatenated-solutions}, \eqref{eq:mixed_ideal_final},
	and \eqref{eq:mixed_final_stability} give
	\begin{align*}
		&\left\|
		\mathcal S_{T}(z_0,0,f+\eta)
		-
		z_0
		+
		\mathcal Q(\widehat\rho)
		-
		\mathcal M(\widehat\zeta,\widehat\rho)
		\right\|_{\mathbb V^1}
		\\
		&\le
		\left\|
		\mathcal S_{T}(z_0,0,f+\eta)
		-
		\mathcal S_{\tau_5}
		\left(
		\mathcal S_{\tau_4}
		\left(
		\mathcal S_{\tau_3}
		\bigl(
		z_0+\alpha\widehat\zeta-\beta\widehat\rho,
		0,
		f
		\bigr),
		0,
		f_{\tau_3}+\eta^4
		\right),
		0,
		f_{\tau_3+\tau_4}+\eta^5
		\right)
		\right\|_{\mathbb V^1}
		\\
		&\quad+
		\left\|
		\mathcal S_{\tau_5}
		\left(
		\mathcal S_{\tau_4}
		\left(
		\mathcal S_{\tau_3}
		\bigl(
		z_0+\alpha\widehat\zeta-\beta\widehat\rho,
		0,
		f
		\bigr),
		0,
		f_{\tau_3}+\eta^4
		\right),
		0,
		f_{\tau_3+\tau_4}+\eta^5
		\right)
		-
		z_0
		+
		\mathcal Q(\widehat\rho)
		-
		\mathcal M(\widehat\zeta,\widehat\rho)
		\right\|_{\mathbb V^1}
		\\
		&<
		\frac{\epsilon}{2}+\frac{\epsilon}{2}
		=\epsilon.
	\end{align*}
	Set $z(t):=\mathcal S_t(z_0,0,f+\eta)$ for $0\le t\le T$.
	By \eqref{eq:concatenated-regularity}, $z\in C([0,T];\mathbb V^1)
	\cap L^2(0,T;\mathbf H^2\times\dot H^4)$.
	Set $e_T:=\|z(T)-z_0+\mathcal Q(\widehat\rho)
	-\mathcal M(\widehat\zeta,\widehat\rho)\|_{\mathbb V^1}<\epsilon$.
	Keeping $\eta$ fixed, by continuity and \Cref{lem:regular-times},
	we may choose $\tau\in(\tau_1+\tau_2+\tau_3+\tau_4,T)$,
	sufficiently close to $T$, such that
	\begin{equation*}
	 z(\tau)\in\mathbf H^2\times\dot H^4 \quad \text{ and } \quad
	 \|z(\tau)-z(T)\|_{\mathbb V^1}<\frac{\epsilon-e_T}{2}.
	\end{equation*}
	Hence we further obtain the following
	\begin{align*}
	 &\|z(\tau)-z_0+\mathcal Q(\widehat\rho)
	     -\mathcal M(\widehat\zeta,\widehat\rho)\|_{\mathbb V^1}\\
	 &\quad\leq \|z(\tau)-z(T)\|_{\mathbb V^1}
	     +\|z(T)-z_0+\mathcal Q(\widehat\rho)
	     -\mathcal M(\widehat\zeta,\widehat\rho)\|_{\mathbb V^1}\\
	 &\quad<\frac{\epsilon-e_T}{2}+e_T<\epsilon.
	\end{align*}
	Set $\eta^\tau:=\eta|_{(0,\tau)}\in L^\infty(0,\tau;\mathcal H)$.
	By \Cref{cor:unique-strong-common}, it follows that
	$\mathcal S_\tau(z_0,0,f+\eta^\tau)=z(\tau).$
	This proves the required approximation and
	\eqref{eq:mixed_generated_regularity} with control $\eta^\tau$.
	This completes the proof.
\end{proof}
The previous two constructions allow us to isolate the transport term.
\begin{lemma}\label{lem:transport_generation}
	Let $N\ge1$ and assume that $\mathcal P_{N-1}$ holds. Let
	$\zeta\in\mathscr G_1$ and $\rho\in\mathcal H_2(N-1)$.
	Let $f$ satisfy the hypotheses of \Cref{small1}.
	Then, for any
	$\epsilon,\sigma>0$ and any $z_0\in{\mathbf H^2\times\dot H^4}$, there exist
	$\tau\in(0,\sigma)$ and $ \eta\in L^\infty(0,\tau;\mathcal H)$ such that
	\begin{equation}\label{eq:transport_generated_final}
		\left\|
		\mathcal S_\tau(z_0,0,f+\eta)
		-
		z_0
		-
		\mathcal M(\widehat\zeta,\widehat\rho)
		\right\|_{\mathbb V^1}
		<\epsilon,
	\end{equation}
	and
	\begin{equation}\label{eq:transport_generated_regularity}
	 \mathcal S_\tau(z_0,0,f+\eta)\in\mathbf H^2\times\dot H^4.
	\end{equation}
\end{lemma}
\begin{proof}
	Let $\zeta\in\mathscr G_1$ and $\rho\in\mathcal H_2(N-1)$. We show that
	the direction
$\mathcal M(\widehat\zeta,\widehat\rho)$
	can be generated in arbitrarily small time. Let $\epsilon,\sigma>0$ and
	$z_0\in{\mathbf H^2\times\dot H^4}$ be given.
	Set $z_1:=z_0+\mathcal Q(\widehat\rho)$.
	{By \Cref{rem:trigonometric-saturation},
	$\rho$ is a trigonometric polynomial. Hence
	$z_1\in\mathbf H^2\times\dot H^4$.}
	Moreover,
	\begin{equation*}
		z_1-\mathcal Q(\widehat\rho)
		+\mathcal M(\widehat\zeta,\widehat\rho)
		=
		z_0+\mathcal M(\widehat\zeta,\widehat\rho).
	\end{equation*}
	Hence, by \Cref{lem:mixed_generation}
	with prescribed force $f$, there exist
	$\tau_2\in(0,\sigma/2)$ and
	$ \eta^2\in L^\infty(0,\tau_2;\mathcal H)$ such that
	\begin{equation}\label{eq:transport_mixed_part}
		\left\|
		\mathcal S_{\tau_2}(z_1,0,f+\eta^2)
		-
		z_0
		-
		\mathcal M(\widehat\zeta,\widehat\rho)
		\right\|_{\mathbb V^1}
		<\frac{\epsilon}{2}.
	\end{equation}

	We now keep $\tau_2$ and $\eta^2$ fixed. By
	{\Cref{cont} and \eqref{eq:cont_stability}, applied to
	the solution from $z_1$ with force $f+\eta^2$ on $[0,\tau_2]$},
	there exist $r>0$ and $C\ge1$ such that, for $v_0\in\mathbb V^1$ and $s\ge0$,
	\begin{equation}\label{eq:transport_stability_condition}
	 \|v_0-z_1\|_{\mathbb V^1}
	 +\|f_s-f\|_{L^2(0,\tau_2;\mathbb V^0)}<r
	\end{equation}
	implies that $\mathcal S_t(v_0,0,f_s+\eta^2)$ exists on
	$[0,\tau_2]$ and satisfies
	\begin{align}\label{eq:transport_stability_estimate}
	 &\sup_{t\in[0,\tau_2]}
	 \|\mathcal S_t(v_0,0,f_s+\eta^2)
	      -\mathcal S_t(z_1,0,f+\eta^2)\|_{\mathbb V^1} \le C\Bigl(
	 \|v_0-z_1\|_{\mathbb V^1}
	 +\|f_s-f\|_{L^2(0,\tau_2;\mathbb V^0)}\Bigr).
	\end{align}
	Choose $0<\kappa<
		\min\left\{
		\frac r2,\frac{\epsilon}{4C}
		\right\}.$
By continuity of translations in $L^2$, choose $\sigma_0\in(0,\sigma/2)$
	such that
	\begin{equation*}
	 \|f_s-f\|_{L^2(0,\tau_2;\mathbb V^0)}<\kappa
	 \qquad\text{for every }0<s<\sigma_0.
	\end{equation*}
By \Cref{lem:cubic_generation}, applied to $z_0$ and $\rho$
	with prescribed force $f$, there exist
	$ \tau_1\in(0,\sigma_0)$ and
	$ \eta^1\in L^\infty(0,\tau_1;\mathcal H)$ such that
	\begin{equation}\label{eq:transport_prepare_Q}
		\left\|
		\mathcal S_{\tau_1}(z_0,0,f+\eta^1)
		-
		z_1
		\right\|_{\mathbb V^1}
		<\kappa.
	\end{equation}
	{Set $T:=\tau_1+\tau_2<\sigma$ and define
	\begin{equation*}
		\eta(t)=
		\begin{cases}
			\eta^1(t), & t\in(0,\tau_1),\\
			\eta^2(t-\tau_1), & t\in(\tau_1,T),\\
			0, & t\ge T.
		\end{cases}
	\end{equation*}}
By \eqref{eq:transport_prepare_Q} and the choice of $\sigma_0$,
	\begin{align*}
	 &\|\mathcal S_{\tau_1}(z_0,0,f+\eta^1)-z_1\|_{\mathbb V^1}+\|f_{\tau_1}-f\|_{L^2(0,\tau_2;\mathbb V^0)}
	 <2\kappa<r.
	\end{align*}
	Thus \eqref{eq:transport_stability_condition} holds with
	$v_0=\mathcal S_{\tau_1}(z_0,0,f+\eta^1)$ and $s=\tau_1$.
	{By \Cref{lem:concatenation} with controls
	$\eta^1,\eta^2$, we have $\mathcal T(z_0,0,f+\eta)>T$.
	Then \eqref{eq:concatenated-solutions},
	\eqref{eq:transport_stability_estimate},
	and \eqref{eq:transport_mixed_part} give}

	\begin{align*}
		&\left\|
		\mathcal S_{T}(z_0,0,f+\eta)
		-
		z_0
		-
		\mathcal M(\widehat\zeta,\widehat\rho)
		\right\|_{\mathbb V^1}
		\\
		&\le
		\left\|
		\mathcal S_{\tau_2}
		\bigl(
		\mathcal S_{\tau_1}(z_0,0,f+\eta^1),
		0,f_{\tau_1}+\eta^2
		\bigr)
		-
		\mathcal S_{\tau_2}(z_1,0,f+\eta^2)
		\right\|_{\mathbb V^1}
		\\
		&\quad+
		\left\|
		\mathcal S_{\tau_2}(z_1,0,f+\eta^2)
		-
		z_0
		-
		\mathcal M(\widehat\zeta,\widehat\rho)
		\right\|_{\mathbb V^1}
		\\
		&<
		C\Bigl(\kappa+
		\|f_{\tau_1}-f\|_{L^2(0,\tau_2;\mathbb V^0)}\Bigr)
		+\frac{\epsilon}{2}
		<2C\kappa+\frac{\epsilon}{2}
		<\epsilon.
	\end{align*}
	{Set $z(t):=\mathcal S_t(z_0,0,f+\eta)$ for $0\le t\le T$.
	By \eqref{eq:concatenated-regularity}, $z\in C([0,T];\mathbb V^1)
	\cap L^2(0,T;\mathbf H^2\times\dot H^4)$.
	Set $e_T:=\|z(T)-z_0-\mathcal M(\widehat\zeta,\widehat\rho)\|_{\mathbb V^1}
	<\epsilon$.
	Keeping $\eta$ fixed, by continuity and \Cref{lem:regular-times},
	we may choose $\tau\in(\tau_1,T)$, sufficiently close to $T$, such that
	\begin{equation*}
	 z(\tau)\in\mathbf H^2\times\dot H^4,\qquad
	 \|z(\tau)-z(T)\|_{\mathbb V^1}<\frac{\epsilon-e_T}{2}.
	\end{equation*}
	Hence,
	\begin{align*}
	 \|z(\tau)-z_0-\mathcal M(\widehat\zeta,\widehat\rho)\|_{\mathbb V^1}
	 &\leq \|z(\tau)-z(T)\|_{\mathbb V^1}
	       +\|z(T)-z_0-\mathcal M(\widehat\zeta,\widehat\rho)\|_{\mathbb V^1}\\
	 &<\frac{\epsilon-e_T}{2}+e_T<\epsilon.
	\end{align*}
	Set $\eta^\tau:=\eta|_{(0,\tau)}\in L^\infty(0,\tau;\mathcal H)$.
	By \Cref{cor:unique-strong-common},
%	\begin{equation*}
	$ \mathcal S_\tau(z_0,0,f+\eta^\tau)=z(\tau).$
%	\end{equation*}
	Thus \eqref{eq:transport_generated_final} and
	\eqref{eq:transport_generated_regularity} hold with control
	$\eta^\tau$.}
\end{proof}
We are in position to prove the main result of this section.
\begin{proof}[Proof of \Cref{small1}]
	{We prove \eqref{propPN} by induction on $N$.
	We use the notation \eqref{eq:control-force-translations} and note that every $f_s$ satisfies the same hypotheses as $f$.}

	\smallskip
	\noindent
	\textbf{Step 1. Base case.}
	Let $N=0$ and $\xi\in\mathcal H_2(0)=\mathcal H_2$.
	{Let $f$ satisfy the hypotheses of \Cref{small1}, and fix
	$\epsilon,\sigma>0$ and $z_0\in\mathbf H^2\times\dot H^4$.}
	Then $\widehat\xi=(0,\xi)\in\mathcal H$. Hence, by
	\Cref{lem:direct_construction}, there exist
	$ {\delta}\in(0,\sigma)$ and
	$ \eta\in L^\infty(0,{\delta};\mathcal H)$ such that
	\begin{equation*}
		\left\|
		\mathcal S_{\delta}(z_0,0,f+\eta)
		-
		z_0
		-
		\widehat\xi
		\right\|_{\mathbb V^1}
		<{\frac{\epsilon}{2}}.
	\end{equation*}
	{By \Cref{wp1}, this solution belongs to
	$C([0,\delta];\mathbb V^1)\cap
	L^2(0,\delta;\mathbf H^2\times\dot H^4)$.
	Keeping $\eta$ fixed, write $z(t):=\mathcal S_t(z_0,0,f+\eta)$
	for $0\leq t\leq\delta$.
	By \Cref{lem:regular-times}, we may choose $\tau\in(0,\delta)$,
	sufficiently close to $\delta$, such that
	\begin{equation*}
	 z(\tau)\in\mathbf H^2\times\dot H^4,\qquad
	 \|z(\tau)-z(\delta)\|_{\mathbb V^1}<\frac{\epsilon}{2}.
	\end{equation*}
	Consequently,
	\begin{align*}
	 \|z(\tau)-z_0-\widehat\xi\|_{\mathbb V^1}
	 &\leq \|z(\tau)-z(\delta)\|_{\mathbb V^1}
	       +\|z(\delta)-z_0-\widehat\xi\|_{\mathbb V^1}\\
	 &<\frac{\epsilon}{2}+\frac{\epsilon}{2}=\epsilon.
	\end{align*}
	Set $\eta^\tau:=\eta|_{(0,\tau)}$.
	By \Cref{cor:unique-strong-common}, 
we obtain
	$ \mathcal S_\tau(z_0,0,f+\eta^\tau)
	 =\mathcal S_\tau(z_0,0,f+\eta)$.
 Thus $\mathcal P_0$ holds.

	\smallskip
	\noindent
	\textbf{Step 2. Inductive step.}
	Assume that $\mathcal P_{N-1}$ holds and let
	$\xi\in\mathcal H_2(N)$. By the definition of $\mathcal H_2(N)$, we
	can write
	\begin{equation}\label{eq:phase_decomp_transport}
		\xi
		=
		\xi_0
		+
		\sum_{j=1}^{m}B_1(\zeta_j,\xi_j),
	\end{equation}
	where $\xi_0,\xi_j\in\mathcal H_2(N-1), 
		\zeta_j\in\mathscr G_1,
		 j=1,\ldots,m.$
	{An argument similar to that in the proof of
	\Cref{lem:velocity_representation} gives
	\eqref{eq:phase_decomp_transport}. Here we use only the linearity
	of $B_1$ in its second argument: the scalar coefficients are absorbed
	into the phase factors, while the $\zeta_j$ remain in the fixed
	family $\mathscr G_1$.}
%	Notice that the velocity factors
%	$\zeta_j$ are not rescaled and they remain in the fixed family
%	$\mathscr G_1$.

{Let $f$ satisfy the hypotheses of \Cref{small1}, and let}
	$\epsilon,\sigma>0$ and
	$z_0\in{\mathbf H^2\times\dot H^4}$ be given, and set
	$	\widehat \epsilon:=\frac{\epsilon}{m+1}.$
	Since $\xi_0\in\mathcal H_2(N-1)$, the induction hypothesis
	$\mathcal P_{N-1}$ gives $\tau_0\in(0,\sigma/(m+1))$ and
	$\eta^0\in L^\infty(0,\tau_0;\mathcal H)$ such that, setting
	$\vartheta_0:=\mathcal S_{\tau_0}(z_0,0,f+\eta^0),\,
		\varPsi_0:=z_0+\widehat\xi_0,$
	we have
	\begin{equation}\label{eq:phase_first_endpoint}
		\|\vartheta_0-\varPsi_0\|_{\mathbb V^1}<\widehat \epsilon.
	\end{equation}
The induction hypothesis also gives
	$\vartheta_0\in\mathbf H^2\times\dot H^4$.

	For $j=1,\ldots,m$, we construct successively $\tau_j$, $\eta^j$, and
	$\vartheta_j$. Suppose that
	$\vartheta_{j-1}\in{\mathbf H^2\times\dot H^4}$ has already been
	constructed. Set
	$s_{j-1}:=\tau_0+\cdots+\tau_{j-1}$.
	Since $\zeta_j\in\mathscr G_1$ and
	$\xi_j\in\mathcal H_2(N-1)$, \Cref{lem:transport_generation}, applied
	with initial condition $\vartheta_{j-1}$
	and prescribed force $f_{s_{j-1}}$, gives
$	\tau_j\in\left(0,\frac{\sigma}{m+1}\right),
		\,
		\eta^j\in L^\infty(0,\tau_j;\mathcal H),$
	such that, setting
$	\vartheta_j:=\mathcal S_{\tau_j}(\vartheta_{j-1},0,f_{s_{j-1}}+\eta^j),$
	we have
	\begin{equation}\label{eq:phase_exact_iteration}
		\left\|
		\vartheta_j
		-
		\vartheta_{j-1}
		-
		\mathcal M(\widehat\zeta_j,\widehat\xi_j)
		\right\|_{\mathbb V^1}
		<\widehat \epsilon.
	\end{equation}
{By \eqref{eq:transport_generated_regularity},
	$\vartheta_j\in\mathbf H^2\times\dot H^4$.} Define
$\varPsi_j
		:=
		\varPsi_{j-1}
		+
		\mathcal M(\widehat\zeta_j,\widehat\xi_j),
		\, j=1,\ldots,m.$
	Using \eqref{eq:phase_exact_iteration}, we obtain
	\begin{align*}
		\|\vartheta_j-\varPsi_j\|_{\mathbb V^1}
		&\le
		\left\|
		\vartheta_j
		-
		\vartheta_{j-1}
		-
		\mathcal M(\widehat\zeta_j,\widehat\xi_j)
		\right\|_{\mathbb V^1}
		+
		\|\vartheta_{j-1}-\varPsi_{j-1}\|_{\mathbb V^1}\\
		&<
		\widehat \epsilon
		+
		\|\vartheta_{j-1}-\varPsi_{j-1}\|_{\mathbb V^1}.
	\end{align*}
	Hence, by induction and \eqref{eq:phase_first_endpoint},
	\begin{equation}\label{eq:phase_accumulated_error}
		\|\vartheta_j-\varPsi_j\|_{\mathbb V^1}
		<
		(j+1)\widehat \epsilon,
		\qquad j=0,\ldots,m.
	\end{equation}
	{Set $\tau:=\sum_{j=0}^{m}\tau_j<\sigma$ and define
	\begin{align*}
	 \eta(r)&=\eta^0(r), \quad 0<r<\tau_0,\\
	 \eta(s_{j-1}+r)&=\eta^j(r), \quad 0<r<\tau_j,\quad j=1,\ldots,m,
	\end{align*}
	and $\eta(t)=0$ for $t\ge\tau$.
	Applying \Cref{lem:concatenation} to the
	$m+1$ solutions defining $\vartheta_0,\ldots,\vartheta_m$ gives
	$\mathcal T(z_0,0,f+\eta)>\tau$ and
	\begin{equation*}
	 \mathcal S_\tau(z_0,0,f+\eta)=\vartheta_m
	 \in\mathbf H^2\times\dot H^4.
	\end{equation*}}
	On the other hand, by \eqref{eq:phase_decomp_transport},
	\begin{align*}
		\varPsi_m
		&=
		z_0+\widehat\xi_0
		+
		\sum_{j=1}^{m}
		\mathcal M(\widehat\zeta_j,\widehat\xi_j)=
		z_0
		+
		\left(
		0,
		\xi_0+\sum_{j=1}^{m}B_1(\zeta_j,\xi_j)
		\right)
		=
		z_0+\widehat\xi.
	\end{align*}
	Therefore, by \eqref{eq:phase_accumulated_error},
	\begin{equation*}
		\left\|
		\mathcal S_\tau(z_0,0,f+\eta)
		-
		z_0
		-
		\widehat\xi
		\right\|_{\mathbb V^1}
		=
		\|\vartheta_m-\varPsi_m\|_{\mathbb V^1}
		<
		(m+1)\widehat \epsilon
		=
		\epsilon.
	\end{equation*}
	Thus $\mathcal P_N$ holds. By induction, $\mathcal P_N$ holds for every
	$N\ge0$.

	\smallskip
	\noindent
	\textbf{Step 3. Initial data in $\mathbb V^1$.}
	Let now $z_0\in\mathbb V^1$ and
	$\xi\in\mathcal H_2(\infty)$. By the definition of
	$\mathcal H_2(\infty)$, there exists $N\ge0$ such that
	$\xi\in\mathcal H_2(N)$.

	Let $z(t):=\mathcal S_t(z_0,0,f)$ be the strong solution with zero control.
	By continuity in $\mathbb V^1$, there exists
	$\sigma_0\in(0,\min\{\sigma/2,\mathcal T(z_0,0,f)\})$ such that

	\begin{equation}\label{eq:initial_short_free}
		\|z(t)-z_0\|_{\mathbb V^1}<\frac{\epsilon}{2},
		\qquad t\in[0,\sigma_0].
	\end{equation}
	{By \Cref{wp1},
	$z\in L^2(0,\sigma_0;\mathbf H^2\times\dot H^4)$.
	Applying \Cref{lem:regular-times} at time $0$ on $[0,\sigma_0]$,
	we obtain $t_0\in(0,\sigma_0)$ such that}
	\begin{equation*}
		\widetilde z:=\mathcal S_{t_0}(z_0,0,f)
		\in{\mathbf H^2\times\dot H^4}.
	\end{equation*}
	Applying $\mathcal P_N$ with initial condition $\widetilde z$, accuracy
	$\epsilon/2$, time bound $\sigma/2$,
	and prescribed force $f_{t_0}$, we obtain
	$\tau_1\in(0,\sigma/2)$ and
	$\eta^1\in L^\infty(0,\tau_1;\mathcal H)$ such that
	\begin{equation}\label{eq:phase_after_regularization}
		\left\|
		\mathcal S_{\tau_1}(\widetilde z,0,f_{t_0}+\eta^1)
		-
		\widetilde z
		-
		\widehat\xi
		\right\|_{\mathbb V^1}
		<\frac{\epsilon}{2}.
	\end{equation}
	Set $\tau:=t_0+\tau_1<\sigma$ and define
	\begin{equation*}
		\eta(t)=
		\begin{cases}
			0, & t\in(0,t_0),\\
			\eta^1(t-t_0), & t\in(t_0,\tau),\\
			0, & t\ge\tau.
		\end{cases}
	\end{equation*}
	By \Cref{lem:concatenation} with controls
	$0,\eta^1$, we have $\mathcal T(z_0,0,f+\eta)>\tau$.
	Then \eqref{eq:concatenated-solutions},
	\eqref{eq:initial_short_free}, and
	\eqref{eq:phase_after_regularization} give}
	\begin{align*}
		\left\|
		\mathcal S_\tau(z_0,0,f+\eta)
		-
		z_0
		-
		\widehat\xi
		\right\|_{\mathbb V^1}
		&\le
		\left\|
		\mathcal S_{\tau_1}(\widetilde z,0,f_{t_0}+\eta^1)
		-
		\widetilde z
		-
		\widehat\xi
		\right\|_{\mathbb V^1}
		+
		\|\widetilde z-z_0\|_{\mathbb V^1}\\
		&<
		\frac{\epsilon}{2}
		+
		\frac{\epsilon}{2}
		=
		\epsilon.
	\end{align*}
	This proves \eqref{limit_phi1}.
\end{proof}

\subsection{Proof of small-time accessibility (\Cref{thm:global-small-time-control})}
\label{sec:unit-ball-accessibility-proof}
\begin{proof}
Let $z_0\in\mathbb V^1$, let $f$ satisfy the hypotheses of
\Cref{small1}, and fix $\sigma>0$.
By \Cref{prop_sat}, there exists
$\rho=(\rho_u,\rho_\phi)\in\mathcal H(\infty)$ such that
\begin{equation}\label{eq:unit-ball-saturation-approximation}
 \|z_0+\rho\|_{\mathbb V^1}<\frac13.
\end{equation}
Since $\rho_u\in\mathcal H_1(\infty)$, \Cref{second_limit} gives
$\tau_1\in(0,\sigma/2)$ and
$\eta^1\in L^\infty(0,\tau_1;\mathcal H)$ such that, setting
$z^1:=\mathcal S_{\tau_1}(z_0,0,f+\eta^1)$, we have
\begin{equation}\label{eq:unit-ball-velocity-step}
 \|z^1-z_0-(\rho_u,0)\|_{\mathbb V^1}<\frac13.
\end{equation}
The force $f_{\tau_1}$ satisfies the same hypotheses as $f$.
Applying \Cref{small1} with initial state $z^1\in\mathbb V^1$
and phase direction $\rho_\phi\in\mathcal H_2(\infty)$, we obtain
$\tau_2\in(0,\sigma/2)$ and
$\eta^2\in L^\infty(0,\tau_2;\mathcal H)$ such that, setting
$z^2:=\mathcal S_{\tau_2}(z^1,0,f_{\tau_1}+\eta^2)$, we have
\begin{equation}\label{eq:unit-ball-phase-step}
 \|z^2-z^1-(0,\rho_\phi)\|_{\mathbb V^1}<\frac13.
\end{equation}
{Set $\tau:=\tau_1+\tau_2<\sigma$ and define
\begin{equation*}
 \eta(t)=
 \begin{cases}
  \eta^1(t),&t\in(0,\tau_1),\\
  \eta^2(t-\tau_1),&t\in(\tau_1,\tau),\\
  0,&t\ge\tau.
 \end{cases}
\end{equation*}}
Then $\eta\in L^\infty(0,\tau;\mathcal H)$.
{By \Cref{lem:concatenation} with controls
$\eta^1,\eta^2$,}
the strong solution
$z(t):=\mathcal S_t(z_0,0,f+\eta)$ of \eqref{eq:chn-accessibility}
is defined on $[0,\tau]$ and satisfies
\begin{align*}
 z&\in C([0,\tau];\mathbb V^1)
       \cap L^2(0,\tau;\mathbf H^2\times\dot H^4),\\
 z(\tau)&=\mathcal S_{\tau_2}(z^1,0,f_{\tau_1}+\eta^2)=z^2.
\end{align*}
Uniqueness follows from \Cref{cor:unique-strong-common}.
Finally, \eqref{eq:unit-ball-saturation-approximation},
\eqref{eq:unit-ball-velocity-step}, and
\eqref{eq:unit-ball-phase-step} give
\begin{align*}
 \|z(\tau)\|_{\mathbb V^1}
 &\le \|z^2-z^1-(0,\rho_\phi)\|_{\mathbb V^1}
      +\|z^1-z_0-(\rho_u,0)\|_{\mathbb V^1}\\
 &\qquad+\|z_0+\rho\|_{\mathbb V^1}
 <\frac13+\frac13+\frac13=1.
\end{align*}
This gives the desired result.
\end{proof}

\section{Further results, remarks and open problems}\label{sec:further-results}
\subsection{Scope of the methodology}\label{sec:scope-construction}
We believe that the strategy developed here to obtain global
strong solutions can be extended to other equations for which
global regularity is still open. The proof in
\Cref{sec:global-prescribed-force} uses only two ingredients:
a robust local theory and small-time accessibility of a fixed
ball. We first present a consequence for the Navier--Stokes equations.

In the proof of \Cref{second_limit}, we may take $\eta_\phi=0$
at every step of the induction.
If $\phi_0=f_\phi=0$, the phase then remains zero. Combining
this construction with \Cref{prop_sat_phase} gives small-time
accessibility of the unit ball in $\mathbf H^1$.
We may then repeat the proof of \Cref{mnthm}, to obtain the following result.
\begin{proposition}\label{prop:ns-scalar-continuation}
Let $\nu>0$, $u_0\in\mathbf H^1$, and
$f_u\in L^\infty_{\mathrm{loc}}
(\mathbb R_+;\mathbf L^2_{\mathrm{div},0})$.
For every $\rho>0$, there exists a scalar correction
$\eta_u\in C^\infty([0,\infty);\mathcal E_{\mathcal K_u})$,
where $\mathcal K_u=\{e_2,e_1+e_3,e_2+e_3\}$, such that
\begin{equation*}
 \partial_tu+\nu A_0u+B_0(u,u)=f_u+\Pi(\eta_ue_3),
 \qquad u(0)=u_0,
\end{equation*}
has a unique global strong solution
\begin{equation*}
 u\in C([0,\infty);\mathbf H^1)
       \cap L^2_{\mathrm{loc}}([0,\infty);\mathbf H^2).
\end{equation*}
Moreover, $\eta_u(t)=0$ outside a union $I$ of intervals with
$|I|<\rho$.
\end{proposition}

To our knowledge, \Cref{prop:ns-scalar-continuation}
is new with a scalar correction supported on these three frequencies. If corrections are allowed in all components of the force and
sufficiently many Stokes modes are used, an analogous result can
be deduced by combining the controllability results of
\cite[Theorem~2.5 and Corollary~2.3]{S07} with the continuation
and smoothing argument in \Cref{sec:global-prescribed-force}.

We also have in mind the periodic Kuramoto--Sivashinsky equation
in two and three dimensions. Global regularity for general smooth
initial data is still open, already in two dimensions; see
\cite{LariosRahmanYamazaki2022}. We ask whether an analogue of
\Cref{prop:ns-scalar-continuation} holds, with forcing directions
fixed independently of the initial data.

For the three-dimensional Boussinesq system with viscosity
and thermal diffusion, it would be particularly interesting
to act on the temperature alone. The two-dimensional results of
\cite{NR25b,rissel2025controllability} suggest ways to use the
temperature to move the fluid, but adapting these ideas
to three dimensions requires further work.

For the three-dimensional MHD equations with viscosity
and magnetic diffusion, we may ask the same question, allowing
corrections in both equations. Global weak solutions and local
strong solutions are known \cite{SermangeTemam1983}, but
global strong existence for general data remains open. Indeed,
setting the magnetic field and its force equal to zero recovers
the three-dimensional Navier--Stokes equations.
Some of the algebra needed for saturation has already been
studied in \cite[Propositions~4.2--4.3]{Yamazaki2019} through
Lie-bracket computations for a three-dimensional Galerkin approximation.
It would be interesting to use these calculations to identify a
saturating space for the full system and justify the associated
small-time limits.

For the ideal MHD, without viscosity or magnetic diffusion, the local
estimates around prescribed velocity fields in the recent work
\cite{KukavicaOzanski2026} may help in studying the large profiles
used in scaling arguments. Their controllability result concerns
boundary controls. In our setting, one would seek small-time
accessibility with additive forcing. Such a result would extend to ideal MHD the
prevention of blow-up for Euler flows discussed in
\cite[Remark~5.33]{GHHM18}, which motivated the present work.
\subsection{Controllability to arbitrary targets}\label{sec:arbitrary-targets}
Apart from establishing the existence of global strong solutions, our analysis can also be used to study the approximate controllability of the Cahn–Hilliard–Navier–Stokes system.
{We first extend \Cref{def:small-time-accessibility} to balls
with arbitrary centres and radii. The proof in
\Cref{sec:unit-ball-accessibility-proof} gives this property with
only a change in the target and the error bounds.}

\begin{defi}\label{def:target-accessibility}
Let $z_1\in\mathbb V^1$ and $\epsilon>0$. We say that the ball
$B_{\mathbb V^1}(z_1,\epsilon)$ is accessible in small time with
controls in $\mathcal H$ if, for every $z_0\in\mathbb V^1$,
$f\in L^\infty_{\mathrm{loc}}(\mathbb R_+;\mathbb V^0)$, and
$\sigma>0$, there exist $\tau\in(0,\sigma)$ and
$\eta\in L^\infty(0,\tau;\mathcal H)$ such that
\eqref{eq:chn-accessibility} has a strong solution $z$ on
$[0,\tau]$, with initial state $z_0$, satisfying
\begin{equation}\label{eq:arbitrary-target-accessibility}
 \|z(\tau)-z_1\|_{\mathbb V^1}<\epsilon.
\end{equation}
\end{defi}

\begin{proposition}\label{t0}
For $\mathcal H$ defined by \eqref{eq:projected_modes_intro},
every ball $B_{\mathbb V^1}(z_1,\epsilon)$, with
$z_1\in\mathbb V^1$ and $\epsilon>0$, is accessible in small time
with controls in $\mathcal H$.
\end{proposition}
\begin{proof}
Fix $z_0$, $f$, and $\sigma>0$ as in
\Cref{def:target-accessibility}. We follow the proof in
\Cref{sec:unit-ball-accessibility-proof}.
By \Cref{prop_sat}, choose
$\rho=(\rho_u,\rho_\phi)\in\mathcal H(\infty)$ such that
\begin{equation*}
 \|z_0+\rho-z_1\|_{\mathbb V^1}<\frac{\epsilon}{3}
\end{equation*}
in place of \eqref{eq:unit-ball-saturation-approximation}.
Apply \Cref{second_limit} with displacement $\rho_u$,
precision $\epsilon/3$, and time horizon $\sigma/2$, and write
$z^1:=\mathcal S_{\tau_1}(z_0,0,f+\eta^1)$ for the resulting
state. Next apply \Cref{small1} to $z^1$ with displacement
$\rho_\phi$ and force $f_{\tau_1}$, using the same precision
and time horizon, and set
$z^2:=\mathcal S_{\tau_2}(z^1,0,f_{\tau_1}+\eta^2)$.
Here $\tau_j\in(0,\sigma/2)$ and
$\eta^j\in L^\infty(0,\tau_j;\mathcal H)$, $j=1,2$, and
\begin{equation*}
 \|z^1-z_0-(\rho_u,0)\|_{\mathbb V^1}<\frac{\epsilon}{3},
 \qquad
 \|z^2-z^1-(0,\rho_\phi)\|_{\mathbb V^1}<\frac{\epsilon}{3}.
\end{equation*}
Set $\tau:=\tau_1+\tau_2<\sigma$.
By \Cref{lem:concatenation}, the control obtained from
$\eta^1$ and $\eta^2$, as in
\Cref{sec:unit-ball-accessibility-proof}, gives a strong
solution $z$ on $[0,\tau]$ with $z(\tau)=z^2$. Therefore,
\begin{align*}
 \|z(\tau)-z_1\|_{\mathbb V^1}
 &\le \|z^2-z^1-(0,\rho_\phi)\|_{\mathbb V^1}
      +\|z^1-z_0-(\rho_u,0)\|_{\mathbb V^1}+\|z_0+\rho-z_1\|_{\mathbb V^1}\\
 &<\frac{\epsilon}{3}+\frac{\epsilon}{3}
      +\frac{\epsilon}{3}=\epsilon.
\end{align*}
This ends the proof.
\end{proof}
In \Cref{t0}, the arrival time is chosen as part of the
construction. We now require the approximation at a prescribed
time $T>0$. This motivates the following definition.
\begin{defi}[Approximate controllability]\label{def:approximate-controllability}
Let $T>0$ and
$f\in L^\infty(0,T;\mathbb V^0)$.
We say that \eqref{eq:chn-accessibility} is approximately
controllable at time $T$ with controls in $\mathcal H$ if, for
every $z_0,z_1\in\mathbb V^1$ and $\epsilon>0$, there exists
$\eta\in L^\infty(0,T;\mathcal H)$ such that the system has a
strong solution $z$ on $[0,T]$, with initial state $z_0$,
satisfying
\begin{equation}\label{eq:global-small-time-access1}
 \|z(T)-z_1\|_{\mathbb V^1}<\epsilon.
\end{equation}
\end{defi}
\begin{theo}\label{thm:global-small-time-control1}
For every $T>0$ and
$f\in L^\infty(0,T;\mathbb V^0)$,
system \eqref{eq:chn-accessibility} is approximately controllable
at time $T$ with controls in $\mathcal H$.
\end{theo}
The proof follows the construction in
\Cref{sec:global-prescribed-force}, now with returns to a ball
around $z_1$ and a prescribed final time $T$.\footnote{{To apply
the results formulated on $\mathbb R_+$, we extend forces and
controls by zero beyond their intervals of definition. By
\Cref{cor:unique-strong-common}, the solution on the original
interval is independent of this extension.}}

\begin{proof}[Proof of \Cref{thm:global-small-time-control1}]
Let $T>0$ and $\epsilon>0$. Given
$z_0=(u_0,\phi_0)$ and $z_1=(u_1,\phi_1)\in\mathbb V^1$,
we aim to show the existence of a control
$\eta\in L^\infty(0,T;\mathcal H)$ such that
\begin{equation*}
 \|\mathcal S_T(z_0,0,{f+}\eta)-z_1\|_{\mathbb V^1}<\epsilon.
\end{equation*}
First we find $r,\tau>0$ such that, from any state in
$B_{\mathbb V^1}(z_1,r)$, evolution under $f_s$ stays within
distance $\epsilon$ of $z_1$ for
$0\le t\le\min\{\tau,T-s\}$, uniformly for $0\le s<T$.
By \Cref{wp1}, choose $\tau_0\in(0,T)$ such that
$y(t):=\mathcal S_t(z_1,0,0)$ is defined on $[0,\tau_0]$.
Apply \Cref{cont} to this solution on $[0,\tau_0]$.
Let $\delta_*>0$ be the threshold in
\eqref{eq:shifted-perturbation-smallness} and $C\ge1$ the
constant in \eqref{eq:cont_stability} (increasing it if necessary).
Set $M_T:=\|f\|_{L^\infty(0,T;\mathbb V^0)}$.
Choose $r>0$ and $\tau\in(0,\tau_0)$ so small that
\begin{equation}\label{eq:free_target}
 r+M_T\sqrt{\tau}<\min\left\{\delta_*,\frac{\epsilon}{2C}\right\},
 \qquad
 \sup_{0\le t\le\tau}\|y(t)-z_1\|_{\mathbb V^1}
 <\frac{\epsilon}{2}.
\end{equation}
For $0\le s<T$, put
\begin{equation*}
 F_s(t):=
 \begin{cases}
  f(s+t),&0<t<\min\{\tau,T-s\},\\
  0,&\min\{\tau,T-s\}\le t<\tau_0.
 \end{cases}
\end{equation*}
Since
$\|F_s\|_{L^2(0,\tau_0;\mathbb V^0)}\le M_T\sqrt{\tau}$,
\Cref{cont} gives, for every $v\in B_{\mathbb V^1}(z_1,r)$,
a strong solution $\mathcal S(v,0,F_s)$ on $[0,\tau_0]$ satisfying
\begin{equation}\label{eq:uniform_stability_target}
 \|\mathcal S_\cdot(v,0,F_s)-y\|_{C([0,\tau_0];\mathbb V^1)}
 \le C\left(\|v-z_1\|_{\mathbb V^1}
       +\|F_s\|_{L^2(0,\tau_0;\mathbb V^0)}\right)
 <\frac{\epsilon}{2}.
\end{equation}
Since $F_s=f_s$ on $(0,\min\{\tau,T-s\})$,
\Cref{cor:unique-strong-common,prop:maximal-strong} imply that
$\mathcal S(v,0,f_s)$ is also defined on
$[0,\min\{\tau,T-s\}]$ and agrees
there with $\mathcal S(v,0,F_s)$. Indeed, if its maximal time
were at most $\min\{\tau,T-s\}$, uniqueness and
\eqref{eq:uniform_stability_target} would bound its
$\mathbb V^1$-norm up to that time, contradicting
\eqref{eq:maximal-strong-blowup}.
Thus \eqref{eq:free_target} and
\eqref{eq:uniform_stability_target} give
\begin{equation}\label{cont_free_estimate}
 \begin{aligned}
 \|\mathcal S_t(v,0,f_s)-z_1\|_{\mathbb V^1}
 &\le \|\mathcal S_t(v,0,F_s)-y(t)\|_{\mathbb V^1}
       +\|y(t)-z_1\|_{\mathbb V^1}\\
 &<\frac{\epsilon}{2}+\frac{\epsilon}{2}=\epsilon,
 \end{aligned}
\end{equation}
for $0\le t\le\min\{\tau,T-s\}$, $0\le s<T$, and
$v\in B_{\mathbb V^1}(z_1,r)$.

With $r$ and $\tau$ fixed, we construct the control recursively,
as in the proof of \Cref{mnthm}, starting with
$a_0=0$ and $Z_0=z_0$.
Suppose that $a_n<T$ and $Z_n\in\mathbb V^1$ have been chosen.
By \Cref{t0}, applied with force $f_{a_n}$, there exist
$\delta_n\in(0,T-a_n)$ and
$\eta^n\in L^\infty(0,\delta_n;\mathcal H)$ such that
\begin{equation*}
 \begin{gathered}
 v^{(n)}(t):=\mathcal S_t(Z_n,0,f_{a_n}+\eta^n),
 \qquad 0\le t\le\delta_n,\\
 \widehat Z_n:=v^{(n)}(\delta_n)
 \in B_{\mathbb V^1}(z_1,r).
 \end{gathered}
\end{equation*}
Set
\begin{equation*}
 \begin{gathered}
 b_n:=a_n+\delta_n,\qquad
 \ell_n:=\min\{\tau,T-b_n\},\qquad
 a_{n+1}:=b_n+\ell_n,\\
 Z_{n+1}:=\mathcal S_{\ell_n}(\widehat Z_n,0,f_{b_n}).
 \end{gathered}
\end{equation*}
By \eqref{cont_free_estimate},
$\|Z_{n+1}-z_1\|_{\mathbb V^1}<\epsilon$.
If $a_{n+1}=T$, we stop. If $a_{n+1}<T$, then $\ell_n=\tau$,
and we repeat the construction from $(a_{n+1},Z_{n+1})$.

The procedure terminates after finitely many steps. Indeed, each
completed step contains an evolution under $f$
of length $\tau$, except possibly the last one.
Hence, after $n\geq1$ completed steps with $a_n<T$,
the elapsed time is strictly larger than $n\tau$.
Consequently, the procedure must reach $T$ after at most
$\lfloor T/\tau\rfloor+1$ steps.

Let $N$ be the final index, so that $a_N=T$, and define
\begin{equation*}
 \eta(t)=
 \begin{cases}
  \eta^n(t-a_n),&a_n<t<b_n,\quad n=0,\ldots,N-1,\\
  0,&b_n<t<a_{n+1},\quad n=0,\ldots,N-1.
 \end{cases}
\end{equation*}
The resulting control belongs to $L^\infty(0,T;\mathcal H)$,
since it is obtained by concatenating finitely many
$\mathcal H$-valued $L^\infty$ controls.

Applying \Cref{lem:concatenation} to the solutions $v^{(n)}$
and $t\mapsto\mathcal S_t(\widehat Z_n,0,f_{b_n})$ on
$[0,\ell_n]$, we obtain a strong solution on $[0,T]$ with
$\mathcal S_T(z_0,0,f+\eta)=Z_N$.
Therefore \eqref{eq:global-small-time-access1} holds.
\end{proof}

\subsection{An open problem: forcing the mixture, not the flow}
In the saturation argument of \Cref{sec:saturation}, the fluid generates the missing phase
directions through transport, while the directions in the fluid are obtained
from self-advection, without using the phase. Note, however, that
capillarity also allows the phase to act on the fluid. We would like
to understand whether this interaction can replace some, or even all,
of the direct forcing in the fluid equation.
In particular, we ask
\begin{center}
\begin{minipage}{0.9\textwidth}
\centering\itshape
Does the conclusion of \Cref{mnthm} hold if only the phase source
is modified, at finitely many fixed Fourier frequencies?
\end{minipage}
\end{center}
A positive answer would give global strong solutions by acting only
on the composition of the mixture, leaving the prescribed fluid
force unchanged. The use of temperature controls to move the fluid
in the primitive equations and the two-dimensional Boussinesq system
in \cite{BGN23,NR25b,rissel2025controllability} motivates this question.
However, in those models, temperature enters the fluid equation
linearly. For CHNS, the capillary force is quadratic in the phase,
while the phase equation contains a cubic term. As discussed in
\Cref{sec:main-ideas}, this difference makes the balance of amplitudes
more delicate.

Moreover, the projected capillary force vanishes on each eigenspace of the Laplacian. Indeed, if
$\Delta\psi=-\lambda\psi$, then
\begin{equation*}
 R_0(\Delta\psi,\psi)
 =-\frac{\lambda}{2}\Pi\nabla(\psi^2)=0.
\end{equation*}
This cancellation removes the fluid
displacement we would like to generate. Thus capillarity can
move the fluid only through interactions between phase frequencies
of different lengths. For instance, the frequencies of
$\psi=\cos x_3+\cos(x_2+x_3)$ have lengths $1$ and $\sqrt2$, and
\begin{equation*}
 R_0(\Delta\psi,\psi)
 =\frac12\Bigl((P_{e_2+2e_3}e_3)\sin(x_2+2x_3)-e_3\sin x_2\Bigr)\neq0.
\end{equation*}
This is, however, only an algebraic property.

The difficulty is to turn this algebraic interaction into a displacement in arbitrarily short time. With the phase amplitude
$\delta^{-1/3}$ used in \Cref{thm:asymp-mixed}, the capillary
contribution over time $\delta$ is of order $\delta^{1/3}$ and
disappears. The amplitude $\delta^{-1/2}$ would formally give a capillary contribution of order one, but the cubic phase term would then
contribute at order $\delta^{-1/2}$. At present, we lack a mechanism
playing the role of the cancellation $B_0(\zeta_u,\zeta_u)=0$ in
\Cref{thm:asymp-mixed}. Its direct analogue, $\Delta(\psi^3)=0$, holds
only for $\psi=0$ among real zero-average periodic profiles. Thus a different cancellation, or a way of compensating the phase evolution without losing the fluid displacement, would be needed.

It would already be interesting to know whether a correction
supported at a single frequency in the fluid equation suffices
if we allow three frequencies in the phase source. For example,
we could consider the $1+3$ choice
\begin{equation*}
 \mathcal K_u=\{e_1\},\qquad
 \mathcal K_\phi=\{e_2,e_3,e_2+e_3\}.
\end{equation*}
The profiles $e_3\cos x_1$ and $e_3\sin x_1$ are shear
flows with zero self-advection. They can introduce dependence
on $x_1$ into the phase through transport, but any new directions
for the velocity would have to be generated through capillarity.

With two frequencies in each equation, a $2+2$ choice would be
\begin{equation*}
 \mathcal K_u=\{e_1,e_1+e_3\},\qquad
 \mathcal K_\phi=\{e_3,e_2+2e_3\}.
\end{equation*}
The profiles for the velocity depend only on $(x_1,x_3)$ and take values in
$\operatorname{span}\{e_1,e_3\}$. Self-advection preserves this
class of two-dimensional flows. The question here is whether
capillarity can produce displacements outside this class in
arbitrarily short time. The cross terms in the capillary force
suggest that this may be possible. For $p=e_3$ and $q=e_2+2e_3$,
linear combinations of the cross terms in
$R_0(\Delta\psi,\psi)$, with $\psi\in\mathcal C_p+\mathcal C_q$,
give the profiles
\begin{equation*}
 (P_{e_2+e_3}e_3)\cos(x_2+x_3),\qquad
 (P_{e_2+e_3}e_3)\sin(x_2+x_3).
\end{equation*}
The coefficients are nonzero because $|q|^2-|p|^2=4$ and
$P_{q-p}p\ne0$. If these displacements could be realized,
we would recover the forcing space used in
\Cref{prop_sat_phase}, up to exchanging $e_1$ and $e_2$. As in the previous example, this computation is only algebraic and we do not yet know how to realize these displacements while keeping
the phase evolution under control.

For either of these choices, or with corrections only in the phase
source, it would suffice to establish small-time accessibility
of the unit ball. The argument of
\Cref{sec:global-prescribed-force} would then give global strong
solutions.

\appendix

\section{Energy estimates and existence of solutions}\label{app:auxiliary-results}
We begin by collecting some auxiliary estimates and then give the proofs of
\Cref{thm:existence-energy-m0,thm:existence-energy-m1}.
\subsection{Preliminary results}
We collect below the interpolation and embedding estimates that will be used throughout the paper.
\begin{lemma}\label{lem:interp-estimates}
Let $\phi\in \dot H^4(\mathbb T^3)$. Then there exists $C>0$ such that
\begin{align}
\|\Delta \phi\|_{L^2}
&\le C\,\|\nabla \phi\|_{L^2}^{2/3}\,\|\Delta^2\phi\|_{L^2}^{1/3},
\label{eq:laplace-interp}\\
\|\nabla \phi\|_{L^\infty}
&\le C\,\|\nabla \phi\|_{L^2}^{1/3}\,\|\Delta^2\phi\|_{L^2}^{2/3},
\label{eq:grad-infty}\\
\|\nabla \phi\|_{L^4}
&\le C\,\|\nabla \phi\|_{L^2}^{3/4}\,\|\Delta^2\phi\|_{L^2}^{1/4},
\label{eq:grad-4}\\
\|\phi\|_{L^\infty}
&\le C\,\|\nabla \phi\|_{L^2}^{1/2}\,\|\Delta\phi\|_{L^2}^{1/2}
\le C\,\|\Delta\phi\|_{L^2}.
\label{eq:phi-infty}
\end{align}
\end{lemma}

\begin{proof}
We begin with \eqref{eq:laplace-interp}. Interpolation between $\dot H^1$ and $\dot H^4$ gives
\[
\|\phi\|_{\dot H^2}
\le C\,\|\phi\|_{\dot H^1}^{2/3}\,\|\phi\|_{\dot H^4}^{1/3},
\]
hence $\|\Delta \phi\|_{L^2}
\le C\,\|\nabla \phi\|_{L^2}^{2/3}\,\|\Delta^2\phi\|_{L^2}^{1/3}$. For \eqref{eq:grad-infty}, 
by the Sobolev embedding $H^2(\mathbb{T}^3) \hookrightarrow L^\infty(\mathbb{T}^3)$, standard $H^k$ interpolation and \eqref{eq:periodic-norm-equivalences} applied to $\phi$ and $\Delta\phi$, we have
\[
\|\nabla\phi\|_{L^\infty}
\le C\,\|\nabla\phi\|_{H^2}
\le C\,\|\nabla\phi\|_{H^1}^{1/2}\,\|\nabla\phi\|_{H^3}^{1/2}
\le C\,\|\Delta\phi\|_{L^2}^{1/2}\,\|\Delta^2\phi\|_{L^2}^{1/2},
\]
and \eqref{eq:laplace-interp} gives \eqref{eq:grad-infty}. For \eqref{eq:grad-4}, interpolation and the Sobolev embedding $H^1(\mathbb T^3)\hookrightarrow L^6(\mathbb T^3)$ give
\begin{equation}\label{eq:paso_interp}
\|\nabla\phi\|_{L^4}
\le C\,\|\nabla\phi\|_{L^2}^{1/4}\,\|\nabla\phi\|_{L^6}^{3/4}
\le C\,\|\nabla\phi\|_{L^2}^{1/4}\,\|\Delta\phi\|_{L^2}^{3/4},
\end{equation}
and we conclude again using \eqref{eq:laplace-interp}. Finally, \eqref{eq:phi-infty} follows from Agmon's inequality 
\[
\|\phi\|_{L^\infty}
\le C\,\|\phi\|_{H^1}^{1/2}\,\|\phi\|_{H^2}^{1/2}
\le C\,\|\nabla\phi\|_{L^2}^{1/2}\,\|\Delta\phi\|_{L^2}^{1/2}.
\]
\end{proof}

\begin{lemma}\label{lem:u-L6-grad-u-L3}
There exists a constant $C>0$ such that the following estimates hold:

\smallskip
\noindent
(i) For every $u\in\mathbf H^1$,
\begin{equation}\label{eq:u-L6}
\|u\|_{L^6} \le C\|\nabla u\|_{L^2}.
\end{equation}

\smallskip
\noindent
(ii) For every $u\in\mathbf H^2$,
\begin{equation}\label{eq:grad-u-L3}
\|\nabla u\|_{L^3}
\le C\|\nabla u\|_{L^2}^{1/2}\|\Delta u\|_{L^2}^{1/2}.
\end{equation}
\end{lemma}
\begin{proof}
Both results follow by applying standard results component by component. Estimate \eqref{eq:u-L6} follows from the Sobolev embedding $H^1(\mathbb T^3)\hookrightarrow L^6(\mathbb T^3)$ and Poincar\'e inequality. Estimate \eqref{eq:grad-u-L3} follows by interpolating between $L^2$ and $L^6$, using the embedding $H^1(\mathbb T^3)\hookrightarrow L^6(\mathbb T^3)$ and \eqref{eq:periodic-norm-equivalences}.
\end{proof}
\begin{rmk}\label{rem:scalar-L6-grad-L3}
Estimates \eqref{eq:u-L6} and \eqref{eq:grad-u-L3} also hold for
scalar functions in $\dot H^1(\mathbb T^3)$ and
$\dot H^2(\mathbb T^3)$, respectively, by the same proof.
\end{rmk}
\begin{lemma}\label{lem:quartic-potential}
For every $r\in\mathbb R$, the potential in \eqref{reg potential}
satisfies
\begin{equation}\label{eq:quartic-potential-bounds}
 \frac18r^4-\frac14\leq F(r)\leq\frac14(r^4+1).
\end{equation}
\end{lemma}
\begin{proof}
The first inequality follows from
$F(r)-\frac18r^4+\frac14=\frac18(r^2-2)^2\geq0$,
and the second from
$\frac14(r^4+1)-F(r)=\frac12r^2\geq0$.
\end{proof}

\subsection{Energy estimates for weak solutions}
Existence of weak solutions can be established by a standard Galerkin
approximation, where the energy estimates are obtained before passing
to the limit. Here we only derive these estimates formally from
\eqref{chn-abstract}, omitting the approximation and compactness arguments.

\begin{proof}[Proof of \Cref{thm:existence-energy-m0}]
Multiplying the velocity equation by $u$ and the
second equation of \eqref{chn-abstract} by $\bar\mu:=\mu-\langle\mu\rangle$, and
integrating over $\mathbb T^3$, we obtain
\begin{align}\label{ener_u}
\frac{1}{2}\frac{\d}{\d t}\|u(t)\|_{L^2}^2
+ \nu\|\nabla u(t)\|_{L^2}^2
&= -(R_0(\Delta \phi,\phi),u)_{L^2}
+ (f_u,u)_{L^2},
\\[0.3ex] \label{ener_phi}
(\partial_t\phi,\bar\mu)_{L^2}
+ (B_1(u,\phi),\bar\mu)_{L^2}
&= (\Delta \mu,\bar\mu)_{L^2}
+ (f_\phi,\bar\mu)_{L^2},
\end{align}
where we have used the cancellation property in
\Cref{lem-cancellation-B0} for the first identity.

Integrating the second equation of \eqref{chn-abstract} over $\mathbb T^3$ and
using periodic boundary conditions, the divergence-free condition
$\nabla\cdot u=0$, and the fact that $f_\phi(\cdot)\in \dot H^1$, we infer
that $\frac{\d}{\d t}\langle \phi(t)\rangle = 0$.
Since $\phi_0\in \dot H^1$, it follows that
$\langle \phi(t)\rangle = 0$ for all $t\geq 0$. In particular, $\langle \partial_t\phi(t)\rangle = 0$, and therefore
\begin{equation}\label{iden_phi_t}
(\partial_t\phi,\bar\mu)_{L^2} = (\partial_t\phi,\mu)_{L^2}.
\end{equation}

Now, using the third equation of \eqref{chn-abstract} and recalling
\eqref{reg potential}, we deduce that
\begin{equation}\label{iden_td_phi}
(\partial_t\phi,\mu)_{L^2}
=
\frac{1}{2}\frac{\d}{\d t}\|\nabla \phi(t)\|_{L^2}^2
+\frac{\d}{\d t}\int_{\TT}F(\phi).
\end{equation}
On the other hand, using the definition of $B_1$ in \eqref{def_b1} and
since $\bar\mu=\mu-\langle \mu\rangle$, we have
\[
(B_1(u,\phi),\bar\mu)_{L^2}
=
(u\cdot\nabla\phi,\mu-\langle \mu\rangle)_{L^2}.
\]
Observe that
\[
(u\cdot\nabla\phi,\langle \mu\rangle)_{L^2}
=
\langle \mu\rangle \int_{\TT} u\cdot\nabla\phi\,dx
=
\langle \mu\rangle \int_{\TT} \nabla\cdot(u\phi)\,dx
=0,
\]
and therefore
\[
(B_1(u,\phi),\bar\mu)_{L^2}
=
(u\cdot\nabla\phi,\mu)_{L^2}.
\]
Using \eqref{iden_chem_pot} and performing integration by parts,
together with the divergence-free condition, we obtain
\begin{equation}\label{eq:cancel_B1_R0}
(B_1(u,\phi),\bar\mu)_{L^2}
=
(u\cdot\nabla\phi,-\Delta\phi+F'(\phi))_{L^2}
=
-(R_0(\Delta\phi,\phi),u)_{L^2}.
\end{equation}
Moreover, by integration by parts and since $\Delta \langle\mu\rangle=0$, we have 
\begin{equation}\label{eq:iden_lap}
(\Delta \mu, \bar\mu)_{L^2}
=
(\Delta \bar\mu, \bar\mu)_{L^2}
=
- \|\nabla \bar\mu\|_{L^2}^2.
\end{equation}
Adding \eqref{ener_u} and \eqref{ener_phi}, and using \eqref{iden_phi_t}--\eqref{eq:iden_lap}, we arrive at the identity
\begin{equation}\label{energy_total}
\frac{1}{2}\frac{\d}{\d t} \Bigl( \|u(t)\|_{L^2}^2 + \|\nabla \phi(t)\|_{L^2}^2 + 2\int_{\TT} F(\phi) \Bigr)
+ \nu\|\nabla u(t)\|_{L^2}^2 + \|\nabla \bar\mu\|_{L^2}^2
= (f_u,u)_{L^2} + (f_\phi, \bar\mu)_{L^2}.
\end{equation}

Since $\langle \bar\mu \rangle=0$, we can use Poincaré inequality together with Cauchy--Schwarz and Young inequalities to deduce that 
\begin{align*}
\frac{1}{2}\frac{\d}{\d t} & \Bigl( \|u(t)\|_{L^2}^2 + \|\nabla \phi(t)\|_{L^2}^2 + 2\int_{\TT} F(\phi) \Bigr)
+ \nu\|\nabla u(t)\|_{L^2}^2 + \|\nabla \bar\mu\|_{L^2}^2
\\
& \leq \|f_u\|_{L^2}\|u\|_{L^2}
+ \|f_\phi\|_{L^2} \|\bar\mu\|_{L^2}
\\
& \leq \frac{1}{2}\|f_u\|_{L^2}^2
+ \frac{1}{2}\|u\|^2_{L^2}
+ C_\epsilon\|f_\phi\|_{L^2}^2
+ \epsilon \|\nabla \bar\mu\|_{L^2}^2,
\end{align*}
for any $\epsilon\in (0,1)$ and some $C_\epsilon>0$. Since $f_\phi(\cdot)\in \dot{H}^1$, we can use once again Poincaré inequality and take $\epsilon$ small enough to obtain
\begin{align}\notag
\frac{\d}{\d t} & \Bigl( \|u(t)\|_{L^2}^2 + \|\nabla \phi(t)\|_{L^2}^2 + 2\int_{\TT} F(\phi) \Bigr)
+ {c\bigl(}\|\nabla u(t)\|_{L^2}^2 + \|\nabla \mu\|_{L^2}^2{\bigr)}
\\ \label{iden_energy_f}
& \leq  C \|u(t)\|^2_{L^2}
+ C\|f\|^2_{\mathbb V^0},
\end{align}
for some constants $c,C>0$.

Integrating \eqref{iden_energy_f} and dropping some of the positive terms, we obtain {for almost every $t\in[0,T]$}
\begin{align}\notag
\|u(t)\|^2_{L^2}
&+\|\phi(t)\|_{\dot{H}^1}^2
+{c}\int_0^{t}\left(\|\nabla u(s)\|^2_{L^2}+\|\nabla \mu(s)\|^2_{L^2}\right)\,ds
\\ \label{est_ff}
&\leq \kappa
+ C\int_0^t\|u(s)\|^2_{L^2}\,ds
+ CT\|f\|^2_{L^\infty(0,T;\mathbb V^0)},
\end{align}
where $\mathcal \kappa$ is a positive constant depending only on $u_0$ and $\phi_0$.
Gr\"onwall inequality allows us to further deduce
\begin{equation}\label{eq:est_linfty_v0}
\|u(t)\|^2_{L^2}
+\|\phi(t)\|^2_{\dot{H}^1}
\leq C\left(1+\|f\|^2_{L^\infty(0,T;\mathbb V^0)}\right),
\quad {\text{for almost every }}t\in[0,T],
\end{equation}
where $C>0$ is a constant depending only on $T$, $u_0$ and $\phi_0$.

By Poincaré inequality, we can improve the norm in the third term on the left-hand side of \eqref{est_ff} to $\int_0^t \|u(s)\|_{H^1}^2ds$, so it remains to incorporate the term $\int_0^t \|\phi(s)\|_{\dot H^2}^2ds$ on the left-hand side. To this end, using the third equation of \eqref{chn-abstract} and integration by parts, we obtain
\begin{align}\notag 
\|\phi\|^2_{\dot{H}^2}
&\leq(\Delta\phi,\Delta \phi)
=(\Delta\phi,\phi^3-\phi-\mu) \\ \notag
& = -3\int_{\TT}\phi^2|\nabla \phi|^2
+\|\nabla \phi\|^2_{L^2}
+(\nabla \phi,\nabla\mu)_{L^2} \\ \label{est_phi_H^2}
& \leq  {\frac{1}{2}}\|\nabla \mu\|^2_{L^2}
+ C \|\nabla \phi\|_{L^2}^2,
\end{align}
where we have dropped the negative term and used the Cauchy--Schwarz and Young's inequalities in the last line.

Integrating \eqref{est_phi_H^2} in time and combining with
\eqref{est_ff} and \eqref{eq:est_linfty_v0}, we deduce that
\begin{align}\label{est_0T_m0}
\|u(t)\|^2_{L^2}
+\|\phi(t)\|^2_{\dot{H}^1}
+\int_0^t\left(\|u(s)\|_{H^1}^2+\|\phi(s)\|^2_{\dot{H}^2}\right)\,ds
\leq C\left(1+\|f\|^2_{L^\infty(0,T;\mathbb V^0)}\right),
\end{align}
for almost every $t\in[0,T]$, where $C>0$ depends only on $T$, $u_0$ and $\phi_0$. This yields \eqref{energy-estimate-m0}.
\end{proof}

\subsection{Local existence of strong solutions}
\label{app:strong-existence}

As for weak solutions, the following estimates hold uniformly for Galerkin approximations. We give the formal energy computations and omit the standard approximation and compactness arguments.

\begin{proof}[Proof of \Cref{thm:existence-energy-m1}]
We begin by multiplying the first equation of system \eqref{chn-abstract} by $A_0u$. Integrating over $\TT$ and integrating by parts, we obtain
\begin{align}\label{ener_u_m1}
\frac{1}{2}\frac{\d}{\d t}\|\nabla u(t)\|_{L^2}^2
+ \nu\|\Delta u\|_{L^2}^2
&= b_0(u,u,\Delta u)
+ (R_0(\Delta \phi,\phi),\Delta u)_{L^2}
- (f_u,\Delta u)_{L^2},
\end{align}

We now estimate the terms on the right-hand side of \eqref{ener_u_m1}. For the convective term, using H\"older inequality together with
\eqref{eq:u-L6}--\eqref{eq:grad-u-L3} in \Cref{lem:u-L6-grad-u-L3}, we obtain
\begin{align}\notag
|b_0(u,u,\Delta u)|
&\le C\|u\|_{L^6}\|\nabla u\|_{L^3}\|\Delta u\|_{L^2}
\\ \notag
&\le C\|\nabla u\|_{L^2}\,
\|\nabla u\|_{L^2}^{1/2}\|\Delta u\|_{L^2}^{1/2}\,
\|\Delta u\|_{L^2}
\\ \notag
&\le C\|\nabla u\|_{L^2}^{3/2}\|\Delta u\|_{L^2}^{3/2}
\\ \label{est_B0}
&\le \gamma\|\Delta u\|_{L^2}^2 + C_\gamma\|\nabla u\|_{L^2}^6
\end{align}
for any $\gamma\in(0,1)$.

For the coupling term, by estimate \eqref{eq:grad-infty} in \Cref{lem:interp-estimates}, we get
\begin{align}\notag 
\left|(R_0(\Delta \phi,\phi),\Delta u)_{L^2}\right|
&\le \|\Delta\phi\|_{L^2}\|\nabla \phi\|_{L^\infty}\|\Delta u\|_{L^2}
\\ \notag
& \le C\|\Delta\phi\|_{L^2}\|\nabla \phi\|_{L^2}^{1/3}\|\Delta^2\phi\|_{L^2}^{2/3}\|\Delta u\|_{L^2}.
\end{align}
Using Young inequalities (with conjugate exponents $(3/2,3)$ and $(4/3,4)$), we deduce that
\begin{align}\notag
\left|(R_0(\Delta \phi,\phi),\Delta u)_{L^2}\right|
&\le \delta \|\Delta^2\phi\|_{L^2}^2
+ C_\delta \|\Delta\phi\|_{L^2}^{3/2}\|\nabla \phi\|_{L^2}^{1/2}\|\Delta u\|_{L^2}^{3/2}
\\ \label{est_R0_Delta-u}
&\le \delta \|\Delta^2\phi\|_{L^2}^2
+ \gamma \|\Delta u\|^2_{L^2}
+ C_{\delta,\gamma}\|\nabla \phi\|^{2}_{L^2}\|\Delta\phi\|_{L^2}^{6}
\end{align}
for any $\delta,\gamma\in (0,1)$.

Finally, for the forcing term we readily have 
\begin{equation}\label{eta1_Deltau}
\left|(f_u,\Delta u)_{L^2}\right|
\le \gamma\|\Delta u\|_{L^2}^2+C_\gamma\|f_u\|_{L^2}^2,
\end{equation}
for any $\gamma\in(0,1)$.

Combining \eqref{ener_u_m1} with \eqref{est_B0}--\eqref{eta1_Deltau}, and fixing $\gamma\in(0,1)$ small enough, we deduce that
\begin{align}\label{eq:est_nabla-ut}
\frac{\d}{\d t}\|\nabla u(t)\|_{L^2}^2
+ {\nu}\|\Delta u\|_{L^2}^2
\le\,
\delta \|\Delta^2\phi\|^2_{L^2}
+ C\|\nabla u\|_{L^2}^6
+ C\|f_u\|_{L^2}^2
+ C_{\delta}\|\nabla \phi\|^{2}_{L^2}\|\Delta\phi\|_{L^2}^{6}.
\end{align}

On the other hand, combining the second and third equations in \eqref{chn-abstract} yields
\begin{equation*}
\partial_t\phi+B_1(u,\phi)=-\Delta^2 \phi+\Delta(\phi^3)-\Delta\phi+f_\phi.
\end{equation*}
Multiplying the above equation by $\Delta^2\phi$ and integrating over $\TT$, we obtain after integration by parts
\begin{align}\notag
\frac{1}{2}\frac{\d}{\d t}\|\Delta\phi(t)\|_{L^2}^2+\|\Delta^2\phi\|^2_{L^2}
&=-\int_{\TT}u\cdot\nabla\phi\, \Delta^2\phi
+\int_{\TT}\Delta(\phi^3)\,\Delta^2\phi
-\int_{\TT}\Delta\phi\,\Delta^2\phi
+\int_{\TT}f_\phi\,\Delta^2\phi
\\ \label{iden_delta_phi}
&=: I_1+I_2+I_3+I_4.
\end{align}
Let us estimate each term $I_i$, $i\in\inter{1,4}$. For the first one, we have
\begin{equation*}
|I_1|
=
\left|\int_{\TT}u\cdot\nabla\phi \,\Delta^2\phi\right|
\le \|u\cdot\nabla \phi\|_{L^2}\|\Delta^2\phi\|_{L^2}
\le \|u\|_{L^2}\|\nabla \phi\|_{L^\infty}\|\Delta^2\phi\|_{L^2}.
\end{equation*}
Using estimate \eqref{eq:grad-infty} in \Cref{lem:interp-estimates} and Young inequality, we obtain
\begin{align}\notag
|I_1|
&\le C\|u\|_{L^2}\|\nabla \phi\|_{L^2}^{1/3}\|\Delta^2\phi\|_{L^2}^{5/3}
\\ \notag
&\le \epsilon \|\Delta^2\phi\|_{L^2}^2
+ C_{\epsilon}\|u\|_{L^2}^6\|\nabla \phi\|_{L^2}^2
\\ \label{est_I1}
&\le \epsilon \|\Delta^2\phi\|_{L^2}^2
+ C_{\epsilon}\|u\|_{L^2}^6\|\Delta \phi\|_{L^2}^2,
\end{align}
for any $\epsilon\in(0,1)$, where in the last line we have used Poincar\'e inequality.

For $I_2$, we note that $\Delta(\phi^3)=6\phi|\nabla\phi|^2+3\phi^2\Delta\phi$, so that, by Cauchy--Schwarz and Young inequalities, we deduce
\begin{align}\label{est_I2_ini}
|I_2|
=\left|\int_{\TT}\Delta(\phi^3)\Delta^2\phi\right|
\leq C_{\epsilon}\left(\left\|\phi|\nabla\phi|^2\right\|_{L^2}^2 + \|\phi^2\Delta\phi\|^2_{L^2}\right)+\epsilon\|\Delta^2\phi\|^2_{L^2}.
\end{align}
Let us estimate the two terms on the right-hand side of \eqref{est_I2_ini}. For the first one, using estimates \eqref{eq:grad-4} and \eqref{eq:phi-infty} in \Cref{lem:interp-estimates}, together with Young inequality, we get
\begin{align}\notag 
\left\|\phi|\nabla\phi|^2\right\|_{L^2}^2
&\leq \|\phi\|_{L^\infty}^2\|\nabla \phi\|_{L^4}^4
\\ \notag
&\leq C\left(\|\nabla \phi\|_{L^2}\|\Delta\phi\|_{L^2}\right)
\left(\|\nabla \phi\|_{L^2}^{3}\|\Delta^2\phi\|_{L^2}\right)
\\ \label{est_I21}
&\leq \epsilon \|\Delta^2\phi\|_{L^2}^2+C_{\epsilon}\|\nabla\phi\|_{L^2}^8\|\Delta\phi\|^2_{L^2}
\end{align}
for any $\epsilon\in(0,1)$.

Similarly, for the second one, using \eqref{eq:phi-infty} and \eqref{eq:laplace-interp}, we obtain
\begin{align}\notag
\|\phi^2\Delta\phi\|_{L^2}^2
&\le \|\phi\|_{L^\infty}^4\|\Delta\phi\|_{L^2}^2
\le C \|\nabla \phi\|_{L^2}^2\|\Delta\phi\|_{L^2}^4
\\ \notag
&\le C \|\nabla \phi\|_{L^2}^4\|\Delta\phi\|_{L^2}\|\Delta^2\phi\|_{L^2}
\\ \label{est_I22}
&\le \epsilon \|\Delta^2\phi\|_{L^2}^2
+ C_{\epsilon}\|\nabla\phi\|_{L^2}^8\|\Delta\phi\|_{L^2}^2,
\end{align}
for any $\epsilon\in(0,1)$. Therefore, combining \eqref{est_I21} and \eqref{est_I22} with \eqref{est_I2_ini}, we deduce that
\begin{equation}\label{est_final_I2}
|I_2|\leq 3\epsilon \|\Delta^2\phi\|^2_{L^2}+C_{\epsilon}\|\nabla\phi\|_{L^2}^8\|\Delta \phi\|_{L^2}^2
\end{equation}
for any $\epsilon\in(0,1)$. 

For the term $I_3$, we readily have
\begin{equation}\label{est_I3}
|I_3|
=\left|\int_{\TT}\Delta\phi \,\Delta^2\phi\right|
\le \epsilon \|\Delta^2\phi\|_{L^2}^2
+ C_{\epsilon}\|\Delta\phi\|_{L^2}^2.
\end{equation}
Similarly, we estimate $I_4$ as
\begin{equation}\label{est_I4}
|I_4|
=\left|\int_{\TT}f_\phi\,\Delta^2\phi\right|
\le \epsilon \|\Delta^2\phi\|_{L^2}^2
+ C_\epsilon\|f_\phi\|_{L^2}^2.
\end{equation}

Putting together \eqref{iden_delta_phi} with \eqref{est_I1}, \eqref{est_final_I2}--\eqref{est_I4}, and taking $\epsilon\in(0,1)$ small enough, we obtain
\begin{equation}\label{est_delta_phi}
\frac{\d}{\d t}\|\Delta\phi(t)\|_{L^2}^2
+\|\Delta^2\phi\|_{L^2}^2
\le
C\Bigl(1+\|u(t)\|_{L^2}^6+\|\nabla\phi(t)\|_{L^2}^8\Bigr)\|\Delta\phi(t)\|_{L^2}^2
+ C \|f_\phi(t)\|_{L^2}^2,
\end{equation}
for some constant $C>0$.

Let us define
\begin{equation*}
\mathcal Y(t):=\|\nabla u(t)\|_{L^2}^2+{\|\Delta\phi(t)\|_{L^2}^2}.
\end{equation*}
Adding estimates \eqref{eq:est_nabla-ut} and \eqref{est_delta_phi}, and taking
$\delta\in(0,1)$ small enough, we obtain {for some $c>0$}
\begin{align}\label{dif_hr}
\frac{\d}{\dt}\mathcal Y(t)
+{c\bigl(}\|\Delta u(t)\|_{L^2}^2+{\|\Delta^2\phi(t)\|_{L^2}^2\bigr)}
&\le C\|\nabla u(t)\|_{L^2}^6
+ C\|u(t)\|_{L^2}^6\|\Delta\phi(t)\|_{L^2}^2
\\ \notag
&\quad
+ C\|\nabla \phi(t)\|_{L^2}^8\|\Delta \phi(t)\|^2_{L^2}
+ C\|\nabla \phi(t)\|^{2}_{L^2}\|\Delta\phi(t)\|_{L^2}^{6}
\\ \notag
&\quad {+ C\|\Delta\phi(t)\|_{L^2}^2}
+ C\|f(t)\|_{\mathbb V^0}^2.
\end{align}

We now use the estimate obtained in the case $m=0$ (see \Cref{thm:existence-energy-m0}). More precisely, from
\eqref{est_0T_m0} we know that, for every $T>0$, there exists a constant $K_T>0$, depending on $u_0,\phi_0,T$ and $f$, such that
\begin{equation}\label{eq:est-kT}
\|u(t)\|_{L^2}^2+\|\phi(t)\|_{\dot H^{1}}^2
+\int_0^t\Big(\|u(s)\|_{H^{1}}^2+\|\phi(s)\|_{\dot H^{2}}^2\Big)\d{s}
\le K_T,
\qquad t\in[0,T].
\end{equation}

Using this bound, we estimate the terms on the right-hand side of \eqref{dif_hr}.
First, by definition of $\mathcal Y$,
\begin{equation*}
\|\nabla u(t)\|_{L^2}^6 \le \mathcal Y(t)^3.
\end{equation*}
Next, using that $\|u(t)\|_{L^2}^2\le K_T$ and $\|\nabla\phi(t)\|_{L^2}^2\le K_T$, we obtain
\begin{equation*}
\|u(t)\|_{L^2}^6\|\Delta\phi(t)\|_{L^2}^2
\le K_T^3\,\mathcal Y(t),
\qquad
\|\nabla \phi(t)\|_{L^2}^8\|\Delta \phi(t)\|_{L^2}^2
\le K_T^4\,\mathcal Y(t),
\end{equation*}
and
\begin{equation*}
\|\nabla \phi(t)\|^{2}_{L^2}\|\Delta\phi(t)\|_{L^2}^{6}
\le K_T\,\mathcal Y(t)^3.
\end{equation*}
Therefore, using also that $f\in L^\infty(0,T;\mathbb V^0)$, we deduce from \eqref{dif_hr} that
\begin{equation*}%\label{ode_Y}
\frac{\d}{\dt}\mathcal Y(t)
\le C_T\bigl(1+\mathcal Y(t)+\mathcal Y(t)^3\bigr),
\qquad t\in[0,T],
\end{equation*}
where $C_T>0$ depends on $u_0,\phi_0,T$ and $f$. A standard comparison argument for ordinary differential equations then yields the existence of a time $T_*\in(0,T]$, depending only on $\mathcal Y(0)$ and $C_T$, such that
\begin{equation*}
\sup_{t\in[0,T_*)}\mathcal Y(t)
\le
C\bigl(\mathcal Y(0),\|f\|_{L^\infty(0,T;\mathbb V^0)},T_*\bigr).
\end{equation*}
Finally, integrating \eqref{dif_hr} over $(0,t)$ with $t\in[0,T_*)$, and using the previous bound, we deduce that
\begin{equation*}%\label{eq:energy-m1-local}
\|z(t)\|_{\mathbb V^1}^2+{\int_0^t \|z(s)\|_{\mathbf H^2\times\dot H^4}^2\,ds}
\le
C\bigl(\|z_0\|_{\mathbb V^1},\|f\|_{L^\infty(0,T;\mathbb V^0)},T_*\bigr),
\qquad \forall\, t\in[0,T_*).
\end{equation*}
This yields the desired estimate. Fix $S\in(0,T_*)$. The bound above gives
\begin{equation}\label{eq:strong-spatial-regularity}
\begin{aligned}
 u&\in L^\infty(0,S;\mathbf H^1)\cap L^2(0,S;\mathbf H^2),\\
 \phi&\in L^\infty(0,S;\dot H^2)\cap L^2(0,S;\dot H^4).
\end{aligned}
\end{equation}
Using $H^1\hookrightarrow L^6$, $H^2\hookrightarrow L^\infty$,
and the algebra property of $H^2$, we have
\begin{align}
 \|B_0(u,u)\|_{L^2}
 &\leq \|(u\cdot\nabla)u\|_{L^2}
 \leq \|u\|_{L^\infty}\|\nabla u\|_{L^2}\leq C\|u\|_{H^2}\|u\|_{H^1},\label{eq:strong-convection-L2}\\
 \|R_0(\Delta\phi,\phi)\|_{L^2}
 &\leq \|\Delta\phi\|_{L^2}\|\nabla\phi\|_{L^\infty}
 \leq C\|\phi\|_{H^2}\|\phi\|_{H^3},\notag\\
 \|B_1(u,\phi)\|_{L^2}
 &\leq \|u\|_{L^6}\|\nabla\phi\|_{L^3}
 \leq C\|u\|_{H^1}\|\phi\|_{H^2},\notag\\
 \|\Delta(\phi^3)\|_{L^2}
 &\leq C\|\phi^3\|_{H^2}
 \leq C\|\phi\|_{H^2}^3.\notag
\end{align}
By \eqref{eq:strong-spatial-regularity}, all four terms belong to
$L^2(0,S;L^2)$. The equations in \eqref{chn-abstract} therefore yield
\begin{equation}\label{eq:strong-time-regularity}
 \partial_tu,\ \partial_t\phi\in L^2(0,S;L^2).
\end{equation}
Applying \cite[Theorem~II.5.14]{BoyerFabrie2013} with the spaces
$(H^2,L^2)$ for $u$ and $(H^4,L^2)$ for $\phi$,
\eqref{eq:strong-spatial-regularity} and \eqref{eq:strong-time-regularity}
give $z\in C([0,S];\mathbb V^1)$.
Since $S<T_*$ is arbitrary, this completes the proof.
\end{proof}
\begin{rmk}\label{rem:pressure-recovery}
For a weak solution $(u,\phi)$, de Rham's theorem gives a pressure $p$
such that the momentum equation in \eqref{chn} holds in the sense of distributions;
see e.g. \cite[Section~IV.2]{BoyerFabrie2013}.
The pressure is determined up to an additive function of
time and is uniquely determined by imposing
$\langle p(t)\rangle=0$.
%Adding a function of time to $p$ leaves the equation unchanged so we can choose $p$ uniquely such that $\langle p(t)\rangle=0$.
For a strong solution, \eqref{eq:phi-infty} gives
$\|\mu(t)\|_{L^2}\leq C(\|\phi(t)\|_{H^2}+\|\phi(t)\|_{H^2}^3)$,
so $\mu\in L^\infty(0,T;L^2)$.
Estimate \eqref{eq:grad-infty} in \Cref{lem:interp-estimates}, together
with $\phi\in L^\infty(0,T;H^2)\cap L^2(0,T;H^4)$, also gives
$\nabla\phi\in L^2(0,T;L^\infty)$, and therefore
$\|\mu\nabla\phi\|_{L^2(0,T;L^2)}<\infty$.
The bound for $(u\cdot\nabla)u$ in \eqref{eq:strong-convection-L2}
also gives $(u\cdot\nabla)u\in L^2(0,T;L^2)$.
Applying $I-\Pi$ to the first equation of \eqref{chn} gives
$\nabla p\in L^2(0,T;L^2)$ and by normalization,
$p\in L^2(0,T;\dot H^1)$.
Thus $(u,\phi,\mu,p)$ satisfies the original system \eqref{chn}
in distributions, and almost everywhere when $(u,\phi)$ is a strong solution.
\end{rmk}
\section{Estimates for the asymptotic limit}\label{app:short-time-estimates}
We collect some identities and estimates used in the proof of \Cref{thm:asymp-mixed}. We begin by presenting two results concerning differentiation and integration by parts.
\begin{lemma}\label{lem:pointwise_leibniz_identities}
Let $a,b\in C^\infty(\mathbb T^3;\mathbb R^3)$ and
$q\in C^\infty(\mathbb T^3)$. Let
\begin{equation*}
 (\nabla a)_{\ell j}=\partial_\ell a_j,
 \qquad \ell,j\in\{1,2,3\}.
\end{equation*}
For matrices $A,B\in\mathbb R^{3\times3}$, we denote
\begin{equation*}
 A:B:=\sum_{\ell,j=1}^3 A_{\ell j}B_{\ell j},
 \qquad |A|^2:=A:A.
\end{equation*}
Then
\begin{equation}\label{eq:leibniz_transport_vector}
 \nabla\big((a\cdot\nabla)b\big)
 ={\nabla a\cdot\nabla b}+(a\cdot\nabla)\nabla b,
\end{equation}
and
\begin{equation}\label{eq:leibniz_laplacian_transport}
 \Delta\big((a\cdot\nabla)q\big)
 =\Delta a\cdot\nabla q
  +2\nabla a:\nabla\nabla q
  +a\cdot\nabla\Delta q,
\end{equation}
where $\nabla a:\nabla\nabla q
 =\sum_{\ell=1}^3
        \partial_\ell a\cdot\partial_\ell\nabla q
 =\sum_{\ell,j=1}^3
        \partial_\ell a_j\,\partial_{\ell j}q$.

Finally, for $G\in C^\infty(\mathbb T^3;\mathbb R^{3\times3})$,
we have
\begin{equation}\label{eq:pointwise_transport_energy}
 (a\cdot\nabla)G:G=\frac12a\cdot\nabla|G|^2.
\end{equation}
The same identity holds for scalar and vector functions, with
ordinary multiplication and the scalar product, respectively.
\end{lemma}
\begin{proof}
For \eqref{eq:leibniz_transport_vector}, we compute componentwise:
\begin{align*}
\partial_\ell\big((a\cdot\nabla)b_k\big)
&=
\partial_\ell
\left(
\sum_{j=1}^3 a_j\partial_j b_k
\right) =
\sum_{j=1}^3
\partial_\ell a_j\,\partial_j b_k
+
\sum_{j=1}^3
a_j\,\partial_j\partial_\ell b_k.
\end{align*}
This is precisely $\partial_\ell\big((a\cdot\nabla)b_k\big)
=
(\nabla a\cdot\nabla b)_{\ell k}
+
\big((a\cdot\nabla)\nabla b\big)_{\ell k}$. Hence \eqref{eq:leibniz_transport_vector} follows.

For \eqref{eq:leibniz_laplacian_transport}, we compute
\begin{align*}
\Delta\big((a\cdot\nabla)q\big)
&=
\sum_{\ell=1}^3
\partial_{\ell\ell}
\left(
\sum_{j=1}^3 a_j\partial_j q
\right) =
\sum_{j=1}^3
(\Delta a_j)\partial_j q
+
2\sum_{\ell,j=1}^3
\partial_\ell a_j\,\partial_{\ell j}q
+
\sum_{j=1}^3
a_j\partial_j\Delta q \\
&=
\Delta a\cdot\nabla q
+
2\nabla a:\nabla\nabla q
+
a\cdot\nabla\Delta q.
\end{align*}
This proves \eqref{eq:leibniz_laplacian_transport}.

Finally, for every smooth scalar function $g$, chain rule gives $(a\cdot\nabla g)g=\frac12a\cdot\nabla(g^2)$. Applying this identity to the entries of $G$, we obtain
\begin{align*}
 (a\cdot\nabla)G:G
 =\sum_{\ell,k=1}^3(a\cdot\nabla G_{\ell k})G_{\ell k}=\frac12a\cdot\nabla
       \left(\sum_{\ell,k=1}^3G_{\ell k}^2\right)
 =\frac12a\cdot\nabla|G|^2.
\end{align*}
The proof is complete.
\end{proof}
The following integral identities use periodicity and, for the
transport cancellation, the divergence-free condition.
\begin{lemma}\label{lem:transport_integration_by_parts}
Let $a,b,c\in C^\infty(\mathbb T^3;\mathbb R^3)$. Then
\begin{equation}\label{eq:ibp_vector_transport}
\int_{\mathbb T^3}
(a\cdot\nabla)b\cdot\Delta c\,dx
=
-
\int_{\mathbb T^3}
\nabla\big((a\cdot\nabla)b\big):\nabla c\,dx.
\end{equation}
Moreover, if $\diver a=0$, then, for every
$G\in C^\infty(\mathbb T^3;\mathbb R^{3\times3})$,
\begin{equation}\label{eq:transport_skew_integral}
\int_{\mathbb T^3}
(a\cdot\nabla)G:G\,dx
=
0.
\end{equation}
The same identity holds for scalar and vector functions, with
ordinary multiplication and the scalar product, respectively.
\end{lemma}
\begin{proof}
For \eqref{eq:ibp_vector_transport}, using integration by parts on the torus,
we obtain
\begin{align*}
\int_{\mathbb T^3}
(a\cdot\nabla)b\cdot\Delta c\,dx
&=
\sum_{k=1}^3
\int_{\mathbb T^3}
(a\cdot\nabla)b_k\,\Delta c_k\,dx =
-
\sum_{\ell,k=1}^3
\int_{\mathbb T^3}
\partial_\ell\big((a\cdot\nabla)b_k\big)\,
\partial_\ell c_k\,dx \\
&=
-
\int_{\mathbb T^3}
\nabla\big((a\cdot\nabla)b\big):\nabla c\,dx.
\end{align*}
For \eqref{eq:transport_skew_integral}, we use
\eqref{eq:pointwise_transport_energy} from
\Cref{lem:pointwise_leibniz_identities} and integrate by parts:
\begin{align*}
 \int_{\mathbb T^3}(a\cdot\nabla)G:G\,dx
 &=\frac12\int_{\mathbb T^3}a\cdot\nabla|G|^2\,dx = -\frac12\int_{\mathbb T^3}(\diver a)|G|^2\,dx
 =0.
\end{align*}
The proof is complete.
\end{proof}

The following results are used in the proof of \Cref{thm:asymp-mixed}. 

\begin{lemma}\label{lem:transport_fluid_first_argument}
	Let $v\in \mathbf{H}^2$ and $\zeta_u\in \mathbf{H}^3$. For any $\delta\in(0,1)$ the following inequality holds
	\begin{equation*}
		\delta^{1/3} \left| b_0(\zeta_u,v,\Delta v)\right| \leq C\delta^{1/3}\|\zeta_u\|_{H^3}\|\nabla v\|_{L^2}^2,
	\end{equation*}
where $C>0$ is independent of $\delta$, $\zeta_u$, and $v$.
\end{lemma}
\begin{proof}
By \eqref{eq:ibp_vector_transport} in \Cref{lem:transport_integration_by_parts} and identity \eqref{eq:leibniz_transport_vector} in \Cref{lem:pointwise_leibniz_identities}, we have
\begin{equation*}
\begin{split}
	\delta^{1/3}\left|b_0(\zeta_u,v,\Delta v)\right| &= \delta^{1/3}\left| \int_{\TT} \nabla\left( (\zeta_u \cdot \nabla )v\right) : \nabla v\right| \\
	&\leq \delta^{1/3}\left|\int_{\TT} (\nabla\zeta_u \nabla v): \nabla v\right| + \delta^{1/3}\left|\int_{\TT}(\zeta_u\cdot\nabla)\nabla v:\nabla v\right|.
\end{split}
\end{equation*}
The second term on the right-hand side of the above equation is equal to zero by \Cref{eq:transport_skew_integral} since $\nabla\cdot\zeta_u=0$. Thus, we only need to estimate the first term. Using H\"older inequality and the fact that $H^2(\TT) \hookrightarrow L^\infty(\TT)$ we obtain
\begin{equation}\label{eq:lemma1first}
	\delta^{1/3}\left|\int_{\TT} (\nabla\zeta_u\nabla v): \nabla v\right| \leq \delta^{1/3}\|\nabla\zeta_u\|_{L^\infty}\|\nabla v\|_{L^2}^2 \leq \delta^{1/3}\|\zeta_u\|_{H^3}\|\nabla v\|_{L^2}^2 .
\end{equation}
The result then follows. 
\end{proof}
When the profile occurs in the second argument, we use the following
estimate.
\begin{lemma}\label{lem:transport_fluid_second_argument}
	Let $v\in \mathbf{H}^2$ and $\zeta_u\in \mathbf{H}^4$. For any $\delta\in(0,1)$ the following inequality holds
	\begin{equation*}
		\delta^{1/3} \left| b_0(v,\zeta_u,\Delta v)\right| \leq C\delta^{1/3}\|\zeta_u\|_{H^4}\|\nabla v\|_{L^2}^2,
	\end{equation*}
where $C>0$ is independent of $\delta$, $\zeta_u$, and $v$.
\end{lemma}
\begin{proof}
By \eqref{eq:ibp_vector_transport} in \Cref{lem:transport_integration_by_parts} and identity \eqref{eq:leibniz_transport_vector} in \Cref{lem:pointwise_leibniz_identities}, we have
\begin{equation*}
	\begin{split}
		\delta^{1/3}\left|b_0(v,\zeta_u,\Delta v)\right| &= \delta^{1/3}\left| \int_{\TT} \nabla\left( (v \cdot \nabla )\zeta_u\right) : \nabla v\right| \\
		&\leq \delta^{1/3}\left|\int_{\TT} (\nabla v \nabla \zeta_u): \nabla v\right| + \delta^{1/3}\left|\int_{\TT}(v\cdot\nabla)\nabla \zeta_u:\nabla v\right|.
	\end{split}
\end{equation*}
Note that the first term can be estimated as \eqref{eq:lemma1first}. For the second term, we use H\"older inequality, Poincaré inequality and the fact that $H^2(\TT)\hookrightarrow L^\infty(\TT)$ to obtain
\begin{equation}\label{eq:lemma2second}
	\delta^{1/3}\left|\int_{\TT}(v\cdot\nabla)\nabla \zeta_u:\nabla v\right| \leq \delta^{1/3}\|v\|_{L^2}\|\nabla^2 \zeta_u\|_{L^\infty}\|\nabla v\|_{L^2} \leq \delta^{1/3}\| \zeta_u\|_{H^4}\|\nabla v\|_{L^2}^2.
\end{equation}
From here, we can deduce the desired result. 
\end{proof}
For the phase equation, we estimate the transport term tested against
$\Delta^2\varphi$.
\begin{lemma}\label{lem:transport_phase_high_order}
Let $\zeta_u\in\mathbf H^4$ and $\varphi\in\dot H^4$. Then, for every
$\delta\in(0,1)$, we have
\begin{equation*}
\delta^{1/3}
\left|
b_1(\zeta_u,\varphi,\Delta^2\varphi)
\right|
\leq
C\delta^{1/3}
\norm{\zeta_u}_{H^4}
\norm{\varphi}_{\dot H^2}^2,
\end{equation*}
where $C>0$ is independent of $\delta$, $\zeta_u$, and $\varphi$.
\end{lemma}
\begin{proof}
Recall that $b_1(\zeta_u,\varphi,\Delta^2\varphi)
=
\int_{\mathbb T^3}
(\zeta_u\cdot\nabla)\varphi\,\Delta^2\varphi\,dx$. Since $\Delta$ is self-adjoint, we get
\begin{align*}
b_1(\zeta_u,\varphi,\Delta^2\varphi)
=
\left(
(\zeta_u\cdot\nabla\varphi),
\Delta^2\varphi
\right)_{L^2}=
\left(
\Delta\big((\zeta_u\cdot\nabla)\varphi\big),
\Delta\varphi
\right)_{L^2}.
\end{align*}

By formula \eqref{eq:leibniz_laplacian_transport}, we have
\begin{align}\label{eq:iden_z1_deltaphi}
b_1(\zeta_u,\varphi,\Delta^2\varphi)
=
\left(
\Delta\zeta_u\cdot\nabla\varphi,
\Delta\varphi
\right)_{L^2}
+
2\left(
\nabla\zeta_u:\nabla\nabla\varphi,
\Delta\varphi
\right)_{L^2} +
\left(
\zeta_u\cdot\nabla\Delta\varphi,
\Delta\varphi
\right)_{L^2}.
\end{align}
The last term in \eqref{eq:iden_z1_deltaphi} vanishes since $\zeta_u$ is divergence-free and \eqref{eq:transport_skew_integral}. For the other two, direct applications of H\"older inequality, yield
\begin{align*}
\left|
b_1(\zeta_u,\varphi,\Delta^2\varphi)
\right|
&\leq
\norm{\Delta\zeta_u}_{L^\infty}
\norm{\nabla\varphi}_{L^2}
\norm{\Delta\varphi}_{L^2} + C\norm{\nabla\zeta_u}_{L^\infty}
\norm{\nabla^2\varphi}_{L^2}
\norm{\Delta\varphi}_{L^2}.
\end{align*}
Using that $\zeta_u\in\mathbf H^4$, $H^2(\TT)\hookrightarrow L^\infty(\TT)$, Poincaré inequality, and \eqref{eq:periodic-norm-equivalences} we have
\begin{equation*}
\left|
b_1(\zeta_u,\varphi,\Delta^2\varphi)
\right|
\leq
C\norm{\zeta_u}_{H^4}
\norm{\varphi}_{\dot H^2}^2.
\end{equation*}
Multiplying by $\delta^{1/3}$ the result then follows. 
\end{proof}

%\printbibliography
\bibliographystyle{alpha}
\small{\bibliography{bib_nonlinear_coupled}}

\medskip
\medskip
\medskip
\medskip

\begin{flushleft}

\bigskip

\textbf{Víctor Hernández-Santamaría \orcidlink{0000-0002-2294-1376}}\\
Departamento de Matem\'aticas, Facultad de Ciencias\\
Universidad Nacional Autónoma de México \\
Circuito Exterior, C.U.\\
04510, Coyoacán, CDMX, Mexico\\
\texttt{victor.santamaria@ciencias.unam.mx}

\bigskip
\bigskip
	
	\textbf{Subrata Majumdar \orcidlink{0000-0001-6724-6943}}\\
Centre For Applicable Mathematics \\	
Tata Institute of Fundamental Research\\
	Post Bag No 6503, Sharada Nagar \\
	Bengaluru 560065, India\\
	\texttt{subrata26@tifrbng.res.in}\\

\bigskip
\bigskip
	
		\textbf{Yingying Wu-Zhang \orcidlink{0009-0001-8606-4536}}\\
	Departamento de Matemáticas \\
	Universidad de Sonora \\
	Blvd. Luis Encinas y Rosales S/N \\
	83000, Hermosillo, Sonora, México \\
	\texttt{yingying.wu@unison.mx}\\
	
\end{flushleft}

\end{document}